%% file: confspace.tex
\documentclass[reqno,12pt]{amsbook}
\usepackage{amsmath,amsthm,amscd,amssymb}
\numberwithin{equation}{section}
\makeatletter

\def\th@plain{%
\let\thm@indent\noindent
\thm@headfont{\bfseries}
\let\thmhead\thmhead@plain
\let\swappedhead\swappedhead@plain
\thm@preskip.5\baselineskip\@plus.2\baselineskip\@minus.2\baselineskip
\thm@postskip\thm@preskip
\slshape
}

\def\th@remark{%
\let\thm@indent\noindent
\thm@headfont{\bfseries}
\let\thmhead\thmhead@plain
\let\swappedhead\swappedhead@plain
\thm@preskip.5\baselineskip\@plus.2\baselineskip\@minus.2\baselineskip
\thm@postskip\thm@preskip
\upshape
}

\makeatother

\theoremstyle{plain}
\newtheorem{Theorem}{Theorem}[section]
\newtheorem{Corollary}[Theorem]{Corollary}
\newtheorem{Lemma}[Theorem]{Lemma}
\newtheorem{Lemma-Definition}[Theorem]{Lemma-Definition}
\newtheorem{Proposition}[Theorem]{Proposition}

\theoremstyle{remark}
\newtheorem{Remark}[Theorem]{\bfit Remark}

\newtheorem{Definition}[Theorem]{\bfit Definition}
\newtheorem{Definition-Notation}[Theorem]{\bfit Definition-Notation}
\newtheorem{Example}[Theorem]{\bfit Example}

\newtheorem{Example-Definition}[Theorem]{\bfit Example-Definition}
\newtheorem{Notation}[Theorem]{\bfit Notation}

\newcommand{\1}{\ensuremath{1\mspace{-4.85mu}{\mathrm I}}}
\newcommand{\Def}{\stackrel{\mathrm{def}}{=\!\!=}}

\input cd
\LaTeXfeatures

\DeclareMathOperator{\Aff}{{\mathbf{Aff}}}      

\DeclareMathOperator{\Aut}{{Aut}} 
  
\DeclareMathOperator{\const}{const}

\DeclareMathOperator{\Hol}{\text{\bfit{Hol}}}
\DeclareMathOperator{\hol}{\text{\bfit{hol}}}
\DeclareMathOperator{\Hom}{Hom}

\DeclareMathOperator{\Id}{Id}
\DeclareMathOperator{\id}{id}
\DeclareMathOperator{\Img}{Im}

\DeclareMathOperator{\Ker}{Ker}

\DeclareMathOperator{\modl}{mod}

\DeclareMathOperator{\rank}{rank}

\DeclareMathOperator{\supp}{supp}
\DeclareMathOperator{\Sym}{Sym}
\DeclareMathOperator{\thh}{{\!}^{th}}
\DeclareMathOperator{\Tors}{Tors}

\font\caps=cmcsc10
\font\bcaps=cmcsc10 scaled \magstep1

\font\bfit=cmbxti10 scaled \magstep1

\font\twelvegtc=eufm10 scaled 1200

\font\ninegtc=eufm9
\font\sevengtc=eufm7

\newfam\gtcfam

\textfont\gtcfam=\twelvegtc
\scriptfont\gtcfam=\ninegtc
\scriptscriptfont\gtcfam=\sevengtc

\def\modo#1{\left| #1 \right|}

\def\eqnum#1{\eqno (#1)}

\begin{document}

\title{\vskip40pt
Configuration spaces of ${\mathbb C}$ and ${\mathbb{CP}}^1$:\\
some analytic properties
}

\author{Vladimir Lin\\
\vskip130pt

\begin{small} 
\hskip-68pt
Department of Mathematics 
                      \hskip140pt Max-Planck-Institut\hskip40pt\\
\hskip-88pt
Technion-Israel Institute of \hskip145pt f{\" u}r Mathematik \hskip141pt\\

\hskip-100pt
Technology         \hskip221pt Vivatsgasse 7  \hskip127pt\\

\hskip-104pt
32000 Haifa         \hskip218pt 53111 Bonn  \hskip133pt\\

\vskip10pt

\hskip-118pt
Israel         \hskip249pt Germany  \hskip117pt\\

\vskip200pt

MPI03-98
\end{small}
}


\thanks{This work was started in 2000 and completed in
2003, both times during my stay at the Max-Planck-Institut f\"ur Mathematik 
in Bonn. I am deeply grateful to MPIM for hospitality and
especially for giving me a fortunate opportunity to meet in 2000  
Marat Gizatullin, discuss with him some examples and learn from him
about the Eisenstein automorphism of the space of non-degenerate binary
cubic forms. This renovated my old interest to configuration
spaces and stimulated the present work.}

\vskip-20pt

\begin{abstract}{\ \ \vskip0.1cm

\hfill {\em I mean... \ \ You know...}
\vskip0.1cm
\noindent 
We study certain analytic properties of the ordered
and unordered configuration spaces
${\mathcal C}_o^n(X)=\{(q_1,...,q_n)\in X^n\,|\ 
q_i\ne q_j \ \forall\,i\ne j\}$ and
${\mathcal C}^n(X)=\{Q\subset X\,|\ \#Q=n\}$ of simply connected
algebraic curves $X={\mathbb C}$ and $X={\mathbb{CP}}^1$.  
\vskip0.1cm

\noindent For $n\ge 3$ the braid group $B_n(X)=\pi_1({\mathcal C}^n(X))$
of the curve $X$ is non-abelian. A map $F\colon{\mathcal C}^n(X)
\to{\mathcal C}^k(X)$ is called {\em cyclic} if the subgroup
$F_*(\pi_1({\mathcal C}^n(X)))\subset\pi_1({\mathcal C}^k(X))$
is cyclic; otherwise $F$ is called {\em non-cyclic}.
\vskip0.1cm

\noindent The diagonal $\Aut X$ action in $X^n$ induces $\Aut X$ actions
in ${\mathcal C}_o^n(X)$ and ${\mathcal C}^n(X)$.
A morphism (holomorphic or regular)
$F\colon{\mathcal C}^n(X)\to{\mathcal C}^k(X)$ is called {\em orbit-like}
if $F({\mathcal C}^n(X))\subset(\Aut X)Q^*$ for some
$Q^*\in{\mathcal C}^k(X)$. An endomorphism $F$ of ${\mathcal C}^n(X)$
is called {\em tame} if there is a morphism
$T\colon{\mathcal C}^n(X)\to\Aut X$ such that $F(Q)=T(Q)Q$ for all
$Q\in{\mathcal C}^n(X)$. Tame endomorphisms preserve
each $\Aut X$ orbit in ${\mathcal C}^n(X)$. We prove 
that {\sl for $n>4$ and $k\ge t(X):=\dim_{\mathbb C}\Aut X$
an endomorphism $F$
of ${\mathcal C}^n(X)$ is tame if and only if it is non-cyclic, and
a morphism $F\colon{\mathcal C}^n(X)\to{\mathcal C}^k(X)$
is orbit-like if and only if it is cyclic}
(``Tame Map Theorem" and ``Cyclic Map Theorem"). 
\vskip0.1cm

\noindent To study non-cyclic maps we establish first certain algebraic
properties of the braid groups $B_n(X)$, including the stability
of the pure braid group $PB_n(X)\subset B_n(X)$ under all non-cyclic
endomorphisms of $B_n(X)$. This implies that any non-cyclic endomorphism
of ${\mathcal C}^n(X)$ lifts to an ${\mathbf S}(n)$ equivariant endomorphism
of ${\mathcal C}_o^n(X)$. To study the latter one, we define a new
invariant of affine varieties $Z$. We constract a finite simplicial complex
$L_\Delta(Z)$ whose set of vertices $L(Z)$ consists
of all non-constant holomorphic functions $Z\to{\mathbb C}\setminus\{0,1\}$.
In many interesting cases the correspondence $Z\mapsto L_\Delta(Z)$ is a
contravariant functor (this is surely the case when restricting to dominant
morphisms). We describe the simplices of the complex
$L_\Delta({\mathcal C}_o^n(X))$ explicitly,
compute $\dim_{\mathbb R}L_\Delta({\mathcal C}_o^n(X))$
and show that any ${\mathbf S}(n)$
equivariant endomorphism of ${\mathcal C}_o^n(X)$ produces the corresponding
simplicial self-map of $L_\Delta({\mathcal C}_o^n(X))$. We determine the
${\mathbf S}(n)$ orbits of simplices 
$\Delta\in L_\Delta({\mathcal C}_o^n(X))$
and find their normal forms with respect to this action. Eventually, this
provides a complete description of ${\mathbf S}(n)$ equivariant
endomorphisms of ${\mathcal C}_o^n(X)$ and proves Tame Map Theorem.
This theorem shows that
{\sl ${\mathcal C}^n(X)/\Aut{\mathcal C}^n(X)={\mathcal C}^n(X)/\Aut X$
and the holomorphic homotopy classes of non-cyclic endomorphisms
of ${\mathcal C}^n(X)$ are in a natural $1-1$ correspondence with
the homotopy classes of all continuous maps
${\mathcal C}^n(X)\to{\mathcal K}(X)$,
where ${\mathcal K}(X)$ is a maximal compact subgroup of $\Aut X$.}
It also implies that {\sl $\dim_{\mathbb C}F({\mathcal C}^n(X))
\ge n-t(X)+1$ for any $n>4$ and any non-cyclic endomorphism $F$
of ${\mathcal C}^n(X)$.}
\vskip0.1cm

\noindent The study of cyclic maps involves the Liouville property
of ${\mathcal C}^n(X)$ and the holomorphic direct decomposition
${\mathcal C}_o^k(X)\cong\Aut X
\times{\mathcal C}_o^m({\mathbb C}\setminus\{0,1\})$, where
$m=k-t(X)$ and ${\mathcal C}_o^m({\mathbb C}\setminus\{0,1\})$ is 
the ordered configuration space of ${\mathbb C}\setminus\{0,1\}$.
\vskip0.1cm

\noindent We prove that {\sl for $n>t(X)+1$ and any non-cyclic morphism
$F\colon{\mathcal C}^n(X)\to{\mathcal C}^k(X)$
there is a point $Q\in{\mathcal C}^n(X)$ such that
$Q\cap F(Q)\ne\varnothing$} (``Linked Map Theorem").
The proof involves an ``explicit" description of the holomorphic universal
covering map $\Pi\colon\widetilde{\Aut X}
\times {\mathbf T}(0,m)\to {\mathcal C}^n(X)$,
where $m=n+3-t(X)$, $\widetilde{\Aut X}$ is the universal cover of $\Aut X$
and ${\mathbf T}(0,m)$ is the Teichm{\" u}ller space of type $(0,m)$.
The construction of $\Pi$ involves the universal
Teichm{\" u}ller family ${\mathbf{V'}}(0,m)\to{\mathbf T}(0,m)$. 
To prove Linked Map Theorem, we use the Hubbard-Earl-Kra theorem,
which says that the universal Teichm{\" u}ller family of type $(g,m)$
has no holomorphic sections if $\dim_{\mathbb C}{\mathbf T}(g,m)>1$. 
\vskip0.1cm

\noindent The above results apply to a version of the
$13$th Hilbert problem. 
We prove that for $n>\max\{4,m+2\}$ the algebraic function
$[x:y]=U_n(z)$ on ${\mathbb{CP}}^n$ defined by the ''universal" equation
$z_0x^n+z_1x^{n-1}y+...+z_ny^n=0$ cannot be represented as a composition
of algebraic functions of $m$ variables in such a way that 
a representing function-composition $F\supseteq U_n$ and
$U_n$ itself have the same branch loci.
}
\end{abstract}

\maketitle

\tableofcontents

\markboth{
\bcaps {Vladimir Lin}}
{\textsc {Configuration spaces 
of ${\mathbb C}$ and ${\mathbb{CP}}^1$}}




\section{Introduction}
\label{Sec: Introduction}

\subsection{Configuration spaces}\label{Ss: Configuration spaces}
The (unordered) $n\thh$ {\em configuration space} 
\index{Configuration spaces\hfill|phantom\hfill}
\index{Configuration spaces!ordered\hfill}
\index{Configuration spaces!unordered\hfill}
${\mathcal C}^n(X)$
of a space $X$ consists of all $n$ point subsets $Q\subseteq X$. 
It may be also viewed as the {\em regular orbit space} of 
the natural ${\mathbf S}(n)$ action in $X^n$. The standard
left action of the symmetric group ${\mathbf S}(n)$ on $X^n$
defined by 
\begin{equation}\label{eq: S(n) action}
(\sigma,q)\mapsto\sigma q
=(q_{\sigma^{-1}(1)},...,q_{\sigma^{-1}(n)})\quad
\forall \, \sigma\in{\mathbf S}(n) \ \text{and} \ 
\forall\,q=(q_1,...,q_n)\in X^n
\end{equation}
is free on the {\em ordered configuration space}
\begin{equation}\label{eq: ordered configuration space} 
{\mathcal C}_o^n(X)=\{q=(q_1,...,q_n) \in X^n\,|\ q_i\ne q_j \ 
\forall\,i\ne j\}\,; 
\end{equation}
this provides the unbranched ${\mathbf S}(n)$ Galois covering
\begin{equation}\label{eq: Co to C covering}
\aligned
&p\colon{\mathcal C}_o^n(X)\to{\mathcal C}^n(X)
={\mathcal C}_o^n(X)/{\mathbf S}(n)\,,\\
&p(q)=p(q_1,...,q_n)=\{q_1,...,q_n\}=Q\in{\mathcal C}^n(X)\quad
\forall\,q=(q_1,...,q_n)\in{\mathcal C}_o^n(X)\,.
\endaligned
\end{equation}
If $X$ is a complex or complex algebraic manifold then
both ${\mathcal C}_o^n(X)$ and ${\mathcal C}^n(X)$ are so.
\vskip0.2cm

\noindent This paper is devoted to the configuration spaces 
${\mathcal C}^n({\mathbb C})$ and ${\mathcal C}^n({\mathbb{CP}}^1)$,
which are irreducible complex affine algebraic manifolds of
dimension $n$; ${\mathcal C}^n(\mathbb C)$ may be viewed
as a Zariski open subset of $\mathcal C^n(\mathbb{CP}^1)$.
We are especially interested in the morphisms of these spaces
in the analytic or algebraic category.
\vskip0.2cm

\noindent Throughout the paper, unless otherwise specified,
we assume that $n\ge 3$ and $X$ is one of the curves 
${\mathbb C}$, ${\mathbb {CP}}^1$. 

\subsection{The fundamental groups $\pi_1({\mathcal C}^n(X))$
and $\pi_1({\mathcal C}_o^n(X))$}
\label{Ss: The fundamental groups of Cn(X) and Con(X)}
The fundamental group $\pi_1({\mathcal C}^n(X))$
and its normal subgroup $\pi_1({\mathcal C}_o^n(X))$ 
are quite the classical objects known as the {\em braid group}
\index{Braid groups\hfill|phantom\hfill}
$B_n(X)$ and the {\em pure braid group} $PB_n(X)$ respectively. 
They fit into the exact sequence
\begin{equation}\label{eq: main exact sequence}
1\to PB_n(X)\overset{p_*}{\longrightarrow}B_n(X)
\overset{\mu}{\longrightarrow}{\mathbf S}(n)\to 1
\end{equation}
related to the Galois covering (\ref{eq: Co to C covering});
the epimorphism $\mu\colon B_n(X)\to{\mathbf S}(n)$ is 
called {\em standard}.
In more details, $\pi_1({\mathcal C}^n({\mathbb C}))=B_n$ and 
$\pi_1({\mathcal C}^n_o({\mathbb C}))=PB_n$ are the classical Artin
braid group and the pure Artin braid group respectively.
\index{Braid groups!Artin braid group\hfill}
\index{Braid groups!Artin braid group!pure Artin braid group\hfill}
Similarly, $\pi_1({\mathcal C}^n({\mathbb{PC}^1}))=B_n(S^2)$ and 
$\pi_1({\mathcal C}^n_o({\mathbb{PC}^1}))=PB_n(S^2)$ are the 
braid group and the pure braid group of the two dimensional
sphere $S^2$. 
\index{Braid groups!sphere braid group\hfill}
\index{Braid groups!sphere braid group!pure sphere braid group\hfill}
Notice that for $n\ge 3$
the braid groups $B_n$ and $B_n(S^2)$ are non-abelian.
Their abelianizations are cyclic: $B_n/B'_n\cong{\mathbb Z}$ and
$B_n(S^2)/B'_n(S^2)\cong{\mathbb Z}/(2n-2){\mathbb Z}$.

\subsection{Tame, degenerate tame and orbit-like morphisms}
\label{Ss: Tame and degenerate tame morphisms}
Let $\Aut X$ be the automorphism group of $X$; thus,
$\Aut X$ is either the group $\Aff({\mathbb C})$ 
of all affine transformations of ${\mathbb C}$
or the group ${\mathbf{PSL}}(2,{\mathbb C})$
of all M{\"o}bius transformations of the Riemann sphere.
The diagonal $\Aut X$ action in ${\mathcal C}^n(X)$
\begin{equation}\label{eq: diagonal Aaut X action}
{\mathcal C}^n(X)\ni Q=\{q_1,...,q_n\}\mapsto\{Aq_1,...,Aq_n\}\Def AQ
\in{\mathcal C}^n(X)\quad\forall\, A\in\Aut X
\end{equation}
gives rise to the simplest examples of automorphisms 
of ${\mathcal C}^n(X)$: every $A\in\Aut X$ induces the automorphism 
${\mathcal C^n}(X)\ni Q\mapsto AQ\in{\mathcal C^n}(X)$.
These examples may be modified by letting $A=A(Q)$ depend
on a point $Q\in{\mathcal C}^n(X)$ holomorphically 
and applying $A(Q)$ either to $Q$ itself or to a chosen and marked 
point $Q^*\in{\mathcal C}^n(X)$. To be more precise, we give the
following definition.

\begin{Definition}\label{Def: tame, degenerate tame and orbit-like maps}
For any morphism
$T\colon{\mathcal C^n}(X)\to\Aut X$ we define the corresponding
{\em tame endomorphism}
\index{Tame endomorphism\hfill}
\index{Endomorphism!tame\hfill}
$F_T\colon{\mathcal C^n}(X)\to{\mathcal C^n}(X)$ by
\begin{equation}\label{eq: tame map}
F_T(Q)=T(Q)Q=\{T(Q)q_1,....,T(Q)q_n\}\ \
\forall\, Q=\{q_1,....,q_n\}\in {\mathcal C}^n(X)\,.
\end{equation}
Choosing a point $Q^*=\{q_1^*,....,q_k^*\}\in{\mathcal C}^k(X)$,
we define also the {\em degenerate tame}
\index{Morphism!degenerate tame\hfill}
morphism
$F_{T,Q^*}\colon{\mathcal C^n}(X)\to{\mathcal C^k}(X)$,
\begin{equation}\label{eq: degenerate tame morphism}
F_{T,Q^*}(Q)=T(Q)Q^*=\{T(Q)q_1^*,....,T(Q)q_k^*\}\ \
\forall\,Q=\{q_1,....,q_n\}\in {\mathcal C}^n(X)\,;
\end{equation}
if $k=n$, we obtain a degenerate tame endomorphism
of ${\mathcal C^n}(X)$. 
Furthermore, a morphism $F\colon{\mathcal C^n}(X)\to{\mathcal C^k}(X)$
is said to be {\em orbit-like}
\index{Morphism!orbit-like\hfill}
if its image $F({\mathcal C^n}(X))$
is contained in an orbit of $\Aut X$ action in ${\mathcal C^k}(X)$.
\hfill $\bigcirc$
\end{Definition}

\begin{Remark}\label{Rmk: tame maps keep each AutX orbit}
{\bfit a}$)$ Of course, any degenerate tame morphism is orbit-like. 
Since the stabilizer $St_{Q^*}=\{A\in\Aut X\,|\ AQ^*=Q^*\}$
of a generic point $Q^*\in{\mathcal C}^k(X)$
is trivial, an orbit-like morphism whose image is contained
in the orbit of such a point $Q^*$ is degenerate tame. 
\vskip0.2cm

{\bfit b}$)$ We may extend the above definition to continuous maps
by saying that a continuous map $T\colon{\mathcal C^n}(X)\to\Aut X$
provides the tame continuous map 
$F_T\colon{\mathcal C}^n(X)\to{\mathcal C}^n(X)$ 
and the degenerate tame continuous map
$F_{T,Q^*}\colon{\mathcal C}^n(X)\to{\mathcal C}^k(X)$
defined by (\ref{eq: tame map}) and (\ref{eq: degenerate tame morphism})
respectively. Any tame map $F$ keeps
each orbit of the $\Aut X$ action in ${\mathcal C}^n(X)$ stable, 
that is, $F((\Aut X)Q)\subseteq (\Aut X)Q$ for all $Q\in{\mathcal C}^n(X)$.
Any degenerate tame map is orbit-like, that is, it 
carries the whole space ${\mathcal C}^n(X)$ into a single $\Aut X$
orbit in ${\mathcal C}^k(X)$. These facts are of a little moment since
tame and degenerate tame maps are the rarities in the class of all
continuous maps of configuration spaces. However, due to Tame Map Theorem
formulated below, a ''generic" endomorphism is tame and hence it
keeps each $\Aut X$ orbit stable. Furthermore, by a quite elementary reason,  
any morphism ${\mathcal C}^n(X)\to{\mathcal C}^k(X)$ must
carry $\Aut X$ orbits to $\Aut X$ orbits; see Section
\ref{Ss: morphisms and Aut X orbits}, 
Proposition \ref{Prp: morphisms respect orbits}.
\hfill $\bigcirc$
\end{Remark}

\subsection{Cyclic homomorphisms and cyclic maps}
\label{Ss: Cyclic homomorphisms and cyclic maps}
It will be convenient to introduce the following notions.

\begin{Definition}\label{Def: cyclic and abelian homomorphisms and maps}
A group homomorphism $\varphi\colon\, G\to H$ is said to be
{\em cyclic}, {\em abelian} or {\em solvable} if its image
$\varphi(G)$ is a cyclic, abelian or solvable subgroup of the group $H$
respectively.
\index{Homomorphism!cyclic\hfill}
\index{Homomorphism!abelian\hfill}
\index{Homomorphism!solvable\hfill}
A homomorphism $\varphi\colon\,G\to H$ 
is abelian if and only if its restriction 
to the commutator subgroup $G'=[G,G]$ of the group $G$ is trivial.
If the abelianization $G/G'$ of $G$ 
is a cyclic group then any abelian homomorphism of $G$ is cyclic.

A continuous map $F\colon \mathcal X\to \mathcal Y$ of 
arcwise connected spaces is called {\em cyclic},
{\em abelian} or {\em solvable} if the induced homomorphism
$F_*\colon\pi_1(\mathcal X)\to\pi_1(\mathcal Y)$ 
of the fundamental groups is cyclic, abelian or solvable respectively. 
\index{Map!abelian\hfill}
\index{Map!cyclic\hfill}
\index{Map!solvable\hfill}
\hfill $\bigcirc$
\end{Definition}

\subsection{Main results}\label{Ss: Main results}
\noindent Any degenerate tame morphism is cyclic since
$\pi_1(\Aut{\mathbb C})\cong{\mathbb Z}$ and 
$\pi_1(\Aut{\mathbb{CP}}^1)\cong{\mathbb Z}/2{\mathbb Z}$.
In fact, we will show in Section 
\ref{Sec: Morphisms Cn(X) to Ck(X) and Aut X orbits.
Cyclic and orbit-like maps}
that for $n>4$ a morphism ${\mathcal C}^n(X)\to{\mathcal C}^k(X)$
is cyclic if and only if it is orbit-like.
\vskip0.2cm

\noindent In contrast to orbit-like morphisms, for $n\ge 3$
any tame endomorphism of ${\mathcal C}^n(X)$ is non-cyclic
(see Remark \ref{Rmk: equivariance}$(${\bfit b}$)$).
Our first main result tells that the for $n>4$ the converse
is also true.

\begin{Theorem}[{\caps Tame Map Theorem}]
\label{Thm: Tame Map Thm} 
\index{Theorem!Tame Map Theorem\hfill}
For $n>4$ every non-cyclic endomorphism $F$ of ${\mathcal C}^n(X)$
is tame. That is, there exists a unique morphism 
$T\colon{\mathcal C}^n(X)\to\Aut X$ such that 
$$
F(Q)=F_T(Q)\Def T(Q)Q\quad\text{\rm for all} \ 
Q\in{\mathcal C}^n(X)\,.
$$
In particular, for $n>4$ every automorphism of ${\mathcal C}^n(X)$
is tame. 
\hfill $\bigcirc$
\end{Theorem}

\noindent In view of the above notes about cyclic morphisms,
this theorem shows that for endomorphisms of ${\mathcal C}^n(X)$
the properties to be tame or orbit-like form a dichotomy, which 
actually coincides with the dichotomy {\em ``non-cyclic -- cyclic"}:

\begin{Corollary}\label{Crl: dichotomy}
For $n>4$ the classes of non-cyclic and cyclic
endomorphisms of ${\mathcal C}^n(X)$ coincide with the
classes of tame and orbit-like endomorphisms respectively.
\hfill $\bigcirc$
\end{Corollary}

\noindent Tame Map Theorem has also the following two immediate
consequences.

\begin{Corollary}\label{Crl: homotopy classification} 
For $n>4$ the holomorphic homotopy classes of non-cyclic holomorphic
maps\footnote{That is, the classes of holomorphic maps homotopic
to each other within the space of all holomorphic maps.}
${\mathcal C}^n(X)\to{\mathcal C}^n(X)$ are in the natural
one-to-one correspondence with the holomorphic homotopy classes
of holomorphic  maps ${\mathcal C}^n(X)\to\Aut X$.
By the famous Grauert theorem, the latter classes are in the
natural one-to-one correspondence with the usual homotopy classes
of all continuous maps  ${\mathcal C}^n(X)\to\Aut X$, or, which is the
same, with the homotopy classes
of all continuous maps ${\mathcal C}^n(X)\to{\mathcal K}(X)$,
where ${\mathcal K}(X)$ is a maximal compact subgroup
of the group $\Aut X$.
\hfill $\bigcirc$
\end{Corollary}

\begin{Corollary}\label{Crl: orbit space}
Let $n>4$ and $G=\Aut{\mathcal C}^n(X)$\,.
Then the orbits of the natural $G$ action in ${\mathcal C}^n(X)$
coincide with the orbits of the diagonal $\Aut X$ action 
in ${\mathcal C}^n(X)$.
\vskip0.2cm

\noindent In particular, the orbit space 
${\mathcal C}^n({\mathbb{CP}}^1)/\Aut{\mathcal C}^n({\mathbb{CP}}^1)$
may be identified to the orbit space
${\mathcal C}^n({\mathbb{CP}}^1)/{\mathbf{PSL}}(2,{\mathbb C})$,
which, in turn, may be viewed as the moduli space $M(0,n)$ 
\index{Moduli space\hfill|phantom\hfill}
of the Riemann sphere with $n$ unordered punctures. 
\hfill $\bigcirc$
\end{Corollary}

\begin{Remark}\label{Rmk: holomorphic homotopy classes}
For $X={\mathbb C}$ one may take ${\mathcal K}({\mathbb C})$
to be the subgroup ${\mathbb T}\subset\Aff{\mathbb C}$
consisting of all rotations $z\mapsto e^{it}z$,
$t\in{\mathbb R}$. Thereby, the holomorphic homotopy classes of non-cyclic
holomorphic maps ${\mathcal C}^n({\mathbb C})\to{\mathcal C}^n({\mathbb C})$
are in the one-to-one correspondence with the homotopy classes
of all continuous maps ${\mathcal C}^n({\mathbb C})\to{\mathbb T}$,
which, in turn, may be identified with the elements of the group
$H^1({\mathcal C}^n({\mathbb C}),{\mathbb Z})\cong\Hom(B_n,{\mathbb Z})
\cong{\mathbb Z}$.
\vskip0.2cm

\noindent For $X={\mathbb{CP}}^1$ one may take 
${\mathcal K}({\mathbb{CP}}^1)\subset{\mathbf{PSL}}(2,{\mathbb C})$
to be the subgroup isomorphic to ${\mathbf{SO}}(3)$ and
consisting of all M{\"o}bius transformations of the form
$$
z\mapsto\frac{az - \bar c}{cz + \bar a}\,,\quad |a|^2 + |c|^2 = 1\,.
$$
Thus, the holomorphic homotopy classes of non-cyclic holomorphic
endomorphisms 
of ${\mathcal C}^n({\mathbb{CP}}^1)$ 
are in the one-to-one correspondence with the homotopy classes
of all continuous maps  ${\mathcal C}^n({\mathbb{CP}}^1)
\to{\mathbf{SO}}(3)$.
\hfill $\bigcirc$
\end{Remark}

%


\begin{Remark}\label{Rmk: non-degenerate forms}
Consider the space of all non-zero complex 
binary forms of degree $n$
\index{Binary form\hfill\phantom\hfill}
\index{Binary form!discriminant of\hfill}
\index{Binary form!non-degenerate\hfill}
\index{Binary form!non-degenerate!projective\hfill}
\begin{equation}\label{binary form}
\phi(x,y;z)=z_0x^n+z_1x^{n-1}y+\cdots+z_{n-1}xy^{n-1}+z_n y^n\,,
\end{equation}
where $z$ is regarded either as a non-zero point
$(z_0,\cdots ,z_n)\in{\mathbb C}^{n+1}$ or as a point
$[z_0:\cdots :z_n]\in{\mathbb{CP}}^n$. In the latter case
the form $\phi(x,y;z)$ is called {\em projective}, 
meaning that proportional forms are identified.
The discriminant $D_n(z)$ of a form $\phi(x,y;z)$
is a homogeneous polynomial in $z_0,...,z_n$ of degree $2(n-1)$:
$$
D_n(z)\!=\!\left|\!\!\begin{array}{cccccccccc}
1 & z_1 & z_2 & ...  & z_{n-2} & z_{n-1} & z_n & 0  & ... & 0\\
0 & z_0 & z_1 & ...  & ...     & ... & z_{n-1} & z_n  & .... & 0\\
... & ...& ... & ... & ... & ... & ... & ... & ... & ... \\
0 & 0 & 0 & ... & z_0 & z_1 & z_2 & .... & z_{n-1} & z_n\\
n & (n-1)z_1 & (n-2)z_2 & ...  & 2z_{n-2} & z_{n-1} & 0 & 0  & ... & 0\\
0 & nz_0 & (n-1)z_1 & ...  & ...     & 2z_{n-2} & z_{n-1} & 0  & .... & 0\\
... & ...& ... & ... & ... & ... & ... & ... & ... & ... \\
0 & 0 & 0 & ... & 0 & nz_0 & (n-1)z_1 & ... & 2z_{n-2} & z_{n-1} 
\end{array}
\!\!\right|;
$$ 
$\phi(x,y;z)$ is said to be {\em non-degenerate}
if $D_n(z)\ne 0$. The set 
\begin{equation}\label{eq: Fn}
{\mathcal F}^n=\{\phi=\phi(x,y;z)\,|\  
z=[z_0:...:z_n]\in{\mathbb{CP}}^n\,, \ \ D_n(z)\ne 0\}
\end{equation}
is called {\em the space of all 
non-degenerate complex projective binary forms of degree} $n$.  
For every $Q=\{q_1,...,q_n\}\in{\mathcal C}^n(\mathbb {CP}^1)$ 
there is a unique form $\phi=\phi(x,y;z)\in{\mathcal F}^n$
whose zero set $Z_\phi=\{[x:y]\in{\mathbb{CP}^1}\,|\ 
\phi(x,y;z)=0\}$ coincides with $Q$.
This provides the natural identification
${\mathcal C}^n(\mathbb {CP}^1)={\mathcal F}^n$.
\vskip0.2cm

\noindent Under this identification, the subspace 
${\mathcal C}^n(\mathbb C)\subset{\mathcal C}^n(\mathbb {CP}^1)$ 
consists of all $\phi([x:y],z)\in{\mathcal F}^n$ with $z_0\ne 0$; 
it may be viewed as the space ${\mathbf G}_n$ of all complex polynomials 
$p_n(t;w)=t^n+w_1t^{n-1}+...+w_n$ with simple roots.
We also regard ${\mathbf G}_n$ as a domain
in ${\mathbb C}^n$ defined by
\begin{equation}\label{eq: domain Gn}
{\mathbf G}_n=\{w=(w_1,...,w_n)\in{\mathbb C}^n\,|\ d_n(w)\ne 0\}\,,
\end{equation}
where $d_n(w)=d_n(w_1,...,w_n)=D_n(1,w_1,...,w_n)$ is the
discriminant of $p_n(t;w)$.
The discriminant map 
\begin{equation}\label{eq: discriminant map}
d_n\colon{\mathcal C}^n(\mathbb C)
={\mathbf G}_n\ni w\mapsto d_n(w)\in{\mathbb C}^*
\end{equation}
is a holomorphic locally trivial fiber bundle;
the induced homomorphism of the fundamental
groups $d_{n*}$ coincides with the abelianization homomorphism
$\chi\colon B_n\to{\mathbb Z}$.  
\vskip0.2cm

\noindent In terms of binary forms, a tame endomorphism of 
${\mathcal C}^n({\mathbb {CP}}^1)={\mathcal F}^n$ is just a 
projective change of variables $[x:y]\mapsto T(\phi)[x:y]$,
which depends holomorphically on a form $\phi$ itself.
A tame endomorphism of the space 
${\mathcal C}^n({\mathbb C})={\mathbf G}_n$ 
is the transformation of polynomials $p_n(t;w)$
corresponding to an affine change 
of the variables $T\colon t\mapsto t'=a(w)t+b(w)$ holomorphically 
depending on $w\in{\mathbf G}_n$.
\hfill $\bigcirc$
\end{Remark}

\noindent Our second main result concerns morphisms
of configuration spaces which do not necessarily coincide to each other. 
It will be convenient to make the following definition and notation.

\begin{Definition}\label{Def: disjoint and linked maps of Cn}
Let $F\colon{\mathcal C}^n(X)\to{\mathcal C}^k(X)$
be a map of configuration spaces. Then
either $Q\cap F(Q)=\varnothing$ for all $Q\in{\mathcal C}^n(X)$
or there is a point $Q_0\in{\mathcal C}^n(X)$ such that
$Q_0\cap F(Q_0)\ne\varnothing$. In the first case $F$
is called {\em disjoint}; in the second case we say that
$F$ is {\em linked}.
\index{Map!disjoint\hfill}
\index{Map!linked\hfill}
\end{Definition}

\begin{Notation}\label{Not: t(X)}
We set $t(X)=\dim_{\mathbb C}\Aut X$, i. e.,
$t(X)=2$ for $X={\mathbb C}$ and $t(X)=3$ for $X={\mathbb{CP}}^1$.
\index{$t(X)$\hfill}
\hfill $\bigcirc$
\end{Notation}

\noindent In all dimensions there are smooth maps
${\mathcal C}^n(X)\to{\mathcal C}^k(X)$ of both kinds.
For holomorphic maps the situation changes drastically.

\begin{Theorem}[{\caps Linked Map Theorem}]
\label{LinkMapThm}
\index{Theorem!Linked Map Theorem\hfill}
For $n>t(X)+1$ and $k\ge 1$ every holomorphic map 
$F\colon{\mathcal C}^n(X)\to{\mathcal C}^k(X)$
is linked.
\hfill $\bigcirc$
\end{Theorem}

\noindent This theorem is actually a reflection of the fact that 
the universal Teichm{\" u}ller family of type $(0,n+3)$ does not
possess holomorphic sections. In Section 
\ref{Sec: Holomorphic universal covers} 
we will relate configuration spaces to Teichm{\" u}ller spaces. The 
proof of Linked Map Theorem  and some its applications 
to algebraic functions will be postponed till Section 
\ref{Sec: Linked Map Theorem}. Examples of disjoint holomorphic
maps in small dimension will be presented in Section
\ref{Sec: Configuration spaces of small dimension and some other
examples}.
In Section \ref{Sec: Some applications to algebraic functions}
we will discuss some applications of Linked Map Theorem
to the 13th Hilbert problem for algebraic functions.

\subsection{Some historical remarks}
\label{Ss: Some historical remarks}
My interest to holomorphic mappings ${\mathbf G}_n\to{\mathbf G}_k$
of polynomials with simple roots was motivated 
by their relations to a version
of the 13th Hilbert problem for algebraic functions.
Endomorphisms of ${\mathcal C}^n(\mathbb C)={\mathbf G}_n$
were completely described in the 1970's 
(see my papers \cite{Lin72b,Lin72c,Lin79}, where some results
about morphsms ${\mathbf G}_n\to{\mathbf G}_k$ 
were obtained as well).
However, the complete proofs were never published, since 
they contained some very long and unpleasant combinatorial computations.  
\vskip0.2cm

\noindent In the 1990's I found that certain
morphisms ${\mathbf G}_n\to{\mathbf G}_k$ are related to the
universal Teichm{\" u}ller family of type $(0,n+1)$, which
led to the proof of Linked Map Theorem for $X={\mathbb C}$ 
in the form presented in Remark
\ref{Rmk: linked maps}$(${\bfit a}$)$ (see \cite{Lin96b}).
\vskip0.2cm

\noindent On the other hand, in \cite{Lin96b}
(see also Section \ref{Ss: G3 to G6} of the present paper) a disjoint
polynomial map $F\colon{\mathcal C}^3({\mathcal C})
\to{\mathcal C}^6({\mathcal C})$ has been constructed.
In view of the result mentioned above, the latter map is so special that
its novelty must look dubious. During my stay at Max-Planck-Institut
f{\" u}r Mathematik in Bonn in 2000, I met Marat Gizatullin
(whose knowledge of and love to the classical mathematical literature 
are exceptional) and asked him whether he have seen such 
a map somewhere. After a few days he came back and told that he did not
meet such a monster but instead he found out that in $1844$
Gotthold Eisenstein
\index{Eisenstein\hfill|phantom\hfill}
has discovered a very exciting
automorphism of non-degenerate binary cubic forms (see \cite{Eis1844}
and Section \ref{Sec: Eisenstein automorphism} of the present paper).
\vskip0.2cm

\noindent This Eisenstein's example renovated
my interest to endomorphisms of configuration spaces. 
In the present work my old approach is improved and simplified 
and now it applies to both ${\mathcal C}^n(\mathbb C)$ and 
${\mathcal C}^n(\mathbb {CP}^1)$.
\vskip0.2cm

\noindent It should be also mentioned that endomorphisms
of the configuration spaces ${\mathcal C}^n({\mathbb C}^*)$
have been described by V. Zinde \cite{Zin77}.
\vskip0.2cm

\noindent Recently my student Y. Feler proved that for $n>4$ all
automorphisms of the configuration space
${\mathcal C}^n({\mathbb T}^2)$ of a torus
${\mathbb T}^2={\mathbb C}/(a{\mathbb Z}+b{\mathbb Z})$
are tame. Together with my old
results quoted above, results of V. Zinde, \cite{Zin77},
and the ones of the present paper, this gives a realistic
hope to complete the description of endomorphisms of configuration 
spaces for all non-hyperbolic Riemann surfaces.
Furthermore, Feler studied certain configuration spaces of different
kind, say the spaces ${\mathcal C}^n_{aff}({\mathbb C}^2)$ and 
${\mathcal C}^n_{prj}({\mathbb{CP}}^2)$ consisting of
subsets $Q=\{q_1,...,q_n\}$ in ${\mathbb C}^2$ 
or in ${\mathbb{CP}}^2$ whose points $q_1,...,q_n$ 
are in general affine and general projective 
position respectively. (My attention to such spaces 
was attracted by A. M. Vershik.) 


\section{Strategy of the proof of Tame Map Theorem}
\label{Sec: Strategy of the proof of Tame Map Theorem}
\index{Theorem!Tame Map Theorem!strategy of the proof\hfill}

\noindent The question about the nature of endomorphisms 
of configuration spaces belongs to analytic or algebraic
geometry and seemingly is a rather complicated one. 
We will reduce it, step by step, to a series 
of algebraic or combinatorial tasks,
which look and as a matter of fact are simpler. 
\vskip0.2cm

\noindent In this section we try first to find out what the inevitable
attributes of tame endomorphisms and those of certain accompanying objects
ought to be. Then we formulate the statements which should be proved
in order to ensure that non-cyclic endomorphisms and the corresponding
objects of their ``suite" actually possess such properties.
The accomplishment of this approach in the consequent
sections eventually will lead to the desired result.
\vskip0.2cm

\noindent Some proofs named {\bfit ``Explanation"} also are placed
here but only when they are completely trivial. 
We start with an obvious observation.
\vskip0.2cm

\noindent{\bfit Lifting and invariance properties.}
Tame and degenerate tame endomorphisms $F$ of
${\mathcal C^n}(X)$ have in common the following 
\vskip0.1cm

\noindent{\bfit Lifting property:}
\index{Lifting property\hfill}
{\sl there exists an endomorphism
$f\colon{\mathcal C}_o^n(X)\to{\mathcal C}_o^n(X)$ which fits into
the commutative diagram}
\begin{equation}\label{CD: lifting diagram}
\CD
{{\mathcal C}_o^n(X)} @ > {f} >> {{\mathcal C}_o^n(X)}\\
@V{p}VV @VV{p}V\\ 
{{\mathcal C}^n(X)} @ >> {F}>{{\mathcal C}^n(X)}
\endCD
\end{equation}

\noindent Such an $f$ may be defined as follows:
\begin{equation}\label{eq: lifting of tame map}
f(q)=f_T(q)\Def T(Q)(q)=(T(Q)q_1,...,T(Q)q_n) \ \
\text{whenever} \ F=F_T\,,
\end{equation}
and
\begin{equation}\label{eq: covering degenerate tame map}
f(q)=f_{T,Q^*}(q)\Def T(Q)(q^*)=(T(Q)q_1^*,...,T(Q)q_n^*) \ \
\text{whenever} \ F=F_{T,Q^*}\,, 
\end{equation}
where $q=(q_1,...,q_n)\in{\mathcal C}_o^n(X)$, \
$Q=\{q_1,...,q_n\}=p(q)\in{\mathcal C}^n(X)$. 
\hfill $\bigcirc$
\vskip0.2cm

\noindent By the covering map theorem,
a continuous map $F\colon{\mathcal C}^n(X)\to{\mathcal C}^n(X)$
possesses a continuous covering map $f$ that fits into
diagram (\ref{CD: lifting diagram}) if and only if
the following invariance condition is fulfilled:
\vskip0.2cm

\noindent{\bfit Invariance property:}
\index{Invariance property\hfill}
{\sl the induced endomorphism 
$F_*\colon\pi_1({\mathcal C}^n(X))\to\pi_1({\mathcal C}^n(X))$
of the fundamental group respects the normal subgroup 
$\pi_1({\mathcal C}_o^n(X))\vartriangleleft\pi_1({\mathcal C}^n(X))$,
that is,}
\begin{equation}\label{eq: conservation property}
F_*(\pi_1({\mathcal C}_o^n(X)))
\subseteq\pi_1({\mathcal C}_o^n(X))\,.
\end{equation}

\subsection{Step 1: Lifting of non-cyclic continuous maps}
\label{Ss: Step 1}
Taking into account that tame endomorphisms
of ${\mathcal C}^n(X)$ {\bfit must} lift to the covering 
${\mathcal C}_o^n(X)$, we are going to prove the following

\begin{Theorem}[{\caps Lifting Theorem}]
\label{LiftingThm}
\index{Theorem!Lifting Theorem\hfill}
Let $n>4$. Then every non-cyclic continuous map 
$F\colon{\mathcal C}^n(X)\to{\mathcal C}^n(X)$ lifts to
a continuous map 
$f\colon{\mathcal C}_o^n(X)\to{\mathcal C}_o^n(X)$ so that
$F$ and $f$ fit into commutative diagram $(\ref{CD: lifting diagram})$. 
Clearly $f$ is an endomorphism if $F$ is so.
\hfill $\bigcirc$
\end{Theorem}

\noindent In view of the covering mapping theorem, this
topological result is equivalent to the following algebraic

\begin{Theorem}[{\caps Invariance Theorem}]
\label{InvarThm}
\index{Theorem!Invariance Theorem\hfill}
Let $n>4$. Then any non-cyclic endomorphism $\varphi$ of the group
$B_n(X)$ respects the normal subgroup
$PB_n(X)$, i. e., $\varphi(PB_n(X))\subseteq PB_n(X)$. In
fact $PB_n(X)=\varphi^{-1}(PB_n(X))$, which certainly implies
the previous inclusion.
\hfill $\bigcirc$
\end{Theorem}


\noindent Invariance Theorem is, in turn, 
an immediate consequence of the following two results:

\begin{Theorem}[{\caps Composition Theorem}]
\label{CompositionThm}
\index{Theorem!Composition Theorem\hfill}
Let $n>4$ and let $\varphi\colon B_n(X)
\to B_k(X)$ be a non-cyclic homomorphism. 
Then the composition  
\begin{equation}\label{eq: composition with mu}
\psi=\mu\circ\varphi\colon B_n(X)\overset{\varphi}
\longrightarrow B_k(X)\overset{\mu}\longrightarrow{\mathbf S}(k) 
\end{equation}
of $\varphi$ and the standard epimorphism
$\mu\colon B_n(X)\to{\mathbf S}(n)$ is non-cyclic as well.
\hfill $\bigcirc$ 
\end{Theorem}


\begin{Theorem}[{\caps Kernel Theorem}]
\label{KerThm}
\index{Theorem!Kernel Theorem\hfill}
For $n>4$ the kernel $\Ker\psi$ of any non-cyclic homomorphism 
$\psi\colon B_n(X)\to{\mathbf S}(n)$ coincides with the pure braid 
group $PB_n(X)$.
\hfill $\bigcirc$
\end{Theorem}


\noindent{\bfit Explanation.} For a non-cyclic
\index{Theorem!Invariance Theorem!explanation of\hfill}
$\varphi\colon B_n(X)\to B_n(X)$ Composition Theorem shows that 
$\psi=\mu\circ\varphi\colon B_n(X)\to{\mathbf S}(n)$ is non-cyclic and,
by Kernel Theorem, $PB_n(X)=\Ker\psi=\Ker(\mu\circ\varphi)
=\varphi^{-1}(\Ker\mu)=\varphi^{-1}(PB_n(X))$,
which is exactly the statement of Invariance Theorem.
\hfill $\square$

\subsection{Step 2: Commutator subgroup, homomorphisms
$B_n(X)\to {\mathbf S}(n)$}\label{Ss: Step 2} 
To prove Composition and Kernel Theorems, we establish first
some other important properties of braid groups $B_n(X)$.
\vskip0.2cm

\noindent Recall that a group $G$ is said to be 
{\em perfect} if it coincides with its commutator subgroup 
$G'=[G,G]$, or, which is the same,
$G$ does not possess non-trivial homomorphisms to abelian groups. 
\vskip0.2cm

\noindent We extensively make use of the following theorem
proved in Section \ref{Ss: Proof of Commutator Theorem}:

\begin{Theorem}[{\caps Commutator Theorem}]
\label{TheoremA}
\index{Theorem!Commutator Theorem\hfill}
The commutator subgroup $B'_n(X)$ of the 
braid group $B_n(X)$ is a perfect group whenever $n>4$.
\hfill $\bigcirc$
\end{Theorem}

\noindent In the next section we show   
(Theorem \ref{Thm: perfect groups and PBk(X)}) that
for any $n\ge 2$ perfect groups do not admit non-trivial homomorphisms
into the pure braid group $PB_k(X)$. 
Together with Commutator Theorem, this implies:

\begin{Theorem}[{\caps homomorphisms $B'_n(X)\to PB_k(X)$}]
\label{TheoremB}
For $n>4$ the commutator subgroup $B'_n(X)$ of the 
braid group $B_n(X)$ does not possess non-trivial homomorphisms
to the pure braid group $PB_k(X)$.
\hfill $\bigcirc$
\end{Theorem}

\noindent Furthermore, in Section \ref{Ss: Proof of Surjectivity Theorem}
we prove the following theorem:

\begin{Theorem}[{\caps Surjectivity Theorem}]
\label{Thm: Surjectivity Theorem}
\index{Theorem!Surjectivity Theorem\hfill}
For $n\ne 4$ every non-cyclic homomorphism 
$\psi\colon B_n(X)\to{\mathbf S}(n)$ is surjective.
\hfill $\bigcirc$
\end{Theorem}


\noindent These results, together with the classical Artin
theorem about transitive representations $B_n\to{\mathbf S}(n)$
(see \cite{Art47b, Lin96b}), lead to the proof of Composition and
Kernel Theorems and thereby Invariance and Lifting Theorems as well 
(see Section \ref{Sec: Proof of Invariance Theorem}
for the proofs). 

\subsection{Step 3: Strictly equivariant maps}\label{Ss: Step 3} 
With Lifting Theorem at our disposal, 
we proceed as follows. Since (\ref{eq: Co to C covering}) 
is a Galois ${\mathbf S}(n)$ covering, any continuous lifting 
\begin{equation}\label{eq: covering map f}
f=(f_1,...,f_n)\colon{\mathcal C}_o^n(X)\to{\mathcal C}_o^n(X)
\end{equation}
of a continuous map 
\begin{equation}\label{eq: base map}
F\colon{\mathcal C}^n(X)\to{\mathcal C}^n(X)
\end{equation}
is {\em equivariant} in the following sense:
\vskip0.3cm

\noindent{\bfit Equivariance condition:} 
\index{Equivariance condition\hfill}
{\sl there exists an endomorphism 
$\alpha\colon{\mathbf S}(n)\to{\mathbf S}(n)$ such that}
\begin{equation}\label{eq: equivariance condition}
f(\sigma q)=\alpha(\sigma)f(q) \ \
\text{for all} \ q\in{\mathcal C}_o^n(X) \
\text{and} \ \sigma\in{\mathbf S}(n)\,,
\end{equation}
{\sl or, in terms of the components of $f$,}
$$
f_m(\sigma q)=f_{[\alpha(\sigma)]^{-1}(m)}(q) \ \
\text{for all} \ m=1,...,n\,, \ q\in{\mathcal C}_o^n(X)\,, \
\sigma\in{\mathbf S}(n)\,.
$$
{\sl Up to conjugation, the endomorphism $\alpha$ is uniquely defined
by the original map $F$. That is, if endomorphisms $\alpha'$, $\alpha''$
correspond to two liftings $f'$, $f''$ of $F$ then there is a permutation
$s\in{\mathbf S}(n)$ such that $\alpha''(\sigma)=s^{-1}\alpha'(\sigma)s$
for all $\sigma\in{\mathbf S}(n)$.}
\vskip0.2cm

\noindent{\bfit Explanation.} Given $F$ which lifts to 
${\mathcal C}_o^n(X)$, we may construct $\alpha$ as follows.
Taking a continuous lifting $f$ and choosing
basepoints $q^\circ\in{\mathcal C}_o^n(X)$ and
$Q^\circ=p(q^\circ)\in{\mathcal C}^n(X)$ 
we obtain the commutative diagram
\begin{equation}\label{CD: topolofical-algebraic lifting diagram}
\CD
{\pi_1({\mathcal C}_o^n(X),q^\circ)} @ > {f_*} >> 
                          {\pi_1({\mathcal C}_o^n(X),f(q^\circ))}\\
@V{p_*}VV @VV{p_*}V\\ 
{\pi_1({\mathcal C}^n(X),Q^\circ)} @ >> {F_*}>
                          {\pi_1({\mathcal C}^n(X),F(Q^\circ))}\\
@V{\mu'}VV @VV{\mu''}V\\
{\pi_1({\mathcal C}^n(X),Q^\circ)
                  /\pi_1({\mathcal C}_o^n(X),q^\circ)} @ > {\alpha} >> 
 {\pi_1({\mathcal C}^n(X),F(Q^\circ))
                  /\pi_1({\mathcal C}_o^n(X),f(q^\circ))},\\
@VVV @VVV\\
{1} @ .{1} \\
\endCD
\end{equation}
\vskip0.3cm

\noindent where $\mu'$ and $\mu''$ are the natural epimorphisms 
onto the quotient groups 
and $\alpha$ is a unique homomorphism which makes the lower rectangle 
commutative. The latter quotient groups are identified with
the symmetric groups of the ordered sets $q^\circ=(q_1^\circ,...,q_n^\circ)$
and $f(q^\circ)=(f_1(q^\circ),...,f_n(q^\circ))$, that is,
with ${\mathbf S}(n)$; thus, $\alpha$ becomes an endomorphism
of ${\mathbf S}(n)$. Another choice of $f$ and/or a 
basepoint $q^\circ\in{\mathcal C}_o^n(X)$ may lead to another 
endomorphism which must be conjugate to $\alpha$. Furthermore,
connecting $f(q^\circ)$ to $q^\circ$ and $F(Q^\circ)$ to
$Q^\circ$ with paths $\gamma\colon [0,1]\to{\mathcal C}_o^n(X)$
and $p\circ\gamma\colon [0,1]\to{\mathcal C}^n(X)$
we come to identifications of pairs 
$$
\aligned
\ & \{\pi_1({\mathcal C}_o^n(X),q^\circ)\vartriangleleft
\pi_1({\mathcal C}^n(X),Q^\circ)\}
\cong \{PB_n(X)\vartriangleleft B_n(X)\} \ \
\text{and}\\ 
\ & \{\pi_1({\mathcal C}_o^n(X),f(q^\circ))\vartriangleleft
\pi_1({\mathcal C}^n(X),F(Q^\circ))\}
\cong \{PB_n(X)\vartriangleleft B_n(X)\}
\endaligned
$$
that conform to each other. Thereby
(\ref{CD: topolofical-algebraic lifting diagram}) 
leads to the commutative diagram
\vskip0.1cm
\begin{equation}\label{CD: algebraic lifting diagram}
\CD
1 @>>> PB_n(X) @>{p_*}>> B_n(X) @>{\mu}>>{\mathbf S}(n) @>>> 1\\
@. @ V {f_*} VV@V{F_*}VV @VV{\alpha}V\\ 
1 @>>> PB_n(X) @>{p_*}>> B_n(X) @>{\mu}>>{\mathbf S}(n) @>>> 1\\
\endCD
\end{equation}
\vskip0.3cm

\noindent whose horizontal lines are exact sequences
(\ref{eq: main exact sequence}) of the covering
$p\colon{\mathcal C}_o^n(X)\to{\mathcal C}^n(X)$
corresponding to the two pairs of basepoints $\{Q^\circ,q^\circ\}$
and $\{F(Q^\circ),f(q^\circ)\}$.
\hfill $\square$
\vskip0.2cm

\begin{Definition}\label{Def: strict equivariance}
\index{Equivariance condition!strict\hfill}
A continuous map $f\colon{\mathcal C}_o^n(X)\to{\mathcal C}_o^n(X)$
is said to be {\em strictly equivariant} if there exists
an {\bfit automorphism} $\alpha$ of the group ${\mathbf S}(n)$
such that the above mentioned condition
$$
f(\sigma q)=\alpha(\sigma)f(q) \ \
\text{for all} \ q\in{\mathcal C}_o^n(X) \
\text{and} \ \sigma\in{\mathbf S}(n)
\eqnum{2.9}
$$
is fulfilled.
\hfill $\bigcirc$
\end{Definition}  

\begin{Remark}\label{Rmk: equivariance}
{\bfit a}$)$ For $n\ne 6$ any automorphism of ${\mathbf S}(n)$
is inner, and therefore (\ref{eq: equivariance condition})
takes the form
\begin{equation}\label{eq: strong form of strict equivariance}
\exists\, s\in{\mathbf S}(n) \ \ \text{such that}\ \
f(\sigma q)=s\sigma s^{-1}f(q) \ \ \forall\, q\in{\mathcal C}_o^n(X)
\ \ \text{and} \ \ \forall\, \sigma\in{\mathbf S}(n)\,.
\end{equation}
In fact, any strictly equivariant {\em endomorphism} of 
${\mathcal C}_o^n(X)$ satisfies 
(\ref{eq: strong form of strict equivariance}) for all $n\ge 3$
(see {\bfit Explanation} in Section \ref{Ss: Step 4}).
\vskip0.2cm

{\bfit b}$)$ It is easily seen from diagram 
(\ref{CD: algebraic lifting diagram}) that a tame 
endomorphism $F=F_T$ of ${\mathcal C}^n(X)$ must be non-cyclic. Indeed,
let $f=f_T$ be the lifting of $F$ defined according to
(\ref{eq: lifting of tame map}), that is,
$f(q)=T(p(q))q=(T(p(q))q_1,...,T(p(q))q_n)$ for $q=(q_1,...,q_n)$.
As we know, $f$ is equivariant; let $\alpha$ be the corresponding
endomorphism of ${\mathbf S}(n)$. Then for any $\sigma\in{\mathbf S}(n)$ and all $q\in{\mathcal C}^n_o(X)$
we have $\alpha(\sigma)f(q)=f(\sigma q)=T(p(\sigma q))\sigma q
=\sigma T(p(q))q=\sigma f(q)$; thereby $\alpha(\sigma)=\sigma$ and 
$\alpha=\id_{{\mathbf S}(n)}$. In particular, $\alpha$ is non-cyclic,
and (\ref{CD: algebraic lifting diagram}) implies that $F_*$ is
non-cyclic. 
\hfill $\bigcirc$
\end{Remark}

\noindent The following result is an easy consequence of 
Composition Theorem.

\begin{Theorem}[{\caps Equivariance Theorem}]
\label{EquivThm}
\index{Theorem!Equivariance Theorem\hfill}
Let $n>4$. Then any continuous lifting
$f\colon{\mathcal C}_o^n(X)\to{\mathcal C}_o^n(X)$
of a non-cyclic continuous map 
$F\colon{\mathcal C}^n(X)\to{\mathcal C}^n(X)$ is
strictly equivariant.
\hfill $\bigcirc$
\end{Theorem}


\noindent{\bfit Explanation.} 
\index{Theorem!Equivariance Theorem!explanation of\hfill}
Since $F$ is non-cyclic, the induced homomorphism $F_*$ in
(\ref{CD: algebraic lifting diagram}) is non-cyclic and,
by Composition Theorem, the composition $\alpha\circ\mu=\mu\circ F_*$
is non-cyclic as well. Thus, $\alpha$ is non-cyclic, that is,
its image $\alpha({\mathbf S}(n))\cong{\mathbf S}(n)/\Ker\alpha$
is a non-cyclic group. For $n>4$ the group ${\mathbf S}(n)$
possesses a single proper normal subgroup, the alternating
group ${\mathbf A}(n)$, and the quotient group 
${\mathbf S}(n)/{\mathbf A}(n)\cong{\mathbb Z}/2{\mathbb Z}$ 
is cyclic. Therefore $\Ker\alpha=\{1\}$ and $\alpha$
is an automorphism.
\hfill $\square$ 

\begin{Definition}\label{Def: tame endomorphism of Co(X)}
\index{Morphism!invariant\hfill}
\index{Tame endomorphism!of ${\mathcal C}_o^n(X)$\hfill}
A morphism $\tau\colon{\mathcal C}_o^n(X)\to\Aut X$ is called
${\mathbf S}(n)$ {\em invariant} if it is constant on each
${\mathbf S}(n)$ orbit in ${\mathcal C}_o^n(X)$, that is,
$$
\tau(\sigma q)=\tau(q)\ \ \text{for all} \ q\in{\mathcal C}_o^n(X) \ 
\text{and} \ \sigma\in{\mathbf S}(n)\,.
$$
An endomorphism
$f\colon{\mathcal C}_o^n(X)\to{\mathcal C}_o^n(X)$
is said to be {\em tame} if there exists an ${\mathbf S}(n)$
invariant morphism
$\tau\colon{\mathcal C}_o^n(X)\to\Aut X$
and a permutation $s\in{\mathbf S}(n)$ such that
$$
f(q)=f^s_{\tau}(q)\Def s\tau(q)q
=(\tau(q)q_{s^{-1}(1)},...,\tau(q)q_{s^{-1}(n)}) \ \
\forall\, q=(q_1,...,q_n)\in{\mathcal C}_o^n(X)\,.
$$
\vskip0.2cm

\noindent Similarly, we say that $f$ is {\em degenerate tame}
\index{Tame endomorphism!of ${\mathcal C}_o^n(X)$!degenerate\hfill}
if there are
$q^\circ\in{\mathcal C}_o^n(X)$, $s\in{\mathbf S}(n)$
and an ${\mathbf S}(n)$
invariant morphism $\tau\colon{\mathcal C}_o^n(X)\to\Aut X$
such that 
\vskip0.2cm

\hskip83pt 
$f(q)=f^s_{\tau,q^\circ}(q)\Def s\tau(q)q^\circ
\ \ \forall\, q\in{\mathcal C}_o^n(X)\,.\hskip120pt\bigcirc$
\end{Definition}

\noindent Any tame endomorphism $f=f^s_{\tau}$ of ${\mathcal C}_o^n(X)$
may be pulled down to a unique tame endomorphism $F=F_T$ 
of ${\mathcal C}^n(X)$ such that $f$ and $F$ fit into
(\ref{CD: lifting diagram}); $T$ and $F$ are defined by
\begin{equation}\label{eq: push down}
T(Q)=\tau(q)\ \ \text{and}\ \ F(Q)=F_T(Q)=T(Q)Q=pf^s_{\tau}(q)\,,
\end{equation}
where $Q\in{\mathcal C}^n(X)$ and $q$ is any point of the fiber 
$p^{-1}(Q)\subset{\mathcal C}_o^n(X)$. 
\vskip0.2cm

\noindent Lifting and Equivariance Theorems show
that in order to prove Tame Map Theorem it suffices
to prove the following

\begin{Theorem}[{\caps Equivariant Map Theorem}]
\label{EquivMapThm}
\index{Theorem!Equivariant Map Theorem\hfill}
For $n>6-t(X)$ every strictly equivariant endomorphism
of ${\mathcal C}_o^n(X)$ is tame. 
\hfill $\bigcirc$
\end{Theorem}

\subsection{Step 4: Coherence Theorem}\label{Ss: Step 4}
The group $\Aut X$ is exactly $t(X)$ times transitive in $X$
(see Notation \ref{Not: t(X)}). Let $n>t(X)$.
Two ordered $n$ point subsets 
$x=(x_1,...,x_n),\,y=(y_1,...,y_n)\subset X$ 
may be carried to each other by an automorphism $A\in\Aut X$
if and only if they satisfy  
$$
\frac{y_1-y_m}{y_1-y_2}=\frac{x_1-x_m}{x_1-x_2}\quad
\text{for all} \ m=3,...,n \quad \text{when} \ \ X={\mathbb C}
$$
and
$$
\frac{y_m-y_1}{y_m-y_2}:\frac{y_1-y_3}{y_2-y_3}
=\frac{x_m-x_1}{x_m-x_2}:\frac{x_1-x_3}{x_2-x_3}\quad
\text{for all} \ m=4,...,n \quad \text{when} \ \ X={\mathbb{CP}}^1\,.
$$
When the corresponding condition is fulfilled a required $A\in\Aut X$
is uniquely defined and regular in $x,y\in{\mathcal C}_o^n(X)$.   
\vskip0.2cm

\noindent Let $n>t(X)$ and let $f=(f_1,...,f_n)$
be an endomorphism of ${\mathcal C}_o^n(X)$. If we are looking for
a representation $f=f_\tau$ with a morphism 
$\tau\colon{\mathcal C}_o^n(X)\to\Aut X$, we must
for every $q=(q_1,...,q_n)\in{\mathcal C}_o^n(X)$
find $\tau(q)\in\Aut X$, which carries $(q_1,...,q_n)$
to $(f_1(q),...,f_n(q))$ and nicely depends on $q$. 
This is possible if and only if the ordered sets
$x=(q_1,...,q_n)\subset X$ and $y=(f_1(q),...,f_n(q))\subset X$  
satisfy the necessary and sufficient condition mentioned above.
Taking into account that permutations are admitted,
we can conclude that Equivariant Map Theorem
is an immediate consequence of the following result.
\index{Theorem!Equivariant Map Theorem!proof of\hfill}

\begin{Theorem}[{\caps Coherence Theorem}]
\label{CohThm}
\index{Theorem!Coherence Theorem\hfill}
Let $n>6-t(X)$ and let $f=(f_1,...,f_n)$ be a strictly equivariant
endomorphism of ${\mathcal C}_o^n(X)$. Then there exists
a permutation $\sigma\in{\mathbf S}(n)$ such that
for all $q\in{\mathcal C}_o^n(X)$ and all \ $r=t(X)+1,...,n$

\begin{itemize}
\item[(SR)] \ $\displaystyle 
\frac{f_1(\sigma q)-f_r(\sigma q)}{f_1(\sigma q)-f_2(\sigma q)}
=\frac{q_1-q_r}{q_1-q_2}$ \ {\rm if} $X={\mathbb C}$\,,
\vskip0.2cm

\item[(CR)] \ $\displaystyle
\frac{f_r(\sigma q)-f_1(\sigma q)}{f_r(\sigma q)-f_2(\sigma q)}:
\frac{f_1(\sigma q)-f_3(\sigma q)}{f_2(\sigma q)-f_3(\sigma q)}
=\frac{q_r-q_1}{q_r-q_2}:\frac{q_1-q_3}{q_2-q_3}$ \
{\rm if} $X={\mathbb{CP}}^1$\,.\hfill $\bigcirc$
\end{itemize}
\end{Theorem}

\noindent{\bfit Explanation.}
\index{Theorem!Coherence Theorem!explanation of\hfill}
For the sake of neatness, let us explain
how to deduce Equivariant Map Theorem from Coherence Theorem.
We start with the following simple remark, which will be useful
in what follows:

\begin{Remark}
\label{Rmk: equality for generic points implies equality everywhere}
For $n>t(X)$ there exist non-empty Zariski open subsets
$U\subset{\mathcal C}_o^n(X)$ and $V\subset{\mathcal C}^n(X)$
such that 
{\bfit a}$)$ if $Aq=sq$ for some $q\in U$, 
$A\in\Aut X$ and $s\in{\mathbf S}(n)$ then $A=\id$
and $s=1$; \
{\bfit b}$)$ if $AQ=Q$ for some $Q\in V$ and $A\in\Aut X$
then $A=\id$.
Indeed, the set of all points
$q=(q_1,...,q_n)\in{\mathcal C}_o^n(X)$
such that some its non-trivial permutation 
$sq=(q_{s^{-1}(1)},...,q_{s^{-1}(n)})$
may be also written as $Aq=(Aq_1,...,Aq_n)$ with some $A\in\Aut X$
is a proper Zariski closed subset of ${\mathcal C}_o^n(X)$; this
implies both $(${\bfit a}$)$ and $(${\bfit b}$)$.
\hfill $\bigcirc$
\end{Remark}

\noindent For each $q=(q_1,...,q_n)\in{\mathcal C}_o^n(X)$ 
we define the transformation $\tau(q)\in\Aut X$
by one of the following equations:
\begin{equation}\label{eq: tau(q) for X=C}
\hskip-170pt\frac{f_1(\sigma q)-\tau(q)x}{f_1(\sigma q)-f_2(\sigma q)}
=\frac{q_1-x}{q_1-q_2} \ \ {\rm if} \ x\in X={\mathbb C}\,,
\end{equation}
\begin{equation}\label{eq: tau(q) for X=CP1}
\frac{\tau(q)x-f_1(\sigma q)}{\tau(q)x-f_2(\sigma q)}:
\frac{f_1(\sigma q)-f_3(\sigma q)}{f_2(\sigma q)-f_3(\sigma q)}
=\frac{x-q_1}{x-q_2}:\frac{q_1-q_3}{q_2-q_3} \ \
{\rm if} \ x\in X={\mathbb{CP}}^1.
\end{equation}
\vskip0.2cm

\noindent Clearly $\tau\colon{\mathcal C}_o^n(X)\to\Aut X$
is a morphism. Conditions (SR) and (CR) respectively imply that
$\tau(q)q_m=f_m(\sigma q)$ for each $m=1,...,n$ and all
$q=(q_1,...,q_n)\in{\mathcal C}_o^n(X)$; thereby 
$\tau(q)q=f(\sigma q)=\alpha(\sigma)f(q)$, or, which is the same,
$f(q)=\alpha(\sigma^{-1})\tau(q)q$ for all $q\in{\mathcal C}_o^n(X)$. 
To complete the proof, we must check that the morphism $\tau$ is 
${\mathbf S}(n)$ invariant. For every $s\in{\mathbf S}(n)$ and all 
$q\in{\mathcal C}_o^n(X)$ we have
$$
\aligned
s\tau(sq)q&=\tau(sq)sq=f(\sigma sq)=\alpha(\sigma s)f(q)\\
&=\alpha(\sigma s)f(\sigma\sigma^{-1}q)
=\alpha(\sigma s)\tau(\sigma^{-1}q)\sigma^{-1}q
=\alpha(\sigma s)\sigma^{-1}\tau(\sigma^{-1}q)q\,,
\endaligned
$$
which can be rewritten as
\begin{equation}\label{eq: a}
[(\tau(sq))^{-1}\cdot\tau(\sigma^{-1}q)]q
=\sigma\alpha(s^{-1}\sigma^{-1})sq\,,
\end{equation}
where $(\tau(sq))^{-1}\cdot\tau(\sigma^{-1}q)\in\Aut X$ is
the product in the group $\Aut X$. In view of Remark
\ref{Rmk: equality for generic points implies equality everywhere},
this implies that $\sigma\alpha(s^{-1}\sigma^{-1})s=1$
and $\tau(sq)=\tau(\sigma^{-1}q)$ for all $s\in{\mathbf S}(n)$ and
all $q$ in a non-empty Zariski open subset $U\subset{\mathcal C}_o^n(X)$. 
Since $\tau$ is continuous, the latter relation holds true for all
$q\in {\mathcal C}_o^n(X)$ and all $s\in{\mathbf S}(n)$,
which shows that the morphism
$\tau\colon{\mathcal C}_o^n(X)\to\Aut X$ is ${\mathbf S}(n)$
invariant.

Moreover, it is easily seen that the second identity implies
$\alpha(s)=\sigma^{-1}s\sigma$ for all $s\in{\mathbf S}(n)$.
Thus, Coherence Theorem implies that {\sl for all $n>4$ the automorphism
$\alpha$ related to a strictly equivariant endomorphism $f$ of
${\mathcal C}_o^n(X)$ must be inner}; in view of
Remark \ref{Rmk: equivariance}, this $\alpha$ must be inner also
for $n=3$ or $4$.
\hfill $\square$  
 
\subsection{Step 5: Simplicial complex of functions omitting
two values}\label{Ss: Step 5}
\index{Holomorphic functions omitting two values!simplicial complex of\hfill}
Let $f=(f_1,...,f_n)\colon{\mathcal C}_o^n(X)
\to{\mathcal C}_o^n(X)$ be a morphism and let $i,j,k,l$
be distinct. Then 
$$
sr_{ijk}[f]=(f_k-f_i)/(f_k-f_j) \ \ \text{and} \ \
cr_{ijkl}[f]=[(f_l-f_i)/(f_l-f_j)]:[(f_i-f_k)/(f_j-f_k)]
$$
are meromorphic function on ${\mathcal C}_o^n(X)$
and do not take the values $0$ and $1$. In fact,
if $X={\mathbb C}$ then $sr_{ijk}[f]$ and $cr_{ijkl}[f]$
are both holomorphic; if $X={\mathbb{CP}^1}$ this is certainly 
the case for $cr_{ijkl}[f]$.
\vskip0.2cm

\noindent In Section \ref{Sec: Coherence Theorem} 
(Lemma \ref{Lm: non-constant simple and cross ratios})
we show that if $f$ is a strictly
equivariant endomorphism of ${\mathcal C}_o^n(X)$ then 
for $X={\mathbb C}$ all $sr_{ijk}[f]$ and $cr_{ijkl}[f]$
are non-constant and for $X={\mathbb{CP}^1}$ this is certainly 
the case for all $cr_{ijkl}[f]$.
\vskip0.2cm

\noindent Notice that the functions standing in the left hand sides 
of relations (SR) and (CR) up to a permutation of the variables
$q_1,...,q_n$ coincide with $sr_{r,2,1}[f]$ and $cr_{1,2,3,r}[f]$
respectively. Thus, what is written above applies to these
functions.
\vskip0.2cm

\noindent It is well known that the set $L(Z)$ of all non-constant
holomorphic functions $\lambda\colon Z\to{\mathbb C}\setminus\{0;1\}$
on an irreducible quasi-projective variety $Z$ is finite. In particular,
the set $L({\mathcal C}_o^n(X))$ is finite and all functions
$sr_{ijk}[f]$ and $cr_{ijkl}[f]$ in the case $X={\mathbb C}$
(or all $cr_{ijkl}[f]$ in the case $X={\mathbb{CP}^1}$)
must be its elements. 
\vskip0.2cm

\noindent This observation suggests the following way to prove 
Coherence Theorem. First, compute $L({\mathcal C}_o^n(X))$, that is,
determine explicitly all non-constant holomorphic functions
on ${\mathcal C}_o^n(X)$ omitting the values $0$ and $1$.
Second, try to identify (perhaps after a suitable
permutation of $q_1,...,q_n$) the functions 
$sr_{r,2,1}[f]$ ($r=3,...,n$) and/or
$cr_{1,2,3,r}[f]$ ($r=4,...,n$) with certain elements
of $L({\mathcal C}_o^n(X))$, and then... if we are lucky... who knows... 
To say the least, this looks as a finite combinatorial problem.
\vskip0.2cm

\noindent To execute the first part of this program,
we prove the following theorem. 

\begin{Theorem}[{\caps non-constant holomorphic functions
${\mathcal C}_o^n(X)\to{\mathbb C}\setminus\{0;\,1\}$}]
\label{Thm: Non-constant holomorphic functions on Con(X)}
\index{Holomorphic functions omitting two values!on ${\mathcal C}_o^n(X)$\hfill}
{\bfit a}$)$ Each non-constant holomorphic function 
$\mu\colon{\mathcal C}_o^n({\mathbb C})
\to{\mathbb C}\setminus\{0;\,1\}$ is either 

a simple ratio
$sr_{ijk}=(q_k-q_i)/(q_k-q_j)$
or 

a cross ratio 
$cr_{ijkl}=[(q_l-q_i)/(q_l-q_j)]:[(q_i-q_k)/(q_j-q_k)]$

\noindent $($for $n=3$ only simple ratios survive$)$.
\vskip0.2cm

{\bfit b}$)$ Each non-constant holomorphic function 
$\mu\colon{\mathcal C}_o^n({\mathbb{CP}}^1)
\to{\mathbb C}\setminus\{0;\,1\}$ is 
a cross ratio $cr_{ijkl}=[(q_l-q_i)/(q_l-q_j)]:[(q_i-q_k)/(q_j-q_k)]$
$($for $n=3$ no such function exists$)$.
\hfill $\bigcirc$
\end{Theorem}

\noindent To avoid unpleasant combinatorics in the course of
the realization of the second part of our plan, we must try to
find the {\em ``genuine"} reason forcing the sequence of the functions
in the left hand side of (SR) or (CR) to be of the form
prescribed by the right hand side (perhaps after a suitable
rearrangement). 
\vskip0.2cm

\noindent To this end, we notice that the set $L(Z)$
of all non-constant holomorphic functions
$\lambda\colon Z\to{\mathbb C}\setminus\{0;1\}$
\index{Holomorphic functions omitting two values!simplicial complex of\hfill}
on an irreducible quasi-projective variety $Z$
carries a natural structure of a finite simplicial complex
$L_{\vartriangle}(Z)$,
which turns out to be amazingly useful for our purposes.\footnote{I came
to this almost obvious structure just very recently when the first version
of this paper was already written. In this renewed version the proof does
not become much shorter but its main ideas are now utterly clear.
Strangely enough I paid no attention to this structure 
30 years ago. But seemingly nobody did so...}
\vskip0.2cm

\noindent A subset $\Delta^m=\{\mu_0,...,\mu_m\}$ consisting 
of $m+1$ distinct elements of $L(Z)$ is declared to be an $m$ 
{\em simplex} if the quotients $\mu_i:\mu_j\in L(z)$ for all $i\ne j$.
In other words, 
\begin{equation}\label{eq: m simplices are maps to M(m+1)}
\Delta^m\colon Z\ni z\mapsto (\mu_0(z),...,\mu_m(z))
\in{\mathbb C}^{m+1}
\end{equation}
must be a holomorphic map such that all its compon
$\mu_i$ and all their pairwise quotients $\mu_i:\mu_j$ ($i\ne j$) 
are non-constant and its image $\Delta^m(Z)$
is contained in the domain
\begin{equation}\label{eq: M(m+1)}
{\mathcal M}^{m+1}=\{w=(w_0,...,w_m)\in{\mathbb C}^{m+1}\,|\ 
           w_i\ne 0,1 \ \forall\,i \ \text{\rm and} \ w_i\ne w_j \
                                              \forall\, i\ne j\}\,,
\end{equation}
which may also be viewed as the ordered configuration space
${\mathcal C}_o^{m+1}({\mathbb C}\setminus\{0,1\})$.
\vskip0.3cm

\noindent For instance, {\sl the sets of functions 
$$
\aligned
\ & \{(q_1-q_r)/(q_1-q_2)\,|\ r=3,...,n\} \ \ \text{\rm and}\\ 
\ & \{[(q_r-q_1)/(q_r-q_2)]:[(q_1-q_3)/(q_2-q_3)]\,|\ r=4,...,n\}
\endaligned
$$
standing in the right hand sides of {\rm (SR)} and {\rm (CR)}
are simplices} of dimension $n-3$ and $n-4$ respectively.
\vskip0.2cm

\noindent For a holomorphic map $f\colon Y\to Z$ 
of irreducible quasi-projective varieties the correspondence 
$$
\lambda\mapsto f^*(\lambda)
\Def\lambda\circ f\colon Y\overset{f}\longrightarrow
Z\overset{\lambda}\longrightarrow {\mathbb C}\setminus\{0;1\}
$$ 
turns out to be a {\em simplicial map}
$f^*\colon L_{\vartriangle}(Z)\to L_{\vartriangle}(Y)$, 
provided that $\lambda\circ f\ne\const$ for each $\lambda\in L(Z)$    
(this is certainly the case if $f$ is dominant
but the latter condition is not a necessary one).
Moreover, the map $f^*$ preserves dimension
of simplices.
\vskip0.2cm

\begin{Remark}\label{Rmk: invariants of L(Z)}
The simplicial complex $L_{\vartriangle}(Z)$ itself and all
its standard topological invariants are invariants of the underlying
complex manifold $Z$ with respect to biholomorphic isomorphisms.
In Section \ref{Ss: Complexes of simple and cross ratios} we shaw
that $\dim L_{\vartriangle}({\mathcal C}_o^n({\mathbb C}))
=n-3$ for $n\ge 3$ and 
$\dim L_{\vartriangle}({\mathcal C}_o^n({\mathbb{CP}}^1))
=n-4$ for $n\ge 4$, which plays important part in the proof of
Coherence Theorem. It is not difficult to show that for $n\ge 4$
the Euler characteristic 
$\chi(L_{\vartriangle}({\mathcal C}_o^n({\mathbb{CP}}^1)))
=n(n-1)(n-2)(13-3n)/4$; for $n>4$ the complex 
$L_{\vartriangle}({\mathcal C}_o^n({\mathbb{CP}}^1))$ is connected and
its highest homology group $H_{n-4}(L_{\vartriangle}
({\mathcal C}_o^n({\mathbb{CP}}^1),{\mathbb Z}))$ is trivial.
Moreover, $\pi_1(L_{\vartriangle}({\mathcal C}_o^5({\mathbb{CP}}^1)))
={\mathbb F}_{31}$ and $\rank H_1(L_{\vartriangle}
({\mathcal C}_o^6({\mathbb{CP}}^1),{\mathbb Z}))=151$.
Does it mean anything?

The complex $L_{\vartriangle}(Z)$ always has the following rather
special property: if $\{\mu_0,\mu_1\}$, $\{\mu_1,\mu_2\}$ and 
$\{\mu_0,\mu_2\}$ are simplices then the ordered set 
$\Delta^2=[\mu_0,\mu_1,\mu_2]$ is an ordered $2$ simplex and therefore
the chain $c=[\mu_1,\mu_2]-[\mu_0,\mu_2]+[\mu_0,\mu_1]
=\partial\Delta^2$ is a $1$ cycle homologic to $0$.
This does not mean, however,
that $H_1(L_{\vartriangle}(Z),{\mathbb Z})=0$, as some 
``longer" $1$ cycle may be not a boundary.
\hfill $\bigcirc$
\end{Remark}

\noindent We show (see Section \ref{Sec: Coherence Theorem}, 
Lemma \ref{Lm: non-constant simple and cross ratios})
that the above construction applies to every strictly
equivariant endomorphism $f$ of ${\mathcal C}^n_o(X)$.
It follows that {\sl the sets of functions standing
in the left hand sides of {\rm (SR)} and {\rm (CR)} 
also are simplices} (of dimension $n-3$ and $n-4$ respectively).
\vskip0.2cm

\noindent The ${\mathbf S}(n)$ action in ${\mathcal C}^n_o(X)$ induces
an ${\mathbf S}(n)$ action on the set of all simplices
$\Delta^m\subseteq L({\mathcal C}^n_o(X))$ of any fixed dimension $m$. 
In view of what was explained above, 
Coherence Theorem just tells us that
{\em the simplices of functions standing in the left hand sides 
of {\rm (SR)} and {\rm (CR)} belong to the ${\mathbf S}(n)$ orbits
of the simplices of functions standing in the right hand sides
of {\rm (SR)} and {\rm (CR)} respectively}. 
\vskip0.2cm

\noindent Using the explicit form of functions
$\lambda\in L({\mathcal C}_o^n(X))$ provided by Theorem
\ref{Thm: Non-constant holomorphic functions on Con(X)},
we will determine {\em all} the orbits of the latter
${\mathbf S}(n)$ action on the set of $m$ simplices. 
In fact, this will prove Coherence Theorem in the form 
presented in the previous paragraph and
complete the proof of Tame Map Theorem.


\section{Holomorphic universal covers and universal Teichm{\" u}ller families}
\label{Sec: Holomorphic universal covers}

\noindent In this section we exhibit holomorphic 
universal coverings of certain configuration spaces, which
will play the key part in the proof of Linked Map Theorem in Sec.
\ref{Sec: Linked Map Theorem}. To my best knowledge,
the holomorphic universal coverings of ${\mathcal C}^n({\mathbb C})$
and some related spaces were for the first time descibed
by S. Kaliman \cite{Kal75,Kal76b,Kal93}; in particular,
it was stated in \cite{Kal75,Kal76b} that the universal cover
of ${\mathcal C}^n({\mathbb C})$ is isomorphic to 
${\mathbb C}^2\times{\mathbf T}(0,n+1)$.\footnote{H. Shiga in his review
(MR1254030 (95b:32032)) of the paper \cite{Kal93}, which is the English
translation of \cite{Kal76b},  noted: ``However, this has been essentially
shown already in a paper by L. Bers and H. L. Royden 
[Acta Math. 157 (1986), no. 3-4, 259--286; MR88i:30034 ($\S3$, Lemma)]".
It is a pity that being aware of the importance of such an historical
remark Reviewer payed no attention to the far earler papers
\cite{Kal75,Kal76b}, MR0364688 (51 \#942) and  MR0590056 (58 \#28654).}
We do this in a different manner more
convenient for our purposes.
\vskip0.2cm

\noindent At the end of the section we describe the standard normal
series with free factors in the pure braid group $PB_n$ and construct
a similar series in $PB_n(S^2)$. In particular, this leads to the proof of
Theorem \ref{TheoremB} formulated in Section \ref{Ss: Step 2}.
We also present some information about higher homotopy groups
of ${\mathcal C}^n(X)$.  
\vskip0.2cm

\noindent Some of the results discussed below are known,
especially those dealing with homotopy groups; see, for instance,
\cite{FoxNeu62,FadNeu62,FadHus01}. In the quoted sources complex structure of
the configuration spaces of the plane or sphere played no role,
unlike in works of V. Arnold, where it was used extensively;
in fact, the complex point of view goes back to
A. Hurwitz and O. Zariski. The direct decomposition
(\ref{eq: j(X,n): Con(X) to AutX x D(n-t(X))}),
which is certainly very well known, and projections
(\ref{eq: universal coverings pi and Pi})
are ``invisible" from the real point of view. 
S. Kojima \cite{Koj02} deals with the decomposition 
(\ref{eq: j(X,n): Con(X) to AutX x D(n-t(X))}) for
$X={\mathbb{CP}^1}$. E. M. Feichtner and M. Ziegler \cite{FeiZie99}
also pointed it out; what they denoted by $M_{0,n}$
and referred to as {\em moduli spaces of $n$-punctured
complex projective line} in the present paper is denoted by
${\mathcal D}^{n-3}({\mathbb{CP}^1})$,  
${\mathcal C}_o^{n-3}({\mathbb C}\setminus\{0,1\})$
or $M_o(0,n)$; see the next section and Sec.
\ref{Ss: Automorphisms of certain moduli spaces}.

\subsection{Holomorphic universal coverings of 
certain configuration spaces}
\label{Ss: Holomorphic universal covering}
\index{Configuration spaces!holomorphic universal covering of\hfill}
We need now to deal with configuration spaces of
the line ${\mathbb C}$ punctured at several distinct points,
say such as ${\mathcal C}_o^n({\mathbb C}
\setminus\{0,1,...,k-1\})$. To shorten notation, we set 
\index{${\mathbb C}^{\bigcirc\hskip-6pt k}$\hfill}
\begin{equation}\label{eq: k punctured C}
{\mathbb C}^{\bigcirc\hskip-6pt k}
={\mathbb C}\setminus\{0,1,...,k-1\}\quad (k\ge 1)\,. 
\end{equation}
For instance, for $m>1$
\index{${\mathcal C}_o^m({\mathbb C}^{\bigcirc\hskip-6pt 2}\,)$\hfill}
\begin{equation}\label{eq: Com(C**)}
\aligned
{\mathcal C}_o^m({\mathbb C}^{\bigcirc\hskip-6pt 2}\,)
&={\mathcal C}_o^{m}({\mathbb C}\setminus\{0,1\})\\
&=\{z=(z_1,...,z_m)\in{\mathbb C}^m\,|\ 
                z_i\ne 0,1 \ \forall\,i \ \text{\rm and} \ z_i\ne z_j \
                                              \forall\, i\ne j\}
\endaligned
\end{equation}
(compare to (\ref{eq: M(m+1)})). 
To avoid the doubling in formulations related to
two cases $X={\mathbb C}$ and $X={\mathbb{CP}}^1$ under consideration,
we will assume that $n\ge t(X)$ (see Notation \ref{Not: t(X)})
and use the notation 
$$
{\mathcal D}^{n-t(X)}(X)\Def\left\{\aligned
\ &\{q=(q_1,...,q_n)\in{\mathcal C}_o^n({\mathbb C})\,|\ 
                        q_{n-1}=0,\ q_n=1\} \hskip24pt 
                            \text{if} \ \ X={\mathbb C}\,, \\
\ &\{q=(q_1,...,q_n)\in{\mathcal C}_o^n({\mathbb{CP}}^1)\,|\ 
                    q_{n-2}=0,\ q_{n-1}=1,\ q_n=\infty\}\\
&\hskip7.6cm\text{if} \ \ X={\mathbb{CP}}^1\,.
\endaligned
\right.
$$
For $n>t(X)$ we identify ${\mathcal D}^{n-t(X)}(X)$
and 
${\mathcal C}_o^{n-t(X)}({\mathbb C}\setminus\{0,1\})$
via the correspondence
\begin{equation}\label{eq: Co(n-t(X))=D(n-t(X))(X)}
\aligned
&{\mathcal D}^{n-2}({\mathbb C})\ni (q_1,...,q_{n-2},0,1)
\overset{\cong}\longleftrightarrow (q_1,...,q_{n-2})\in
{\mathcal C}_o^{n-2}({\mathbb C}\setminus\{0,1\})\,, \\
&{\mathcal D}^{n-3}({\mathbb{CP}}^1)\ni (q_1,...,q_{n-3},0,1,\infty)
\overset{\cong}\longleftrightarrow (q_1,...,q_{n-3})\in
{\mathcal C}_o^{n-3}({\mathbb C}\setminus\{0,1\})\,.
\endaligned
\end{equation}

\noindent It is well known (see, for instance, \cite{Kal75,Kal76b,Kal93})
that the holomorphic universal cover
$\widetilde{{\mathcal C}_o^m({\mathbb C}\setminus\{0,1\})}$
of 
${\mathcal C}_o^m({\mathbb C}\setminus\{0,1\})$
may be viewed as the Teichm{\" u}ller space ${\mathbf T}(0,m+3)$
\index{Teichm{\" u}ller space\hfill}
of the Riemann sphere with $m+3$ punctures.\footnote{Notice that
${\mathbf T}(0,m+3)$ is homeomorphic to ${\mathbb R}^{2m}$;
moreover, it is biholomorphic to a
holomorphically convex Bergmann domain in ${\mathbb C}^m$.}
The construction presented in the quoted papers makes use of 
the classical techniques related to quasi-conformal mappings,
which is, in a sense, inevitable when dealing with 
Teichm{\" u}ller spaces. We present here a more direct construction
of holomorphic covering mappings\footnote{$\widetilde{\Aut X}$
is the holomorphic universal covering of the complex
Lie group ${\Aut X}$; see (\ref{eq: universal covering of Aut X}).}
\begin{equation}\label{eq: T(0,m+3) to Com(C**)}
\tau\colon{\mathbf T}(0,m+3)
\to
{\mathcal C}_o^m({\mathbb C}\setminus\{0,1\})\,,
\end{equation}
\begin{equation}\label{eq: univ cov maps}
\aligned
&\pi\colon\widetilde{\Aut X}\times {\mathbf T}(0,n+3-t(X))
\to {\mathcal C}_o^n(X)\quad\text{and}\\ 
&\Pi\colon\widetilde{\Aut X}
\times {\mathbf T}(0,n+3-t(X))\to {\mathcal C}^n(X)
\endaligned
\end{equation}
based on the concept of the universal Teichm{\" u}ller family, 
or, which is the same, the ``universal Teichm{\" u}ller curve"
\index{Universal Teichm{\" u}ller family\hfill}
\begin{equation}\label{eq: projection V'(0,l) to T(0,l)}
{\mathcal T}(0,l):= \ \, \rho\colon{\mathbf{V'}}(0,l)
\to {\mathbf T}(0,l)\,,
\end{equation}
which may be roughly described as follows. 
\vskip0.2cm

\noindent Any point ${\mathbf t}$ of the Teichm{\" u}ller
space ${\mathbf T}(0,l)$ ``is" a curve $\Gamma_{\mathbf t}$
of type $(0,l)$ (i. e., the Riemann sphere
${\mathbb {CP}}^{1}\cong\overline{\mathbb C}={\mathbb C}\cup\{\infty\}$
punctured at certain $l$ points and endowed with a marked free basis
of its fundamental group, which we prefer to forget about).
The natural projection
\begin{equation}\label{eq: natural projection CP1 x T(0,l) to T(0,l)}
\rho\colon{\mathbb {CP}}^{1}\times {\mathbf T}(0,l)
\to{\mathbf T}(0,l)
\end{equation}
possesses $l$ holomorphic sections
${\mathbf s}_i\colon {\mathbf T}(0,l)\ni {\mathbf t}
\longmapsto (s_i({\mathbf t}),{\mathbf t})\in
{\mathbb {CP}}^1\times {\mathbf T}(0,l)$
with the following properties: all the maps
$s_1,...,s_l\colon{\mathbf T}(0,l)\to{\mathbb {CP}}^1$
are holomorphic and take distinct values 
at each point ${\mathbf t}\in {\mathbf T}(0,l)$,
\begin{equation}\label{eq: normalization of sections}
s_{l-2}({\mathbf t})\equiv 0\,, \ \
s_{l-1}({\mathbf t})\equiv 1\,, \ \ s_l({\mathbf t})
\equiv \infty \,\footnote{This
normalization is not a necessary but a very convenient one.}
\end{equation}
and, moreover, for every ${\mathbf t}\in {\mathbf T}(0,l)$
the curve $\Gamma({\mathbf t})\Def{\mathbb {CP}}^1\setminus
\{s_1({\mathbf t}),...,s_l({\mathbf t})\}$
is ``the same" curve $\Gamma_{\mathbf t}$ whose ``shape" in
${\mathbf T}(0,l)$ is the point ${\mathbf t}$.
We set 
\index{${\mathbf V}^{\prime}(0,l)$\hfill}
\begin{equation}\label{eq: V'(0,l)}
{\mathbf V}^{\prime}(0,l)=[{\mathbb {CP}}^{1}\times {\mathbf T}(0,l)]
\setminus[{\mathbf s}_1({\mathbf T}(0,l))\cup...
\cup{\mathbf s}_l({\mathbf T}(0,l))]\,;
\end{equation}
restricting (\ref{eq: natural projection CP1 x T(0,l) to T(0,l)})
to ${\mathbf V}^{\prime}(0,l)$ we obtain
(\ref{eq: projection V'(0,l) to T(0,l)}).
Because of normalization (\ref{eq: normalization of sections}),
the mappings $s_1,...,s_{l-3}$ are just holomorphic functions
on $s_i\colon{\mathbf T}(0,l)\to{\mathbb C}\setminus\{0,1\}$
distinct at each point ${\mathbf t}\in{\mathbf T}(0,l)$.
\vskip0.2cm

\noindent The family ${\mathcal T}(0,l):= \ \, 
                               \rho\colon{\mathbf{V'}}(0,l)
\to {\mathbf T}(0,l)$ is {\em universal} in the following sense:
\vskip0.2cm

\noindent{\bfit Universality:} 
{\sl Let ${\mathcal F}:= \ \, \pi\colon X\to Y$ be a holomorphic family
of curves of type $(0,l)$ over a simply connected base $Y$,
i. e., $\pi$ is a holomorphic map of complex manifolds,
$Y$ is simply connected and for each $y\in Y$ the fiber
$\pi^{-1}(y)$ is a curve of type $(0,l)$. Then there exists
a holomorphic map $\phi\colon Y\to{\mathbf T}(0,l)$ such
that the induced family $\phi^*({\mathcal T}(0,l))$
is equivalent to the family ${\mathcal F}$.}
\vskip0.2cm

\noindent With the above construction at our disposal,
we take $l=m+3$ and define the required holomorphic covering map
(\ref{eq: T(0,m+3) to Com(C**)})
as follows:
\begin{equation}
\label{eq: holomorphic universal covering T(0,m+3) to Com(C**)}
\tau\colon{\mathbf T}(0,m+3)\ni{\mathbf t}
\mapsto (s_1({\mathbf t}),...,s_m({\mathbf t}))\in
{\mathcal C}_o^m({\mathbb C}\setminus\{0,1\})\,,
\end{equation}
where ${\mathbf s}_i=(s_i,\cdot)\colon{\mathbf T}(0,m+3)
\to{\mathbb {CP}}^1 \qquad (i=1,...,m+3)$
are the corresponding normalized sections of the projection
${\mathbf{V'}}(0,m+3)\to{\mathbf T}(0,m+3)$.
\vskip0.2cm

\noindent Furthermore, let $\Theta\colon{\mathcal C}_o^n(X)\to\Aut X$ 
denote the morphism whose value at any point
$q=(q_1,...,q_n)\in{\mathcal C}_o^n(X)$ is a unique element
$\Theta(q)\in\Aut X$ that carries $(0,1)$ to $(q_{n-1},q_n)$
when $X={\mathbb C}$ and $(0,1,\infty)$ to $(q_{n-2},q_{n-1},q_n)$ 
when $X={\mathbb{CP}}^1$. The diagonal $\Aut X$ action in
${\mathcal C}_o^n(X)$ is free and biregular, and each its orbit
$(\Aut X)q$ intersects the submanifold ${\mathcal D}^{n-t(X)}(X)$
in a single point that is exactly $(\Theta(q))^{-1}q$. 
In view of (\ref{eq: Co(n-t(X))=D(n-t(X))(X)}),
this provides us with the biregular isomorphisms  
\begin{equation}\label{eq: j(X,n): Con(X) to AutX x D(n-t(X))}
j_{X,n}\colon{\mathcal C}_o^n(X)\overset{\cong}{\longrightarrow}
\Aut X\times {\mathcal D}^{n-t(X)}(X)=\Aut X
\times
{\mathcal C}_o^{n-t(X)}({\mathbb C}\setminus\{0,1\})\,,
\end{equation}
where $j_{X,n}$ and its inverse $j_{X,n}^{-1}$ are defined by
\begin{equation}\label{eq: j(X,n)}
\aligned
j_{X,n}(q)=(\Theta(q),\ (&\Theta(q))^{-1}q)\,,\qquad
j_{X,n}^{-1}(A,z)=Az \\
&(q\in{\mathcal C}_o^n(X)\,,\ \
A\in\Aut X\,, \ \ z\in{\mathcal D}^{n-t(X)}(X))\,. 
\endaligned
\end{equation}
\vskip0.1cm

\noindent Let ${\mathbb C}\leftthreetimes{\mathbb C}$
be the semi-direct product with the group multiplication 
$(b,\zeta)\cdot (b',\zeta')=(b+e^\zeta b', \zeta+\zeta')$
and the complex structure of ${\mathbb C}\times {\mathbb C}$.
For $X={\mathbb C}$ or $X={\mathbb{CP}}^1$
the holomorphic universal covering 
$\pi_X\colon\widetilde{\Aut X}\to\Aut X$ is given by
\begin{equation}\label{eq: universal covering of Aut X}
\aligned
\ &\pi_{\mathbb C}\colon\hskip8pt{\mathbb C}\leftthreetimes{\mathbb C}
\ni (b,\zeta)\mapsto\pi_{\mathbb C}(b,\zeta)
=A\in \Aff {\mathbb C}\,,\quad
Az=e^\zeta z + b \ \ \forall\, z\in{\mathbb C}\,,\\
\ &\pi_{{\mathbb{CP}}^1}\colon{\mathbf{SL}}(2,{\mathbb C})
\to{\mathbf{SL}}(2,{\mathbb C})/\{\pm I\}
={\mathbf{PSL}}(2,{\mathbb C})=\Aut{\mathbb{CP}}^1\,.
\endaligned
\end{equation}
Finally, we define the required holomorphic covering maps
(\ref{eq: univ cov maps})
\index{Configuration spaces!holomorphic universal covering of\hfill}
as follows:\footnote{We invite the reader to write down these
mappings for $X={\mathbb C}$ and $X={\mathbb{CP}}^1$ separately.}
\begin{equation}\label{eq: universal coverings pi and Pi}
\aligned
\pi\colon \widetilde{\Aut X}\times{\mathbf T}(0,n+3-t(X))
&\ni(A,{\mathbf t})\mapsto \pi(A,{\mathbf t})\\
&=(\pi_X(A)s_1({\mathbf t}),...,\pi_X(A)s_n({\mathbf t}))
\in{\mathcal C}_o^n(X)\,,\\
\Pi\colon \widetilde{\Aut X}\times{\mathbf T}(0,n+3-t(X))
&\ni(A,{\mathbf t})\mapsto\Pi(A,{\mathbf t})\\
&=\{\pi_X(A)s_1({\mathbf t}),...,\pi_X(A)s_n({\mathbf t})\}
\in{\mathcal C}^n(X)\,,
\endaligned
\end{equation}
where $\pi_X\colon\widetilde{\Aut X}\to\Aut X$
is defined by (\ref{eq: universal covering of Aut X}).

\subsection{Homotopy groups and the standard normal series for $PB_n(X)$}
\label{Ss: Homotopy groups and PBn(X)}
\index{Standard normal series for $PB_n(X)$\hfill}
According to (\ref{eq: j(X,n): Con(X) to AutX x D(n-t(X))}),
${\mathcal C}_o^n(X)\cong\Aut X
\times
{\mathcal C}_o^{n-t(X)}({\mathbb C}\setminus\{0,1\})$.
Hence for $n\ge 3$

\begin{equation}\label{eq: pi1(Con)}
\aligned
\pi_1({\mathcal C}_o^n(X))\ =\pi_1&(\Aut X)\times 
\pi_1(
       {\mathcal C}_o^{n-t(X)}({\mathbb C}\setminus\{0,1\}))\\
\ &=\left\{\aligned
   &{\mathbb Z}\times
        \pi_1(
              {\mathcal C}_o^{n-2}({\mathbb C}\setminus\{0,1\}))                 
                   \hskip37pt\text{\rm if} \ X={\mathbb C}\,,\\
   &({\mathbb Z}/2{\mathbb Z})\times
         \pi_1(
               {\mathcal C}_o^{n-3}({\mathbb C}\setminus\{0,1\})) 
                   \hskip8pt\text{\rm if} \ X={\mathbb{CP}}^1\,.
        \endaligned
 \right.
\endaligned
\end{equation}
The fundamental group
$\pi_1({\mathcal C}_o^{n-2}({\mathbb C}\setminus\{0,1\}))$
is one of the important objects of the classical braid theory, 
which usually appears by the following way. The projection 
$$
{\mathcal C}_o^n({\mathbb C})\ni q
=(q_1,...,q_n)\mapsto q_n\in{\mathbb C}
$$
is a smooth fiber bundle with the fiber 
$$
{\mathcal C}_o^{n-1}({\mathbb C}^*)
=\{(q_1,...,q_{n-1})\in{\mathbb C}\,|\
q_i\ne 0 \ \forall i, \ q_i\ne q_j\ \forall i\ne j\}\,.
$$
Since the base ${\mathbb C}$ is contractible, the exact homotopy
sequence of this fibration shows that
$$
\pi_k({\mathcal C}_o^n({\mathbb C}))
\cong\pi_k({\mathcal C}_o^{n-1}({\mathbb C}^*)) \ \
\forall\,k\,.
$$
The next projection
$$
{\mathcal C}_o^{n-1}({\mathbb C}\setminus\{0,1\})
\ni (q_1,...,q_{n-1})\mapsto q_{n-1}\in{\mathbb C}^*
={\mathbb C}^{\bigcirc\hskip-6pt 1}
$$
is a smooth fiber bundle as well and the fiber is 
$$
{\mathcal C}_o^{n-2}({\mathbb C}^{\bigcirc\hskip-6pt 2}\,) 
={\mathcal C}_o^{n-2}({\mathbb C}\setminus\{0,1\})
=\{(q_1,...,q_{n-2})\in{\mathbb C}\,|\ 
q_i\ne 0,1 \ \forall i, \ q_i\ne q_j\ \forall i\ne j\}\,.
$$
Since the base ${\mathbb C}^{\bigcirc\hskip-6pt 1}={\mathbb C}^*$
is aspherical, we have
$$
\pi_k({\mathcal C}_o^{n-1}({\mathbb C}^{\bigcirc\hskip-6pt 1}\,))  
\cong\pi_k({\mathcal C}_o^{n-2}({\mathbb C}^{\bigcirc\hskip-6pt 2}\,))
\ \ \forall\,k\ge 1\,.
$$
Moreover, the final segment
of the corresponding exact homotopy sequence is
$$
1\to\pi_1({\mathcal C}_o^{n-2}({\mathbb C}^{\bigcirc\hskip-6pt 2}\,))  
\to\pi_1({\mathcal C}_o^{n-1}({\mathbb C}^{\bigcirc\hskip-6pt 1}\,))
\to \pi_1({\mathbb C}^*)\to 0
$$
may be written as
$$
1\to PB_n^{(2)}\to PB_n^{(1)}\to{\mathbb Z}\to 0\,,
$$
where $PB_n^{(2)}
\Def\pi_1({\mathcal C}_o^{n-2}({\mathbb C}^{\bigcirc\hskip-6pt 2}\,))$
and $PB_n^{(1)}\Def\pi_1({\mathcal C}_o^n({\mathbb C}))=PB_n$. 
This projecting process may be continued; say the next step
leads to a fiber bundle with the fiber
$$
{\mathcal C}_o^{n-3}({\mathbb C}^{\bigcirc\hskip-6pt 3}\,)
=\{(q_1,...,q_{n-3})\in{\mathbb C}\,|\
q_i\ne 0,1,2 \ \forall i, \ q_i\ne q_j\ \forall i\ne j\}
$$
over the aspherical base 
${\mathbb C}^{\bigcirc\hskip-6pt 2}={\mathbb C}\setminus\{0,1\}$,
the relations 
$$\pi_k({\mathcal C}_o^{n-3}({\mathbb C}^{\bigcirc\hskip-6pt 3}\,))
\cong\pi_k({\mathcal C}_o^{n-2}({\mathbb C}^{\bigcirc\hskip-6pt 2}\,))\ \ 
\forall\, k>1
$$
and the exact sequence
$$
1\to PB_n^{(3)}\to PB_n^{(2)}\to{\mathbb F}_2\to 1\,,
$$
where $PB_n^{(3)}
\Def\pi_1({\mathcal C}_o^{n-3}({\mathbb C}^{\bigcirc\hskip-6pt 3}\,))$
and ${\mathbb F}_2$ stands for a free group of rank $2$. 
This inductive procedure eventually
leads to the following important result going back to 
E. Artin \cite{Art25}, A. Markov \cite{Mar45}
(at least part $(${\bfit c}$)$ and the determination of the fundamental
groups in part $(${\bfit a}$)$) and also to 
E. Fadell, R. H. Fox and L. Neuwirth
\cite{FadNeu62,FoxNeu62} (asphericity in parts $(${\bfit {a,b}}$))$.   
\vskip0.2cm

\begin{Theorem}[{\caps Artin-Markov-Fadell-Neuwirth}]
\label{Thm AMFFN}
\index{Theorem!Artin-Markov-Fadell-Neuwirth\hfill}
{\bfit a}$)$ The spaces ${\mathcal C}_o^n({\mathbb C})$ and 
${\mathcal C}^n({\mathbb C})$ 
are Eilenberg-MacLane $K(\pi,1)$ spaces for the pure braid group
$PB_n$ and the braid group $B_n$ respectively, i. e.,
these spaces are aspherical and have the fundamental groups
that are pointed out.\footnote{This certainly follows from
(\ref{eq: j(X,n): Con(X) to AutX x D(n-t(X))}),
(\ref{eq: universal covering of Aut X}) and 
(\ref{eq: universal coverings pi and Pi});
but the standard argument exhibited above is more preferable since it
is quite elementary.}
\vskip0.2cm

{\bfit b}$)$ All configuration
spaces\footnote{See notation (\ref{eq: k punctured C})}  
${\mathcal C}_o^m({\mathbb C}^{\bigcirc\hskip-6pt k}\,)$
are aspherical. 
\vskip0.2cm

{\bfit c}$)$ The pure braid group $PB_n$
possesses a finite series of normal subgroups
\begin{equation}\label{eq: normal series for PBn}
\{1\}=PB_n^{(n)}\subset PB_n^{(n-1)}\subset\cdots\subset 
PB_n^{(2)}\subset PB_n^{(1)}=PB_n 
\end{equation}
such that $PB_n^{(r)}/PB_n^{(r+1)}\cong {\mathbb F}_r$
{\rm (the free group of rank $r$)}, \ $1\le r\le n-1$.
\hfill $\square$
\end{Theorem}

\noindent The normal series (\ref{eq: normal series for PBn})
will be referred to as the {\em standard} one.
\vskip0.2cm  

\noindent Replacing $n$ in 
(\ref{eq: normal series for PBn}) with $n-1$, we see that 
the fundamental group 
$\pi_1({\mathcal C}_o^{n-3}(
                            {\mathbb C}\setminus\{0,1\}))$,
which participate in the formula for
$\pi_1({\mathcal C}_o^n({\mathbb{CP}}^1))=PB_n(S^2)$
given by (\ref{eq: pi1(Con)}), coincides with second term
$PB_{n-1}^{(2)}\vartriangleleft PB_{n-1}$
of the standard normal series for $PB_{n-1}$.  
Together with (\ref{eq: j(X,n): Con(X) to AutX x D(n-t(X))}),
(\ref{eq: universal covering of Aut X}) and asphericity of
${\mathcal C}_o^{n-3}(
                      {\mathbb C}\setminus\{0,1\})$,
this implies the following result.

\begin{Theorem}\label{Thm: structure of PBn(S2)}
{\bfit a}$)$ The pure sphere braid group 
$PB_3(S^2)\cong{\mathbb Z}/2{\mathbb Z}$. 
For $n>3$ the group $PB_n(S^2)\cong({\mathbb Z}/2{\mathbb Z})
\times PB_{n-1}^{(2)}$, where $PB_{n-1}^{(2)}$
possesses a series of normal subgroups
\begin{equation}\label{eq: normal series for PB(n-1)(2)}
\{1\}=PB_{n-1}^{(n-1)}\subset PB_{n-1}^{(n-2)}\subset\cdots\subset 
PB_{n-1}^{(3)}\subset PB_{n-1}^{(2)}
\end{equation}
such that
$PB_{n-1}^{(r)}/PB_{n-1}^{(r+1)}\cong {\mathbb F}_r$
$($the free group of rank $r)$, \ $1\le r\le n-2$.
\vskip0.2cm

{\bfit b}$)$ The higher homotopy groups of 
${\mathcal C}^n({\mathbb{CP}^1})$ and 
${\mathcal C}_o^n({\mathbb{CP}^1})$ are as follows:
\begin{equation}\label{eq: homotopy groups of Con(CP1)}
\aligned
\pi_k({\mathcal C}^n({\mathbb{CP}^1}))
&=\pi_k({\mathcal C}_o^n({\mathbb{CP}^1}))
=\pi_k({\mathbf{PSL}}(2,{\mathbb C}))\\
&=\pi_k({\mathbf{SL}}(2,{\mathbb C}))=\pi_k(S^3)=
               \left\{\aligned
                   & 0\hskip31pt\text{\rm if} \ k=2\,,\\
                   & \pi_k(S^2)\ \ \text{\rm if} \ k>2\,.
                     \endaligned
               \right.
\endaligned
\end{equation}
\end{Theorem}

\noindent The following theorem, which contains also Theorem
\ref{TheoremB} formulated in Section \ref{Ss: Step 2}, 
is an immediate consequence of part $(${\bfit c}$)$
of Artin-Markov-Fadell-Fox-Neuwirth Theorem, 
part $(${\bfit a}$)$ of Theorem \ref{Thm: structure of PBn(S2)}
and Commutator Theorem, which we still need to prove.

\begin{Theorem}\label{Thm: perfect groups and PBk(X)}
Perfect groups do not admit non-trivial homomorphisms
into the pure braid group $PB_k(X)$.
In particular, for $n>4$ the commutator subgroup $B_n'(X)$
of the braid group $B_n(X)$ does not possess non-trivial
homomorphisms to the pure braid group $PB_k(X)$.
\hfill $\bigcirc$
\end{Theorem}



\section{Proof of Invariance Theorem}
\label{Sec: Proof of Invariance Theorem}

\noindent The main aim of this section is to prove 
Invariance Theorem; we follow the way
outlined in Section \ref{Sec: Strategy of the proof of Tame Map Theorem}.

\subsection{Canonical presentations of $B_n$ and $B_n(S^2)$}
\label{Ss: Canonical presentations of Bn and Bn(S2)}
\index{Braid groups!canonical presentations\hfill}
The inclusion 
${\mathbb C}\hookrightarrow{\mathbb C}\cup\{\infty\}
={\mathbb{CP}}^1$  provides the inclusions
$i_o\colon{\mathcal C}_o^n({\mathbb C})\hookrightarrow
{\mathcal C}_o^n({\mathbb{CP}}^1)$ and
$i\colon{\mathcal C}^n({\mathbb C})\hookrightarrow
{\mathcal C}^n({\mathbb{CP}}^1)$, which, in turn, 
induce the epimorphisms of the fundamental groups 
\begin{equation}\label{eq: epimorphism io*}
i_{o*}\colon PB_n=\pi_1({\mathcal C}_o^n({\mathbb C}))
\to\pi_1({\mathcal C}_o^n({\mathbb{CP}}^1))=PB_n(S^2)
\end{equation}
and 
\begin{equation}\label{eq: epimorphism i*}
i_*\colon B_n=\pi_1({\mathcal C}^n({\mathbb C}))
\to\pi_1({\mathcal C}^n({\mathbb{CP}}^1))=B_n(S^2)\,.
\end{equation}
Recall that the canonical presentation of the 
Artin braid group $B_n=\pi_1({\mathcal C}^n({\mathbb C}))$ involves
$n-1$ generators $\sigma _{1},...,\sigma _{n-1}$
and the defining system of relations
\begin{eqnarray}
&&\sigma_i\sigma_j
=\sigma_j\sigma_i\qquad\qquad\qquad \ \ (\modo{i-j}\ge 2)\,,
                            \label{eq: commutation relation in Bn}\\
&&\sigma_i\sigma_{i+1}\sigma_i
=\sigma_{i+1}\sigma_i\sigma_{i+1}\qquad (1\le i\le n-2)\,.
                                             \label{eq: braid relation}
\end{eqnarray}
It is well known that the kernel of the epimorphism
$i_*\colon B_n\to B_n(S^2)$ is the normal subgroup
$S\vartriangleleft B_n$ generated (as a normal subgroup)
by the single element
\begin{equation}\label{eq: normal generator of the kernel}
s=\sigma_1\sigma_2\cdots\sigma_{n-1}
\sigma_{n-1}\cdots\sigma_2\sigma_1\in PB_n\,,
\end{equation}
so that $S\vartriangleleft PB_n$. Thus, involving also the monomorphisms 
$p_*$ that correspond to the coverings 
$p\colon{\mathcal C}_o^n({\mathbb C})\to{\mathcal C}^n({\mathbb C})$
and $p\colon{\mathcal C}_o^n({\mathbb{CP}}^1)
\to{\mathcal C}^n({\mathbb{CP}}^1)$,
we obtain the following commutative diagram 
with exact lines and rows: 

\begin{equation}\label{CD: embeddings and coverings diagram}
\CD
@.      1 @.    1             @.            1               \\
@.     @VVV    @VVV                         @VVV            \\
1  @>>> S @>>> PB_n           @>>i_{o*}> PB_n(S^2) @ >>> 1  \\
@.     @VVV   @V{p_*}VV                     @VV{p_*}V       \\ 
1  @>>> S @>>> B_n            @>>i_*>    B_n(S^2)@ >>> 1    \\
@.     @VVV   @V{\mu}VV                     @VV{\mu}V       \\
@.      1 @>>> {\mathbf S}(n) @>>\id>  {\mathbf S}(n) @>>> 1\\
@. @.          @VVV                         @VVV            \\
@. @.          1              @.            1               \\
\endCD
\end{equation}
\vskip0.1cm

\noindent(compare to (\ref{eq: main exact sequence})).
It follows that the sphere braid group 
$B_n(S^2)=\pi_1({\mathcal C}^n({\mathbb{CP}}^1))$ 
admits a presentation (we call it the {\em canonical} one) 
with $n-1$ generators $\sigma _{1},..,\sigma _{n-1}$
and the defining system of relations
\index{Sphere braid group $B_n(S^2)$\hfill}
\begin{equation}\label{eq: defining relations in Bn(S2)}
\aligned
&\sigma_i\sigma_j
=\sigma_j\sigma_i\qquad\qquad\qquad \ \ (\modo{i-j}\ge 2)\,,\\
&\sigma_i\sigma_{i+1}\sigma_i
=\sigma_{i+1}\sigma_i\sigma_{i+1}\qquad (1\le i\le n-2)\,,\\
&\sigma_1\sigma_2\cdots\sigma_{n-1}
\sigma_{n-1}\cdots\sigma_2\sigma_1=1
\endaligned
\end{equation} 
(see, for instance, \cite{FadBusk61,FadBusk62}).

\subsection{Proof of Commutator Theorem}
\label{Ss: Proof of Commutator Theorem}
\index{Theorem!Commutator Theorem!proof of\hfill}
A finite presentation of the commutator subgroup $B_n'$ found 
in \cite{GorLin69} implies immediately that for $n>4$
the abelianization $B_n'/[B_n',B_n']$ of $B_n'$
is trivial; that is, $B_n'$ is perfect.
Since $B_n(S^2)$ is a quotient group of $B_n$, its commutator subgroup
$B_n'(S^2)$ is a quotient group of $B_n'$. A quotient group of 
a perfect group is obviously perfect; hence for $n>4$
the group $B_n'(S^2)$ is perfect as well.
\hfill $\square$

\begin{Remark}\label{Rmk: Gorin formula}
This remark allows to avoid the quotation of
\cite{GorLin69} in the above proof of Commutator Theorem.

In 1967 I found certain identities in $B_n$, $n>6$,
which led to the first proof of the fact that for such $n$
the group $B'_n$ is perfect; during a long time,
the only known proof for smaller $n$ was based
on the presentation of $B'_n$ found in \cite{GorLin69}. 

In the middle 1980's E. A. Gorin discovered the following
beautiful relation, which holds for any $n\ge 4$:
\begin{equation}\label{eq: Gorin relation}
\index{Gorin's relation\hfill}
\sigma_3\sigma_1^{-1}
= (\sigma_1\sigma_2)^{-1}\cdot
\left[\sigma_3\sigma_1^{-1},\sigma_1\sigma_2^{-1}\right]
\cdot (\sigma_1\sigma_2)\,,
\end{equation}
where
$\left[\sigma_3\sigma_1^{-1},\sigma_1\sigma_2^{-1}\right]
=\left(\sigma_3\sigma_1^{-1}\right)^{-1}
\cdot\left(\sigma_1\sigma_2^{-1}\right)^{-1}
\cdot\left(\sigma_3\sigma_1^{-1}\right)
\cdot\left(\sigma_1\sigma_2^{-1}\right)$ is the
commutator of the elements $g_1=\sigma_3\sigma_1^{-1}$ and
$g_2=\sigma_1\sigma_2^{-1}$. Clearly, $g_1,g_2\in B_n'$,
and (\ref{eq: Gorin relation}) shows that 
$\sigma_3\sigma_1^{-1}\in[B_n',B_n']$. Hence
the normal subgroup $N\vartriangleleft B_n$ generated
(as a normal subgroup) by the element $\sigma_3\sigma_1^{-1}$
is contained in $[B_n',B_n']$.
For $n>4$ it follows readily from
(\ref{eq: commutation relation in Bn}) and (\ref{eq: braid relation}) 
that $N$ contains the whole
commutator subgroup $B_n'$. Indeed, the presentation of $B_n/N$
involves the generators $\sigma_1,...,\sigma_{n-1}$ 
and the defining system of relations consisting of 
(\ref{eq: commutation relation in Bn}), 
(\ref{eq: braid relation}) and the additional relation
$\sigma_3\sigma_1^{-1}=1$; since $n>4$ 
relations (\ref{eq: commutation relation in Bn}) and 
$\sigma_3\sigma_1^{-1}=1$ imply 
$\sigma_3\sigma_4=\sigma_4\sigma_3$; in view of 
(\ref{eq: braid relation}), this shows that $\sigma_3=\sigma_4$
and, finally, $\sigma_1=\sigma_2=...=\sigma_{n-1}$.
Thus, $B_n/N\cong{\mathbb Z}$ and $N\supseteq B_n'$.
Thereby $B_n'\subseteq N\subseteq [B_n',B_n']$ 
and $B'_n$ is perfect. 
\hfill $\bigcirc$
\end{Remark} 

\subsection{Proof of Surjectivity Theorem \ref{Thm: Surjectivity Theorem}}
\label{Ss: Proof of Surjectivity Theorem}
\index{Theorem!Surjectivity Theorem!proof of\hfill}
Since $B_n(S^2)$ is a quotient group of $B_n$, 
it suffices to prove the theorem for $B_n$. This 
was done in \cite{Lin96b} (Lemma 2.7). It turns out that this theorem
is a very convenient technical tool. Since I still did 
not take care of publication of \cite{Lin96b}, 
I decided to outline the proof here.
\vskip0.3cm

\noindent First, let me recall (in a convenient form) 
some classical results of E. Artin \cite{Art47b}.

\begin{Definition}\label{Def: conjugate homomorphisms}
Two group homomorphisms 
$\phi,\psi\colon G\rightarrow H$
are said to be {\em conjugate} if there is an element 
$h\in H$ such that $\psi(g)=h\phi(g)h^{-1}$ for all $g\in G$.
If this is the case, we write $\phi\sim\psi$; otherwise, we
write $\phi\nsim\psi$. 
\end{Definition}

\begin{Definition}\label{Def: transitive homomorphism}
\index{Homomorphism!transitive\hfill}
A group homomorphism $\psi\colon G\to{\mathbf S}(n)$
is said to be {\em transitive} if its image $\psi(G)$ is
a transitive subgroup of the symmetric ${\mathbf S}(n)$
(meaning the standard left action of ${\mathbf S}(n)$ on
$\boldsymbol\Delta_n\Def\{1,...,n\}$). 
\end{Definition}

\noindent Recall that the elements $\sigma_1$
and $\alpha=\sigma_1\cdots\sigma_{n-1}$
generate the whole Artin braid group $B_n$;
this follows from the relations 
\begin{equation}\label{eq: conjugation of sigma by alpha}
\aligned
\alpha\sigma_j&=\sigma_1\cdots\sigma_{j-1}\cdot
(\sigma_j\sigma_{j+1}\sigma_j)\cdot\sigma_{j+2}\cdots\sigma_{n-1}\\
&=\sigma_1\cdots\sigma_{j-1}\cdot
(\sigma_{j+1}\sigma_j\sigma_{j+1})\cdot\sigma_{j+2}\cdots\sigma_{n-1}
=\sigma_{j+1}\alpha\quad (1\le j\le n-2)
\endaligned
\end{equation}
(the same relations hold true in $B_n(S^2)$).
Thus, a homomorphism 
$\psi\colon B_n\to H$ is uniquely defined by
its values $\psi(\sigma_1)$, $\psi(\alpha)$.
\vskip0.3cm  

\noindent {\bfit Artin Theorem.}
\index{Theorem!Artin Theorem\hfill}
{\sl Let $\psi\colon B_n\to{\mathbf S}(n)$ be a non-cyclic
transitive homomorphism.
\vskip0.2cm

{\bfit a}$)$ If $n\ne 4$ and $n\ne 6$ then
the homomorphism $\psi$ is conjugate to the standard epimorphism
$\mu\colon B_n\to{\mathbf S}(n)$.
\vskip0.2cm

{\bfit b}$)$ If $n = 6$ and $\psi\nsim\mu$ then $\psi$ is conjugate 
to the homomorphism $\nu_6$ defined by
$$
\nu_6(\sigma_1) = (1,2)(3,4)(5,6),
\qquad \nu_6(\alpha) = (1,2,3)(4,5).
$$

{\bfit c}$)$ If $n = 4$ and $\psi\nsim\mu$ then $\psi$ is conjugate to
one of the following three homomorphisms 
$\nu_{4,1}$, $\nu_{4,2}$, $\nu_{4,3}$:
$$
\aligned
\nu_{4,1}(\sigma_1) &= (1,2,3,4),\qquad \\
\nu_{4,2}(\sigma_1) &= (1,3,2,4),\qquad \\
\nu_{4,3}(\sigma_1) &= (1,2,3), \qquad  \\
\endaligned
\aligned
\nu_{4,1}(\alpha) &= (1,2);\qquad \\
\nu_{4,2}(\alpha) &= (1,2,3,4);\qquad \\
\nu_{4,3}(\alpha) &= (1,2)(3,4);\qquad \\
\endaligned
\aligned
[\nu_{4,1}(\sigma_3) &= \nu_{4,1}(\sigma_1)]\\
[\nu_{4,2}(\sigma_3) &= \nu_{4,2}(\sigma_1^{-1})]\\
[\nu_{4,3}(\sigma_3) &= \nu_{4,3}(\sigma_1)].
\endaligned
$$

{\bfit d}$)$ Except of the case $n=4$ and
$\psi\sim\nu_{4,3}$, the homomorphism $\psi$ is surjective.
In the exceptional case when $\psi\sim \nu_{4,3}$
the image of $\psi$ coincides with the alternating subgroup
${\mathbf A}(4)\subset{\mathbf S}(4)$.}
\hfill $\bigcirc$
\vskip0.3cm

\noindent The following lemma is the heart 
of Artin's methods developed in \cite{Art47b}.\footnote{In fact, 
this lemma is not formulated explicitly in \cite{Art47b}; 
moreover, Artin treats only the special case $n=m$. 
However, the corresponding part of the proof of 
Lemma 6 in \cite{Art47b}, after an evident minor modification, 
leads to the desired result. See \cite{Lin96b} for more details.}
\vskip0.3cm

\noindent{\bfit Artin Lemma.}
\index{Artin Lemma\hfill}
{\sl Let $n>4$ and $m$
be natural numbers. Suppose that there is a prime $p>2$ such that
\begin{equation}\label{eq: p between m/2 and n-2}
m/2 < p\le n-2\,.
\end{equation}
Then for every non-cyclic transitive homomorphism
$\psi\colon B_n\to{\mathbf S}(m)$
the permutation $\widehat\sigma_1=\psi(\sigma_1)$ 
has at least $n-2$ fixed points; in particular, 
such a homomorphism cannot exist if $m<n$.}
\hfill $\bigcirc$
\vskip0.2cm

\noindent We use Artin Lemma to prove the following
useful addition to Artin Theorem (see \cite{Lin72a,Lin74,Lin79}, 
and also \cite{Lin96b}, Theorem 2.1).

\begin{Theorem}\label{Thm: Bn(X) to S(k) for n>k}
For $n\ne 4$ and $n>k$ any homomorphism $\psi\colon B_n(X)\to{\mathbf S}(k)$
is cyclic.
\end{Theorem}

\begin{proof}
As $B_n(S^2)$ is a quotient group of $B_n$ it suffices
to prove the theorem for homomorphisms $\psi\colon B_n\to{\mathbf S}(k)$.

For $n\le 3$ all is trivial, since 
${\mathbf S}(2)\cong{\mathbb Z}/2{\mathbb Z}$. 
Suppose that $n>\max\{4,k\}$ and there is a non-cyclic homomorphism
$\psi\colon B_n\to{\mathbf S}(k)$.
Let $H=\Img\psi\subseteq{\mathbf S}(k)$.
For each orbit $O\subseteq\boldsymbol\Delta_k=\{1,...,k\}$ of $H$,
define the {\em reduction} $\psi_O\colon B_n\to{\mathbf S}(O)$
of $\psi$ to $O$ as follows: $\psi_O(g)=\psi(g)|O$ for all $g\in B_n$, 
where $\psi(g)|O$ denotes the restriction of the permutation
$\psi(g)\in H\subseteq{\mathbf S}(k)$ to the $H$-orbit $O$. 
Clearly, all such reductions are transitive (meaning that
$\Img\psi_O$ is a transitive subgroup of the symmetric group
${\mathbf S}(O)$), and for at least one orbit $O_\circ$
the corresponding reduction $\psi_{O_\circ}\colon B_n
\to{\mathbf S}({O_\circ})\cong{\mathbf S}(m)$ ($m=\#O_\circ$)
must be non-cyclic (otherwise, $\psi$ itself would be abelian and
hence cyclic). Since $n>4$ and $m\le k<n$, it follows 
from the famous Chebyshev Theorem that there is
a prime $p>2$ such that $m/2<p\le n-2$.
In view of the above inequalities $m\le k<n$, 
Artin Lemma implies that any transitive homomorphism 
$B_n\to{\mathbf S}(m)$ must be cyclic, which contradicts
our choice of the orbit $O_\circ$.
\end{proof}

\noindent We are now in a position to prove Theorem \ref{Thm: Surjectivity Theorem}.
\vskip0.2cm

\begin{proof}
Let $n\ne 4$ and let $\psi\colon B_n(X)\to{\mathbf S}(n)$ be a 
non-cyclic homomorphism. We must prove that $\psi$ is surjective.

We have already noticed that it suffices to deal with
the Artin braid group $B_n$. If $\psi$ is transitive the statement
follows directly from Artin Theorem. Assume 
that $\psi$ is intransitive. Then $\#O<n$
for any $(\Img\psi)$-orbit $O\subset\boldsymbol\Delta_n$ and
Theorem \ref{Thm: Bn(X) to S(k) for n>k} implies that the reduction of 
$\psi$ to any such orbit is cyclic. Hence $\psi$ itself is cyclic,
contradicting our assumption.
\end{proof}

\subsection{Proof of Composition Theorem}
\label{Ss: Proof of Composition Theorem}
\index{Theorem!Composition Theorem!proof of\hfill}
Suppose, on the contrary, that the composition
$\psi=\mu\circ\varphi\colon B_n(X)\overset{\varphi}
\longrightarrow B_k(X)\overset{\mu}\longrightarrow{\mathbf S}(k)$
is cyclic. Then $\mu(\varphi(B_n'(X)))=\{1\}$ and hence
$\varphi(B_n'(X))\subseteq\Ker\mu=PB_k(X)$. By Commutator Theorem,
the commutator subgroup $B_n'(X)$ is perfect; thus, 
the inclusion $\varphi(B_n'(X))\subseteq PB_k(X)$
and Theorem \ref{Thm: perfect groups and PBk(X)}
imply $\varphi(B_n'(X))=\{1\}$, which contradicts the assumption
that $\varphi$ is non-cyclic.
\hfill $\square$ 

\subsection{Homomorphisms $B_n(X)\to B_k(X)$ for $n>k$} 
\label{Ss: Homomorphisms Bn(X) to Bk(X) for n>k}
Theorem \ref{Thm: Bn(X) to S(k) for n>k} and Composition Theorem
imply the following result.

\begin{Theorem}\label{Thm: Bn(X) to Bk(X) for n>k}
For $n\ne 4$ and $n>k$ any homomorphism
$\varphi\colon B_n(X)\to B_k(X)$
is cyclic.
\end{Theorem}

\begin{proof}
For $n=3$ all is trivial since $B_2(X)$ is cyclic.
Let $n>4$. Were $\varphi$ non-cyclic, then, by Composition Theorem,
the map
$\psi=\mu\circ\varphi\colon B_n(X)\overset{\varphi}{\longrightarrow}
B_k(X)\overset{\mu}{\longrightarrow} {\mathbf S}(k)$  
would be non-cyclic,
which contradicts Theorem \ref{Thm: Bn(X) to S(k) for n>k}.
\end{proof}

\subsection{Proof of Equivariance Theorem}
\index{Theorem!Equivariance Theorem!proof of\hfill}
As we mentioned in Section 
\ref{Sec: Strategy of the proof of Tame Map Theorem},
Equivariance Theorem follows immediately from Composition Theorem.
\hfill $\square$

\subsection{Proof of Kernel Theorem.}
\label{Ss: Proof of Kernel Theorem}
\index{Theorem!Kernel Theorem!proof of\hfill}
Let $n>4$ and let $\psi\colon B_n(X)\to{\mathbf S}(n)$
be a non-cyclic homomorphism. By Theorem \ref{Thm: Surjectivity Theorem}, $\psi$
is surjective and, all the more, 
{\em transitive}\footnote{That is, its image 
$\psi(B_n(X))$ is a transitive permutation group.}.
\vskip0.2cm

\noindent All non-cyclic transitive homomorphisms 
$\psi\colon B_n\to{\mathbf S}(n)$ were described
by E. Artin \cite{Art47b} (see also \cite{Lin96b}, Section 0.5);
it follows from this description that $\Ker\psi=PB_n$.
\vskip0.2cm

\noindent Finally, for a non-cyclic homomorphism 
$\psi\colon B_n(S^2)\to{\mathbf S}(n)$ of the sphere braid group,
we apply the above result to the composition 
$$
\widetilde\psi=\psi\circ i_*\colon B_n\overset{i_*}\longrightarrow B_n(S^2)
\overset{\psi}\longrightarrow{\mathbf S}(n)
$$
and conclude that 
$(i_*)^{-1}(\Ker\psi)=\Ker(\psi\circ i_*)=\Ker\widetilde\psi=PB_n$.
Since the homomorphism $i_*\colon B_n\to B_n(S^2)$ is surjective,
$i_*((i_*)^{-1}(H))=H$ for any $H\subseteq B_n(S^2)$;
hence
$\Ker\psi=i_*((i_*)^{-1}(\Ker\psi))=i_*(PB_n)=PB_n(S^2)$
(see (\ref{CD: embeddings and coverings diagram})).
\hfill $\square$

\subsection{Proof of Invariance Theorem}
\label{Proof of Invariance Theorem}
\index{Theorem!Invariance Theorem!proof of\hfill}
We prove the following stronger statement: 
\vskip0.3cm

\noindent {\sl Let $n>4$. Then $\varphi^{-1}(PB_n(X))=PB_n(X)$
for every non-cyclic endomorphism $\varphi\colon B_n(X)\to B_n(X)$.}
\vskip0.3cm

\noindent By Composition Theorem, the homomorphism 
$\psi=\mu\circ\varphi\colon B_n(X)\overset{\varphi}
\longrightarrow B_n(X)\overset{\mu}\longrightarrow{\mathbf S}(n)$
is non-cyclic, and Kernel Theorem shows that $\Ker\psi=PB_n(X)$.
Thus,
$$
PB_n(X)=\Ker\psi=\varphi^{-1}(\Ker\mu)=\varphi^{-1}(PB_n(X))\,.
$$
Of course, this implies $\varphi(PB_n(X))
=\varphi(\varphi^{-1}(PB_n(X)))\subseteq PB_n(X)$.
\hfill $\square$
\vskip0.3cm

\noindent As we mentioned in Section \ref{Ss: Step 1}, 
Lifting Theorem follows immediately from Invariance Theorem.
\index{Theorem!Lifting Theorem!proof of\hfill}
\hfill $\square$ 


\section{Holomorphic functions omitting two values}
\label{Sec: Holomorphic functions omitting two values}

\noindent In this section we prove Theorem
\ref{Thm: Non-constant holomorphic functions on Con(X)},
that is, determine explicitly all non-constant holomorphic
functions 
${\mathcal C}_o^n({\mathbb C})\to{\mathbb C}\setminus\{0;\,1\}$
and
${\mathcal C}_o^n({\mathbb{CP}}^1)\to{\mathbb C}\setminus\{0;\,1\}$.

Furthermore, we establish that the set $L(Z)$ of all non-constant
holomorphic functions $Z\to\setminus\{0;\,1\}$ on a 
quasi-projective algebraic variety $Z$ carries a natural
structure of a finite simplicial complex.
For $Z={\mathcal C}_o^n({\mathbb C})$ and 
$Z={\mathcal C}_o^n({\mathbb{CP}}^1)$ we study this
structure in more details, which provides us with a very useful tool 
for the proof of Coherence Theorem.

\subsection{Proof of Theorem
\ref{Thm: Non-constant holomorphic functions on Con(X)}}
\label{Ss: Proof of Thm: Non-constant holomorphic functions on Con(X)}
For non-zero complex functions 
$f,g$ of the same variables, we write $f\approx g$ whenever
$f=\gamma g$ for some $\gamma\in{\mathbb C}^*$.
The following lemma is, in a sense, a multivariate 
analog of so-called {\em "abc-lemma``} for polynomials of
one variables.

\begin{Lemma}\label{Lm: 3 polynomials}
Let $P,Q,R\in{\mathbb C}[z_1,...,z_n]$ be
pairwise coprime polynomials that do not vanish
outside of the "big diagonal``
${\mathcal D}=\bigcup_{i\ne j}\{z_i=z_j\}$.
Assume that at least one of them is non-constant and
\begin{equation}\label{eq: zero sum condition}
P+Q+R=0\,.
\end{equation}
Then either
\begin{equation}\label{eq: simple PQR}
P\approx z_q-z_r\,,\qquad
Q\approx z_r-z_p\,,\qquad R\approx z_p-z_q
\end{equation}
or
\begin{equation}\label{eq: double PQR}
P\approx (z_q-z_r)(z_s-z_p)\,, \ \
Q\approx (z_r-z_p)(z_s-z_q)\,, \ \
R\approx (z_p-z_q)(z_s-z_r)\,.
\end{equation}
\end{Lemma}
\vskip0.2cm

\begin{proof} Let $a_i, b_i, c_i$ denote the degrees of $P,Q,R$
with respect to the variable $z_i$ and let 
$d_i=\max \{a_i,b_i,c_i\}$;
the condition (\ref{eq: zero sum condition}) implies that at least 
two of the three numbers $a_i,b_i,c_i$ must coincide with $d_i$.
\vskip0.2cm

Every polynomial non-vanishing outside of ${\mathcal D}$
is of the form 
$\gamma\prod_{i\ne j} (z_i-z_j)^{s_{ij}}$,
where $\gamma\in{\mathbb C}^*$ and all $s_{ij}\in{\mathbb Z}_+$. 
Hence for every $i_\circ$ we have representations
\begin{equation}\label{eq: product representations of PQR}
P=P_{i_\circ}
\prod_{j\ne i_\circ} (z_{i_\circ}-z_j)^{a_{i_\circ j}}\,, 
\ \
Q=Q_{i_\circ}
\prod_{j\ne i_\circ} (z_{i_\circ}-z_j)^{b_{i_\circ j}}\,, \ \
R=R_{i_\circ}
\prod_{j\ne i_\circ} (z_{i_\circ}-z_j)^{c_{i_\circ j}}\,,
\end{equation}
where
\begin{equation}\label{eq: Pio Qio Rio}
P_{i_\circ}\approx\prod\limits_{\stackrel{i\ne j}{i,j\ne i_\circ}}
(z_i-z_j)^{a_{i_\circ; i,j}}\,, \ \
Q_{i_\circ}\approx\prod\limits_{\stackrel{i\ne j}{i,j\ne i_\circ}}
(z_i-z_j)^{b_{i_\circ; i,j}}\,, \ \
R_{i_\circ}\approx\prod\limits_{\stackrel{i\ne j}{i,j\ne i_\circ}}
(z_i-z_j)^{c_{i_\circ; i,j}}
\end{equation}
are pairwise coprime polynomials not depending on 
$z_{i_\circ}$, and the products
$$
\prod_{j\ne i_\circ}
   (z_{i_\circ}-z_j)^{a_{i_\circ j}}\,, \ \
\prod_{j\ne i_\circ}
   (z_{i_\circ}-z_j)^{b_{i_\circ j}}\,, \ \
\prod_{j\ne i_\circ}
   (z_{i_\circ}-z_j)^{c_{i_\circ j}}
$$
are pairwise coprime polynomials in the variable $z_{i_\circ}$ 
of degrees
$a_{i_\circ},b_{i_\circ},c_{i_\circ}$, respectively.
\vskip0.2cm

We consider the following two mutually complementary cases:

\begin{itemize}

\item[{\bfit a}$)$] there exists a value $i=p$ such that the triple
$\{a_p,b_p,c_p\}$
contains at least two distinct numbers;

\item[{\bfit b}$)$] $a_i=b_i=c_i$ for every $i$.

\end{itemize}

{\bfit a}$)$ We may assume that $a_p<b_p=c_p$. Taking into account 
(\ref{eq: product representations of PQR})
(with $i_\circ = p$) and comparing the leading terms with respect
to the variable $z_p$ in identity 
(\ref{eq: zero sum condition}), we obtain
$Q_p + R_p = 0$. As $Q_p$ and $R_p$ are coprime, 
they are non-zero constants. Hence
$$
Q\approx\prod_{j\ne p} (z_p-z_j)^{b_{p j}}\,,\qquad
R\approx\prod_{j\ne p} (z_p-z_j)^{c_{p j}}
$$
and identity (\ref{eq: zero sum condition}) may be written as
\begin{equation}\label{eq: first specification of P+Q+R=0}
P_p
\prod_{j\ne p} (z_p-z_j)^{a_{p j}} +
B \prod_{j\ne p} (z_p-z_j)^{b_{p j}} +
C \prod_{j\ne p} (z_p-z_j)^{c_{p j}}=0\,,
\end{equation}
where $B,C\in{\mathbb C}^*$.

Each exponent $a_{p j}$ in (\ref{eq: first specification of P+Q+R=0})
must be zero. Indeed, each factor
$z_p-z_j$ that actually occurs in the first summand
of (\ref{eq: first specification of P+Q+R=0}) 
cannot occur in the second and the third ones;
were $a_{p j_\circ}>0$ for some $j_\circ$, then
the left hand side would be a non-trivial polynomial 
in the variable $z_{j_\circ}$, 
contradicting (\ref{eq: first specification of P+Q+R=0}).

This means that $a_p=\sum_{j\ne p} a_{p j}=0$,
the polynomial $P$ does not depend on
$z_p$, \ $P=P_p$, and identity 
(\ref{eq: first specification of P+Q+R=0}) takes the form
$$
P + B \prod_{j\ne p} (z_p-z_j)^{b_{p j}} +
C \prod_{j\ne p} (z_p-z_j)^{c_{p j}}=0\,.
$$
For the particular value $z_p=0$ the latter identity 
shows that
$$
P + B \prod_{j\ne p} (-z_j)^{b_{p j}} +
C \prod_{j\ne p} (-z_j)^{c_{p j}}=0\,,
$$
where the polynomial $P=P_p$ is given by (\ref{eq: Pio Qio Rio}) 
(with $i_\circ=p$), and the products
$$
\prod_{j\ne p} (-z_j)^{b_{p j}}\,, \qquad
\prod_{j\ne p} (-z_j)^{c_{p j}}
$$
are coprime. One can readily see that this may only happen if
$P,Q,R$ satisfy (\ref{eq: simple PQR}) with some distinct
$q,r\ne p$, which complete the proof in case $(${\bfit a}$)$.
\vskip0.2cm

{\bfit b}$)$ Since $a_i=b_i=c_i$ for any $i$, and at least 
one of the polynomials
$P,Q,R$ is non-constant, there is a value $i=s$ such that
$z_s$ actually occurs in $P$, $Q$, and $R$.
Taking into account (\ref{eq: product representations of PQR}) 
and (\ref{eq: Pio Qio Rio}) (with $i_\circ=s$),
we can write (\ref{eq: zero sum condition}) as
\begin{equation}\label{eq: second specification of P+Q+R=0}
P_s\prod_{j\ne s} (z_s-z_j)^{a_{s j}}
+ Q_s\prod_{j\ne s} (z_s-z_j)^{b_{s j}}
+ R_s\prod_{j\ne s} (z_s-z_j)^{c_{s j}} = 0\,.
\end{equation}
By our choice of $s$,
$$
\sum_{j\ne s} a_{s j}=\sum_{j\ne s} b_{s j}=\sum_{j\ne s} c_{s j}>0\,;
$$
hence there is $p\ne s$ such that $a_{s p}>0$.
As $P$ is coprime with $Q$, $R$ and contains the factor 
$(z_s-z_p)^{a_{s p}}$, we have $b_{s p}=c_{s p}=0$.
Comparing the leading terms with respect to $z_s$ 
in (\ref{eq: second specification of P+Q+R=0}), we see that
\begin{equation}\label{eq: third specification of P+Q+R=0}
P_s + Q_s + R_s = 0.
\end{equation}
Let $a_{s p}^\prime$, $b_{s p}^\prime$, $c_{s p}^\prime$
be the degrees of the polynomials $P_s,Q_s,R_s$
with respect to $z_p$. Clearly,
$$
a_{s p}^\prime=a_p-a_{s p}<a_p\,, \ \ \
b_{s p}^\prime=b_p-b_{s p}=b_p=a_p\,, \ \ \
c_{s p}^\prime=c_p-c_{s p}=c_p=a_p\,;
$$
hence
$a_{s p}^\prime < b_{s p}^\prime = c_{s p}^\prime$; in particular,
$Q_s$ and $R_s$ are non-constant. 
In view of (\ref{eq: third specification of P+Q+R=0}),
this means that the triple $\{P_s,Q_s,R_s\}$ is of the type considered 
in case $(${\bfit a}$)$. Therefore, 
$$
P_s\approx(z_q-z_r)\,, \ \ \
Q_s\approx(z_r-z_p)\,, \ \ \
R_s\approx(z_p-z_q)
$$
for certain distinct $q,r$ such that
$q\ne p,s$ and $r\ne p,s$. Consequently, 
identity (\ref{eq: second specification of P+Q+R=0}) may be written as
\begin{equation}\label{eq: fourth specification of P+Q+R=0}
\aligned
A(z_q-z_r)&(z_s-z_p)^{a_{s p}}
\prod_{j\ne s,p} (z_s-z_j)^{a_{s j}} \\
&+ B(z_r-z_p)\prod_{j\ne s,p} (z_s-z_j)^{b_{s j}}
+ C(z_p-z_q)\prod_{j\ne s,p} (z_s-z_j)^{c_{s j}}
= 0
\endaligned
\end{equation}
(with certain $A,B,C\in{\mathbb C}^*$).
In particular, on the hyperplane $z_p=z_q$ we have the identity
$$
A(z_q-z_r)(z_s-z_q)^{a_{s p}}
\prod_{j\ne s,p} (z_s-z_j)^{a_{s j}}
+ B(z_r-z_q)\prod_{j\ne s,p} (z_s-z_j)^{b_{s j}}
= 0\,,
$$
or
\begin{equation}
\label{eq: fourth specification of P+Q+R=0 on hyperplain}
A(z_s-z_q)^{a_{s p}}\prod_{j\ne s,p}(z_s-z_j)^{a_{s j}}
- B\prod_{j\ne s,p} (z_s-z_j)^{b_{s j}} = 0\,.
\end{equation}
The products $\prod_{j\ne s,p} (z_s-z_j)^{a_{s j}}$ and
$\prod_{j\ne s,p} (z_s-z_j)^{b_{s j}}$ must be coprime;
hence identity 
(\ref{eq: fourth specification of P+Q+R=0 on hyperplain}) 
can hold true only if
$$
A=B,\qquad
a_{s j}=0 \ \ \text{for all} \ \ j\ne s,p\,,\ \ \
b_{s j}=0  \ \ \text{for all} \ \ j\ne s,p,q\,,\ \ \
\text{and} \ \ a_{s p}=b_{s q}\,.
$$
Under these conditions we have
\begin{equation}\label{eq: intermediate PQR}
P\approx (z_q-z_r)(z_s-z_p)^{a_{s p}},\
Q\approx (z_r-z_p)(z_s-z_q)^{a_{s p}},\
R\approx (z_p-z_q)\prod_{j\ne s,p,q} (z_s-z_j)^{c_{s j}},
\end{equation}
and identity (\ref{eq: fourth specification of P+Q+R=0}) 
takes the form
\begin{equation}\label{eq: fifth specification of P+Q+R=0}
A(z_q-z_r)(z_s-z_p)^{a_{s p}}
+ A(z_r-z_p)(z_s-z_q)^{a_{s p}} + C R = 0\,.
\end{equation}
Taking into account the conditions $a_q=b_q=c_q$ and $a_r=b_r=c_r$,
we obtain from (\ref{eq: fifth specification of P+Q+R=0}) 
and (\ref{eq: intermediate PQR}) that
$a_{s p}=1$, \ $c_{s j}=0$ for all $j\ne r$, and $c_{s r}=1$;
hence
$P\approx (z_q-z_r)(z_s-z_p)$\,, \
$Q\approx (z_r-z_p)(z_s-z_q)$\,, \
$R\approx (z_p-z_q)(z_s-z_r)$\,. 
\end{proof}

\noindent Now we are ready to prove Theorem \ref{Thm: Non-constant holomorphic functions on Con(X)}.
\vskip0.2cm

\begin{proof}
{\bfit a}$)$ Let $\ell\subset{\mathbb C}^n$ be any com lex straight 
line such that $\ell\not\subseteq{\mathcal D}$; 
by Great Picard Theorem,
the restriction $\mu|_{{\mathcal C}_o^n({\mathbb C})\cap\ell}$
is a regular function. Thus, the function $\mu$ 
itself is rational on ${\mathbb C}^n$ and regular on 
${\mathcal C}_o^n({\mathbb C})$. That is, $\mu=-P/R$, where
polynomials $P,R\in{\mathbb C}[q_1,...,q_n]$ 
are coprime, do not vanish on 
${\mathcal C}_o^n({\mathbb C})$,
and at least one of them is non-constant.
The polynomial $Q=-P-R$ is coprime to $P$ and $R$; 
it has no zeros in ${\mathcal C}_o^n({\mathbb C})$,
as $-Q/R=1-\mu$ and the function $\mu$ does not assume the value $1$. 
Lemma \ref{Lm: 3 polynomials} implies the desired result.
\vskip0.3cm

{\bfit b}$)$ By $(${\bfit a}$)$, the restriction of $\mu$ to the 
Zariski open subset ${\mathcal C}_o^n({\mathbb C})
\subset{\mathcal C}_o^n({\mathbb{CP}}^1)$
must be either a simple ratio $sr_{ijk}$ or
a cross ratio $cr_{ijkl}$. However, the function
$sr_{ijk}=(q_k-q_i)/(q_k-q_j)$ has a pole
in ${\mathcal C}_o^n({\mathbb{CP}}^1)$
along the divisor $q_i=\infty$,
and only the cross ratios survive when passing from
${\mathcal C}_o^n({\mathbb C})$
to ${\mathcal C}_o^n({\mathbb{CP}}^1)$. 
\end{proof}

\begin{Remark}\label{Rmk: function omitting 0,1 on Com(C**) etc}
We have seen in Section \ref{Ss: Holomorphic universal covering}
that for any $m\ge 1$ there are a natural identification
\begin{equation}\label{eq: Com=Dm(CP1)}
{\mathcal D}^m({\mathbb{CP}}^1)\ni (q_1,...,q_m,0,1,\infty)
\overset{\cong}\longleftrightarrow (q_1,...,q_m)\in
{\mathcal C}_o^m(
                 {\mathbb C}\setminus\{0,1\})
\end{equation}
and a biholomorphic isomorphism
\begin{equation}
\label{eq: j(CP1,m+3): Co(m+3)(CP1) to PSL(2,C) x Dm}
j_{{\mathbb{CP}}^1,m+3}\colon{\mathcal C}_o^{m+3}({\mathbb{CP}}^1)
\overset{\cong}{\longrightarrow}
{\mathbf{PSL}}(2,{\mathbb C})
\times {\mathcal D}^m({\mathbb{CP}}^1)
={\mathbf{PSL}}(2,{\mathbb C})
\times{\mathcal C}_o^m(
                       {\mathbb C}\setminus\{0,1\})\,;
\end{equation}
$j_{{\mathbb{CP}}^1,m+3}$ and its inverse
$j_{{\mathbb{CP}}^1,m+3}^{-1}$ are defined by
\begin{equation}\label{eq: j(X,m+3)}
\aligned
j_{{\mathbb{CP}}^1,m+3}(q)=(\Theta(q),\ (&\Theta(q))^{-1}q)\,,\qquad
j_{{\mathbb{CP}}^1,m+3}^{-1}(A,z)=Az \\
&(q\in{\mathcal C}_o^{m+3}({\mathbb{CP}}^1)\,,\ \
A\in{\mathbf{PSL}}(2,{\mathbb C})\,, \ \ 
z\in{\mathcal D}^m({\mathbb{CP}}^1))\,, 
\endaligned
\end{equation}
where $\Theta\colon{\mathcal C}_o^{m+3}({\mathbb{CP}}^1)
\to{\mathbf{PSL}}(2,{\mathbb C})$ 
denote the morphism whose value at any point
$q=(q_1,...,q_{m+3})\in{\mathcal C}_o^{m+3}({\mathbb{CP}}^1)$
is a unique element
$\Theta(q)\in{\mathbf{PSL}}(2,{\mathbb C})$ 
that carries $(0,1,\infty)$ to $(q_{m+1},q_{m+2},q_{m+3})$.
By the Picard theorem, any holomorphic function
${\mathcal C}_o^{m+3}({\mathbb{CP}}^1)
\to{\mathbb C}\setminus\{0;\,1\}$ is constant on each fiber
of the projection ${\mathcal C}_o^{m+3}({\mathbb{CP}}^1)
\to{\mathcal C}_o^m(
                      {\mathbb C}\setminus\{0,1\})$;
therefore there is a one-to-one correspondence between
the set $L({\mathcal C}_o^{m}(
                              {\mathbb C}\setminus\{0,1\}))$
of all non-constant holomorphic functions
${\mathcal C}_o^m(
                 {\mathbb C}\setminus\{0,1\})
\to{\mathbb C}\setminus\{0;\,1\}$ and 
the set $L({\mathcal C}_o^{m+3}({\mathbb{CP}}^1))$
of all non-constant holomorphic functions
${\mathcal C}_o^{m+3}({\mathbb{CP}}^1)
\to{\mathbb C}\setminus\{0;\,1\}$. Due to Theorem \ref{Thm: Non-constant holomorphic functions on Con(X)} we know
all functions of the latter set; this provides us with the following explicit
description of all functions
$h\in L({\mathcal C}_o^m(
                         {\mathbb C}\setminus\{0,1\}))$, which
was also obtained by S. Kaliman \cite{Kal76b,Kal93} in a different way. 
Let us set $q_{m+1}=0$, $q_{m+2}=1$, $q_{m+3}=\infty$.
Then every function 
$h\in L({\mathcal C}_o^m(
                         {\mathbb C}\setminus\{0,1\}))$
is of the form
\begin{equation}\label{eq: function omitting 0 and one on Com(C-0,1)}
h(q_1,...,q_m)=h_{i,j,k,l}(q_1,...,q_m)
=\frac{q_i-q_j}{q_i-q_k}:\frac{q_j-q_l}{q_k-q_l}\,,
\end{equation}
where $i,j,k,l\in\{1,...,m+3\}$ are pairwise distinct. We shall use this fact
in Sec. \ref{Ss: Automorphisms of certain moduli spaces}
\vskip0.3cm

\noindent Similarly, using an isomorphism
${\mathcal C}_o^m({\mathbb C}^*))\cong{\mathbb C}^*
\times{\mathcal C}_o^{m-1}(
                           {\mathbb C}\setminus\{0,1\})$
we may identify the set
$L({\mathcal C}_o^m({\mathbb C}^*))$ first with
$L({\mathcal C}_o^{m-1}(
                        {\mathbb C}\setminus\{0,1\}))$
and then with $L({\mathcal C}_o^{m+2}({\mathbb{CP}}^1))$,
which provides an explicit description of all
non-constant holomorphic functions 
${\mathcal C}_o^m({\mathbb C}^*)\to{\mathbb C}\setminus\{0;\,1\}$
obtained in \cite{Zin77} in a different way.
\hfill $\bigcirc$
\end{Remark}

\subsection{Simplicial complex of holomorphic functions 
omitting two values}
\label{Ss: Simplicial complex of f: Z to C-0,1}
\index{Holomorphic functions omitting two values!simplicial complex of\hfill}
\index{Simplicial complex of holomorphic functions omitting two values\hfill}
Here I will describe an almost obvious
combinatorial structure, which occurs amazingly useful
for the proof of Coherence Theorem. I came to it
very recently when the first version of this paper was 
already written. In this renewed version the proof does not
become much shorter but its main ideas are now absolutely clear.
Strangely enough I paid no attention to this structure 
30 years ago. But seemingly nobody did so...
\vskip0.2cm

\noindent For an irreducible quasi-projective variety $Z$
we denote by $L(Z)$ the set of all non-constant holomorphic functions
$\lambda\colon Z\to{\mathbb C}\setminus\{0;1\}$.
It is well known that all $\lambda\in L(Z)$
are regular and $L(Z)$ is finite.\footnote{In what follows 
we are interested only in the case when $L(Z)\ne\varnothing$.}  
Clearly $\lambda^{-1}\in L(Z)$ whenever $\lambda\in L(Z)$
and $\iota\colon L(Z)\ni\lambda
\mapsto\iota(\lambda)\Def\lambda^{-1}\in L(Z)$ is an involution.
The set $L(Z)$ carries a natural structure of a finite simplicial
complex $L_{\vartriangle}(Z)$ described in the following definition.

\begin{Definition}\label{Def: divisibility and simplices}
Let $\mu,\nu\in L(Z)$; we say that $\nu$
is {\em a proper divisor of} $\mu$ and write $\nu\mid\mu$
if the quotient $\lambda=\mu:\nu\in L(Z)$, i. e.,
$\lambda\neq\const$ and $\lambda(z)\ne 1$ for
all $z\in Z$; otherwise, we write $\nu\nmid\mu$.
Clearly, $\nu\mid\mu$ is equivalent to $\mu\mid\nu$.

A non-empty subset (respectively a non-empty ordered subset) 
$\Delta\subseteq L(Z)$ is said to be a 
{\em simplex} (respectively an {\em ordered simplex}) 
of dimension $m=\#\Delta - 1$ if $\nu\mid\mu$
for each pair of distinct $\mu,\nu\in\Delta$;
the functions $\mu\in \Delta$ are said to be the {\em vertices}
of $\Delta$. Any non-empty subset of an (ordered) simplex is an (ordered)
simplex. We define the {\em boundary} $\partial\Delta$
of an ordered simplex $\Delta=[\mu_0,...,\mu_m]\subseteq L(Z)$
by the standard formula
\begin{equation}\label{eq: boundary of simplex}
\partial \Delta\Def \sum_{i=0}^m (-1)^i 
[\mu_0,...,\widehat{\mu_i},...,\mu_m]\,.
\end{equation}
Thereby $L(Z)$ becomes a finite simplicial complex
with the simplicial involution $\iota$. 

When we need to emphasize that $L(Z)$ is regarded as a simplicial complex,
i. e., the collection of all simplices $\Delta\subset L(Z)$ rather
than the set of all vertices $\mu\in L(Z)$ we write $L_{\vartriangle}(Z)$
instead of $L(Z)$ and $\Delta\in L_{\vartriangle}(Z)$ 
instead of $\Delta\subset L(Z)$.
\hfill $\bigcirc$
\end{Definition}

\begin{Lemma}\label{Lm: simplicial map f*: L(Y) to L(Z)}
Let $f\colon Z\to Y$ be a holomorphic mapping of
irreducible quasi-projective varieties.
Suppose that for each $\lambda\in L(Y)$ the composition
$$
f^*(\lambda)\Def\lambda\circ f\colon Z
\overset{f}\longrightarrow Y
\overset{\lambda}\longrightarrow{\mathbb C}\setminus\{0;1\}
$$
is non-constant $($this is certainly the case whenever $f$
is dominant$)$. Then 
\begin{equation}\label{eq: simplicial map f*: L(Y) to L(Z)}
f^*\colon L(Y)\ni\lambda\mapsto\lambda\circ f\in L(Z)
\end{equation}
is a simplicial mapping whose restriction to each
simplex $\Delta\subseteq L(Y)$ is injective. Consequently,
$f^*$ preserves dimension of simplices. Moreover,
$f^*(\iota\lambda)=\iota f^*(\lambda)$ for all $\lambda\in L(Y)$. 
\end{Lemma}

\begin{proof}
Let $\mu,\nu$ be two distinct vertices of a 
simplex $\Delta\in L_{\vartriangle}(Y)$ and let $\lambda=\mu:\nu$.
By Definition \ref{Def: divisibility and simplices}, 
$\lambda\in L(Y)$. By our assumption, 
$f^*(\mu):f^*(\nu)=f^*(\lambda)\ne\const$ and
$\mu(f(z)):\nu(f(z))=\lambda(f(z))\ne 1$ for all 
$z\in Z$. Hence $f^*(\nu)\mid f^*(\mu)$.
This shows that $f^*(\Delta)\in L_{\vartriangle}(Z)$ and
the restriction of $f^*$ to the set of all vertices of $\Delta$
is injective. The last statement of Lemma is evident.
\end{proof}  

\subsection{Complexes of simple and cross ratios}
\label{Ss: Complexes of simple and cross ratios}
Let, as usual, $X={\mathbb C}$ or $X={\mathbb{CP}}^1$
and let $Z={\mathcal C}^n_o(X)$.
Due to Theorem \ref{Thm: Non-constant holomorphic functions on Con(X)}, where we met all non-constant holomorphic 
functions $Z\to{\mathbb C}\setminus\{0;1\}$ face to face,
certain properties of the simplicial complex $L(Z)$ may be described 
in a more explicit manner.
\vskip0.2cm

\noindent{\bf Notation.}
\index{Simple and cross ratios\hfill}
It is convenient to start with a little bit
more formal situation. Let ${\mathbf q}=\{q_i\,|\ i\in{\mathbb N}\}$
be a countable set of independent variables.
For any three distinct $i,j,k\in{\mathbb N}$ 
and any four distinct $i,j,k,l\in{\mathbb N}$
the rational functions\footnote{The ground field ${\mathbf k}$ 
is not important as soon as it is commutative and 
of characteristic $0$; for instance, ${\mathbf k}={\mathbb C}$
and ${\mathbf k}={\mathbb Q}$ both are suitable.}
\begin{equation}\label{eq: simple and cross ratios}
\aligned
&sr_{ijk}\ =sr(q_i,q_j,q_k)\ \ \
\Def\frac{q_k-q_i}{q_k-q_j}\quad\text{and}\\
&cr_{ijkl}=cr(q_i,q_j,q_k,q_l)
\Def\frac{q_l-q_i}{q_l-q_j}:\frac{q_i-q_k}{q_j-q_k}
=\frac{(q_l-q_i)(q_j-q_k)}{(q_l-q_j)(q_i-q_k)}
\endaligned
\end{equation}
are referred to as {\em simple} and {\em cross} ratios,
respectively. The {\em unordered} set of variables $q$
that actually participate in such a function $\mu$
is called its {\em support}; we denote it by $\supp\mu$;
the differences of these variables in the numerator 
and the denominator of $\mu$ are said to be the {\em factors} 
of $\mu$. There are $6$ distinct simple ratios with 
support $\{q_i,q_j,q_k\}$ and $6$ distinct cross ratios with 
support $\{q_i,q_j,q_k,q_l\}$. 
(Permutations of $i,j,k,l$ formally lead to
$24$ cross ratios with the same support, but
only $6$ of them are distinct, since
Kleinian permutations  
$1,\,(i,j)(k,l),\,(i,k)(j,l),\,(i,l)(j,k)$
do not change the rational function $cr_{ijkl}$.)
The sets of all simple and all cross ratios are denoted by 
$SR({\mathbf q})$ and $CR({\mathbf q})$, respectively.
We set $L({\mathbf q})=SR({\mathbf q})\cup CR({\mathbf q})$.
The following definition is similar 
to Definition \ref{Def: divisibility and simplices}.

\begin{Definition}\label{Def: proper divisors}
\index{Proper divisors\hfill}
Let $\mu,\nu\in L({\mathbf q})$ and $\mu\ne\nu$. 
We say that $\nu$ is a {\em proper divisor} of $\mu$
and write $\nu\mid\mu$ if the quotient $\mu:\nu\in L({\mathbf q})$;
otherwise, we write $\nu\nmid\mu$. As $\lambda\in L({\mathbf q})$
is equivalent to $\lambda^{-1}\in L({\mathbf q})$, \
$\nu\mid\mu$ is equivalent to $\mu\mid\nu$.
The {\em set $(\lambda,\mu)$ of all common proper divisors}
of elements $\lambda,\mu\in L({\mathbf q})$ consists of
all $\nu\in L({\mathbf q})$ such that $\nu\mid\lambda$ and $\nu\mid\mu$. 

A non-empty finite subset $\Delta\subset L({\mathbf q})$ 
is called a {\em simplex} (respectively, an {\em ordered simplex} when
$\Delta$ is ordered) of dimension $\dim\Delta=\#\Delta -1$  
if $\mu\mid\nu$ for every two distinct $\mu,\nu\in\Delta$.\footnote{The 
latter definition {\em does not} require that the quotients 
$\nu:\mu$, $\mu:\nu$ must be elements of $\Delta$.}
The elements $\mu\in\Delta$ are said to be the {\em vertices}
of $\Delta$ and the union 
$\displaystyle\supp\Delta\Def\bigcup_{\mu\in \Delta}\supp\mu$
is called the {\em support} of $\Delta$. 
Any non-empty subset of an (ordered) simplex is an (ordered)
simplex. The boundary $\partial\Delta$ of an ordered simplex
$\Delta=[\mu_0,...,\mu_m]\subseteq L({\mathbf q})$
is defined by formula (\ref{eq: boundary of simplex}).
Thereby, $L({\mathbf q})$ becomes a simplicial complex
with the simplicial involution $\iota\colon \lambda\mapsto\lambda^{-1}$. 
\index{Simplicial complexes of simple and cross ratios\hfill} 

As in Section \ref{Ss: Simplicial complex of f: Z to C-0,1},
we use also the notation $L_{\vartriangle}({\mathbf q})$
instead of $L({\mathbf q})$ to emphasize that 
we deal with a simplicial complex, i. e., the collection of all 
simplices $\Delta\subset L({\mathbf q})$ rather
than just the set of all vertices $\mu\in L({\mathbf q})$.
There are three types of simplices 
in $L_{\vartriangle}({\mathbf q})$: two types of {\em ``pure"} ones, 
which belong either to $SR_{\vartriangle}({\mathbf q})$
(they are referred to as {\em ``simple"} ones) 
or to $CR_{\vartriangle}({\mathbf q})$
(they are referred to as {\em ``cross"} ones), 
and {\em ``mixed"} ones, which have vertices in both $SR({\mathbf q})$
and $CR({\mathbf q})$.

We denote by $SR({\mathbf q}_n)\subset SR({\mathbf q})$, 
$CR({\mathbf q}_n)\subset CR({\mathbf q})$
and $L({\mathbf q}_n)\subset L({\mathbf q})$ the subsets
consisting of all vertices of the corresponding type 
whose supports are contained in the set 
${\mathbf q}_n=\{q_1,...,q_n\}\subset{\mathbf q}$
of the first $n$ variables.
The subcomplexes $SR_{\vartriangle}({\mathbf q}_n)$, 
$CR_{\vartriangle}({\mathbf q}_n)$ and $L_{\vartriangle}({\mathbf q}_n)$
consist of all simplices
whose vertices are in $SR({\mathbf q}_n)$, $CR({\mathbf q}_n)$
and $L({\mathbf q}_n)$, respectively.\footnote{$L({\mathbf q}_2)
=CR({\mathbf q}_3)=\varnothing$. 
It should be emphasized that $SR_{\vartriangle}({\mathbf q}_n)$, 
$CR_{\vartriangle}({\mathbf q}_n)$ and $L_{\vartriangle}({\mathbf q}_n)$
{\em are not} the $n$-skeletons of the corresponding complexes.}
\hfill $\bigcirc$
\end{Definition}
\vskip0.2cm

\noindent The main purpose of the following three lemmas 
is to estimate dimension of simplices of a certain type 
provided their supports are contained in ${\mathbf q}_n$.
\vskip0.2cm

\begin{Lemma}\label{Lm: simplices of simple ratios}
\index{Simplices of simple ratios\hfill}
{\bfit a}$)$ If $\mu,\nu\in SR({\mathbf q})$ and 
$\nu\mid\mu$ then either $\mu$ and $\nu$ have the same denominator
or they have the same numerator;
in particular, $\#(\supp\mu\cap\supp\nu)=2$.
\vskip0.2cm

{\bfit b}$)$ Each proper divisor $\mu\in SR({\mathbf q})$ of $sr_{ijk}$ 
is obtained from $sr_{ijk}$ by replacing of one of the indices 
$i,j$ with some $l\ne i,j,k$.

In fact, up to a choice of $l\ne i,j,k$ and
the order of multipliers, the formula
\begin{equation}\label{eq: sr prod}
sr_{ijk}=sr_{ilk}\cdot sr_{ljk}\,.
\end{equation}
provides us with a unique possible decompositions of $sr_{ijk}$ 
to a product of two other simple ratios.
\vskip0.2cm

{\bfit c}$)$ For every four distinct $i,j,k,l\in{\mathbb N}$
\begin{equation}\label{eq: common proper divisors of two sr}
\aligned
(sr_{ijk},sr_{ilk})\cap SR({\mathbf q})
=\{sr_{imk}\,|\ m\ne i,j,k,l\}\,,\\
(sr_{ijk},sr_{ljk})\cap SR({\mathbf q})
=\{sr_{mjk}\,|\ m\ne i,j,k,l\}\,.
\endaligned
\end{equation}
\par{\bfit d}$)$ If $\Delta=\{\lambda,\mu,\nu\}\subset SR({\mathbf q})$ 
is a simplex, then either $\lambda,\mu,\nu$ have the same denominator
or they all have the same numerator; in any case $\#\supp \Delta=5$. 
\vskip0.2cm

{\bfit e}$)$ Let $n\ge 3$ and $\Delta\in SR_{\vartriangle}({\mathbf q}_n)$.
Then either all vertices of $\Delta$ have the same denominator
or they all have the same numerator. Moreover, 
$\dim\Delta\le n-3$, i. e., $\dim SR_{\vartriangle}({\mathbf q}_n)\le n-3$.
\end{Lemma}

\begin{proof}
$(${\bfit a}$)$ is evident; the proof of $(${\bfit b}$)$
is straightforward. $(${\bfit c}$)$ is a simple 
consequence of $(${\bfit b}$)$, and $(${\bfit d}$)$ follows
from $(${\bfit c}$)$. The first statement of $(${\bfit e}$)$
follows immediately from $(${\bfit d}$)$, and the second statement
follows from the first one.
\end{proof}

\begin{Lemma}\label{Lm: simplices of cross ratios}
\index{Simplices of cross ratios\hfill}
{\bfit a}$)$ If $\mu,\nu\in CR({\mathbf q})$ and 
$\nu\mid\mu$ then $\#(\supp\mu\cap\supp\nu)=3$.
\vskip0.2cm

{\bfit b}$)$ Each proper divisor $\mu\in CR({\mathbf q})$ of 
$cr_{ijkl}$ is obtained from $cr_{ijkl}$ 
by replacing of one of the indices 
$i,j,k,l$ with some $m\ne i,j,k,l$.

In fact, up to a choice of $m\ne i,j,k,l$ and
the order of multipliers, the formulae
\begin{equation}\label{eq: cr prod}
cr_{ijkl}=cr_{ijkm}\cdot cr_{ijml}=cr_{imkl}\cdot cr_{mjkl}
\end{equation}
provide us with all decompositions of $cr_{ijkl}$ 
to a product of two other cross ratios.\footnote{As $cr_{imkl}=cr_{klim}$ 
and $cr_{mjkl}=cr_{klmj}$, the Kleinian permutation $(i,k)(j,l)$, 
which does not change $cr_{ijkl}$, 
transforms these decompositions to each other.}
\vskip0.2cm

{\bfit c}$)$ For every five distinct $i,j,k,l,m\in{\mathbb N}$
\begin{equation}\label{eq: common proper divisors}
\aligned
(cr_{ijkl},cr_{ijkm})\cap CR({\mathbf q})
&=\{cr_{ijkn}\,|\ n\ne i,j,k,l,m\}\,,\\
(cr_{ijkl},cr_{ijml})\cap CR({\mathbf q})
&=\{cr_{ijnl}\,|\ n\ne i,j,k,l,m\}\,,\\
(cr_{ijkl},cr_{imkl})\cap CR({\mathbf q})
&=\{cr_{inkl}\,|\ n\ne i,j,k,l,m\}\,,\\
(cr_{ijkl},cr_{mjkl})\cap CR({\mathbf q})
&=\{cr_{njkl}\,|\ n\ne i,j,k,l,m\}\,.
\endaligned
\end{equation}
\par{\bfit d}$)$ If $\Delta=\{\lambda,\mu,\nu\}\subset CR({\mathbf q})$ 
is a simplex then $\#\supp\Delta=6$. 
\vskip0.2cm

{\bfit e}$)$ Let $n\ge 4$. Then 
$\dim CR_{\vartriangle}({\mathbf q}_n)\le n-4$,
that is, $\dim\Delta=\#\Delta-1\le n-4$ for every cross simplex  
$\Delta$ whose vertices are supported in $\{q_1,...,q_n\}$. 
\end{Lemma}

\begin{proof}
{\bfit a}$)$ The quotient $\mu:\nu$ formally 
contains $8$ factors; a mutual cancellation of two pairs of them 
cannot happen unless $\#(\supp\mu\cap\supp\nu)\ge 3$.
A straightforward verification shows that 
$\nu\nmid\mu$ whenever $\supp\mu=\supp\nu$.
\vskip0.2cm

{\bfit b}$)$ Let $cr_{ijkl}=\mu'\cdot\mu''$
for some two cross ratios $\mu',\mu''$;
due to $(${\bfit a}$)$,
we may also assume that $\supp\mu'$
is one of the sets $\{q_i,q_j,q_k,q_m\}$, 
$\{q_i,q_j,q_l,q_m\}$,
$\{q_i,q_k,q_l,q_m\}$, $\{q_j,q_k,q_l,q_m\}$. 
There are precisely $4\cdot 6=24$ 
cross ratios $\mu'$ with such a support, 
and one may check them all, one by one; this leads to 
the desired result.
\vskip0.2cm

{\bfit c}$)$ By $(${\bfit b}$)$, to get a common proper
divisor $\nu\in CR({\mathbf q})$ of $cr_{ijkl}$ and $cr_{ijkm}$,
one must replace one of the indices in each of the latter two
cross ratios with some $n\ne i,j,k,l,m$ so that 
the results coincide with each other. Only the cross ratio
$\nu=cr_{ijkn}$ appears in such a way; this proves the
first of relations (\ref{eq: common proper divisors}).
As Kleinian permutations $(i,j)(k,l)$, $(i,k)(j,l)$ and
$(i,l)(j,k)$ do not change the rational function $cr_{ijkl}$, 
and for any $s\ne i,j,k,l$ we have 
$cr_{jils}=cr_{ijsl}$, $cr_{klis}=cr_{iskl}$ 
and $cr_{lkjs}=cr_{sjkl}$, three other relations
follow from the first one.
\vskip0.2cm

{\bfit d}$)$ Let $\lambda=cr_{ijkl}$. By $(${\bfit b}$)$, 
there is $m\ne i,j,k,l$ such that $\mu$ coincides with one 
of the four cross ratios $cr_{ijkm}$, $cr_{ijml}$, $cr_{imkl}$, 
$cr_{mjkl}$. In each of these cases one of relations 
(\ref{eq: common proper divisors}) shows that 
$\#(\supp\lambda\cup\supp\mu\cup\supp\nu)=6$.
\smallskip

{\bfit e}$)$ By $(${\bfit a}$)$, statement ({\bfit e}) 
is valid for $n=4$. Suppose that for some $n\ge 4$ 
we have $\dim CR_{\vartriangle}({\mathbf q}_n)\le n-4$ 
and consider any simplex 
$\Delta\in CR_{\vartriangle}({\mathbf q})$ with support 
$\supp\Delta\subseteq\{q_1,...,q_{n+1}\}$. 
We must prove that $\dim\Delta\le n-3$, i. e., $\#\Delta\le n-2$. 

Take two distinct $\lambda,\mu\in\Delta$. 
According to $(${\bfit a}$)$, $\#(\supp\lambda\cap\supp\mu)=3$;
thus, $\supp\lambda=\{q_i,q_j,q_k,q_l\}$
and $\supp\mu=\{q_i,q_j,q_k,q_m\}$ for some
distinct $i,j,k,l,m\in\{1,...,n+1\}$.

Let $\Delta(\widehat{m})=\{\tau\in\Delta\,|\ q_m\notin\supp\tau\}$; 
clearly, $\lambda\in\Delta(\widehat{m})$ 
and $\mu\notin\Delta(\widehat{m})$. By the induction hypothesis, 
$\dim\Delta=\#\Delta(\widehat{m})-1\le n-4$. Let us show that actually 
$\Delta\setminus\Delta(\widehat{m})=\{\mu\}$; 
this will imply $\#\Delta =\#\Delta(\widehat{m})+1\le (n-3)+1=n-2$ 
and complete the step of induction. 

Suppose to the contrary that $\nu\in\Delta\setminus\Delta(\widehat{m})$
and $\nu\ne\mu$; clearly, $q_m\in\supp\nu$. 
As $\lambda\in\Delta(\widehat{m})$, we have $\nu\ne\lambda$. 
Thereby, $\{\lambda,\mu,\nu\}$ is a simplex; by $(${\bfit d}$)$, 
$\#\supp\Delta=\#(\supp\lambda\cup\supp\mu\cup\supp\nu)=6$.
Since $\supp\lambda\cup\supp\mu=\{q_i,q_j,q_k,q_l,q_m\}$, 
the set $\supp\nu$ must contain at least one additional
variable, say $q_p$, distinct from all $q_i,q_j,q_k,q_l,q_m$.
As $\supp\lambda=\{q_i,q_j,q_k,q_l\}$ and 
$q_m,q_p\in\supp\nu$, we see that 
$\#(\supp\lambda\cap\supp\nu)\le 2$, which contradicts
$(${\bfit a}$)$.
\end{proof}

\begin{Lemma}\label{Lm: mixed simplices}
\index{Mixed simplices\hfill}
{\bfit a}$)$ Each mixed simplex $\Delta$ contains precisely one vertex
$\mu\in SR({\mathbf q})$.

{\bfit b}$)$ $\dim\Delta\le n-3$ for any mixed simplex
with support $\supp\Delta\subseteq\{q_1,...,q_n\}$. 
\end{Lemma}

\begin{proof}
{\bfit a}$)$ As $\Delta$ is mixed, there are 
$\mu\in\Delta\cap SR({\mathbf q})$ 
and $\nu\in\Delta\cap CR({\mathbf q})$. We may assume that
$\mu=sr_{ijk}$; then it is easily seen that 
$\nu=cr_{ijlk}$ for a certain $l\ne i,j,k$.
Up to an order of multipliers, there exist exactly 
two distinct decompositions of $\nu=cr_{ijlk}$
to a product of two simple ratios:
$$
cr_{ijlk}=sr_{ijk}\cdot sr_{jil}=sr_{kli}\cdot sr_{lkj}\,.
$$
Each of the three simple ratios $sr_{jil}$, $sr_{kli}$, $sr_{lkj}$
is not a proper divisor of $\mu=sr_{ijk}$; hence $\Delta$ does not
contain any simple ratio $\lambda\ne\mu$.
\vskip0.2cm

{\bfit b}$)$ By $(${\bfit a}$)$, there is a single vertex
$\mu\in\Delta\cap SR({\mathbf q})$. Thus, 
$\Delta'\Def\Delta\setminus\{\mu\}$ is a simplex contained
in $CR({\mathbf q})$ with support in $\{q_1,...,q_n\}$.
By Lemma \ref{Lm: simplices of cross ratios}$(${\bfit e}$)$
$\dim\Delta'\le n-4$, which implies $\dim\Delta\le n-4$.
\end{proof}

\begin{Example}\label{Ex: simplices} 
{\bfit a}$)$ Let $n\ge 3$.
Then $\Delta=\{sr_{3,2,1},\,sr_{4,2,1},...,sr_{n,2,1}\}
\subset SR({\mathbf q})$ is a simplex of dimension $n-3$ 
with support in $\{q_1,...,q_n\}$.
\vskip0.2cm

{\bfit b}$)$ Let $s\ge 4$. Then 
$\Delta=\{cr_{1,2,3,4},\,cr_{1,2,3,5},...,cr_{1,2,3,n}\}
\subset CR({\mathbf q})$ 
is a simplex of dimension $n-4$ with support in $\{q_1,...,q_n\}$. 
\vskip0.2cm

{\bfit c}$)$ Let $n\ge 3$. 
Then $\Delta=\{sr_{1,2,3};\,cr_{1,2,4,3},\,cr_{1,2,5,3},...,cr_{1,2,n,3}\}$
is a mixed simplex of dimension $n-3$ 
with support in $\{q_1,...,q_n\}$.
\vskip0.2cm
 
These examples show that the upper bounds of $\dim\Delta$ established in 
Lemma \ref{Lm: simplices of simple ratios}$(${\bfit e}$)$,
Lemma \ref{Lm: simplices of cross ratios}$(${\bfit e}$)$ and
Lemma \ref{Lm: mixed simplices}$(${\bfit b}$)$
are the best possible ones.
\hfill $\bigcirc$
\end{Example}

\noindent Due to Theorem \ref{Thm: Non-constant holomorphic functions on Con(X)} the complexes 
$L_{\vartriangle}({\mathcal C}_o^n({\mathbb C}))$ and 
$L_{\vartriangle}({\mathcal C}_o^n({\mathbb{CP}}^1))$
defined according to Section 
\ref{Ss: Simplicial complex of f: Z to C-0,1}
may be identified to the complex $L_{\vartriangle}({\mathbf q}_n)$ 
and to the subcomplex $CR_{\vartriangle}({\mathbf q}_n)
\subset L_{\vartriangle}({\mathbf q}_n)$,
respectively. 
When it is convenient we use also the notation 
$$
\aligned
&SR({\mathcal C}_o^n({\mathbb C}))\,, \ 
SR_{\vartriangle}({\mathcal C}_o^n({\mathbb C}))\,, \
CR({\mathcal C}_o^n({\mathbb C}))
=CR({\mathcal C}_o^n({\mathbb{CP}}^1))\\ 
&\text{and} \ \ 
CR_{\vartriangle}({\mathcal C}_o^n({\mathbb C}))
=CR_{\vartriangle}({\mathcal C}_o^n({\mathbb{CP}}^1))
\endaligned
$$
instead of $SR({\mathbf q}_n)$, $SR_{\vartriangle}({\mathbf q}_n)$, 
$CR({\mathbf q}_n)$ and $CR_{\vartriangle}({\mathbf q}_n)$.
\vskip0.3cm

\noindent Lemma \ref{Lm: simplices of simple ratios}, 
Lemma \ref{Lm: simplices of cross ratios} and 
Example \ref{Ex: simplices} imply the following corollary.

\begin{Corollary}\label{Crl: dimension of L(Con(X))}
{\bfit a}$)$ $\dim SR_{\vartriangle}({\mathbf q}_n)=n-3$ \
and \ $\dim CR_{\vartriangle}({\mathbf q}_n)=n-4$. 
\vskip0.2cm

\begin{itemize}

\item[{\bfit b}$)$]
\index{Dimension of $L_{\vartriangle}({\mathcal C}_o^n({\mathbb C}))$\hfill}
$\dim L_{\vartriangle}({\mathcal C}_o^n({\mathbb C}))
=\dim SR_{\vartriangle}({\mathcal C}_o^n({\mathbb C}))
=\dim SR_{\vartriangle}({\mathbf q}_n)=n-3$ \ and

\item[]
\index{Dimension of $L_{\vartriangle}({\mathcal C}_o^n({\mathbb{CP}^1}))$\hfill}
$\dim L_{\vartriangle}({\mathcal C}_o^n({\mathbb{CP}^1}))
=\dim CR_{\vartriangle}({\mathcal C}_o^n({\mathbb{CP}^1}))
=\dim CR_{\vartriangle}({\mathbf q}_n)=n-4$. 
\hfill $\bigcirc$
\end{itemize}
\end{Corollary}

\subsection{${\mathbf S}(n)$ action in $L({\mathbf q}_n)$}
\label{Ss: S(n) action in Ln(q)}
\index{${\mathbf S}(n)$ action in $L({\mathbf q}_n)$\hfill}
The natural ${\mathbf S}(n)$ action in $\{1,...,n\}$
induces the left ${\mathbf S}(n)$ action in $L({\mathbf q}_n)$;
in more details, for $\mu\in L({\mathbf q}_n)$ and 
$\sigma\in{\mathbf S}(n)$ the rational function
$\sigma\mu$ is defined by $(\sigma\mu)(q)=\mu(\sigma^{-1}q)$, i. e.,
$(\sigma\mu)(q_1,...,q_n)
=\mu(q_{\sigma(1)},...,q_{\sigma(n)})$. 
This action is transitive on the sets of vertices
$SR({\mathbf q}_n)$ and $CR({\mathbf q}_n)$ (for 
${\mathbf S}(n)$ action in $\{1,...,n\}$ is $n$ transitive). 

Furthermore, if $\mu,\nu\in L({\mathbf q}_n)$ and $\nu\mid\mu$ 
then for each $\sigma\in {\mathbf S}(n)$ we have
$(\sigma\mu)(q):(\sigma\nu)(q)=\mu(\sigma^{-1}q):\nu(\sigma^{-1}q)
=(\mu:\nu)(\sigma^{-1}q)$, which shows that 
$\sigma\nu\mid\sigma\mu$. Hence for each simplex 
$\Delta=\{\mu_0,...,\mu_s\}\in L_{\vartriangle}({\mathbf q}_n)$
and each $\sigma\in{\mathbf S}(n)$ the set 
$\sigma\Delta=\{\sigma\mu_0,...,\sigma\mu_s\}$ is the simplex
of the same type and dimension. Thus we have the (left)
dimension preserving ${\mathbf S}(n)$ action in
$SR_{\vartriangle}({\mathbf q}_n)$ and 
$CR_{\vartriangle}({\mathbf q}_n)$.

We denote by $SR_{\vartriangle^m}({\mathbf q}_n)\subseteq
SR_{\vartriangle}({\mathbf q}_n)$ and 
$CR_{\vartriangle^m}({\mathbf q}_n)\subseteq
CR_{\vartriangle}({\mathbf q}_n)$ 
the ${\mathbf S}(n)$ invariant subsets consisting of all
$m$ simplices of the corresponding type. 

\begin{Definition}\label{Def: normal form of simplices}
\index{Normal form of simplices\hfill}
We define the following three {\em normal forms} of $m$ simplices:
\begin{equation}\label{eq: 1st normal form}
\Delta_S^m\Def\{sr_{3,2,1},\,sr_{4,2,1},...,sr_{m+3,2,1}\}
\in SR_{\vartriangle}({\mathbf q}_n) \ \ (0\le m\le n-3)\,,
\end{equation}
\begin{equation}\label{eq: 2nd normal form}
\Delta_S^{-m}\Def\{sr_{2,3,1},\,sr_{2,4,1},...,sr_{2,m+3,1}\}
\in SR_{\vartriangle}({\mathbf q}_n) \ \ (0\le m\le n-3)
\end{equation}
and 
\begin{equation}\label{eq: 3rd normal form}
\Delta_C^m\Def\{cr_{1,2,3,4},\,cr_{1,2,3,5},...,cr_{1,2,3,m+4}\}
\in CR_{\vartriangle}({\mathbf q}_n) \ \ (0\le m\le n-4)\,.
\end{equation}
A simplex in a normal form we call {\em normal}.
\end{Definition}

\noindent The following lemma shows that 
the ${\mathbf S}(n)$ action in $L({\mathbf q}_n)$
allows to carry any simplex to a normal form.

\begin{Lemma}\label{Lm: S(n) orbits in SRn and CRn}
{\bfit a}$)$
\index{${\mathbf S}(n)$ orbits of simplices of simple ratios\hfill}
For $m\ge 1$ the set of all simple $m$ simplices
$SR_{\vartriangle^m}({\mathbf q}_n)$ is the disjoint union
of two ${\mathbf S}(n)$ orbits  
${\mathbf S}(n)\Delta_S^m$ and ${\mathbf S}(n)\Delta_S^{-m}$
of the normal $m$ simplices $\Delta_S^m$ and 
$\Delta_S^{-m}$. That is, any simple $m$ simplex can be carried 
to one of the normal forms $(\ref{eq: 1st normal form})$,
$(\ref{eq: 2nd normal form})$ 
by a permutation $\sigma\in{\mathbf S}(n)$.
\vskip0.2cm

{\bfit b}$)$
\index{${\mathbf S}(n)$ orbits of simplices of cross ratios\hfill}
The ${\mathbf S}(n)$ action on
$CR_{\vartriangle^m}({\mathbf q}_n)$ is transitive,
that is, any cross $m$ simplex can be carried  
to the normal form $(\ref{eq: 3rd normal form})$
by a permutation $\sigma\in{\mathbf S}(n)$.
\end{Lemma}

\begin{proof}
{\bfit a}$)$ Let 
$\Delta=\{\mu_0,...,\mu_m\}\in SR_{\vartriangle^m}({\mathbf q}_n)$.
As permutations are admitted we may assume that $\mu_0=sr_{3,2,1}$.
By Lemma \ref{Lm: simplices of simple ratios}$(${\bfit b}$)$,
$\mu_1$ must be one of the functions 
$sr_{l,2,1}$, $sr_{3,l,1}$ with a certain $l\ne 1,2,3$.
A suitable permutation of $\{4,...,n\}$ 
keeps $\mu_0$ and carries $\mu_1$ to
one of the functions $sr_{4,2,1}$ and $sr_{3,4,1}$.

In the first case $\{\mu_0,\,\mu_1\}=\{sr_{3,2,1},\,sr_{4,2,1}\}$;
in the second case we apply the transposition $(2,3)$ 
to all vertices of $\Delta$ and get
$\{\mu_0,\,\mu_1\}=\{sr_{2,3,1},\,sr_{2,4,1}\}$.

If $m=1$ we are done. For $m\ge 2$ the simplex $\Delta$ has at 
least one more vertex $\mu_2$; so either 

\begin{itemize}

\item[$(i)$] $\Delta=\{\mu_0,\,\mu_1,\,\mu_2,...,\mu_m\}
=\{sr_{3,2,1},\,sr_{4,2,1},\,\mu_2,...,\mu_m\}$ \ or

\item[$(ii)$] $\Delta=\{\mu_0,\,\mu_1,\,\mu_2,...,\mu_m\}
=\{sr_{2,3,1},\,sr_{2,4,1},\,\mu_2,...,\mu_m\}$\,.
\end{itemize}

In case $(i)$, $\mu_0$ and $\mu_1$ 
have the same denominator $q_1-q_2$. 
Lemma \ref{Lm: simplices of simple ratios}$(${\bfit e}$)$ shows
that all $\mu_2,...,\mu_m$ also have the
same denominator $q_1-q_2$.  
By Lemma \ref{Lm: simplices of simple ratios}$(${\bfit c}$)$, 
they must be of the form 
$\mu_2=sr_{k_2,2,1},...,\mu_m=sr_{k_m,2,1}$,
where all $k_2,...,k_m\in \{5,...,n\}$ are distinct.
A suitable permutation of $\{5,...,n\}$ keeps the forms of
$\mu_0$, $\mu_1$ and carries $\Delta$ to the normal form
(\ref{eq: 1st normal form}).  

In case $(ii)$, $\mu_0$ and $\mu_1$ 
have the same numerator $q_1-q_2$ and 
Lemma \ref{Lm: simplices of simple ratios}$(${\bfit e}$)$ shows
that so do all other vertices $\mu_2,...,\mu_m$ of $\Delta$. 
As above, we conclude that they must be of the form 
$\mu_2=sr_{2,k_2,1},...,\mu_m=sr_{2,k_m,1}$,
where all $k_2,...,k_m\in \{5,...,n\}$ are distinct.
A suitable permutation of $\{5,...,n\}$ does not
change $\mu_0$, $\mu_1$ and brings $\Delta$ to the normal form
(\ref{eq: 2nd normal form}). 
\vskip0.2cm

{\bfit b}$)$ Let 
$\Delta=\{\mu_0,...,\mu_m\}\in CR_{\vartriangle^m}({\mathbf q}_n)$.
We may assume that $\mu_0=cr_{1,2,3,4}$.
By Lemma \ref{Lm: simplices of cross ratios}$(${\bfit b}$)$,
$\mu_1$ is one
of the functions $cr_{1,2,3,l}$, \ $cr_{1,2,l,4}$, \ 
$cr_{1,l,3,4}$, \ $cr_{l,2,3,4}$ for a certain $l\ne 1,2,3,4$.
A suitable permutation of $\{5,...,n\}$ 
keeps $\mu_0$ and carries $\mu_1$ to one of the functions 
$cr_{1,2,3,5}$, \ $cr_{1,2,5,4}$, \ $cr_{1,5,3,4}$, \ $cr_{5,2,3,4}$.
The Kleinian permutation $(1,3)(2,4)$ 
does not change the rational function $cr_{1,2,3,4}$ and carries
the pair $cr_{1,2,3,5}$, \ $cr_{1,2,5,4}$ to the pair 
$cr_{1,5,3,4}$, \ $cr_{5,2,3,4}$; hence we may restrict ourselves
to the following two cases: either
$\{\mu_0,\,\mu_1\}=\{cr_{1,2,3,4},\,cr_{1,2,3,5}\}$ or
$\{\mu_0,\,\mu_1\}=\{cr_{1,2,3,4},\,cr_{1,2,5,4}\}$.
In fact, the Kleinian permutation $(1,2)(3,4)$
keeps $cr_{1,2,3,4}$ and brings $cr_{1,2,5,4}$
to $cr_{1,2,3,5}$; thus we are left with the case when
$\Delta=\{\mu_0,\,\mu_1,...,\mu_m\}=
\{cr_{1,2,3,4},\,cr_{1,2,3,5},...,\mu_m\}$.

If $m=1$ we are done. For $m\ge 2$ 
the first of relations (\ref{eq: common proper divisors})
in Lemma \ref{Lm: simplices of cross ratios}$(${\bfit c}$)$ 
shows that $\mu_2,...,\mu_m$ must be of the form 
$\mu_2=cr_{1,2,3,k_2},...,\mu_m=cr_{1,2,3,k_m}$,
where all $k_2,...,k_m\in \{6,...,n\}$ are distinct.
A suitable permutation of $\{6,...,n\}$ keeps  
$\mu_0$, $\mu_1$ and brings $\Delta$ to the normal form
(\ref{eq: 3rd normal form}). 
\end{proof}

\section{Coherence Theorem}
\label{Sec: Coherence Theorem}

\noindent In this section we prove Coherence Theorem
formulated in Section \ref{Ss: Step 4}, which
completes the proof of Tame Map Theorem.
\vskip0.2cm

\noindent It is clear that all components of a
strictly equivariant endomorphism 
$$
f=(f_1,...,f_n)\colon {\mathcal C}_o^n(X)\to {\mathcal C}_o^n(X)
$$
are non-constant holomorphic mappings ${\mathcal C}_o^n(X)\to X$.
The following lemma shows that 
Lemma \ref{Lm: simplicial map f*: L(Y) to L(Z)} applies to
strictly equivariant endomorphisms.

\begin{Lemma}\label{Lm: non-constant simple and cross ratios}
Let $f=(f_1,...,f_n)\colon{\mathcal C}_o^n(X)
\to{\mathcal C}_o^n(X)$ be a strictly equivariant 
endomorphism and $\mu\colon{\mathcal C}_o^n(X)
\to{\mathbb C}\setminus\{0;\,1\}$ be a non-constant
holomorphic function. Then the composition
$\mu\circ f\colon{\mathcal C}_o^n(X)
\overset{f}\longrightarrow{\mathcal C}_o^n(X)
\overset{\mu}\longrightarrow{\mathbb C}
\setminus\{0;\,1\}$ is non-constant.
\end{Lemma}

\begin{proof}
Suppose to the contrary that $\mu\circ f=c=\const\ne 0, 1$.
Then $(\mu\circ f)(\sigma q)\equiv c$
for all $\sigma \in{\mathbf S}(n)$. Since $f$ is strictly equivariant, 
there is $\alpha\in\Aut{\mathbf S}(n)$ such that
$f(\theta q)=\alpha(\theta) f(q)$
for all $\theta\in{\mathbf S}(n)$ and $q\in{\mathcal C}_o^n(X)$,
so that $c\equiv \mu(f(\theta q))=\mu(\alpha(\theta) f(q))$.
By Theorem \ref{Thm: Non-constant holomorphic functions on Con(X)}, either
$\mu\in SR({\mathbf q}_n)$ or $\mu\in CR({\mathbf q}_n)$
(for $X={\mathbb{CP}}^1$ the first case does not occur).
${\mathbf S}(n)$ acts transitively on the sets $SR({\mathbf q}_n)$ and 
$CR({\mathbf q}_n)$, which are invariant under the involutions 
$\lambda\mapsto\lambda^{-1}$ and $\lambda\mapsto 1-\lambda$;
thus $\mu^{-1}=s\mu$ and $1-\mu = t\mu$ for certain
$s,t\in{\mathbf S}(n)$. Let $\sigma=\alpha^{-1}(s^{-1})$
and $\tau=\alpha^{-1}(t^{-1})$. Then 
$$
c^{-1}\equiv[\mu(f(q))]^{-1}\equiv
s\mu(f(q))\equiv\mu(s^{-1}f(q))\equiv\mu(\alpha(\sigma)f(q))
\equiv\mu(f(\sigma q))\equiv c
$$
and 
$$
1-c\equiv 1-\mu(f(q))\equiv t\mu(f(q))\equiv\mu(t^{-1}f(q))
\equiv\mu(\alpha(\tau)f(q))
\equiv\mu(f(\tau q))\equiv c\,,
$$
and the contradiction $-1=c=1/2$ ensues.\footnote{My original proof
was at least twice as long; I am grateful to Yoel Feler who noticed that
it may be simplified by involving the second involution
$\lambda\mapsto 1-\lambda$.}
\end{proof}

\noindent Recall the notation $t(X)=2$ for $X={\mathbb C}$
and $t(X)=3$ for $X={\mathbb{CP}}^1$. The following
result is just an evident combination of
Theorem \ref{Thm: Non-constant holomorphic functions on Con(X)},
Lemma \ref{Lm: non-constant simple and cross ratios}
and Lemma \ref{Lm: simplicial map f*: L(Y) to L(Z)}. 

\begin{Corollary}
\label{Crl: for strictly equivariant f f^* is simplicial}
Let $n>t(X)$\footnote{See Notation \ref{Not: t(X)}.}
and let $f$ be a strictly equivariant endomorphism
of ${\mathcal C}_o^n(X)$.
The induced map $f^*\colon L({\mathcal C}_o^n(X))\ni\mu
\mapsto\mu\circ f\in L({\mathcal C}_o^n(X))$ is simplicial.
\hfill $\square$
\end{Corollary}

\begin{Lemma}\label{Lm: S(n) action and f*}
Let $f$ be a strictly equivariant endomorphism of
${\mathcal C}_o^n(X)$ and let $\alpha$ be the automorphism of 
${\mathbf S}(n)$ related to $f$ 
{\rm (see (\ref{eq: equivariance condition}) and 
Definition \ref{Def: strict equivariance})}.
\vskip0.2cm

{\bfit a}$)$ The map $f^*\colon L({\mathcal C}_o^n(X))\to
L({\mathcal C}_o^n(X))$ is strictly equivariant meaning that
\begin{equation}\label{eq: f* is strictly equivariant}
f^*(\sigma\mu)=\alpha^{-1}(\sigma)f^*(\mu) \ \text{for all} \ 
\sigma\in{\mathbf S}(n) \ \text{and all} 
\ \mu\in L({\mathcal C}_o^n({\mathbb C}))\,.
\end{equation}
Consequently, for any vertex 
$\mu\in L({\mathcal C}_o^n(X))$ the $f^*$-image of its 
${\mathbf S}(n)$ orbit coincides with the ${\mathbf S}(n)$ orbit
of the vertex $f^*(\mu)$. 
\vskip0.2cm

{\bfit b}$)$ If $X={\mathbb C}$ then 
$f^*(SR({\mathcal C}_o^n({\mathbb C})))
=SR({\mathcal C}_o^n({\mathbb C}))$.
\end{Lemma}

\begin{proof} 
{\bfit a}$)$ Let $s=\alpha^{-1}(\sigma)$. 
The proof of (\ref{eq: f* is strictly equivariant}) 
is straightforward:
$$
\aligned
\ [f^*(\sigma\mu)](q)&=(\sigma\mu)(f(q))=(\alpha(s)\mu)(f(q))
=\mu(\alpha(s^{-1})f(q))\\
&=\mu(f(s^{-1} q))=(f^*(\mu))(s^{-1}q)
=[sf^*(\mu)](q)=[\alpha^{-1}(\sigma)f^*(\mu)](q)\,.
\endaligned
$$
The permutation $s$ runs over the whole ${\mathbf S}(n)$ 
whenever $\sigma$ does so, and vice versa. 
Hence (\ref{eq: f* is strictly equivariant})
implies $f^*({\mathbf S}(n)\mu)={\mathbf S}(n)f^*(\mu)$
for each $\mu\in L({\mathcal C}_o^n(X))$.
\vskip0.2cm

{\bfit b}$)$ By Lemma \ref{eq: simplicial map f*: L(Y) to L(Z)}, 
the map $f^*$
preserves dimension of simplices; taking into account
Corollary \ref{Crl: dimension of L(Con(X))}, we have
$\dim f^*(SR_{\vartriangle}({\mathcal C}_o^n({\mathbb C})))
=\dim SR_{\vartriangle}({\mathcal C}_o^n({\mathbb C}))=n-3$
and $\dim CR_{\vartriangle}({\mathcal C}_o^n({\mathbb C}))=n-4$,
which shows that $f^*(SR_{\vartriangle}({\mathcal C}_o^n({\mathbb C})))
\ne CR_{\vartriangle}({\mathcal C}_o^n({\mathbb C}))$. Thus
$f^*(SR_{\vartriangle}({\mathcal C}_o^n({\mathbb C})))
=SR_{\vartriangle}({\mathcal C}_o^n({\mathbb C}))$, which is certainly
equivalent to $f^*(SR({\mathcal C}_o^n({\mathbb C})))
=SR({\mathcal C}_o^n({\mathbb C}))$.
\end{proof}

\begin{Remark}\label{Rmk: f*(CR(Con(C)))=CR(Con(C))}
It follows from Coherence Theorem (which is not proved yet!)  
that $f^*(CR({\mathcal C}_o^n({\mathbb C})))
=CR({\mathcal C}_o^n({\mathbb C}))$; we do not
need this fact at the moment.
\hfill $\bigcirc$
\end{Remark}

\subsection{Proof of Coherence Theorem for $X={\mathbb C}$}
\label{Ss: Proof of Coherence Theorem for X=C}
\index{Theorem!Coherence Theorem!proof of\hfill}
Let $n>4$ and let $f$ be a strictly equivariant endomorphism of
${\mathcal C}_o^n({\mathbb C})$. In view of what has been
done above, in this case Coherence Theorem
is equivalent to the following
\vskip0.3cm

\noindent{\bfit Claim 1.} {\sl Let $\Delta
=\Delta_S^{n-3}\in SR_{\vartriangle}({\mathcal C}_o^n({\mathbb C}))$
be the normal simplex of the form $(\ref{eq: 1st normal form})$,
that is, the simplex with the vertices $\mu_i=sr_{i,2,1}$, 
$i=3,...,n$. Then $f^*(\Delta)\in{\mathbf S}(n)\Delta$.}

\begin{proof}
By Lemma \ref{Lm: S(n) action and f*}$(${\bfit b}$)$,
$f^*(\Delta)\in SR_{\vartriangle^{n-3}}({\mathcal C}_o^n({\mathbb C}))$.
According to Lemma 
\ref{Lm: S(n) orbits in SRn and CRn}$(${\bfit a}$)$,
$SR_{\vartriangle^{n-3}}({\mathcal C}_o^n({\mathbb C}))$ is
the disjoint union of the orbits ${\mathbf S}(n)\Delta_S^{n-3}$ 
and ${\mathbf S}(n)\Delta_S^{-(n-3)}$.
So it suffices to show that 
$f^*(\Delta)\not\in{\mathbf S}(n)\Delta_S^{-(n-3)}$; in view of 
Lemma \ref{Lm: S(n) action and f*}$(${\bfit a}$)$, the latter
is equivalent to $f^*({\mathbf S}(n)\Delta)\not\in\Delta_S^{-(n-3)}$.

Suppose to the contrary that $f^*(s\Delta)=\Delta_S^{-(n-3)}$
for some $s\in{\mathbf S}(n)$. Then 
$$
\frac{f_1(\sigma q)-f_i(\sigma q)}{f_1(\sigma q)-f_2(\sigma q)}
=\frac{q_1-q_2}{q_1-q_i}
$$
for all $i=3,...,n$, where $\sigma=\alpha^{-1}(s^{-1})$
and $\alpha$ is the automorphism of ${\mathbf S}(n)$ related to $f$.
This implies that for all 
distinct $j,k\in\{3,...,n\}$. 
$$
\aligned
\frac{f_k(\sigma q)-f_j(\sigma q)}{f_1(\sigma q)-f_2(\sigma q)}
&=\frac{f_1(\sigma q)-f_j(\sigma q)}{f_1(\sigma q)-f_2(\sigma q)}
-\frac{f_1(\sigma q)-f_k(\sigma q)}{f_1(\sigma q)-f_2(\sigma q)}\\
&=\frac{q_1-q_2}{q_1-q_j}-\frac{q_1-q_2}{q_1-q_k}
=\frac{(q_1-q_2)(q_j-q_k)}{(q_1-q_j)(q_1-q_k)}\,.
\endaligned
$$
Consequently, for distinct $i,j,k\in\{3,...,n\}$ we have
$$
\aligned
\ [f^*(&sr_{ijk})](\sigma q)=sr_{ijk}(f(\sigma q))
=\frac{f_k(\sigma q)-f_i(\sigma q)}{f_k(\sigma q)-f_j(\sigma q)}\\
&=\frac{f_k(\sigma q)-f_i(\sigma q)}{f_1(\sigma q)-f_2(\sigma q)}
:\frac{f_k(\sigma q)-f_j(\sigma q)}{f_1(\sigma q)-f_2(\sigma q)}\\
&=\frac{(q_1-q_2)(q_i-q_k)}{(q_1-q_i)(q_1-q_k)}
:\frac{(q_1-q_2)(q_j-q_k)}{(q_1-q_j)(q_1-q_k)}
=\frac{q_k-q_i}{q_k-q_j}:\frac{q_i-q_1}{q_j-q_1}=cr_{i,j,1,k}(q)\,,
\endaligned
$$
which contradicts Lemma \ref{Lm: S(n) action and f*}$(${\bfit b}$)$.
\end{proof}

\subsection{Proof of Coherence Theorem for $X={\mathbb{CP}}^1$}
\label{Ss: Proof of Coherence Theorem for X=CP1}
\index{Theorem!Coherence Theorem!proof of\hfill}
Let $n\ge 4$ and let $f=(f_1,...,f_n)$ be a strictly equivariant 
endomorphism of ${\mathcal C}_o^n({\mathbb{CP}}^1)$.
In this case Coherence Theorem is equivalent to the following
\vskip0.3cm

\noindent{\bfit Claim 2.} {\sl Let $\Delta=\Delta_C^{n-4}
\in CR_{\vartriangle^{n-4}}({\mathcal C}_o^n({\mathbb{CP}^1}))$
be the normal simplex of the form $(\ref{eq: 3rd normal form})$,
that is, the simplex with the vertices $\mu_i=cr_{1,2,3,i}$, 
$i=4,...,n$. Then $f^*(\Delta)$ is a $n-4$ simplex
contained in the ${\mathbf S}(n)$ orbit ${\mathbf S}(n)\Delta$
of $\Delta$ itself.}
\vskip0.3cm
 
\begin{proof} 
By Lemma \ref{Lm: non-constant simple and cross ratios}
and Lemma \ref{Lm: simplicial map f*: L(Y) to L(Z)},
$f^*(\Delta)$ is a $n-4$ simplex in $L(({\mathcal C}_o^n({\mathbb{CP}}^1))$.
By Theorem \ref{Thm: Non-constant holomorphic functions on Con(X)},
$L(({\mathcal C}_o^n({\mathbb{CP}}^1)))=CR({\mathbf q}_n)$; 
hence $f^*(\Delta)\in CR_{\vartriangle^{n-4}}({\mathbf q}_n)$.
Finally, Lemma \ref{Lm: S(n) orbits in SRn and CRn}$(${\bfit b}$)$
asserts that ${\mathbf S}(n)$ is transitive on
each $CR_{\vartriangle^m}({\mathbf q}_n)$.
This completes the proof of Coherence Theorem and thereby
the proof of Tame Map Theorem too. 
\index{Theorem!Tame Map Theorem!proof of\hfill}
\end{proof}




\section{Linked Map Theorem}
\label{Sec: Linked Map Theorem}
\index{Theorem!Linked Map Theorem\hfill}

\noindent The main aim of this section is to prove
Linked Map Theorem \ref{LinkMapThm}. In fact, there are
several similar results, which we also intend to discuss.
In all cases the key points of the proofs are the 
the ``explicit" description of the holomorphic
universal covering $\tau\colon{\mathbf T}(0,m+3)
\to{\mathcal C}_o^m(
                    {\mathbb C}\setminus\{0,1\})$
presented in Section \ref{Ss: Holomorphic universal covering}  
and the following remarkable property of the universal
Teichm{\" u}ller families established first by
J. H. Hubbard \cite{Hub72,Hub76} 
for compact Riemann surfaces and then generalized to all Riemann
surfaces of finite type by C. J. Earle and I. Kra
\cite{EarKra74,EarKra76}:

\begin{Theorem}[{\caps Hubbard-Earl-Kra Theorem}]
\label{Thm: TheoremHEK}
\index{Theorem!Hubbard-Earl-Kra Theorem\hfill}
The universal punctured Teichm{\" u}ller family 
${\mathcal T}(g,m):= \ \, \rho\colon
{\mathbf{V'}}(g,m)\to{\mathbf T}(g,m)$
has no holomorphic sections ${\mathbf T}(g,m)\to{\mathbf{V'}}(g,m)$
whenever $\dim_{\mathbb C}{\mathbf T}(g,m)>1$.
\hfill $\square$
\end{Theorem}

\noindent We start with the following definition
similar to Definition \ref{Def: disjoint and linked maps of Cn}.

\begin{Definition}
\label{Def: disjoint and linked maps of Con and Con(C**)}
A map $h\colon{\mathcal C}_o^n(X)\to{\mathcal C}^k(X)$
is said to be {\em linked} if there is a point
$q=(q_1,...,q_n)\in{\mathcal C}_o^n(X)$ such that
$\{q_1,...,q_n\}\cap h(q)\ne\varnothing$,
and {\em disjoint} otherwise.
A map $h\colon{\mathcal C}_o^m(
                               {\mathbb C}\setminus\{0,1\})
\to {\mathcal C}^k(X)$ is said to be
{\em linked} if there is $z=(z_1,...,z_m)\in
{\mathcal C}_o^m(
                 {\mathbb C}\setminus\{0,1\})$ such that
$\{z_1,...,z_m,0,1,\infty\}\cap h(z)\ne\varnothing$,
and {\em disjoint} otherwise.
\hfill $\bigcirc$
\end{Definition}

\subsection{Linked Map Theorem for
${\mathcal C}_o^m({\mathbb C}\setminus\{0,1\})$}
\label{Ss: Linked Map Theorem for Com(C**)}
Among all the forms of Linked Map Theorem that we were
able to prove there is the strongest one, which implies all others as simple
consequences; it will be proved in this section.

\begin{Theorem}[{\caps Linked Map Theorem for 
${\mathcal C}_o^m({\mathbb C}\setminus\{0,1\})$}]
\label{LinkMapThmCom(C**)}
\index{Theorem!Linked Map Theorem!for
${\mathcal C}_o^m({\mathbb C}\setminus\{0,1\})$\hfill}
For $m>1$ and $k\ge 1$ every holomorphic map 
$h\colon{\mathcal C}_o^m({\mathbb C}\setminus\{0,1\})
\to{\mathcal C}^k(X)$ is linked.
\hfill $\bigcirc$
\end{Theorem}

\noindent Our proof of this theorem is based on 
Lemma \ref{Lm: Eqivalence Lemma}, which directly involves
the universal punctured Teichm{\" u}ller family 
${\mathbf{V'}}(0,m+3)\to {\mathbf T}(0,m+3)$,
and on the analytic fact stated in Theorem 
\ref{Thm: absence of certain unbranched coverings over Com(C**))}
below. The latter theorem, in turn, is an easy consequence of
Hubbard-Earl-Kra Theorem \ref{Thm: TheoremHEK} formulated above.
\vskip0.2cm

\noindent Define 
\begin{equation}\label{eq: E(C**)}
\aligned
&{\mathcal E}_o^{m+1}({\mathbb C}\setminus\{0,1\})
=\{(\zeta,z)\in{\mathbb C}\times
{\mathcal C}_o^m({\mathbb C}\setminus\{0,1\})\,|\
\zeta\ne z_1,...,z_m,0,1\}\,,\\
&\Pi:= \ \, \pi\colon{\mathcal E}_o^{m+1}({\mathbb C}^{\bigcirc
                                              \hskip-6pt 2}\,)
\ni (\zeta,z)\mapsto\pi(\zeta,z)\Def z
\in{\mathcal C}_o^m({\mathbb C}\setminus\{0,1\})\,.
\endaligned
\end{equation}
Then ${\mathcal E}_o^{m+1}({\mathbb C}\setminus\{0,1\})$
is a complex manifold and $\Pi$ is a smooth locally trivial fiber
bundle with the fiber diffeomorphic to
${\mathbb C}\setminus \{m+2 \ \text{points}\}$.
Of course, the projection $\pi$ is holomorphic and $\Pi$ 
may be also regarded as a holomorphic family of curves of
type $(0,m+3)$. Set
\begin{equation}\label{eq: total space of pullback of tau}
\tau^*({\mathcal E}_o^{m+1}({\mathbb C}\setminus\{0,1\}))
=\{({\mathbf t},(\zeta,z))\in{\mathbf T}(0,m+3)
\times{\mathcal E}_o^{m+1}({\mathbb C}\setminus\{0,1\})
\,|\ \tau({\mathbf t})=z\}\,;
\end{equation}
then the fibration 
\begin{equation}\label{eq: pullback of tau}
\tau^*(\Pi):= \ \, \tau^*(\pi)
\colon\tau^*({\mathcal E}_o^{m+1}({\mathbb C}\setminus\{0,1\}))
\ni ({\mathbf t},(\zeta,z))\mapsto {\mathbf t}\in {\mathbf T}(0,m+3)
\end{equation}
is the pullback of the fibration $\Pi$ along the map  
$\tau\colon{\mathbf T}(0,m+3)
\to{\mathcal C}_o^m({\mathbb C}\setminus\{0,1\})$ and we have the
commutative diagram
\begin{equation}\label{CD: change of base tau}
\CD
\tau^*({\mathcal E}_o^{m+1}({\mathbb C}\setminus\{0,1\}))
          @ >>> {\mathcal E}_o^{m+1}({\mathbb C}\setminus\{0,1\})\\
@V\tau^*(\pi) VV \!\!@VV{\pi}V \\ 
{\mathbf T}(0,m+3) 
   @ >> \tau > {\mathcal C}_o^m({\mathbb C}\setminus\{0,1\})\,,
\endCD
\end{equation}
where the top horizontal map is just the restriction of the projection
of the product ${\mathbf T}(0,m+3)
\times{\mathcal E}_o^{m+1}({\mathbb C}\setminus\{0,1\})$
to the second factor.
\vskip0.2cm

\begin{Lemma}\label{Lm: Eqivalence Lemma}
The family 
$\tau^*(\Pi):= \ \, \tau^*(\pi)
\colon\tau^*({\mathcal E}_o^{m+1}({\mathbb C}\setminus\{0,1\}))
\to{\mathbf T}(0,m+3)$ is isomorphic to the universal
Teichm{\" u}ller family
$\rho\colon{\mathbf{V'}}(0,m+3)\to{\mathbf T}(0,m+3)$,
that is, there is a commutative diagram
\begin{equation}\label{CD: diagram over T(0,m+3) to Com(C**)}
\CD
{\mathbf{V'}}(0,m+3)
@> \lambda >> \tau^*({\mathcal E}_o^{m+1}({\mathbb C}\setminus\{0,1\}))
          @> >> {\mathcal E}_o^{m+1}({\mathbb C}\setminus\{0,1\}) \\
@.@/SE//\rho/ @VV\tau^*(\pi) V     @VV\pi V  \\
@.{\mathbf T}(0,m+3) 
          @> \tau >> {\mathcal C}_o^m({\mathbb C}\setminus\{0,1\})\,,\\
\endCD
\end{equation}
where $\lambda$ is a biholomorphic fiber isomorphism.
\end{Lemma}

\begin{proof} 
Keeping in mind that a point $(\zeta,{\mathbf t})\in
{\mathbb C}\times{\mathbf T}(0,m+3)$ belongs to
${\mathbf{V'}}(0,m+3)$ {\em iff} $\zeta
\ne s_1({\mathbf t}),...,s_m({\mathbf t}),0,1$, while 
a point $({\mathbf t},\zeta,z)
\in{\mathbf T}(0,m+3)\times{\mathbb C}
\times{\mathcal C}_o^m({\mathbb C}\setminus\{0,1\})$
belongs to
$\tau^*({\mathcal E}_o^{m+1}({\mathbb C}\setminus\{0,1\}))$
{\em iff} $\rho({\mathbf t})=z
=(z_1,...,z_m)\in{\mathcal C}_o^m({\mathbb C}\setminus\{0,1\})$
and
$\zeta\ne z_1,...,z_m,0,1$, we define $\lambda$ by
\begin{equation}\label{eq: lambda}
\aligned
\hskip-5pt\lambda\colon{\mathbf{V'}}(0,m+3)\ni {\mathbf v}
&=(\zeta,\rho({\mathbf v}))
\mapsto({\mathbf t},\zeta, s_1({\mathbf t}),...,s_m({\mathbf t}))\\
&\Def(\rho({\mathbf v}), \zeta, s_1(\rho({\mathbf v})),...,
s_m(\rho({\mathbf v})))
\in\tau^*({\mathcal E}_o^{m+1}({\mathbb C}\setminus\{0,1\}))\,.
\endaligned
\end{equation}
The map $\lambda^{-1}$ is defined by
\begin{equation}\label{eq: inverse of lambda}
\lambda^{-1}
\colon\tau^*({\mathcal E}_o^{m+1}({\mathbb C}\setminus\{0,1\}))
\ni ({\mathbf t},(\zeta,z))\mapsto (\zeta,{\mathbf t})
\in{\mathbf{V'}}(0,m+3)\,.
\end{equation}
These $\lambda$ and $\lambda^{-1}$ are
holomorphic fiber mappings inverse to each other.
\end{proof}

\begin{Theorem}
\label{Thm: absence of certain unbranched coverings over Com(C**))}
For $m>1$ the restriction $\pi|_A\colon A
\to{\mathcal C}_o^m({\mathbb C}\setminus\{0,1\})$
of the projection 
$\pi\colon{\mathcal E}_o^{m+1}({\mathbb C}\setminus\{0,1\})
\to{\mathcal C}_o^m({\mathbb C}\setminus\{0,1\})$
to an analytic subset
$A\subset{\mathcal E}_o^{m+1}({\mathbb C}\setminus\{0,1\})$
cannot be an unbranched covering over the whole manifold
${\mathcal C}_o^m({\mathbb C}\setminus\{0,1\})$.
\end{Theorem}

\begin{proof}
Suppose, on the contrary, that an analytic subset 
$A\subset{\mathcal E}_o^{m+1}({\mathbb C}\setminus\{0,1\})$
covers the whole base
${\mathcal C}_o^m({\mathbb C}\setminus\{0,1\})$.
The pullback 
$\tau^*(\pi|_A)\colon\tau^*(A)\to{\mathbf T}(0,m+3)$
of the covering $\pi|_A\colon A
\to{\mathcal C}_o^m({\mathbb C}\setminus\{0,1\})$
along the mapping $\tau$ is a holomorphic unbranched
covering as well. It must be trivial since the Teichm{\" u}ller
space ${\mathbf T}(0,m+3)$ is contractible;
hence there exists a holomorphic section
$s\colon{\mathbf T}(0,m+3)\to\tau^*(A)\subset
\tau^*({\mathcal E}_o^{m+1}({\mathbb C}\setminus\{0,1\}))$.
As $\lambda$ in diagram (\ref{CD: diagram over T(0,m+3) to Com(C**)})
is an isomorphism, the latter section provides a holomorphic section
$S\colon{\mathbf T}(0,m+3)\to {\mathbf{V'}}(0,m+3)$
of the universal Teichm{\" u}ller family 
${\mathcal T}(0,m+3):= \ \, \rho\colon{\mathbf{V'}}(0,m+3)
\to {\mathbf T}(0,m+3)$. Since $\dim_{\mathbb C}{\mathbf T}(0,m+3)=m>1$,
this contradicts Hubbard-Earl-Kra Theorem \ref{Thm: TheoremHEK}.
\end{proof}

\noindent We are now in a position to prove
Linked Map Theorem for
${\mathcal C}_o^m({\mathbb C}\setminus\{0,1\})$
(Theorem \ref{LinkMapThmCom(C**)}).
\vskip0.2cm

\begin{proof}
\index{Theorem!Linked Map Theorem!for ${\mathcal C}_o^m({\mathbb C}\setminus\{0,1\})$\hfill}
Suppose, on the contrary, that for some $m>1$ and $k\ge 1$
there is a disjoint holomorphic map 
$h\colon{\mathcal C}_o^m({\mathbb C}\setminus\{0,1\})
\to{\mathcal C}^k(X)$. 
\vskip0.2cm

\noindent Notice first that for each point
$z=(z_1,...,z_m)
\in{\mathcal C}_o^m({\mathbb C}\setminus\{0,1\})$
its image $h(z)\subset X$  is in fact contained in ${\mathbb C}$,
i. e., $\infty\notin h(z)$.  
In the case $X={\mathbb C}$ this is self evident;
in the case $X={\mathbb{CP}}^1$ this follows from the condition
$\{z_1,...,z_m,0,1,\infty\}\cap h(z)=\varnothing$
(see Definition \ref{Def: disjoint and linked maps of Con and Con(C**)}).
Let $P_k(t,h(z))$ be the polynomial in $t$
with the leading coefficient $1$ and with the set of the roots
$h(z)$; it may clearly be regarded as an analytic function on
${\mathbb C}
\times{\mathcal C}_o^m({\mathbb C}\setminus\{0,1\})$.
Let $A\subset{\mathbb C}
\times{\mathcal C}_o^m({\mathbb C}\setminus\{0,1\})$
be the analytic subset defined by
$$
\aligned
A&=\{(t,z)\in{\mathbb C}\times
{\mathcal C}_o^m({\mathbb C}\setminus\{0,1\})\,|\
t\in h(z)\}\\
&=\{(t,z)\in{\mathbb C}
\times{\mathcal C}_o^m({\mathbb C}\setminus\{0,1\})\,|\ 
P_k(t,h(z))=0\}\,.
\endaligned
$$
As $\{z_1,...,z_m,0,1,\infty\}\cap h(z)=\varnothing$, 
the image $h(z)$ of any point 
$z\in{\mathcal C}_o^m({\mathbb C}\setminus\{0,1\})$
contains no punctures of the fiber 
$\pi^{-1}(z)$ of the family 
$\pi\colon{\mathcal E}_o^{m+1}({\mathbb C}\setminus\{0,1\})
\to{\mathcal C}_o^m({\mathbb C}\setminus\{0,1\})$
defined in (\ref{eq: E(C**)}). That is, all $k$ distinct
roots $t_1,...,t_k\in h(z)$ of the polynomial $P_k(t,h(z))$ are in
$\pi^{-1}(z)$. Hence 
$A\subseteq{\mathcal E}_o^{m+1}({\mathbb C}\setminus\{0,1\})$
and the restriction $\pi|_A\colon A
\to{\mathcal C}_o^m({\mathbb C}\setminus\{0,1\})$
of the projection 
$\pi\colon{\mathcal E}_o^{m+1}({\mathbb C}\setminus\{0,1\})
\to{\mathcal C}_o^m({\mathbb C}\setminus\{0,1\})$
to $A$ is a holomorphic unbranched covering over the whole 
base ${\mathcal C}_o^m({\mathbb C}\setminus\{0,1\})$,
which contradicts Theorem
\ref{Thm: absence of certain unbranched coverings over Com(C**))}.
\end{proof}

\subsection{Linked Map Theorem for ${\mathcal C}_o^n(X)$}
\label{Ss: Linked Map Theorem for Con(X)}
\index{Theorem!Linked Map Theorem!for ${\mathcal C}_o^n(X)$\hfill}
If for a map $f\colon{\mathcal C}_o^n(X)\to{\mathcal C}^k(X)$
and a subset $M\subseteq{\mathcal C}_o^n(X)$
there is $q=(q_1,...,q_n)\in M$ such that
$\{q_1,...,q_n\}\cap f(q)\ne\varnothing$ 
then $f$ is certainly linked.
Hence the following form of Linked Map Theorem,
which is certainly stronger that the original one 
formulated in Section \ref{Ss: Main results} (Theorem \ref{LinkMapThm}),
follows immediately from Linked Map Theorem for
${\mathcal C}_o^m({\mathbb C}\setminus\{0,1\})$
(Theorem \ref{LinkMapThmCom(C**)}) proved above.

\begin{Theorem}[{Linked Map Theorem for ${\mathcal C}_o^n(X)$}]
\label{LinkMapThmCon}
\index{Theorem!Linked Map Theorem!for ${\mathcal C}_o^n(X)$\hfill}
For $n>t(X)+1$ and $k\ge 1$ every holomorphic map 
$h\colon{\mathcal C}_o^n(X)\to{\mathcal C}^k(X)$
is linked.
\hfill $\square$
\end{Theorem}


\begin{Remark}\label{Rmk: linked maps}
Let $n>t(X)$ and $k\ge 1$. The above theorem implies that 
\vskip0.2cm

{\bfit a}$)$ {\sl any holomorphic
map $F\colon{\mathcal C}^n(X)\to{\mathcal C}^k(X)$
is linked, i. e., there is $Q\in{\mathcal C}^n(X)$ such that
$Q\cap F(Q)\ne\varnothing$;}
\vskip0.2cm

{\bfit b}$)$ {\sl any holomorphic
map $f\colon{\mathcal C}_o^n(X)\to{\mathcal C}_o^k(X)$
is linked, i. e., there is $q\in{\mathcal C}_o^n(X)$ such that
$p(q)\cap p(f(q))\ne\varnothing$;}
\vskip0.2cm

{\bfit c}$)$ {\sl any holomorphic
map $g\colon{\mathcal C}^n(X)\to{\mathcal C}_o^k(X)$
is linked, i. e., there is $Q\in{\mathcal C}^n(X)$ such that
$Q\cap p(g(Q))\ne\varnothing$.}
\vskip0.2cm

\noindent To prove these statements,
one just take $h=F\circ p$ in case $(a)$, 
$h=p\circ f$ in case $(b)$ and $h=p\circ g\circ p$
in case $(c)$ and applies Theorem \ref{LinkMapThmCon}
to the map $h$.
\hfill $\bigcirc$
\end{Remark}

\section{Morphisms ${\mathcal C}^n(X)\to{\mathcal C}^k(X)$ and
$\Aut X$ orbits. Cyclic and orbit-like maps}
\label{Sec: Morphisms Cn(X) to Ck(X) and Aut X orbits.
Cyclic and orbit-like maps}

\noindent In this section, using the direct decomposition
(\ref{eq: j(X,n): Con(X) to AutX x D(n-t(X))}), we 
show that any morphism $F\colon{\mathcal C}^n(X)\to{\mathcal C}^k(X)$
respects $\Aut X$ orbits, that is, $F((\Aut X)Q)\subseteq(\Aut X)F(Q)$
for all $Q\in{\mathcal C}^n(X)$. Then we try to 
clarify, as far as possible, the structure of cyclic morphisms.

\subsection{Liouville and Picard spaces}
\label{Ss: Liouville and Picard spaces}
Here we exhibit a certain class of holomorphic maps
$Z\to{\mathcal C}^k(X)$ whose images must always be contained
in an $\Aut X$ orbit. We will use these results  
in Section \ref{Ss: morphisms and Aut X orbits} below.
All complex spaces under consideration are assumed to be reduced.

\begin{Definition}\label{Def: Liouville and Picard spaces}
\index{Liouville and Picard spaces\hfill}
\index{Ultra-Picard spaces\hfill}
A connected complex space $Z$ is said to be {\em Liouville}
if it carries no non-constant bounded holomorphic functions. 
We say that $Z$ is a {\em Picard space} if it carries no non-constant
holomorphic functions omitting the values $0$ and $1$;
$Z$ is called {\em ultra-Picard} whenever 
every its connected unbranched finite cover is Picard.
\hfill $\bigcirc$
\end{Definition}

\begin{Lemma}\label{Lm: finite cover of Liouville space}
\index{Finite covers of Liouville spaces\hfill}
Let $Y$ be a connected complex space and $p\colon Y\to Z$
be a holomorphic finite unbranched covering map onto a Liouville space $Z$.
Then $Y$ is Liouville.\footnote{Of course, this is a very well-known
elementary result. See \cite{Lin88,LinZai98} and the references therein
for more interesting results about the Liouville property.}
\end{Lemma}

\begin{proof}
Let $h$ be a bounded holomorphic function
on $Y$. The coefficients of the polynomial
$P(t,z)=\prod_{y\in p^{-1}(z)} (t-h(y))$
are bounded holomorphic functions on $Z$ and hence constant.
Since $Y$ is connected this implies $h=\const$.
\end{proof}

\begin{Proposition}\label{Prp: holomorphic maps to Cm(X)}
\index{Holomorphic maps!of ultra-Picard space to ${\mathcal C}^m(X)$\hfill}
\index{Holomorphic maps!Solvable!of Liouville space to ${\mathcal C}^m(X)$\hfill}
Let $Z$ be a connected complex space and
$F\colon Z\to{\mathcal C}^m(X)$ be a holomorphic map.
Suppose that either $(i)$ $Z$ is an ultra-Picard space
or $(ii)$ $Z$ is a Liouville space and
the homomorphism $F_*\colon\pi_1(Z)\to\pi_1({\mathcal C}^m(X))$
is solvable. Then there is a point $Q^*\in{\mathcal C}^m(X)$ such that
$F(Z)\subseteq (\Aut X)Q^*$.
\end{Proposition}

\begin{proof}
Let $F^*(p)\colon F^*({\mathcal C}_o^m(X))\to Z$
be the pullback of the covering
$p\colon{\mathcal C}_o^m(X)\to{\mathcal C}^m(X)$ along the mapping
$F$ and let $Y$ be a connected component of $F^*({\mathcal C}_o^m(X))$.
Then we have the commutative diagram
\begin{equation}\label{CD: induced cover}
\CD
Y @>{f}>> {\mathcal C}_o^m(X) @>{j_{X,m}}>{\cong}>
{(\Aut X)\times{\mathcal D}^{m-t(X)}(X)}
@>{\xi}>> {\mathcal D}^{m-t(X)}(X) \\
@V{\nu}VV @VV{p}V\\
Z @>>{F}> {\mathcal C}^m(X)\,,
\endCD
\end{equation}
where $\xi\colon (\Aut X)\times{\mathcal D}^{m-t(X)}(X)
\to{\mathcal D}^{m-t(X)}(X)$ is the projection onto the second factor
and $\nu$ is the restriction of $F^*(p)$ to $Y$ 
(see (\ref{eq: j(X,n): Con(X) to AutX x D(n-t(X))}) and (\ref{eq: j(X,n)})).
The coordinate functions $z_1,...,z_{m-t(X)}$ on
${\mathcal D}^{m-t(X)}(X)
={\mathcal C}_o^{m-t(X)}({\mathbb C}\setminus\{0,1\})$
are holomorphic and do not assume the values $0$ and $1$; 
hence for each $i$, $1\le i\le m-t(X)$, the function
$h_i=z_i\circ\xi\circ j_{X,m}\circ f$ is holomorphic on $Y$
and does not take the values $0$ and $1$.
\vskip0.2cm

\noindent In case $(i)$, $Y$ is a Picard space. Thus, 
$\xi\circ j_{X,m}\circ f=\const$ and $f(Y)$ is contained in the $\Aut X$
orbit of some point $q^*\in{\mathcal C}_o^m(X)$, which implies
$F(Z)\subseteq (\Aut X)p(q^*)$.
\vskip0.2cm

\noindent In case $(ii)$ for each $i=1,...,m-t(X)$
the image $G_i=(h_i)_*(\pi_1(Y))$ of the homomorphism 
$(h_i)_*\colon\pi_1(Y)\to\pi_1({\mathbb C}\setminus\{0,1\})$
is a solvable subgroup of 
$\pi_1({\mathbb C}\setminus\{0,1\})\cong{\mathbb F}_2$.
Thus, either $G_i$ is trivial or $G_i\cong{\mathbb Z}$.
Let $\nu_i\colon\Gamma_i\to{\mathbb C}\setminus\{0,1\}$
be the covering corresponding to the subgroup
$G_i$ 
so that either $\Gamma_i\cong{\mathbb D}
=\{\zeta\in{\mathbb C}\,|\ |\zeta| < 1\}$
or $\Gamma_i\cong K_{r_i} = \left\{\zeta\in{\mathbb C}\,|\
r_i^{-1} < |\zeta| < r_i\right\}$ for some $r_i>1$.
By the covering map theorem,
there is a holomorphic map $\widetilde{h_i}\colon Y\to\Gamma_i$
that fits into the commutative diagram
$$
\CD
@. \Gamma_i\\
@. @/NE/{\,\widetilde{h_i}}// @VV{\nu_i}V \\
Y @>>{h_i}>{\mathbb C}\setminus\{0,1\}\,.\\
\endCD
$$
By Lemma \ref{Lm: finite cover of Liouville space},
the space $Y$ is Liouville; since both ${\mathbb D}$ and 
$K_{r_i}$ are bounded domains, $\widetilde{h_i}=\const$, which
implies $h_i=\const$. Thus, $\xi\circ j_{X,m}\circ f=\const$,
$f(Y)$ is contained in the $\Aut X$
orbit of some point $q^*\in{\mathcal C}_o^m(X)$ and
$F(Z)\subseteq (\Aut X)p(q^*)$.
\end{proof}

\noindent Any connected complex Lie group ${\mathcal L}$
satisfies both $(i)$ and $(ii)$ since the 
union of the images of all holomorphic homomorphisms
${\mathbb C}\to{\mathcal L}$
contains an open neighborhood of the unity and, moreover,
$\pi_1({\mathcal L})$ is abelian. Thus, we have:

\begin{Corollary}\label{Crl: Lie group to Com(X)}
Let ${\mathcal L}$ be a connected complex Lie group. Then for every
holomorphic map $\Phi\colon{\mathcal L}\to {\mathcal C}^m(X)$
there is a point $Q^*\in{\mathcal C}^m(X)$ such that
$\Phi({\mathcal L})\subseteq (\Aut X)Q^*$.
\hfill $\square$
\end{Corollary}

\subsection{Morphisms ${\mathcal C}^n(X)\to{\mathcal C}^k(X)$ and
$\Aut X$ orbits}\label{Ss: morphisms and Aut X orbits}
Here we prove the statement about morphisms of configuration
spaces mentioned in Remark \ref{Rmk: tame maps keep each AutX orbit}.

\begin{Proposition}\label{Prp: morphisms respect orbits}
Every morphism $F\colon{\mathcal C}^n(X)\to{\mathcal C}^k(X)$
respects $\Aut X$ orbits, that is, $F((\Aut X)Q)\subseteq (\Aut X)F(Q)$
for all $Q\in{\mathcal C}^n(X)$.
\end{Proposition}

\begin{proof}
Let $Q\in{\mathcal C}^n(X)$ and let $\Psi\colon\Aut X\to{\mathcal C}^n(X)$
be the holomorphic map defined by $\Psi(A)=AQ$ for each $A\in\Aut X$.
Set $\Phi=F\circ\Psi$; clearly $\Phi(\id_X)=F(Q)$.
By Corollary \ref{Crl: Lie group to Com(X)}, the image
$F((\Aut X)Q)=F(\Psi(\Aut X))=\Phi(\Aut X)$ is contained in
an $\Aut X$ orbit, which certainly is the orbit of the point $F(Q)$.
\end{proof}

\begin{Remark}\label{Rmk: Dm(X) is Kobayashi hyperbolic}
Being a domain in the Kobayashi hyperbolic manifold
$({\mathbb C}\setminus\{0,1\})^m$, the space ${\mathcal D}^m(X)$
\end{Remark}

\subsection{Stabilizers of points in ${\mathcal C}^k(X)$}\label{Ss: stabilizer}

\noindent In what follows we need to use the 
classification of finite subgroups $G\subset{\mathbf{SO}}(3)$;
this is a very classical matter, which goes back at least 
to F. Klein's ``Vorlesungen {\" u}ber das Ikosaeder" \cite{Kle}.
The following list exhibits all such
subgroups\footnote{We take no care about possible
repetitions in this list. See, for instance, \cite{Fo51},
Chapter VI, Sections 51-57, for more details.}:

\begin{itemize}

\item[{\bfit a}$)$] cyclic groups ${\mathbb Z}/m{\mathbb Z}$, $m\ge 2$; 
\vskip0.1cm

\item[{\bfit b}$)$] $({\mathbb Z}/2{\mathbb Z})
\times{\mathbb Z}/(2{\mathbb Z})$;
\vskip0.1cm

\item[{\bfit c}$)$] dihedral groups 
$({\mathbb Z}/m{\mathbb Z})\leftthreetimes ({\mathbb Z}/2{\mathbb Z})$,
$m>2$, having two generators $a,b$ and the defining system of relations
$a^m=1$, $b^2=1$, $bab^{-1}=a^{-1}$;
\vskip0.1cm

\item[{\bfit d}$)$] the rotation groups of Platonic solids, that is,
the symmetric group ${\mathbf S}(4)$ (rotation groups of the cube and
the octahedron), the alternating groups ${\mathbf A}(4)$
(the rotation group of the tetrahedron) and ${\mathbf A}(5)$
(rotation groups of the icosahedron and the dodecahedron), and
their subgroups. 
\end{itemize}
\vskip0.2cm

\noindent Notice that all groups of this list except of ${\mathbf A}(5)$
are solvable.

\begin{Proposition}\label{Prp: stabilizers}
For $k\ge t(X)$ the stabilizer $St_{Q^*}=\{A\in\Aut X\,|\ AQ^*=Q^*\}$
of any point $Q^*\in{\mathcal C}^k(X)$ is a finite group which is either
solvable or isomorphic to the alternating group ${\mathbf A}(5)$.
\end{Proposition}

\begin{proof}
Let $Q^*=\{q_1^*,...,q_k^*\}$. Any element $A\in St_{Q^*}$
is either an affine or a M{\"o}bius transformation permuting the points 
$q_1^*,...,q_k^*$; \ $A$ is uniquely determined by its values
$Aq_i^*$, $1\le i\le t(X)$, which shows that $St_{Q^*}$ is finite.

If $X={\mathbb C}$ then $St_{Q^*}$ is a finite
subgroup of $\Aff{\mathbb C}$. Every such subgroup is conjugate
to a finite subgroup of the maximal compact subgroup
${\mathbb T}\subset\Aff{\mathbb C}$, which consists of all rotations
$z\mapsto e^{it}z$, $t\in{\mathbb R}$. Every finite subgroup of the
latter group is cyclic.

If $X={\mathbb{CP}}^1$ then $St_{Q^*}$ is a finite
subgroup of ${\mathbf{PSL}}(2,{\mathbb C})$.
Every such subgroup is conjugate
to a finite subgroup of the maximal compact subgroup
$K\subset{\mathbf{PSL}}(2,{\mathbb C})$, which is isomorphic
to ${\mathbf{SO}}(3)$ and consists of all M{\"o}bius transformations
of the form
$$
z\mapsto\frac{az - \bar c}{cz + \bar a}\,,\quad |a|^2 + |c|^2 = 1\,.
$$
Thus, in this case $St_{Q^*}$ is isomorphic to one of the groups
exhibited in the above list, which proves the proposition.
\end{proof} 

\subsection{Cyclic Map Theorem}\label{Ss: Cyclic Map Theorem}
Recall that a continuous map
$F\colon{\mathcal C}^n(X)\to{\mathcal C}^k(X)$ is said to be 
orbit-like if its image $F({\mathcal C}^n(X))$ is contained
in some orbit of $\Aut X$ action in ${\mathcal C}^k(X)$,
that is, there is a point $Q^*\in{\mathcal C}^k(X)$
such that $F({\mathcal C}^n(X))\subseteq(\Aut X)Q^*$.
We have already mentioned that the stabilizer $St_{Q^*}$
of a generic point $Q^*\in{\mathcal C}^k(X)$
is trivial; thereby a morphism $F$ whose image
is contained in the orbit of such a point $Q^*$ is degenerate tame,
that is, there is a morphism $T\colon{\mathcal C}^n(X)\to\Aut X$
such that $F(Q)=T(Q)Q^*$ for all $Q\in{\mathcal C}^n(X)$
(see Definition \ref{Def: tame, degenerate tame and orbit-like maps}
and Remark \ref{Rmk: tame maps keep each AutX orbit}).
\vskip0.2cm

\begin{Theorem}[{\caps Cyclic Map Theorem}]
\label{Thm: Cyclic Map Theorem}
Let $n>4$ and $k\ge t(X)$. Then 

{\bfit a}$)$ every cyclic morphism 
$F\colon{\mathcal C}^n(X)\to{\mathcal C}^k(X)$ is orbit-like;
\vskip0.1cm

{\bfit b}$)$ every continuous orbit-like map
$F\colon{\mathcal C}^n(X)\to{\mathcal C}^k(X)$ is cyclic.
\end{Theorem}

\noindent Statement $(${\bfit a}$)$ follows from Proposition
\ref{Prp: holomorphic maps to Cm(X)}$(ii)$ since 
${\mathcal C}^n(X)$ is a connected algebraic manifold
and therefore Liouville, whereas $F$ is supposed to be cyclic.
\vskip0.2cm

\noindent Statement $(${\bfit b}$)$ is an immediate consequence
of the following lemma.

\begin{Lemma}\label{Lm: homomorphisms Bn(X) to pi1((Aut X)Q*)}
Let $n>4$ and $k\ge t(X)$. Then for each point $Q^*\in{\mathcal C}^k(X)$
every homomorphism $\varphi\colon B_n(X)\to\pi_1((\Aut X)Q^*)$
is cyclic.
\end{Lemma}

\begin{proof}
If the stabilizer $St_{Q^*}$ is trivial then 
$(\Aut X)Q^*\cong\Aut X$, $\pi_1((\Aut X)Q^*)\cong\pi_1(\Aut X)$
is cyclic and $\varphi$ is cyclic. 
\vskip0.2cm

\noindent Suppose that $St_{Q^*}$ is non-trivial; then the natural
map $h\colon\Aut X\to (\Aut X)Q^*$ is a finite
regular covering with the fiber $St_{Q^*}$ and we have
the commutative diagram
\begin{equation}\label{CD: diagram for phi and covering}
\CD
@. @.  {B_n(X)} @. @. \\
@. @.  @V{\varphi}VV @/SE//{\,\delta\circ\varphi}/ @. \\ 
1@>>>\pi_1(\Aut X)@>>{h_*}>{\pi_1((\Aut X)Q^*)}@>>{\delta}>
{St_{Q^*}}@>>>1\,,\\
\endCD
\end{equation}
where the horizontal line is the exact sequence
corresponding to the covering $h$.
\vskip0.2cm

\noindent Consider the restriction of the homomorphisms $\varphi$
and $\delta\circ\varphi$ to the commutator subgroup $B_n'(X)$,
which is perfect since $n>4$. 
\vskip0.2cm

\noindent By Proposition \ref{Prp: stabilizers}, 
$St_{Q^*}$ is either solvable or isomorphic to ${\mathbf A}(5)$.
\vskip0.2cm

\noindent If $St_{Q^*}$ is solvable then the restriction
of $\delta\circ\varphi$ to $B_n'(X)$ is trivial and hence 
$$
\varphi(B_n'(X))\subseteq\Ker\delta=h_*(\pi_1(\Aut X))
\cong\pi_1(\Aut X)\,.
$$
The group $\pi_1(\Aut X)$ is cyclic,
which shows that $\varphi(B_n'(X))$ is trivial and $\varphi$ is cyclic.
\vskip0.2cm

\noindent Thus, we are left with the case 
$St_{Q^*}\cong{\mathbf A}(5)\subset{\mathbf S}(5)$.
If $n>5$ then, by Theorem
\ref{Thm: Bn(X) to S(k) for n>k}, the homomorphism $\delta\circ\varphi$
is cyclic, which implies that $(\delta\circ\varphi)(B_n'(X))=\{1\}$ 
and, as above, $\varphi(B_n'(X))\subseteq\Ker\delta=\pi_1(\Aut X)$,
$\varphi(B_n'(X))=\{1\}$ and $\varphi$ is cyclic. 
Finally, if $n=5$ then the homomorphism $\delta\circ\varphi$
maps $B_5(X)$ into the proper
subgroup ${\mathbf A}(5)\subset{\mathbf S}(5)$ and, 
by Surjectivity Theorem \ref{Thm: Surjectivity Theorem}, it is cyclic again;
as above, we conclude that $\varphi$ is cyclic.
This proves the lemma and completes the proof of Cyclic Map Theorem.
\end{proof}

\subsection{Holomorphic maps to the balanced configuration space}
\label{Ss: Holomorphic maps to the balanced configuration space}
We start with the following

\begin{Definition}\label{Def: balanced configuration space of C}
\index{Balanced configuration spaces\hfill}
\index{${\mathcal C}^{m-1}_b({\mathbb C})$\hfill}
The space of {\em balanced} configurations,
or {\em balanced configuration space} ${\mathcal C}^{m-1}_b({\mathbb C})$
is the $(m-1)$-dimensional subspace of ${\mathcal C}^m({\mathbb C})$
consisting of all sets $Q=\{q_1,...,q_n\}\in{\mathcal C}^m({\mathbb C})$
having its mass centers at the origin: 
$$
{\mathcal C}^{m-1}_b({\mathbb C})
=\{Q=\{q_1,...,q_m\}\in{\mathcal C}^m({\mathbb C})\,|\ q_1+...+q_m=0\}\,.
$$
We set also
$$
\aligned
\hskip30pt
&{\mathcal C}^{m-1}_{ob}({\mathbb C})
=\{q=(q_1,...,q_m)\in{\mathcal C}^m_o({\mathbb C})\,|\ q_1+...+q_m=0\}\,,\\
&{\mathcal C}^{m-2}_{ob,1}({\mathbb C})
=\{q=(q_1,...,q_m)\in{\mathcal C}^n_o({\mathbb C})\,|\ 
q_1+...+q_m=0,\ \ q_1-q_2=1\}\,. \hskip30pt \bigcirc
\endaligned
$$
\end{Definition}

\noindent The multiplicative group ${\mathbb C}^*$ acts
in ${\mathcal C}^{m-1}_b({\mathbb C})$ and
${\mathcal C}^{m-1}_{ob}({\mathbb C})$ via 
\begin{equation}
\label{eq: C* actions in C[m-1]b(C) and C[m-1]ob(C)}
{\mathbb C}^*\ni\zeta: \ \{q_1,...,q_m\}
\mapsto \{\zeta q_1,...,\zeta q_m\}\quad\text{and}\quad  
(q_1,...,q_m)\mapsto (\zeta q_1,...,\zeta q_m)
\end{equation}
respectively. The ${\mathbb C}^*$-orbit
of any point $q=(q_1,...,q_m)\in{\mathcal C}^{m-1}_{ob}({\mathbb C})$
intersects the subspace ${\mathcal C}^{m-2}_{ob,1}({\mathbb C})
\subset{\mathcal C}^{m-1}_{ob}({\mathbb C})$ at a single point
${\widetilde q}=(q_1/(q_1-q_2),...,q_m/(q_1-q_2))$;
this provides the biholomorphic direct decomposition
${\mathcal C}^{m-1}_{ob}({\mathbb C})
={\mathbb C}^*\times{\mathcal C}^{m-2}_{ob,1}({\mathbb C})$,
$$
{\mathcal C}^{m-1}_{ob}({\mathbb C})\ni q=(q_1,...,q_m)\mapsto
(q_1-q_2;\, {\widetilde q})
\in{\mathbb C}^*\times{\mathcal C}^{m-2}_{ob,1}({\mathbb C})\,.
$$
On the other hand, the functions $\phi_i(q)=(q_1-q_i)/(q_1-q_2)$,
$i=3,...,n$, are holomorphic on ${\mathcal C}^{m-1}_{ob}({\mathbb C})$,
separate points of ${\mathcal C}^{m-2}_{ob,1}({\mathbb C})$
and do not take the values $0$ and $1$. This implies that

\begin{itemize}

\item[$(i)$] {\sl for every
holomorphic map $\Phi\colon{\mathcal L}\to {\mathcal C}^m({\mathbb C})$
of a connected complex Lie group ${\mathcal L}$
there is a point $Q^*\in{\mathcal C}^{m-1}_b({\mathbb C})$ such that
$\Phi({\mathcal L})\subseteq {\mathbb C}^*Q^*$} and 
\vskip0.1cm

\item[$(ii)$] {\sl any holomorphic map 
$F\colon{\mathcal C}^{n-1}_b({\mathbb C})
\to{\mathcal C}^{k-1}_b({\mathbb C})$ respects
${\mathbb C}^*$-orbits, that is, $F({\mathbb C}^*Q)\subseteq
{\mathbb C}^*F(Q)$ for all $Q\in{\mathcal C}^{n-1}_b({\mathbb C})$.}
\end{itemize}
\vskip0.2cm

\noindent The proof of the following analog of 
Proposition \ref{Prp: holomorphic maps to Cm(X)} is clear:

\begin{Proposition}\label{Prp: holomorphic maps to Cmb(C)}
Let $Z$ be a connected complex
space and $F\colon Z\to{\mathcal C}^{m-1}_b({\mathbb C})$
be a holomorphic map. Suppose that either $Z$ is an ultra-Picard space
or the homomorphism $F_*\colon\pi_1(Z)
\to\pi_1({\mathcal C}^{m-1}_b({\mathbb C}))$
is solvable and $Z$ is Liouville. Then there is a point
$Q^*\in{\mathcal C}^{m-1}_b({\mathbb C})$ such that
$F(Z)\subseteq {\mathbb C}^*Q^*$.
\hfill $\square$
\end{Proposition}

\noindent In what follows, instead of ${\mathcal C}^m({\mathbb C})$
and its balanced subspace ${\mathcal C}^{m-1}_b({\mathbb C})$,
it is convenient to deal with the space of polynomials with simple roots
${\mathbf G}_m$ and its subspace 
\index{${\mathbf G}_m^0$\hfill}
\index{Balanced configuration space ${\mathbf G}_m^0$\hfill}
\index{${\mathcal C}^{m-1}_b({\mathbb C})$\hfill}
\begin{equation}\label{eq: Gm0}
{\mathbf G}_m^0=\{w=(w_1,...,w_m)\in{\mathbf G}_m\,|\ w_1=0\}
\end{equation}
that corresponds to ${\mathcal C}^{m-1}_b({\mathbb C})$ under the
identification ${\mathcal C}^m({\mathbb C})={\mathbf G}_m$ (see Remark
\ref{Rmk: non-degenerate forms}). The ${\mathbb C}^*$ action in
${\mathbf G}_m$ corresponding to the ${\mathbb C}^*$ action
(\ref{eq: C* actions in C[m-1]b(C) and C[m-1]ob(C)}) in 
${\mathcal C}^{m-1}_b({\mathbb C})$ look as follows:
\begin{equation}\label{eq: C* action in Gm0}
{\mathbb C}^*\ni\zeta\colon{\mathbf G}_m^0\ni (0,w_2,...,w_m)
\mapsto{\mathcal U}_\zeta w=(0,\zeta^2 w_2,...,\zeta^m w_m)
\in{\mathbf G}_m^0\,.
\end{equation}

\begin{Notation}\label{Not: NZ(w), Z(w), D(w)}
\index{${\mathcal D}(w)$\hfill}
For $w=(0,w_2,...,w_m)\in{\mathbf G}_m^0$
set ${\mathcal Z}(w)=\left\{j\,|\ j\ge 2, \ w_j=0\right\}$,
${\mathcal{NZ}}(w)=\left\{j\,|\ j\ge 2, \ w_j\ne 0\right\}$
and denote by ${\mathcal D}(w)$ the greatest common divisor of
all $j\in{\mathcal{NZ}}(w)$. Since ${\mathcal Z}(w)$ cannot contain
both $m-1$ and $m$, the number ${\mathcal D}(w)$ divides $m(m-1)$.
As usual, ${\mathcal O}(Z)$ denotes the algebra of all
holomorphic functions on $Z$; ${\mathcal O}^*(Z)$ is the
multiplicative group of all non-vanishing functions
$f\in{\mathcal O}(Z)$ and its subgroup
$\exp{\mathcal O}(Z)\subseteq{\mathcal O}^*(Z)$ consists
of all functions $e^f,\ f\in{\mathcal O}(Z)$.
\end{Notation}

\noindent The following evident lemma provides a complete description
of ${\mathbb C}^*$ orbits in the space of balanced configurations
${\mathcal C}^{m-1}_b({\mathbb C})={\mathbf G}_m^0$; it will be used
in next sections.

\begin{Lemma}\label{Lm: parameterization of C* orbit in Gm0}
\index{Parameterization of ${\mathbb C}^*$ orbits in ${\mathbf G}_m^0$\hfill}
Let $w^\circ=(0,w_2^\circ,...,w_m^\circ)\in{\mathbf G}_m^0$
and ${\mathcal D}={\mathcal D}(w^\circ)$. Then the stabilizer
$St_{w^\circ}\subset{\mathbb C}^*$ of $w^\circ$
is the cyclic subgroup $\{\zeta\in{\mathbb C}^*\,|\ \zeta^{\mathcal D}=1\}
\cong{\mathbb Z}/{\mathcal D}{\mathbb Z}$, and
the ${\mathbb C}^*$-orbit 
${\mathbb C}^*w^\circ\cong{\mathbb C}^*/St_{w^\circ}\cong{\mathbb C}^*$
of $w^\circ$ admits the following biholomorphic parameterization:
\begin{equation}\label{eq: parametrizarion of C* orbit in Gm0}
f\colon{\mathbb C}^*\stackrel{\cong}\to{\mathbb C}^*w^\circ\,,
\quad
f(\zeta)
=(0,\zeta^{2/{\mathcal D}}w_2^\circ,...,\zeta^{m/{\mathcal D}}w_m^\circ)
\in{\mathbb C}^*w^\circ\subset{\mathbf G}_m^0\subset{\mathbf G}_m
\end{equation}
$($notice that the set $\{j/{\mathcal D}\,|\ w_j^\circ\ne 0\}$
consists of mutually co-prime integers$)$. 
\hfill $\square$
\end{Lemma}

\subsection{Solvable holomorphic maps of Liouville spaces
to ${\mathcal C}^m({\mathbb C})={\mathbf G}_m$}
\label{Solvable holomorphic maps of Liouville spaces to Cm(C)=Gm}

\begin{Proposition}\label{Prp: ultra-Picard or Liouville space to Gm}
\index{Holomorphic maps of ultra-Picard space to ${\mathbf G}_m$\hfill}
\index{Solvable holomorphic maps of Liouville space to ${\mathbf G}_m$\hfill}
Let $Z$ be a connected complex space and $F\colon Z\to{\mathbf G}_m$
be a holomorphic map. Suppose that either $(i)$ $Z$ is an ultra-Picard space
or $(ii)$ the homomorphism $F_*\colon\pi_1(Z)\to\pi_1({\mathbf G}_m)$
is solvable and $Z$ is Liouville. 
Then: 
\vskip0.1cm

{\bfit a}$)$ There exist $w^\circ=(0,w^\circ_2,...,w^\circ_m)
\in{\mathbf G}_m^0$, $\zeta\in{\mathcal O}^*(Z)$
and $c,\eta_2,...,\eta_m\in{\mathcal O}(Z)$
such that $\eta_j=0$ for $j\in{\mathcal Z}(w^\circ)$, 
$\eta_j=\zeta^{j/{\mathcal D}(w^\circ)}$ for $j\in{\mathcal{NZ}}(w^\circ)$
and 
\begin{equation}\label{eq: general form of F: ultra-Picard space Z to Gm}
F(z)
=(t+c(z))^m + \sum_{j=2}^m w_j^\circ \eta_j(z)(t+c(z))^{m-j}\,,
\quad z\in Z\,.
\end{equation}
\vskip0.1cm

{\bfit b}$)$ If $w^\circ_m\ne 0$ then 
the map {\rm(\ref{eq: general form of F: ultra-Picard space Z to Gm})}
is homotopic to the model map
\begin{equation}\label{eq: A-map}
{\mathcal A}_a\colon Z\ni z\mapsto (0,...,0,0,a(z))
=t^m + a(z)\in{\mathbf G}_m^0\subset{\mathbf G}_m\,,
\end{equation}
where $a(z)=[\zeta(z)]^{m/{\mathcal D}(w^\circ)}$.
\vskip0.1cm

{\bfit c}$)$ If $w_m^\circ=0$ then $w^\circ_{m-1}\ne 0$
and the map {\rm(\ref{eq: general form of F: ultra-Picard space Z to Gm})}
is homotopic to the model map
\begin{equation}\label{eq: B-map}
{\mathcal B}_b\colon Z\ni z\mapsto (0,...,0,b(z),0)
=t^m + b(z)t\in{\mathbf G}_m^0\subset{\mathbf G}_m\,,
\end{equation}
where $b(z)=[\zeta(z)]^{(m-1)/{\mathcal D}(w^\circ)}$.
\end{Proposition}

\begin{proof}
{\bfit a}$)$ Let $F(z)=(F_1(z),F_2(z),...,F_m(z))$; consider
a holomorphic map $F^0=(0,F_2^0,...,F_m^0)\colon Z\to{\mathbf G}_m^0$
with the components $F_2^0,...,F_m^0$ defined by the condition
\begin{equation}\label{eq: F and F^0}
t^m+\sum_{j=2}^m F_j^0(z)t^{m-j}
\equiv(t-F_1(z)/m)^m+\sum_{j=1}^m F_j(z)(t-F_1(z)/m)^{m-j}\,.
\end{equation}
By Proposition \ref{Prp: holomorphic maps to Cmb(C)},
$F^0(Z)\subseteq{\mathbb C}^*w^\circ$ for a certain point
$w^\circ=(0,w_2^\circ,...,w_m^\circ)\in{\mathbf G}_m^0$. 
It follows from Lemma \ref{Lm: parameterization of C* orbit in Gm0}
that there are $\zeta\in{\mathcal O}^*(Z)$
and $\eta_2,...,\eta_m\in{\mathcal O}(Z)$
such that $\eta_j=0$ for $j\in{\mathcal Z}(w^\circ)$, 
$\eta_j=\zeta^{j/{\mathcal D}(w^\circ)}$ for $j\in{\mathcal{NZ}}(w^\circ)$
and
\begin{equation}\label{eq: general form of F0: Z to Gm0}
F^0(z)=(0,w_2^\circ\eta_2(z),...,w_m^\circ\eta_m(z))
=t^m + \sum_{j=2}^m w_j^\circ\eta_j(z)t^{m-j}\,,\quad z\in Z\,.
\end{equation}
Let $c(z)=F_1(z)/m$; then (\ref{eq: F and F^0}) and
(\ref{eq: general form of F0: Z to Gm0}) show that 
\begin{equation}\label{eq: general form of F: Z to Gm}
\aligned
F(z)&=
t^m+\sum_{j=1}^m F_j(z)t^{m-j}=
(t+c(z))^m+\sum_{j=2}^m F_j^0(z)(t+c(z))^{m-j}\\
&=(t+c(z))^m+\sum_{j=2}^m w_j^\circ\eta_j(z)(t+c(z))^{m-j}\,,
\endaligned
\end{equation}
which proves $(${\bfit a}$)$.
\vskip0.2cm

$(${\bfit{b,c}}$)$ By $(${\bfit a}$)$, $F$ is of the
form (\ref{eq: general form of F: ultra-Picard space Z to Gm}).
The homotopy 
\begin{equation}\label{eq: homotopy F to F0}
F(z,s)
=(t+sc(z))^m + \sum_{j=2}^m w_j^\circ \eta_j(z)(t+sc(z))^{m-j}\,,
\quad z\in Z\,, \ \ 0\le s\le 1\,,
\end{equation}
joins $F=F(\cdot,1)$ to the map
$F(\cdot,0)$, $F(z,0)=t^m + \sum_{j=2}^m w_j^\circ\eta_j(z) t^{m-j}$ 
in the class of holomorphic maps $Z\to{\mathbf G}_m$. At least one of
the numbers $m-1,m$ is not in ${\mathcal Z}$ and the vector subspace
$$
V(w^\circ)\Def\{w=(0,w_2,...,w_m)\in{\mathbb C}^m\,|\
w_j=0 \ \ \text{whenever} \ j\in{\mathcal Z}\}
$$ 
is of positive dimension. The discriminant $d_m(w)$ cannot
vanish identically on $V(w^\circ)$, for
$w^\circ\in {\mathbf G}_m^0\cap V(w^\circ)$;
thus, the intersection
$$
{\mathbf G}_m^0\cap V(w^\circ)
=V(w^\circ)\setminus\{w=(0,w_2,...,w_m)\in{\mathbb C}^m\,|\ d_m(w)=0\}
$$
is arcwise connected. Set $w'=(0,...,0,0,1)=t^m+1$ when $w^\circ_m\ne 0$
and $w'=(0,...,0,1,0)=t^m+t$ otherwise;
notice that in both cases $w'\in{\mathbf G}_m^0\cap V(w^\circ)$.
Let $w(s)=(0,w_2(s),...,w_m(s))$, \ $0\le s\le 1$, be a path in
${\mathbf G}_m^0\cap V(w^\circ)$ joining the point
$w(0)=w^\circ$ to the point $w(1)=w'$. It is easily seen that
$$
d_m(0,w_2(s)\eta_2(z),...,w_m(s)\eta_m(z))
=[\zeta(z)]^{m(m-1)/{\mathcal D}(w^\circ)}d_m(0,w_2(s),...,w_m(s))\,;
$$
therefore the formula
$h_s(\zeta)=t^m + \sum_{j=2}^m w_j(s)\eta_j(z) t^{m-j}$
provides a well-defined homotopy joining the map $F(\cdot,0)$ to
a model map ${\mathcal A}_a$ in case $(${\bfit b}$)$ 
and to a model map ${\mathcal B}_b$ in case $(${\bfit c}$)$.
\end{proof}

\begin{Remark}\label{Rmk: homotopy between model maps}
Let $a,b,a_i,b_i\in{\mathcal O}^*(Z)$, \ 
$i=1,2$.
Then\footnote{Here ``$\sim$" means free homotopy.}
${\mathcal A}_{a_1}\sim {\mathcal A}_{a_2}$ $($respectively
${\mathcal B}_{b_1}\sim {\mathcal B}_{b_2})$ if and only if
$a_1/a_2\in \exp{\mathcal O}(Z)$ $($respectively if and only if
$b_1/b_2\in \exp{\mathcal O}(Z))$. Moreover,
${\mathcal A}_a\sim {\mathcal B}_b$
if and only if $a^{m-1}/b^m\in\exp{\mathcal O}(Z)$.
\hfill $\bigcirc$
\end{Remark}

\begin{Definition}\label{Def: model maps C* to Cm(C)=Gm}
\index{Model holomorphic maps
${\mathbb C}^*\to{\mathcal C}^m({\mathbb C})$\hfill}
Define the points $w_\circ^{(m)}, w_\circ^{(m-1)}\in{\mathbf G}_m$ as 
$$
w_\circ^{(m)}=(0,...,0,0,-1)=t^m-1\quad\text{and}\quad
w_\circ^{(m-1)}=(0,...,0,-1,0)=t^m-t
$$
respectively. For any $r\in{\mathbb Z}$ we define the {\em model
holomorphic maps}
\begin{equation}\label{eq: model maps C* to Cm(C)=Gm}
\aligned
&{\mathcal A}_r\colon{\mathbb C}^*\ni\zeta
\mapsto{\mathcal U}_{\zeta^r} w_\circ^{(m)}\hskip10pt
=(0,...,0,0,-\zeta^r)=t^m-\zeta^r
\in{\mathbf G}_m^0\subset{\mathbf G}_m={\mathcal C}^m({\mathbb C})\,,\\
&{\mathcal B}_r\colon{\mathbb C}^*\ni\zeta
\mapsto{\mathcal U}_{\zeta^r} w_\circ^{(m-1)}=(0,...,0,-\zeta^r,0)
=t^m-\zeta^r t\in{\mathbf G}_m^0\subset{\mathbf G}_m={\mathcal C}^m({\mathbb C})\,.
\endaligned
\hskip-4pt
\end{equation}
Notice that ${\mathcal A}_p\sim{\mathcal B}_q$ if and only if $(m-1)p=mq$.
\hfill $\bigcirc$
\end{Definition}

\begin{Corollary}\label{Crl: holomorphic maps of C* to Cm(C)=Gm}
Any holomorphic map ${\mathbb C}^*\to{\mathcal C}^m({\mathbb C})={\mathbf G}_m$
is homotopic to one of the model maps ${\mathcal A}_r$, ${\mathcal B}_r$ \
$(r\in{\mathbb Z})$.
\hfill $\square$
\end{Corollary}

\subsection{Cyclic morphisms ${\mathcal C}^n({\mathbb C})
\to{\mathcal C}^k({\mathbb C})$}
\label{Cyclic morphisms Cn(C) to Ck(C)}
\index{Cyclic morphisms ${\mathcal C}^n({\mathbb C})
\to{\mathcal C}^k({\mathbb C})$\hfill}
Notice that any non-vanishing holomorphic function $f(Q)$ on
${\mathcal C}^n({\mathbb C})$ is of the form
$f(Q)=[d_n(Q)]^l e^{\varphi(Q)}$, where $l\in{\mathbb Z}$ and
$\varphi\in{\mathcal O}({\mathcal C}^n({\mathbb C}))$.
Since any connected algebraic variety is a Liouville space, Proposition
\ref{Prp: ultra-Picard or Liouville space to Gm} applies to cyclic
morphisms ${\mathcal C}^n({\mathbb C})\to{\mathcal C}^k({\mathbb C})$:
 
\begin{Theorem}\label{Thm: cyclic morphisms Gn to Gk}
\index{Cyclic morphisms ${\mathcal C}^n({\mathbb C})
\to{\mathcal C}^k({\mathbb C})$\hfill}
Let $F\colon{\mathcal C}^n({\mathbb C})={\mathbf G}_n\to{\mathbf G}_k
={\mathcal C}^k({\mathbb C})$ be a cyclic holomorphic map.
Then: 
\vskip0.1cm

{\bfit a}$)$ There exist $w^\circ=(0,w^\circ_2,...,w^\circ_k)
\in{\mathbf G}_k^0={\mathcal C}^{k-1}_b({\mathbb C})$, 
$l\in{\mathbb Z}$ and
$c,\varphi\in{\mathcal O}({\mathcal C}^n({\mathbb C}))$
such that
\begin{equation}\label{eq: general form of F: Cn(C) to Ck(C)}
F(Q)=(t+c(Q))^k 
+ \sum_{j=2}^k w_j^\circ [d_n(Q)]^{jl/{\mathcal D}(w^\circ)}
e^{jl\varphi(Q)}(t+c(Q))^{k-j}\,,
\quad Q\in{\mathcal C}^n({\mathbb C})\,.
\end{equation}
\vskip0.1cm

{\bfit b}$)$ If $w^\circ_k\ne 0$ then
$kl/{\mathcal D}(w^\circ)\in{\mathbb Z}$ and
the map {\rm(\ref{eq: general form of F: Cn(C) to Ck(C)})}
is homotopic to the model map
\begin{equation}\label{eq: Al-map}
{\mathcal A}_l\colon{\mathcal C}^n({\mathbb C})\ni Q
\mapsto (0,...,0,0,[d_n(Q)]^{kl/{\mathcal D}(w^\circ)})
=t^k + [d_n(Q)]^{kl/{\mathcal D}(w^\circ)}\,.
\end{equation}
\vskip0.1cm

{\bfit c}$)$ If $w_k^\circ=0$ then $w^\circ_{k-1}\ne 0$,
$(k-1)l/{\mathcal D}(w^\circ)\in{\mathbb Z}$
and the map {\rm(\ref{eq: general form of F: Cn(C) to Ck(C)})}
is homotopic to the model map
\begin{equation}\label{eq: Bl-map}
{\mathcal B}_l\colon{\mathcal C}^n({\mathbb C})\ni Q
\mapsto (0,...,0,[d_n(Q)]^{(k-1)l/{\mathcal D}(w^\circ)},0)
=t^k + [d_n(Q)]^{(k-1)l/{\mathcal D}(w^\circ)}t\,.
\end{equation}

{\bfit d}$)$ ${\mathcal A}_{l_1}\sim {\mathcal A}_{l_2}$
$($respectively
${\mathcal B}_{l_1}\sim {\mathcal B}_{l_2})$ if and only if
$l_1=l_2$. Moreover, ${\mathcal A}_{l_1}\sim {\mathcal B}_{l_2}$
if and only if $(k-1)l_1=kl_2$.
\vskip0.2cm

\noindent In particular, every cyclic holomorphic map 
$F\colon{\mathcal C}^n({\mathbb C})\to{\mathcal C}^k({\mathbb C})$
homotopically splits, that is, there exists a holomorphic map
$G\colon{\mathbb C}^*\to{\mathcal C}^k({\mathbb C})$
such that $F$ is homotopic to the composition
$$
\hskip127pt 
G\circ d_n\colon {\mathcal C}^n({\mathbb C})
\stackrel{d_n}{\longrightarrow}{\mathbb C}^*
\stackrel{G}{\longrightarrow}{\mathcal C}^k({\mathbb C})\,.
\hskip128pt \square
$$
\end{Theorem}

\section{Morphisms ${\mathcal C}^n({\mathbb C})
\to{\mathcal C}^k({\mathbb C})$ and special homomorphisms $B_n\to B_k$}
\label{Morphisms Cn(C) to Ck(C) and special homomorphisms Bn to Bk} 

\noindent Due to Theorem \ref{Thm: Bn(X) to Bk(X) for n>k}, 
for $n\ne 4$ and $n>k$ every continuous map (and all the more
every morphism) ${\mathcal C}^n(X)\to{\mathcal C}^k(X)$ must be cyclic. 
On the other hand, for $3\le n\le k$ it is easy to construct
non-cyclic continuous maps ${\mathcal C}^n({\mathbb C})
\to{\mathcal C}^k({\mathbb C})$.
It is not so in the case $X={\mathbb{CP}}^1$. To see this, one may
use a theorem of Kunio Murasugi \cite{Mur82}, which provides
a complete description of the torsion in
the sphere braid group $B_m(S^2)$ and in the quotient group
$B_m/CB_m$ of the Artin braid group $B_m$ by its center $CB_m$.
\vskip0.2cm

\noindent Let us set $\alpha=\sigma_1\cdots\sigma_{m-1}$,
$\beta=\sigma_1\cdots\sigma_{m-1}\sigma_1$ (both in 
$B_m$ and $B_m(S^2)$),
and $\gamma=\sigma_1\cdots\sigma_{m-2}\sigma_1$ in $B_m(S^2)$;
furthermore, let us denote by $\hat\alpha$, $\hat\beta$
the images of the elements $\alpha,\beta\in B_m$
in the quotient group $B_m/CB_m$. 
\vskip0.2cm

\noindent{\caps Murasugi Theorem.}
\index{Theorem!Murasugi Theorem\hfill}
\index{Torsion of $B_m/CB_m$ and $B_m(S^2)$\hfill}
{\bfit a}$)$ {\sl Each element of finite order in $B_m/CB_m$
is conjugate to an element of the form
${\hat\alpha}^p$ or ${\hat\beta}^p$, where $p\in{\mathbb Z}$.
The elements $\hat\alpha,\hat\beta\in B_m/CB_m$
are of the order $m$ and $m-1$ respectively.}
\vskip0.1cm

{\bfit b}$)$ {\sl Each element of finite order in $B_m(S^2)$
is conjugate to an element of the form $\alpha^p$, $\beta^p$ or
$\gamma^p$, where $p\in{\mathbb Z}$.
The elements $\alpha,\beta,\gamma\in B_m(S^2)$
are of the order $2m$, $2(m-1)$ and $2(m-2)$ respectively.}  
\vskip0.2cm

\noindent For any homomorphism $\varphi\colon B_n(S^2)\to B_k(S^2)$
we have $\varphi(\Tors B_n(S^2))\subseteq\Tors B_k(S^2)$. Hence
part $(${\bfit b}$)$ of Murasugi Theorem implies that $\varphi$
cannot be non-cyclic unless the number $n(n-1)(n-2)$ divides
$k(k-1)(k-2)$.
\vskip0.2cm

\noindent A homomorphism $\varphi\colon B_n\to B_k$
(cyclic or not) must not induce a homomorphism of quotient groups
$B_n/CB_n\to B_k/CB_k$. Hence Murasugi Theorem puts no restraints
in homomorphisms $B_n\to B_k$. In fact, for $k>n\ge 3$ there is plenty
of non-cyclic homomorphisms $\varphi\colon B_n\to B_k$ and each of
them is induced by a certain non-cyclic continuous map
${\mathcal C}^n({\mathbb C})\to{\mathcal C}^k({\mathbb C})$.
\vskip0.2cm

\noindent For holomorphic maps the situation is completely different.
In fact, I do not know any example of a non-cyclic holomorphic map
${\mathcal C}^n(X)\to{\mathcal C}^k(X)$ for $4<n<k$
(compare to Linked Map Theorem); most likely they do not exist
at all. Being unable to prove this I present
here a rather restricted result, which says that in the case
$X={\mathbb C}$ such a map cannot exist unless the number $n(n-1)$
divides $k(k-1)$. In order to prove this, we will show that
a homomorphism $F_*\colon B_n\to B_k$ induced by a holomorphic
map ${\mathcal C}^n({\mathbb C})\to{\mathcal C}^k({\mathbb C})$
must satisfy a strong algebraic condition. 

\subsection{Special presentations and special homomorphisms}
\label{Sec: Special presentations and special homomorphisms}
Recall that $\alpha=\sigma_1\cdots\sigma_{m-1}$
and $\beta=\alpha\sigma_1$ generate the whole braid group $B_m(X)$
(see \ref{eq: conjugation of sigma by alpha}).

\begin{Definition}\label{Def: special system of generators}
\index{Special system of generators in $B_m$\hfill}
\index{Standard system of generators in $B_m$\hfill}
A pair of elements $a,b\in B_m$ is said to
be a {\em special system of generators} in $B_m$
if there exists an automorphism $\psi$ of $B_m$
such that $\psi(\alpha) = a$ and $\psi(\beta) = b$.
If $\{a,b\}$ is such a system of generators then the elements
\begin{equation*}
s_{i} = \psi(\sigma_{i}) = a^{i-2}ba^{-(i-1)}\qquad (1\le i\le m-1)
\end{equation*}
also form a system of generators of $B_m$ that satisfy
relations (\ref{eq: commutation relation in Bn}),
(\ref{eq: braid relation});
we call such a system of generators {\em standard}.
\hfill $\bigcirc$
\end{Definition}

\begin{Definition}\label{Def: special homomorphism}
\index{Special presentations of $B_m$\hfill}
\index{Special homomorphisms $B_m\to B_k$\hfill}
Given a special system of generators $\{a,b\}$ in $B_m$,
we denote by ${\mathcal H}_m(a,b)$ the set
of all elements of the form $g^{-1}a^p g$ or $g^{-1}b^p g$, 
where $g$ runs over
$B_m$ and $p$ runs over ${\mathbb Z}$.
A homomorphism $\varphi\colon B_n\to B_k$
is said to be {\em special} if $\varphi({\mathcal H}_n(a,b))\subseteq
{\mathcal H}_k(a',b')$ for some choice of special systems of
generators $a,b\in B_n$ and $a',b'\in B_k$.
\hfill $\bigcirc$
\end{Definition}

\noindent Part $(${\bfit a}$)$ of Murasugi Theorem quoted above
shows that {\sl the set ${\mathcal H}_m={\mathcal H}_m(a,b)$
does not depend on a choice of a special system of generators
$a,b\in B_m$, and a homomorphism $\varphi\colon B_n\to B_k$
is special if and only if for any element $g\in B_m$ of finite order
modulo $CB_n$ its image $\varphi(g)\in B_k$ is an element of finite
order modulo $CB_k$} (we do not use this fact below;
it is convenient however to keep it in mind). 

\begin{Remark}\label{Rmk: endomorphism is special}
It is easily seen that for  $n\ge 3$ every {\em endomorphism}
of the group $B_n$ is special. On the other hand, for $k>n\ge 3$
the standard embedding $B_n\hookrightarrow B_k$ can never be
a special homomorphism.
\hfill $\bigcirc$
\end{Remark}

\noindent The following result stated in \cite{Lin71}
and \cite{Lin79} is the crucial one for our purposes.

\begin{Theorem}
\label{Thm: homomorphisms Bn to Bk induced by holomorphic maps}
\index{Theorem!Homomorphisms $B_n\to B_k$!induced by morphisms\hfill}
For every holomorphic map $F\colon{\mathcal C}^n({\mathbb C})
\to{\mathcal C}^k({\mathbb C})$,
every point $Q^\circ\in{\mathcal C}^n({\mathbb C})$ and any choice 
of isomorphisms $B_n\cong\pi_1({\mathcal C}^n({\mathbb C}),Q^\circ)$
and $\pi_1({\mathcal C}^n({\mathbb C}),F(Q^\circ))\cong B_k$
the induced homomorphism
\begin{equation}\label{eq: induced homomorphism in details}
F_*\colon B_n\cong\pi_1({\mathcal C}^n({\mathbb C}),Q^\circ)
\to\pi_1({\mathcal C}^n({\mathbb C}),F(Q^\circ))\cong B_k
\end{equation}
is special.
\hfill $\bigcirc$
\end{Theorem}

\noindent To prove this theorem we will use the description
of holomorphic maps ${\mathbb C}^*\to{\mathcal C}^m({\mathbb C})$
given in Corollary \ref{Crl: holomorphic maps of C* to Cm(C)=Gm}.
But first I want to present a general concept of holomorphic elements
of homotopy groups introduced in \cite{Lin71}.

\subsection{Holomorphic elements of homotopy groups}
\label{Ss: Holomorphic elements of homotopy groups}
\index{Holomorphic elements of homotopy groups\hfill}
The idea of this section is to choose for each $i\in{\mathbb N}$
an appropriate {\em test complex manifold} ${\mathcal M}_i$
homotopy equivalent to the $i$-dimensional sphere $S^i$ and
define corresponding {\em holomorphic elements} in the
homotopy groups $\pi_i(Y)$ of connected complex spaces $Y$
as those admitting holomorphic representatives
$M_i\to Y$. Doing so one may sometimes come to certain non-trivial
algebraic restraints imposed by the condition that a homomorphism
$\pi_i(Y)\to\pi_i(Z)$ can be induced by a holomorphic map $Y\to Z$,
since such a homomorphism must carry holomorphic elements of
$\pi_i(Y)$ to holomorphic elements of $\pi_i(Z)$. 
What are ``test manifolds" depends on the problem
under consideration; for our purpose the following choice is an
appropriate one just because it works.

\begin{Definition}\label{Def: holomorphic part of pii(Y)}
\index{Holomorphic elements of homotopy groups\hfill}
\index{Holomorphic part $\hol\pi_i(Y)$ of the homotopy group $\pi_i(Y)$\hfill}
Let $Y$ be a connected complex space. We define the
{\em holomorphic part} ${\text{\bfit hol}}\,\pi_i(Y,y_\circ)$ 
of the homotopy group $\pi_i(Y,y_\circ)$ as follows. 
In the space ${\mathbb C}^{i+1}$ with coordinates
$\zeta_0,...,\zeta_i$ consider the quadric 
$$
{\mathbb C}S^i = \left\{\zeta=(\zeta_0,...,\zeta_i)
\in{\mathbb C}^{i+1}\,|\ \zeta_0^2 +...+ \zeta_i^2 = 1\right\}\,.
$$
Its intersection
$S^i = {\mathbb C}S^i\cap{\mathbb R}^{i+1}$ with the real
subspace ${\mathbb R}^{i+1}\subset{\mathbb C}^{i+1}$ is the standard
unit sphere in ${\mathbb R}^{i+1}$; moreover, $S^i$ is a deformation
retract of ${\mathbb C}S^i$ (for ${\mathbb C}S^i$ is diffeomorphic
to the total space of the tangent bundle $TS^i$).  
Therefore choosing $\zeta_\circ=(1,0,...,0)\in{\mathbb C}S^i$
as a basepoint we may regard any element
$\gamma\in\pi_i(Y,y_\circ)$ as a homotopy class
of continuous mappings $({\mathbb C}S^i,\zeta_\circ)\to (Y,y_\circ)$.
Let us say that an element $\gamma\in \pi_i(Y,y_\circ)$ is
{\em holomorphic} if there exists a holomorphic mapping 
$h\colon{\mathbb C}S^i\to Y$ freely homotopic
to some (and hence to any) map contained in $\gamma$.
The set of all holomorphic elements in $\pi_i(Y,y_\circ)$ is
denoted by $\hol\pi_i(Y,y_\circ)$; in general, it is not a subgroup;
the subgroup of $\pi_i(Y,y_\circ)$ generated by $\hol\pi_i(Y,y_\circ)$
is denoted by $\Hol\pi_i(Y,y_\circ)$. 
\hfill $\bigcirc$
\end{Definition}

\begin{Proposition}\label{Prp: properties of hol pii(Y)}
The subset $\hol\pi_i(Y,y_\circ)$ and the subgroup
$\Hol\pi_i(Y,y_\circ)$ have the following properties:
\vskip0.2cm

\begin{itemize}

\item[$i)$] If $Y$ is a homogeneous space of a complex
Lie group then $\hol\pi_i(Y,y_\circ)=\pi_i(Y,y_\circ)$
for every natural $i$. In particular, the complex orthogonal group
${\mathbf O}(n+1,{\mathbb C})$ acts transitively on
${\mathbb C}S^m$ and hence
$$
\hol\pi_i({\mathbb C}S^m,\zeta_\circ)=
\pi_i({\mathbb C}S^m,\zeta_\circ)\cong\pi_i(S^m)
$$
for all natural $i$ and $m$.
\vskip0.1cm

\item[$ii)$] If $\gamma\in\hol\pi_i(Y,y_\circ)$ then 
the whole cyclic subgroup generated by $\gamma$ is contained
in $\hol\pi_i(Y,y_\circ)$.
\vskip0.1cm

\item[$iii)$] The isomorphism $\chi_i\colon \pi_i(Y,y_1)\to\pi_i(Y,y_0)$
induced by a path $\chi\colon [0,1]\to Y$
joining points $y_0,y_1\in Y$ carries $\hol\pi_i(Y,y_1)$ into
$\hol\pi_i(Y,y_0)$.
\vskip0.1cm

\item[$iv)$] The homomorphism $f_{*i}\colon\pi_i(Y,y_\circ)
\to\pi_i(Z,z_\circ)$ induced by a holomorphic map $f\colon (Y,y_\circ)
\to (Z,z_\circ)$ carries $\hol\pi_i(Y,y_\circ))$
into $\hol\pi_i(Z,z_\circ)$.
\vskip0.1cm

\item[$v)$] $\Hol\pi_1(Y,y_\circ)$ is a normal subgroup
of the fundamental group $\pi_1(Y,y_\circ)$.
\end{itemize}
\end{Proposition}

\begin{proof}
Property $(i)$ follows from the Cartan--Grauert--Ramspott
theorem (see \cite{Car40,Gra58,Rams65}). It follows from $(i)$ that
$\hol\pi_i({\mathbb C}S^i,\zeta_\circ)
=\pi_i({\mathbb C}S^i,\zeta_\circ)\cong\pi_i(S^i,\zeta_\circ)
\cong{\mathbb Z}$. Hence if $h\colon{\mathbb C}S^i\to X$
is a holomorphic
representative of a class $\gamma\in\hol\pi_i(Y,y_\circ)$ and 
$f^{(q)}\colon{\mathbb C}S^i\to{\mathbb C}S^i$ is a holomorphic
map homotopic to a map $S^i\to S^i$ of degree $q\in{\mathbb Z}$
then the composition 
${\mathbb C}S^i\stackrel{f^{(q)}}{\longrightarrow}{\mathbb C}S^i
\stackrel{h}{\longrightarrow} Y$
is a holomorphic representative of the class $q\gamma$
(for $i=1$ we use the multiplicative notation $\gamma^q$);
this proves $(ii)$. Properties
$(iii)$ and $(iv)$ follow directly from the definitions
and $(v)$ is a consequence of $(iii)$.
\end{proof}


\subsection{Holomorphic elements of $\pi_i({\mathcal C}^m({\mathbb C}))$}
\label{Holomorphic elements of pii(Cm(C))}
\index{Holomorphic elements of $\pi_i({\mathcal C}^m({\mathbb C}))$\hfill}
It follows from
Proposition \ref{Prp: holomorphic maps to Cm(X)} that
the image of any holomorphic map
$f\colon{\mathbb C}S^i\to{\mathcal C}^m({\mathbb C})$
is contained in the $\Aut X$-orbit of some point
$Q^*\in{\mathcal C}^m({\mathbb C})$.
\vskip0.2cm

\noindent First, let $i\ge 2$. Then
$\hol\pi_i({\mathcal C}^m({\mathbb C}))
=\pi_i({\mathcal C}^m({\mathbb C}))=0$ (moreover,
it is easily shown that for such $i$ any holomorphic map
${\mathbb C}S^i\to{\mathcal C}^m({\mathbb C})$ is contractible
in the class of holomorphic maps).
Thus, on this way no difference occurs in the behavior   
of higher homotopy group homomorphisms induced by holomorphic
maps and continuous ones. 
\vskip0.2cm

\noindent Let us compute holomorphic elements of
$\pi_1({\mathcal C}^m({\mathbb C}))$; we use the concept of
model holomorphic maps ${\mathcal A_r},{\mathcal B_r}$ introduced
in Definition \ref{Def: model maps C* to Cm(C)=Gm}. 

\begin{Proposition}\label{Prp: holomorphic elemens of pi1(Cm(C))}
\index{Holomorphic part $\hol B_m$ of $B_m$\hfill}
\index{$\hol B_m$\hfill}
{\bfit a}$)$ For any $Q^*\in{\mathcal C}^m({\mathbb C})={\mathbf G}_m$
there is a special system of generators $\{a,b\}$ in
$\pi_1({\mathcal C}^m({\mathbb C}))\cong B_m$
such that the set $\hol\pi_1({\mathcal C}^m({\mathbb C}),Q^*)$
of all holomorphic elements in $\pi_1({\mathcal C}^m({\mathbb C}),Q^*)$
coincides with ${\mathcal H}_m(a,b)$.
\vskip0.1cm
{\bfit b}$)$ The set $\hol B_m$ of all elements that for some choice
of a basepoint $Q^*\in{\mathcal C}^m({\mathbb C})$ and an isomorphism
$B_m\cong\pi_1({\mathcal C}^m({\mathbb C}),Q^*)$ admit
holomorphic representatives ${\mathbb C}^*\to{\mathcal C}^m({\mathbb C})$  
coincides with the set of all elements of finite order
modulo the center $CB_m$.
\end{Proposition}

\begin{proof}
For any element $h\in\pi_1({\mathcal C}^m({\mathbb C}),Q^*)$
let us denote by $[h]$ the free homotopy class
of continuous mappings ${\mathbb C}^*\to{\mathcal C}^m({\mathbb C})$
containing some (and hence each one) representative of $h$.
A straightforward geometric consideration shows that
a standard system of generators $s_1,...,s_{m-1}$ in the group 
$\pi_1({\mathcal C}^m({\mathbb C}),Q^*)=\pi_1({\mathbf G}_m,w^*)
\cong B_m$ may be chosen in such a way that the model holomorphic map
${\mathcal A}={\mathcal A}_1\colon{\mathbb C}^*\ni\zeta
\mapsto (0,...,0,0,-\zeta)=t^m-\zeta$ belongs to the free homotopy
class $[a]$ corresponding to the element $a=s_1\cdots s_{m-1}$ and
the map ${\mathcal B}={\mathcal B}_1\colon{\mathbb C}^*\ni\zeta
\mapsto (0,...,0,-\zeta,0)=t^m-\zeta t$ belongs to the free homotopy
class $[gbg^{-1}]$ corresponding to the element $gbg^{-1}$, where
$b=as_1=s_1\cdots s_{m-1}s_1$ and $g$ is a suitable element
of $\pi_1({\mathbf G}_m^0,w_*)$. Then it is clear that
${\mathcal A}_r\in [a^r]$ and ${\mathcal b}_r\in g[b^rg^{-1}]$.
On the other hand, by Corollary
\ref{Crl: holomorphic maps of C* to Cm(C)=Gm},
any holomorphic map ${\mathbb C}^*\to{\mathcal C}^m({\mathbb C})
={\mathbf G}_m$
is homotopic to one of the model maps ${\mathcal A}_r$, ${\mathcal B}_r$ \
$(r\in{\mathbb Z})$; this completes the proof of statement $(${\bfit a}$)$.
In view of the Murasugi theorem quoted above,
statement $(${\bfit b}$)$ is just a consequence of $(${\bfit a}$)$
\end{proof}

\noindent Now we are ready to prove Theorem
\ref{Thm: homomorphisms Bn to Bk induced by holomorphic maps}.
\vskip0.2cm

{\bcaps Proof of Theorem} \ref{Thm: homomorphisms Bn to Bk induced
by holomorphic maps}.
\index{Theorem!Homomorphisms $B_n\to B_k$!induced by morphisms\hfill}
Let $F\colon{\mathcal C}^n({\mathbb C})
\to{\mathcal C}^k({\mathbb C})$ be a holomorphic map.
By Proposition \ref{Prp: properties of hol pii(Y)}$(iv)$,
the homomorphism $F_*\colon\pi_1({\mathcal C}^n({\mathbb C}),Q^\circ)
\to\pi_1({\mathcal C}^k({\mathbb C}),F(Q^\circ))$
of the fundamental groups induced by $F$
carries $\hol\pi_1({\mathcal C}^n)$
into $\hol\pi_1({\mathcal C}^k)$. According to Proposition
\ref{Prp: holomorphic elemens of pi1(Cm(C))}$(${\bfit a}$)$, this implies
that $F_*({\mathcal H}_n(a,b))\subseteq{\mathcal H}_k(a,b)$
and $F_*$ is special.
\hfill $\square$

\subsection{Special homomorphisms $B_n\to B_k$
and morphisms ${\mathcal C}^n({\mathbb C})\to{\mathcal C}^k({\mathbb C})$}
\label{Ss: Special homomorphisms Bn to Bk and morphisms Cn(C) to Ck(C)}
The following result proved in \cite{Lin96b}, Theorem 8.1
(see also \cite{Lin04b}, Sec. 9) clarifies
the arithmetical consequences of the assumption that the Artin braid
group $B_n$ admits a non-cyclic special homomorphism to the
Artin braid group $B_k$ and, in view of Theorem
\ref{Thm: homomorphisms Bn to Bk induced by holomorphic maps},
imposes the corresponding restraints on the possible existence
of non-cyclic morphisms
${\mathcal C}^n({\mathbb C})\to{\mathcal C}^k({\mathbb C})$.

\begin{Theorem}\label{Thm: non-cyclic special homomorphisms}
\index{Theorem!on non-cyclic special homomorphisms $B_n\to B_k$\hfill}
If $n\ne 4$ and there is a non-cyclic special homomorphism
$B_n\to B_k$ then $n(n-1)$ divides $k(k-1)$.
\hfill $\bigcirc$
\end{Theorem}

\noindent Theorems \ref{Thm: homomorphisms Bn to Bk induced
by holomorphic maps} and \ref{Thm: non-cyclic special homomorphisms}
imply the following

\begin{Corollary}\label{Crl: morphisms Cn(C) to Ck(C)}
If $n\ne 4$ and $n(n-1)$ does not divide $k(k-1)$ then every holomorphic
map ${\mathcal C}^n({\mathbb C})\to{\mathcal C}^k({\mathbb C})$
is cyclic.
\hfill $\square$
\end{Corollary}

\begin{Remark}\label{Rmk: existence of non-cyclic special homomorphisms} 
The group $B_3$ possesses non-cyclic special homomorphisms
$B_3\to B_k$ for any $k$ that is not forbidden by Theorem 
\ref{Thm: non-cyclic special homomorphisms}.
Moreover, let $\alpha,\beta\in B_3$ and $a,b\in B_4$ 
be special systems of generators; then
the map $\alpha\mapsto b$, $\beta\mapsto a^2$ extends uniquely
to a surjective (and thereby non-cyclic) homomorphism
$B_4\to B_3$; hence, if $6$ divides $k(k-1)$ then
non-cyclic special homomorphisms $B_4\to B_k$ do exist.
\vskip0.2cm

\noindent Suppose that $n>4$ and $k$ is not forbidden by Theorem 
\ref{Thm: non-cyclic special homomorphisms}, i.e.,
$n(n-1)$ divides $k(k-1)$. Then non-cyclic special homomorphisms
$B_n\to B_k$ certainly exist whenever $n$ divides $k$ (see \cite{Lin96b},
Section 8.1, Theorem 8.3); however, I do not know whether such
homomorphisms can exist when $n$ does not divide $k$.
\hfill $\bigcirc$
\end{Remark}


\vskip0.2cm



\section{Dimension of the image}
\label{Sec: Dimension of the image}

\noindent Let $n>4$ and let
$F\colon{\mathcal C}^n(X)\to{\mathcal C}^n(X)$
be a morphism. If $F$ is cyclic then, 
by Cyclic Map Theorem \ref{Thm: Cyclic Map Theorem},
it is orbit-like and hence
$\dim_{\mathbb C}F({\mathcal C}^n(X))\le\dim_{\mathbb C}\Aut X=t(X)$
(see Notation \ref{Not: t(X)}). If $F$ is non-cyclic then it is tame
and hence $F((\Aut X)Q)\subseteq (\Aut X)Q$ for all
$Q\in{\mathcal C}^n(X)$; thus, each $\Aut X$ orbit
in ${\mathcal C}^n(X)$ contains at least one point of its own
image and therefore $\dim_{\mathbb C}F({\mathcal C}^n(X))
\ge n-\dim_{\mathbb C}\Aut X=n-t(X)$. The next theorem strengthens
this estimate.

\begin{Theorem}\label{Thm: Dimension of the image}
Let $n>4$ and let
$F\colon{\mathcal C}^n(X)\to{\mathcal C}^n(X)$
be a non-cyclic morphism. Then $\dim_{\mathbb C}F({\mathcal C}^n(X))
\ge n-t(X)+1$.
\end{Theorem}
 
\noindent This theorem follows immediately from
Tame Map Theorem, the existence of points $Q\in {\mathcal C}^n(X)$
with non-trivial stabilizers $St_Q$
and the following local result.

\begin{Lemma}\label{Lm: local dimension}
Let $T\colon U\to\Aut X$ be a holomorphic map of a non-empty open
$\Aut X$ invariant subset $U\subseteq{\mathcal C}^n(X)$  
and let $F\colon U\to U$ be the map defined by $F(u)=T(u)u$ for all
$u\in U$. Suppose that the stabilizer $St_{u^*}$
of some point $u^*\in U$ is non-trivial. 
There is a non-empty open $\Aut X$ invariant subset $V\subseteq U$ such that
\vskip0.1cm

\begin{itemize}
\item[$(*)$] $F((\Aut X)v)$ contains a subset biholomorphic
to the open disc ${\mathbb D}$
\end{itemize}
\vskip0.1cm

\noindent for each point $v\in V$
Moreover, for such a set $V$ we have 
$\dim_{\mathbb C}F(V)\ge n-t(X)+1$.
\end{Lemma}

\begin{proof}
Notice that $\dim_{\mathbb C}(\Aut X)u=\dim_{\mathbb C}\Aut X=t(X)$
for each $u\in U$; consider the holomorphic map
$$
f_u\colon\Aut X\ni A\mapsto T(Au)Au=F(Au)\in(\Aut X)u\,.
$$
Let us show that the differential
$d_Af_u\colon{\mathbb C}^{t(X)}\to{\mathbb C}^{t(X)}$
of $f_u$ at $A$ cannot be zero identically on $U\times(\Aut X)$.
Since $\Aut X$ is connected, the
assumption that $d_Af_u=0$ on $U\times(\Aut X)$ implies 
$f_u=\const$ for each $u\in U$. I. e.,
$T(Au)Au$ does not depend on $A$. Hence
$T(Au)Au=T(u)u$ and $(T(u))^{-1}T(Au)A\in St_u$
on $U\times(\Aut X)$. Since $St_u$ is finite,
$(T(u))^{-1}T(Au)A$ cannot depend on $A$
and hence $(T(u))^{-1}T(Au)A=(T(u))^{-1}T(u)=1$ for all
$(u,A)\in U\times(\Aut X)$. For $u=u^*$ there is a non-trivial
$A\in St_{u^*}$ and we have 
$1=(T(u^*))^{-1}T(Au^*)A=(T(u^*))^{-1}T(u^*)A=A$,
which is impossible.

Take $(u_0,A_0)\in U\times(\Aut X)$ such that $d_{A_0}f_{u_0}\ne 0$;
then $\rank d_{A_0}f_{u_0}\ge 1$ and, by semi-continuity,
$\rank d_{A}f_u\ge 1$ in some open neighborhood
$U_0\times{\mathcal A}_0$ of $(u_0,A_0)$.
The implicit map theorem implies that  
for each $u\in U_0$ the image $f_u(\Aut X)$ contains
a subset biholomorphic to ${\mathbb D}$.

Clearly $F((\Aut X)u)=f_u(\Aut X)$; hence each point $v$
of the non-empty open $\Aut X$ invariant set
$V=(\Aut X)U_0$ satisfies $(*)$, which implies 
$\dim_{\mathbb C}F(V)\ge n-t(X)+1$.
\end{proof}

\begin{Remark}\label{Rmk: more about the image dimension for X=C}
The statement of the above lemma may fail if
$St_u=\{1\}$ for all $u\in U$. Indeed, let
$u_0\in{\mathcal C}^n(X)$ and $St_{u_0}=\{1\}$; then
$St_u=\{1\}$ for all $u$ in a small open ball $B$ centered
at $u_0$. Let $D$ be a very small open complex disc of
dimension $n-t(X)$ contained in $B$, centered at $u_0$ and transversal
to the orbit $(\Aut X)u_0$ at $u_0$.
Then $[(\Aut X)v]\cap D=\{v\}$ for each $v\in D$
and the correspondence $(A,v)\mapsto Av$ gives rise to the
biholomorphic isomorphism 
$f\colon(\Aut X)\times D\stackrel{\cong}\longrightarrow (\Aut X)D
\Def U$. The composition $\tau=\pi_1\circ f^{-1}\colon U\stackrel{f^{-1}}
\longrightarrow (\Aut X)\times D\stackrel{\pi_1}
\longrightarrow\Aut X$ of the inverse isomorphism $f^{-1}$
and the projection $\pi_1$ to the first factor
is a holomorphic surjection, the map
$T\colon U\ni u\mapsto (\tau(u))^{-1}\in\Aut X$ is holomorphic,
and the map $F\colon U\ni u\mapsto T(u)u$ carries $U$ to
$D$. Indeed, for each $u\in U$, \ 
$u=f(A,v)=Av\in U$, we have $\tau(u)
=(\pi_1\circ f^{-1})(u)=\pi_1(f^{-1}(f(A,v)))
=\pi_1(A,v)=A$, \ $T(u)= (\tau(u))^{-1}=A^{-1}$ and, finally, 
$F(u)=T(u)u=A^{-1}Av=v\in D$ so that $\dim_{\mathbb C}F(U)=n-t(X)$.
\hfill $\bigcirc$
\end{Remark}

\noindent The linear map $L=(L_1,...,L_n)\colon{\mathcal C}^n({\mathbb C})
\to{\mathcal C}^n({\mathbb C})$ defined by
$$
L_i(\{q_1,...,q_n\})=q_i-(q_1+...+q_n)/n\,, \ \ i=1,...,n\,,
$$
carries the whole space ${\mathcal C}^n({\mathbb C})$ to its
$(n-1)$-dimensional submanifold 
${\mathcal C}^{n-1}_b({\mathbb C})$
consisting of all balanced configurations (see Definition
\ref{Def: balanced configuration space of C}).
\index{Balanced configuration spaces\hfill}
\index{${\mathcal C}^{m-1}_b({\mathbb C})$\hfill}
Thus, in the case $X={\mathbb C}$ the lower bound 
$\dim_{\mathbb C}F({\mathcal C}^n({\mathbb C}))
\ge n-t({\mathbb C})+1=n-1$ is the best possible one.
\vskip0.2cm

\noindent Let 
$L=(L_1,...,L_n)\colon{\mathcal C}^n({\mathbb C})
\to{\mathcal C}^{n-1}_b({\mathbb C})$ be the above linear map;
the following theorem provides a more precise result about images
of non-cyclic endomorphisms, which strengthens
Theorem \ref{Thm: Dimension of the image} in the case $X={\mathbb C}$.

\begin{Theorem}\label{Thm: image of non-cyclic endomorphism of Cn(C)}
\index{Image of non-cyclic endomorphism of ${\mathcal C}^n({\mathbb C})$\hfill}
Let $n>4$ and let $F\colon {\mathcal C}^n({\mathbb C})
\to{\mathcal C}^n({\mathbb C})$ be a non-cyclic holomorphic map.
Then the linear map $L$ caries the image $F({\mathcal C}^n({\mathbb C}))$
of $F$ on the whole ${\mathcal C}^{n-1}_b({\mathbb C})$; more precisely,
$$
L(F({\mathcal C}^n({\mathbb C})))
=L(F({\mathcal C}^{n-1}_b({\mathbb C})))
={\mathcal C}^{n-1}_b({\mathbb C})\,.
$$ 
\end{Theorem}

\begin{proof}
By Tame Map Theorem, $F(Q)=T(Q)Q$, where
$T\colon{\mathcal C}^n({\mathbb C})\to\Aff{\mathbb C}$ is a 
holomorphic map, that is, $T(Q)\colon{\mathbb C}\ni z
\mapsto A(Q)z+B(Q)\in{\mathbb C}$ with certain holomorphic functions
$A\colon{\mathcal C}^n({\mathbb C})\to{\mathbb C}^*$ and
$B\colon{\mathcal C}^n({\mathbb C})\to{\mathbb C}$. 
The cohomology class
$d^*_n\in H^1({\mathcal C}^n({\mathbb C}),{\mathbb Z})\cong{\mathbb Z}$
provided by the discriminant map $d_n\colon{\mathcal C}^n({\mathbb C})
\to{\mathbb C}^*$ generates $H^1({\mathcal C}^n({\mathbb C}),{\mathbb Z})$;
hence $A(Q)=(d_n(Q))^s\exp{\varphi(Q)}$, where $s\in{\mathbb Z}$ and
$\varphi$ is a holomorphic function on ${\mathcal C}^n({\mathbb C})$.
Take any
$Q^\circ=\{q_1^\circ,...,q_n^\circ\}
\in{\mathcal C}^{n-1}_b({\mathbb C})$ and consider the holomorphic
function $a\colon{\mathbb C}^*\ni\zeta\mapsto a(\zeta)
=\zeta A(\zeta Q^\circ)\in{\mathbb C}^*$. Then
$a(\zeta)=\zeta A(\zeta Q^\circ)=\zeta\cdot(d_n(\zeta Q^\circ))^s
\exp{\varphi(\zeta Q^\circ)}=\zeta^{sn(n-1)+1}(d_n(Q^\circ))^s
\exp{\varphi(\zeta Q^\circ)}$, which shows that $a\ne\const$.
Since $a$ does not take the value $0$, it follows from the
Picard theorem that $a(\zeta^\circ)=1$ for some $\zeta^\circ\in{\mathbb C}^*$.
A straightforward computation shows that
$G_i(\zeta^\circ Q^\circ)=L_i(F(\zeta^\circ Q^\circ)=
\zeta^\circ A(\zeta^\circ Q^\circ)q_i^\circ=a(\zeta^\circ)q_i^\circ
=q_i^\circ$ for all $i=1,...,n$, that is,
$G(\zeta^\circ Q^\circ)=Q^\circ$, which proves the theorem.
\end{proof}

\begin{Remark}\label{Rmk: no non-cyclic endomorphism F with dim of
image < n is known}
I do not know any example of a non-cyclic holomorphic endomorphism
$F$ of ${\mathcal C}^n({\mathbb{CP}}^1)$ such that 
$\dim_{\mathbb C}F({\mathcal C}^n({\mathbb{CP}}^1))<n$.
\hfill $\bigcirc$
\end{Remark}

\section{Some applications to algebraic functions}
\label{Sec: Some applications to algebraic functions}

\noindent The $n$-valued algebraic function $[x:y]=U_n([z])$ 
defined by the equation
\begin{equation}\label{eq: universal algebraic equation}
\aligned
&\phi_n(x,y;z)=z_0 x^n + z_1 x^{n-1}y +...+ z_n y^n=0\,, \\ 
&z=(z_0,...,z_n)\in{\mathbb C}^{n+1}\,, \ \ 
[z]=[z_0:...:z_n]\in{\mathbb{CP}}^n\,,
\endaligned
\end{equation}
and the $n$-valued {\em entire} algebraic function $t=u_n(w)$ 
defined by the equation
\begin{equation}\label{eq: universal entire algebraic equation}
p_n(t,w)=t^n+w_1t^{n-1}+...+w_n=0\,, \ \
w = (w_1,...,w_n)\in{\mathbb C}^n\,,
\end{equation}
are said to be the {\em universal} ones. It is easily seen
that any $n$-valued algebraic function on ${\mathbb{CP}}^m$
(respectively any $n$-valued entire algebraic function on
${\mathbb C}^m$) is a composition of a rational map
${\mathbb{CP}}^m\to{\mathbb{CP}}^n$ (respectively a polynomial
map ${\mathbb C}^m\to{\mathbb C}^n$) and the corresponding
universal function. The branch loci of the universal functions
are also called universal ones.
\vskip0.2cm 
 
\noindent In this section we study algebraic functions
whose branch loci coincide with a universal one
and apply the obtained results to a version of the 13th Hilbert
problem for algebraic functions.  

\subsection{Algebraic functions with the universal branch locus}
\label{Ss: Algebraic functions with the universal branch locus}
To avoid the doubling related to our wish to consider both ``general"
and ``entire" algebraic functions, we will use the following
notation.

\begin{Notation}
\label{Not: algebraic functions on Cn(X)} 
We consider the following two cases:   

\begin{itemize}
\item[$(${\bfit a}$)$] $X={\mathbb{CP}}^1$,
$X^n/{\mathbf S}(n)={\mathbb{CP}}^n$,
${\mathcal U}_n=U_n$ is the universal $n$-valued algebraic function
on ${\mathbb{CP}}^n$ defined by (\ref{eq: universal algebraic equation})
with the branch locus
$$
\Sigma_{{\mathcal U}_n}=\Sigma_{U_n}=\{[z]\in{\mathbb{CP}}^n\,|\
D_{\phi_n}([z])=0\}\,,
\eqnum{a'}
$$
and 
$(X^n/{\mathbf S}(n))\setminus\Sigma_{{\mathcal U}_n}
={\mathbb{CP}}^n\setminus\Sigma_{{\mathcal U}_n}
\cong{\mathcal C}^n(X)={\mathcal C}^n({\mathbb{CP}}^1)$\,;
\end{itemize}
\vskip0.2cm

\begin{itemize}
\item[$(${\bfit b}$)$] $X={\mathbb C}$,
$X^n/{\mathbf S}(n)={\mathbb C}^n$,
${\mathcal U}_n=u_n$ is the universal $n$-valued entire algebraic
function on ${\mathbb C}^n$ defined by
(\ref{eq: universal entire algebraic equation}) with the branch locus
$$
\Sigma_{{\mathcal U}_n}=\Sigma_{u_n}
=\{w\in{\mathbb C}^n\,|\ d_n(w)=0\}\,,
\eqnum{b'}
$$
and 
$(X^n/{\mathbf S}(n))\setminus\Sigma_{{\mathcal U}_n}
={\mathbb C}^n\setminus\Sigma_{{\mathcal U}_n}
\cong{\mathcal C}^n(X)={\mathcal C}^n({\mathbb C}^1)$.
\end{itemize}
\end{Notation}

\noindent In case $(${\bfit a}$)$ the term ``an $m$-valued
algebraic function ${\mathcal F}$ on $X^n/{\mathbf S}(n)$"
will mean an algebraic function $[x:y]={\mathcal F}([z])$
on ${\mathbb{CP}}^n$ defined by an equation
$$
\varphi(x,y;z)=a_0(z)x^m + a_1(z)x^{m-1}y +...+ a_m(z)y^m=0\,,
$$
where $\varphi(x,y;z)=a_0(z)x^m + a_1(z)x^{m-1}y +...+ a_m(z)y^m$ is
a {\em homogeneous} binary form whose coefficients $a_0,...,a_m$
are coprime {\em homogeneous} polynomials in $z_0,...,z_n$ {\em of the same
positive degree}.\footnote{In this section, all binary forms under
consideration are supposed to be of such kind. It is possible that some
coefficients $a_i=0$.}
The form $\varphi(x,y;z)$ is called the {\em defining form}
of the algebraic function $F$. 
The zero set $\Sigma_\varphi=\{z=(z_0,...,z_n)\in{\mathbb C}^{n+1}\,|\
D_\varphi(z)=0\}$ of the discriminant $D_\varphi(z)$ is said to be 
{\em degeneration locus} of the form $\varphi(x,y;z)$; 
in particular, $\Sigma_\varphi$ contains all common zeros of
the coefficients $a_0,...,a_m$. The quotient space
$\Sigma_F=\Sigma_\varphi/{\mathbb C}^*\subset{\mathbb{CP}}^n$ is 
the {\em branch locus} of the function $F$ (it may happen
that $D_\varphi=0$ and $\Sigma_F={\mathbb{CP}}^n$).
\vskip0.2cm 

\noindent For two algebraic functions $[x:y]=G([z])$
and $[x:y]=F([z])$ with the defining forms $\psi(x,y;z)$ and
$\varphi(x,y;z)$ we say that $G$ is a  (multivalued) {\em branch} of $F$
and write $G\subseteq F$ if there is a form $\omega(x,y;z)$ such that
$\varphi(x,y;z)=\psi(x,y;z)\cdot\omega(x,y;z)$.
We write $F=G$ if $\varphi(x,y;z)=c\psi(x,y;z)$
for some $c\in{\mathbb C}^*$.
\vskip0.2cm

\noindent Let us notice that if $t=f(w)$ is an algebraic function
of $w\in{\mathbb C}^n$ defined by an algebraic
equation $P(t,w)=A_0(w)t^m + A_1(w)t^{m-1} +...+A_m(w)=0$
with coefficients 
$A_i\in{\mathbb C}[w]$ then one may replace $t$ with $x/y$,
each of $w_i$ with $z_i/z_0$, and rewrite the function as
$[x:y]=F([z])$ and the defining equation as $\varphi(x,y;z)=0$,
where $\varphi(x,y;z)$ is a binary form of the type described above.
\vskip0.2cm
  
\noindent In case $(${\bfit b}$)$ the term ``an $m$-valued
algebraic function ${\mathcal F}$ on $X^n/{\mathbf S}(n)$"
will mean an entire algebraic function $t={\mathcal F}(w)$
on ${\mathbb C}^n$ defined by an equation
$$
P(t;w)=t^m + a_1(w)t^{m-1} +...+ a_m(w)=0\,,
\eqnum{b'}
$$
where all $a_i\in{\mathbb C}[w]$.
For two entire algebraic functions $t=g(w)$ and $t=f(w)$
with the defining polynomials $Q(t,w)$ and $P(t,w)$ we
say that $g$ is a (multivalued) {\em branch} of $f$ and write $g\subseteq f$
if there is a polynomial $R(t,w)$ such that $P=Q\cdot R$.
We write $g=f$ if $P=Q$.  
\vskip0.2cm
  
\noindent The following result is in fact equivalent to Tame Map Theorem
for regular endomorphisms.

\begin{Theorem}
\label{Thm: n-valued algebraic functions with universal branch locus}
\index{Theorem!on $n$-valued algebraic functions with the
universal branch locus\hfill}
Let $n>4$ and let ${\mathcal F}$ be an $n$-valued algebraic function on
$X^n/{\mathbf S}(n)$ with the branch locus
$\Sigma_{\mathcal F}=\Sigma_{{\mathcal U}_n}$
and a non-cyclic monodromy representation 
\begin{equation}\label{eq: monodromy representation}
\pi_1((X^n/{\mathbf S}(n))\setminus\Sigma_{{\mathcal U}_n})
\to{\mathbf S}(n)\,.
\end{equation}
Then there exists a rational map $T\colon X^n/{\mathbf S}(n)\to\Aut X$
regular on $(X^n/{\mathbf S}(n))\setminus\Sigma_{{\mathcal U}_n}$
such that ${\mathcal F}(z)=T(z){\mathcal U}_n(z)$
for all $z\in X^n/{\mathbf S}(n)$.
\end{Theorem}

\begin{proof}
The algebraic function ${\mathcal F}$ is the composition
${\mathcal F}(z)={\mathcal U}_n(A(z))$ of the rational map
$A\colon X^n/{\mathbf S}(n)\to X^n/{\mathbf S}(n)$ given by the
coefficients of the equation defining ${\mathcal F}$ and
the universal function ${\mathcal U}_n$. The assumption
$\Sigma_{\mathcal F}=\Sigma_{{\mathcal U}_n}$ implies that  
$A$ is regular on $(X^n/{\mathbf S}(n))\setminus\Sigma_{{\mathcal U}_n}$
and carries the latter set into itself. Thus, under the
identification of the defining equations to the sets of their roots
established in Remark \ref{Rmk: non-degenerate forms}, $A$ may be regarded as
a regular endomorphism $A'$ of the configuration space ${\mathcal C}^n(X)$.
Since the monodromy representation (\ref{eq: monodromy representation})
is non-cyclic, $A'$ is non-cyclic as well.
By Tame Map Theorem, there is a regular map
$T'\colon{\mathcal C}^n(X)\to\Aut X$ such that
$A'(Q)=T'(Q)Q$ for all $Q\in{\mathcal C}^n(X)$.
Denoting by $T$ the regular morphism
$(X^n/{\mathbf S}(n))\setminus\Sigma_{{\mathcal U}_n}\to\Aut X$
corresponding to $T'$, we obtain that
${\mathcal F}(z)\equiv {\mathcal U}_n(A(z))\equiv T(z){\mathcal U}_n(z)$
meaning that each value of ${\mathcal F}$ at
$z\in(X^n/{\mathbf S}(n))\setminus\Sigma_{{\mathcal U}_n}$ is obtained
from the corresponding value of ${\mathcal U}_n$ at $z$ with help of the
M{\" o}bius transformation $T(z)$.
\end{proof}

\noindent The following result is just a simple consequence
of Linked Map Theorem (for $X={\mathbb C}$ it has been proved
in \cite{Lin96b}). 

\begin{Theorem}
\label{Thm: algebraic functions containing the universal one}
\index{Algebraic functions with the universal branch locus\hfill}
\index{Theorem!on algebraic functions containing the universal one\hfill}
Let $n>t(X)+1$ and let ${\mathcal F}$ be an $N$-valued algebraic function
on $X^n/{\mathbf S}(n)$. If ${\mathcal U}_n\subseteq{\mathcal F}$ 
then either $N = n$ and ${\mathcal F}={\mathcal U}_n$ or
$N>n$ and ${\mathcal F}$ has branch points that are not branch
points of ${\mathcal U}_n$.
\end{Theorem}

\begin{proof}
We consider here only the case $X={\mathbb{CP}}^1$. 
Let $\varphi(x,y;z)$ be the defining form of the function $F$.
The condition $U_n\subseteq F$ means that there is a binary form
$\omega(x,y;z)$ of degree $k=N-n\ge 0$
such that $\varphi(x,y;z)=\phi_n(x,y;z)\cdot\omega(x,y;z)$;
in this case $D_{\phi_n}(z)=0$ implies $D_\varphi(z)=0$, which means that
$\Sigma_F\supseteq\Sigma_{U_n}$.
\vskip0.15cm

\noindent Assume that $k>0$ and $\Sigma_F=\Sigma_{U_n}$.
Then for every
$z\in{\mathbb{CP}}^n\setminus\Sigma_{U_n}
={\mathcal C}^n({\mathbb{CP}}^1)$ the form $\varphi(x,y;z)$
is non-degenerate; therefore $\omega(x,y;z)$ is non-degenerate
as well and has no common roots with $\phi_n(x,y;z)$.
This means that the correspondence $z\mapsto\omega(x,y;z)$
defines a disjoint morphism 
${\mathcal C}^n({\mathbb{CP}}^1)\cong{\mathcal F}^n
\to{\mathcal F}^k\cong{\mathcal C}^k({\mathbb{CP}}^1)$,
which contradicts Linked Map Theorem \ref{LinkMapThm}.
Thus, $\Sigma_{U_n}$ must be a proper subset of 
$\Sigma_F$, which shows that $F$ has branch points that are not
branch points of $U_n$.
\vskip0.15cm

\noindent Finally, if $k=0$ then $N=n$, $\omega=\const\ne 0$ and $F=U_n$.  
\end{proof}
  
\noindent The problem about existence of an algebraic function with
the preassigned branch locus and monodromy group goes back
to B. Riemann (in one variable), F. Enriques and O. Zariski
(for functions of two variables; see \cite{Enr23,Zar29}).
It was discovered that in two variables there are
obstacles related to the fundamental group of the complement of the
branch curve. 
\vskip0.2cm

\noindent We deal with very distinguished subvarieties
$\Sigma_{U_n}\subset{\mathbb{CP}}^n$ and 
$\Sigma_{u_n}\subset{\mathbb C}^n$ and we know
in advance that there exist $n$ valued algebraic functions, namely
the universal ones $U_n$ and $u_n$,
with these very branch loci and the standard monodromy representations 
$$
\aligned
&\mu\colon\pi_1({\mathbb{CP}}^n\setminus\Sigma_{U_n})
=\pi_1({\mathcal C}^n({\mathbb{CP}}^1))\cong B_n(S^2)
\to{\mathbf S}(n)\,,\\
&\mu\colon\pi_1({\mathbb C}^n\setminus\Sigma_{u_n})
=\pi_1({\mathbf G}_n)
=\pi_1({\mathcal C}^n({\mathbb C}))\cong B_n\to{\mathbf S}(n)\,.
\endaligned
$$  
Theorem 
\ref{Thm: n-valued algebraic functions with universal branch locus}
provides a reasonable description of all non-cyclic $n$-valued algebraic
functions of $n$ variables $z\in X^n/{\mathbf S}(n)$ with the universal
branch locus $\Sigma_{{\mathcal U}_n}$ in terms of the universal function
${\mathcal U}_n(z)$ and automorphisms of the target domain
$X={\mathbb C}$ or $X={\mathbb{CP}}^1$ depending regularly on a point
$z\in (X^n/{\mathbf S}(n))\setminus\Sigma_{{\mathcal U}_n}$. 
On the other hand, Theorem 
\ref{Thm: algebraic functions containing the universal one}
tells that the universal function cannot be a proper branch of
an algebraic function with the same branch locus.


\subsection{Compositions of multivalued functions}
\label{Ss: Compositions of multivalued functions}
\index{Tame algebraic ensembles\hfill}
\index{Compositions of multivalued functions\hfill}

\noindent One possible form of the question posed in
the 13th Hilbert problem is as follows:
\vskip0.2cm

\noindent What is the minimal natural number $m$ for which the
function $u_n$ {\em can be represented as a composition of algebraic
functions of $m$ variables}?
\vskip0.2cm

\noindent Since $u_n$ is multivalued,
the emphasized expression may have quite different meanings.
It has been proved in \cite{Arn70c,Lin76,Lin96a}
that for certain possible interpretations of the problem the answer
is $m=n-1$. Here we discuss the corresponding problem
for the function $[x:y]=U_n([z])$ and show that in this case
$m\ge n-2$.
\vskip0.2cm

\noindent Following the papers quoted above, we explain first some
formalities related to multivalued maps; we do this in a manner
convenient for working with non-entire algebraic functions.
\vskip0.2cm

\noindent Given two topological spaces $Y$ and $Z$, we use the notation
$f\colon Z\stackrel{n}\longrightarrow Y$ for a map
defined and continuous on an open dense subset $Z_f\subseteq Z$
and carrying it into the space $\Sym^m Y=Y^m/{\mathbf S}(m)$;
such a map is said to be an {\em $m$-valued map of $Z$ to} $Y$.
The points $y_1,...,y_m\in f(z)$ are said to be the
{\em values of $f$ at} $z\in Z_f$. The subset $\Img f\subseteq Y$
consisting of all values of $f$ at all points $z\in Z_f$ is called
the {\em image of} $f$. The subset $\Gamma_f\subseteq Z\times Y$
consisting of all pairs $(z,y)$ such that $y$
is a value of $f$ at $z\in Z_f$ is called the {\em graph} of $g$.
\vskip0.2cm

\noindent Let $f\colon Z\stackrel{m}\longrightarrow Y$ and
$g\colon Y\stackrel{n}\longrightarrow W$
be multivalued maps; suppose that there are open dense 
subsets $Z'\subseteq Z_f$ and $Y'\subseteq Y_g$ such that 
all values of $f$ at any point $z\in Z'$ belong to $Y'$.
The {\em composition}
$$
g\circ f\colon Z\stackrel{mn}\longrightarrow W
$$
is the $(mn)$-valued map whose values at a point $z\in Z'$
are all $mn$ values of $g$ at the $m$ points of $Y'$
which are the values of $f$ at $z$.
The {\em complete composition}
$$
g*f\colon Z\stackrel{mn}\longrightarrow Y\times W
$$
is the $(mn)$-valued map whose values at a point $z\in Z'$
are all pairs $(y,w)$ such that $y$ is one of the values
of $f$ at $z$ and $w$ is one of the values of $g$ at $y$.
\vskip0.2cm

\noindent Let $k,m,M$ be natural numbers. Assume that for each
$i=1,\ldots ,M$ we are given with an $n_i$-valued algebraic function
$[x:y]=g_i([w])$, $[w]=[w_0:...:w_m]\in{\mathbb{CP}}^m$,
defined by an equation $\varphi_i(x,y;[w])
=a_0(w)x^{n_i}+a_1(w)x^{n_i-1}y+...+a_{n_i}(w)y^{n_i}=0$,
where $a_0,...,a_{n_i}$ are co-prime homogeneous polynomials of the
same positive degree.
Assuming that the discriminant $D_{\varphi_i}(w)$ of $\varphi_i$
is not identically zero, we may regard $g_i$ as an $n_i$-valued map 
$
g_i\colon{\mathbb{CP}}^m\stackrel{n_i}\longrightarrow{\mathbb{CP}}^1\,,
$
which is surely defined on the non-empty Zariski open set
$W_i=\{w\in{\mathbb{CP}}^m\,|\ D_{\varphi_i}(w)\ne 0\}$.
\vskip0.2cm

\noindent Furthermore, for each $i=1,\ldots ,M$ let
$R_i\colon{\mathbb{CP}}^k\times({\mathbb{CP}}^1)^{i-1}\to{\mathbb{CP}}^m$ 
be a rational map such that its regularity domain
contains a non-empty Zariski open subset $V_i$ with the image
$R_i(V_i)\subseteq W_i$.
\vskip0.2cm

\noindent Let $F_0\colon{\mathbb{CP}}^k\to{\mathbb{CP}}^k$
denote the identity map. Using the compositions $\circ$ and $*$
of multivalued functions defined above, we construct algebraic functions 
$$
G_i = g_i\circ R_i\colon
{\mathbb{CP}}^k\times({\mathbb{CP}}^1)^{i-1}
\stackrel{R_i}\longrightarrow {\mathbb{CP}}^m
\stackrel{n_i}\longrightarrow{\mathbb{CP}}^1\,, \quad i=1,...,M\,,
$$

$$
\aligned
&
\CD
F_1=G_1*F_0\colon{\mathbb{CP}}^k @>{n_1}>> 
    {\mathbb{CP}}^k\times{\mathbb{CP}}^1\,,
\endCD\\
&
\CD
F_2=G_2*F_1\colon{\mathbb{CP}}^k @>{n_1n_2}>>
           {\mathbb{CP}}^k\times({\mathbb{CP}}^1)^2\,,
\endCD\\
& \ ...................................................................\\
&
\CD
F_i=G_i*F_i\colon{\mathbb{CP}}^k @>{n_1n_2\cdots n_i}>>
                  {\mathbb{CP}}^k\times({\mathbb{CP}}^1)^i\,,
\endCD\\
& \ ...................................................................\\
&
\CD
F_{M-1}=G_{M-1}*F_{M-2}\colon{\mathbb{CP}}^k @>{n_1n_2\cdots n_{M-1}}>>
        {\mathbb{CP}}^k\times({\mathbb{CP}}^1)^{M-1}\,,
\endCD\\
\endaligned
$$
and, finally,
$$
\CD
{F=F_M=G_M\circ F_{M-1}\colon{\mathbb{CP}}^k}
                       @>{n_1n_2\cdots n_M}>>
                               {\mathbb{CP}}^1\,.
\endCD
$$
$F$ is a well-defined $N=n_1\cdots n_M$-valued
algebraic function on ${\mathbb{CP}}^k$ whenever
$\Img F_i\subseteq V_{i+1}$ for all $i=1,..., M-1$.
If this is the case then the collection
${\mathfrak A}=\{V_i,W_i,R_i,g_i\,|\ i=1,..., M\}$ is said
to be a {\em tame algebraic $(k,m)$-ensemble in}  
${\mathbb{CP}}^k$, and the algebraic function
$F_{\mathfrak A}=F$ is called the {\em function generated by the ensemble}
${\mathfrak A}$. Two tame algebraic $(k,m)$-ensembles ${\mathfrak A}$ and
${\mathfrak A}^\prime$ in ${\mathbb{CP}}^k$
are said to be {\em equivalent} if the functions $F_{\mathfrak A}$
and $F_{{\mathfrak A}^\prime}$ generated by them coincide.
\vskip0.2cm

\noindent Some of the algebraic functions $g_i$ appearing
in a tame ensemble ${\mathfrak A}$ may be single-valued; the number
of those ones which are actually multivalued is called
the {\em rank of the ensemble} ${\mathfrak A}$.
 
\begin{Example} For $M=2$ we have
$F(z) = g_2\{R_2[z,g_1(R_1(z))]\}$. If ${\mathfrak A}$ is an
ensemble of rank $1$ then $F_{\mathfrak A}=R_2(z,g(R_1(z))$,
where $[x:y]=g(w)$ is a multivalued algebraic function and
$R_1,R_2$ are rational.
\hfill $\bigcirc$
\end{Example}

\begin{Definition}\label{Def: representation as a composition}
We say that an algebraic function $g$ in ${\mathbb{CP}}^k$
{\em can be represented as a composition 
of algebraic functions of $m$ variables} 
if there exists a tame algebraic $(k,m)$-ensemble ${\mathfrak A}$
such that $g\subseteq F_{\mathfrak A}$, that is, $g$
is a (multivalued) branch of $F_{\mathfrak A}$.
The function $F_{\mathfrak A}$ is called a {\em representing
function-composition for} $g$; clearly
$\Sigma_g\subseteq\Sigma_{F_{\mathfrak A}}$.
If the ensemble ${\mathfrak A}$
can be chosen so that $g = F_{\mathfrak A}$ we say that $g$
can be {\em faithfully} represented as a composition of
algebraic functions of $m$ variables; in this case
$\Sigma_g=\Sigma_{F_{\mathfrak A}}$. In general, the 
conditions $g\subseteq F_{\mathfrak A}$ and
$\Sigma_g=\Sigma_{F_{\mathfrak A}}$ do not imply $g = F_{\mathfrak A}$.
\hfill $\bigcirc$
\end{Definition}

\noindent In the following theorem we deal with a representation
of the universal algebraic function $U_n$ as a composition
$U_n\subseteq F_{\mathfrak A}$, which a priori is not supposed to be
faithful but satisfies the formally weaker condition
$\Sigma_{U_n}=\Sigma_{F_{\mathfrak A}}$. We will see however that
the conditions $U_n\subseteq F_{\mathfrak A}$ and 
$\Sigma_{U_n}=\Sigma_{F_{\mathfrak A}}$ imply $U_n=F_{\mathfrak A}$.

\begin{Theorem}\label{Thm: compositions without superflous branch points}
\index{Theorem!on compositions without superfluous branch points\hfill}
For $n>\max\{4,m+2\}$ the universal $n$-valued algebraic function
$[x:y]=U_n([z])$ cannot be represented as a composition
of algebraic functions of $m$ variables in such
a way that a representing function-composition $F=F_{\mathfrak A}$
and $U_n$ itself have the same branch loci.
\end{Theorem}

\begin{proof}
Suppose to the contrary that there exists a tame algebraic
$(n,m)$-ensemble ${\mathfrak A}=\{U_i,V_i,R_i,g_i\,|\ i=1,..., M\}$ in
${\mathbb{CP}}^n$ such that $U_n\subseteq F=F_{\mathfrak A}$
and the branch locus $\Sigma_F$ of the $N$-valued
algebraic function-composition $F=F_{\mathfrak A}$ defined by
${\mathfrak A}$ coincides with the branch locus $\Sigma_{U_n}$.
Then, by Theorem
\ref{Thm: algebraic functions containing the universal one},
$N=n$ and $U_n=F$. In particular, the monodromy group of
$F$ coincides with the monodromy group
${\mathbf S}(n)$ of the function $U_n$. This implies that the unbranched
$n$-covering
$p\colon\Gamma_F\to{\mathbb{CP}}^n\setminus\Sigma_{U_n}$ defined by
the function $F$ is {\em indecomposable}, that is,
the projection $\pi$ cannot be represented as a composition
$$
p=q_1\circ q_2\colon\Gamma_F\stackrel{q_2}{\longrightarrow}
E\stackrel{q_1}{\longrightarrow}{\mathbb{CP}}^n\setminus\Sigma_{U_n}
$$
of two non-univalent unbranched
coverings $q_1$,$q_2$ (see \cite{Arn70c,Lin76,Lin96a}).
This property implies that the ensemble ${\mathfrak A}$
must be of rank 1, that is, it cannot contain more than one multivalued
algebraic function and thereby the corresponding representation of the function $U_n$
as a superposition must be of the form $U_n(z)=F(z)=R_2(z,g(R_1(z))$, where
$[x:y]=g(w)$ is an $n$-valued algebraic function
on ${\mathbb{CP}}^m$ with a defining binary form
\begin{equation}\label{eq: defining form of g}
\varphi(x,y;w)=a_0(w)x^n+a_1(w)x^{n-1}y+...+a_n(w)y^n\,,
\end{equation}
$R_1\colon{\mathbb{CP}}^n\to{\mathbb{CP}}^m$ is a rational map
regular on ${\mathbb{CP}}^n\setminus\Sigma_{U_n}$ with
the image
$$
R_1({\mathbb{CP}}^n\setminus\Sigma_{U_n})
\subseteq \{w\in{\mathbb{CP}}^m\,|\ D_\varphi(w)\ne 0\}\,,
$$
and $R_2\colon{\mathbb{CP}}^n\times{\mathbb{CP}}^1\to{\mathbb{CP}}^1$
is a rational function regular at each point $(z,[x:y])\in
({\mathbb{CP}}^n\setminus\Sigma_{U_n})\times{\mathbb{CP}}^1$
such that $\varphi(x,y;R_1(z))=0$.
Let $A =(a_1,\ldots ,a_n)\colon{\mathbb{CP}}^m\to{\mathbb{CP}}^n$
be the rational map defined by the coefficients of
the binary form (\ref{eq: defining form of g}).
The identity $U_n(z)\equiv R_2(z,g(R_1(z))$ shows that
$A$ must be regular on $\{w\in{\mathbb{CP}}^m\,|\ D_\varphi(w)\ne 0\}$
and the $A$-image $A(R_1({\mathbb{CP}}^n\setminus\Sigma_{U_n}))$
of the set $R_1({\mathbb{CP}}^n\setminus\Sigma_{U_n})$
must be contained in ${\mathbb{CP}}^n\setminus\Sigma_{U_n}$.
Thus, the composition $f=A\circ R_1$,
$$
f\colon{\mathcal C}^n({\mathbb{CP}}^1)=
{\mathbb{CP}}^n\setminus\Sigma_{U_n}
\stackrel{R_1}\longrightarrow
\{w\in{\mathbb{CP}}^m\,|\ D_\varphi(w)\ne 0\}
\stackrel{A}\longrightarrow{\mathbb{CP}}^n\setminus\Sigma_{U_n}
={\mathcal C}^n({\mathbb{CP}}^1)\,,
$$
is a regular map. This map is certainly non-cyclic, for
otherwise the monodromy group of the algebraic function
$g((R_1(z))$ would be cyclic,
as well as the monodromy group of the function $R_2(z,g(R_1(z))=U_n(z)$.
Since $m<n-2$, we have 
$\dim_{\mathbb C}(f({\mathcal C}^n({\mathbb{CP}}^1)))\le m<n-2$,
which contradicts Theorem \ref{Thm: Dimension of the image}.
\end{proof}

\begin{Remark}\label{Rmk: compositions of entire algebraic functions}
{\bfit a}$)$ The corresponding result for the universal entire
algebraic function $u_n$ on ${\mathbb C}^n$ has been proved in
\cite{Lin96a}:

\begin{Theorem}\label{Thm: compositions of entire algebraic functions}
\index{Theorem!on compositions without superfluous branch points\hfill}
For $n>\max\{4,m+1\}$ the universal $n$-valued entire algebraic function
$t=u_n(z)$ cannot be represented as a composition
of entire algebraic functions of $m$ variables in such
a way that a representing function-composition $F=F_{\mathfrak A}$
and $u_n$ itself have the same branch loci.
\end{Theorem}

\noindent The proof is also based on Theorems \ref{LinkMapThm},
\ref{Thm: Dimension of the image} and 
\ref{Thm: algebraic functions containing the universal one}, which
in the case $X={\mathbb C}$ have been proved in the paper quoted above. 
\vskip0.2cm

{\bfit b}$)$ The compositions of {\em entire} algebraic functions treated
in Theorem \ref{Thm: compositions of entire algebraic functions} 
are defined similarly to what has been done above. Just 
in the definition of an algebraic $(k,m)$-ensemble one should consider
affine spaces instead of projective ones and work with polynomials instead
of rational functions; that is, {\em the operation of division is
in fact prohibited}. V. I. Arnold many times expressed his conviction that 
{\sl whenever rational functions are admitted $u_n$ cannot be
faithfully represented as a composition of algebraic functions of less
than $n-2$ variables}. Seemingly, Theorem
\ref{Thm: compositions without superflous branch points} confirms this
Arnold's opinion (and even for a formally wider class of representations
{\em without superfluous branch points} instead of the {\em faithful}
ones).
\hfill $\bigcirc$
\end{Remark}

\section{Automorphisms}
\label{Sec: Automorphisms}
\index{Automorphisms of ${\mathcal C}^n(X)$\hfill}

\subsection{How to catch automorphisms?}
\label{Ss: How to catch automorphisms}
By Tame Map Theorem, every non-cyclic
endomorphisms $F$ of ${\mathcal C}^n(X)$ is of the form $F=F_T$,
$F_T(Q)=T(Q)Q$, where $T\colon{\mathcal C}^n(X)\to\Aut X$
is a morphism uniquely defined by $F$. It is reasonable
to ask the following question: {\em for which morphism
$T\colon{\mathcal C}^n(X)\to\Aut X$ the corresponding
tame endomorphism $F_T$ is an automorphism?}
It is easy to give the following ``answer":
\vskip0.2cm

\noindent{\sl $F_T$ is an automorphism if and only if there exists
a morphism $S\colon{\mathcal C}^n(X)\to\Aut X$ such that}
$T(S(Q)Q)S(Q)=S(T(Q)Q)T(Q)=I$ (the unity of the group $\Aut X$)
{\sl for all $Q\in{\mathcal C}^n(X)$.}
\vskip0.2cm

\noindent It is unclear how ensure that $T$ satisfies such a condition.
It is naturally to look for a relationship between
the group multiplication in $\Aut X$ and the composition of
tame maps. In general, these operations know nothing about
each other; however in the special case of $\Aut X$ {\em invariant}
maps ${\mathcal C}^n(X)\to\Aut X$ they are in a close connection.  

\begin{Definition}\label{Def: Aut X invariant maps}
\index{Invariant maps to $\Aut X$\hfill}
As usually, the pointwise product
$ST\colon{\mathcal C}^n(X)\to\Aut X$ of two 
maps $S,T\colon{\mathcal C}^n(X)\to\Aut X$ is defined by
$(ST)(Q)=S(Q)T(Q)$, where $S(Q)T(Q)$
denotes the group product of the elements $S(Q),T(Q)\in\Aut X$.
\vskip0.2cm

\noindent A map $S\colon{\mathcal C}^n(X)\to\Aut X$
is said to be $\Aut X$ {\em invariant} if $S(AQ)=S(Q)$ for all
$A\in\Aut X$ and $Q\in{\mathcal C}^n(X)$.
\hfill $\bigcirc$
\end{Definition}

\begin{Lemma}\label{Lm: Aut X invariant maps}
Let $F_S, F_T\colon{\mathcal C}^n(X)\to {\mathcal C}^n(X)$
be the tame maps induced by two maps $S,T\colon {\mathcal C}^n(X)\to \Aut X$. 
If $S$ is $\Aut X$ invariant then $F_{ST} = F_S\circ F_T$.
\end{Lemma}

\begin{proof}
It suffices to know that $S(Q)=S(T(Q)Q)$ $\forall\,Q=\{q_1,...,q_n\}
\in{\mathcal C}^n(X)$; indeed, 
$\displaystyle F_{ST}(Q)
=\left\{\big((ST)(Q)\big)q_1,...,\big((ST)(Q)\big)q_n\right\}
=\left\{S(Q)T(Q)q_1,...,S(Q)T(Q)q_n\right\}
=\left\{S(Q)Aq_1,...,S(Q)Aq_n\right\}
=\left\{S(AQ)Aq_1,...,S(AQ)Aq_n\right\}
=F_S(AQ)=F_S(T(Q)Q))=F_S(F_T(Q))$,
that is, $F_{ST}=F_S\circ F_T$.
\end{proof}

\begin{Corollary}\label{Crl: Aut X invariance and automorphisms}
\index{Theorem!on invariant morphisms to $\Aut X$ and
automorphisms of ${\mathcal C}^n(X)$\hfill}
The tame map $F_T\colon{\mathcal C}^n(X)\to{\mathcal C}^n(X)$ 
induced by an $\Aut X$ invariant morphism
$T\colon {\mathcal C}^n(X)\to\Aut X$
is an automorphism of ${\mathcal C}^n(X)$.
\end{Corollary}

\begin{proof}
Set $S=T^{-1}$, that is, $S(Q)=(T(Q))^{-1}$ for all
$Q\in{\mathcal C}^n(X)$. Then $S$ is an $\Aut X$ invariant morphism and, 
by Lemma \ref{Lm: Aut X invariant maps},
$F_S\circ F_T=F_{ST}=F_I=\Id_{{\mathcal C}^n(X)}=F_{TS}=F_T\circ F_S$.
\end{proof}

\begin{Remark}\label{Rmk: weaker condition}
If $T$ satisfies $T((T(Q))^{-1}Q)=T(Q)$ for all
$Q\in{\mathcal C}^n(X)$ then $F_T$ is an automorphism
whose inverse is $F_{T^{-1}}$.
\hfill $\bigcirc$
\end{Remark}

\subsection{Some examples}\label{Ss: Some examples}
Corollary \ref{Crl: Aut X invariance and automorphisms} gives
rise to some simple examples.

\begin{Example}
\label{Ex: biregular automorphisms generated by invariant morphisms}
Take $c\in{\mathbb C}^*$ and an $\Aff{\mathbb C}$ invariant
regular function $b\colon{\mathcal C}^n({\mathbb C})\to{\mathbb C}$;
then the map
\begin{equation}\label{eq:biregular automorphism of Cn(C)}
F\colon{\mathcal C}^n({\mathbb C})\ni Q\mapsto cQ
+b(Q)\cdot\1\in{\mathcal C}^n({\mathbb C})\,,
\quad \1=\{1,...,1\}\in{\mathbb C}^n\,,
\end{equation}
is a biregular automorphism of ${\mathcal C}^n({\mathbb C})$.
The inverse map is given by
$$
{\mathcal C}^n({\mathbb C})\ni Q\mapsto c'Q+b'(Q)\cdot\1
\in{\mathcal C}^n({\mathbb C})\,,
$$
where $c'=c^{-1}$ and $b'(Q)=-c^{-1}b(Q)$. For instance, one may
take a homogeneous polynomial $P$ of degree $rn(n-1)/2$ in $n(n-1)/2$
variables $\xi=\{\xi_{ij}\,|\ 1\le i<j\le n\}$ and for each
$Q=\{q_1,...,q_n\}\in{\mathcal C}^n({\mathbb C})$
set $\xi_{ij}(Q)=q_i-q_j$, \ $b(Q)=P(\xi(Q))/d_n^r(Q)$.
\vskip0.2cm

\noindent A similar example of a biregular automorphism of
${\mathcal C}^n({\mathbb{CP}}^1)$ looks as follows:
$Q=\{q_1,...,q_n\}\mapsto\{Cq_1+B(Q),...,Cq_n+B(Q)\}$,
where $C\in{\mathbb C}^*$ and $B\colon{\mathcal C}^n({\mathbb{CP}}^1)
\to{\mathbb C}$ is a regular ${\mathbf{PSL}}(2,{\mathbb C})$ invariant
function.
\hfill $\bigcirc$
\end{Example}

\begin{Example}\label{Ex: biregular automorphism of Cn(C) generated
by non-invariant morphism Cn(C) to Aff C}
The map
$$
F_T\colon{\mathcal C}^n({\mathbb C})\ni Q\mapsto cQ
+\lambda a(Q)\cdot\1\in{\mathcal C}^n({\mathbb C})\,,
\quad \1=\{1,...,1\}\in{\mathbb C}^n\,,
$$
where $\displaystyle a(Q)=\sum_{\zeta\in Q}\zeta$,
$c\in{\mathbb C}^*$ and $\lambda\in{\mathbb C}$,
$\lambda\ne -c/n$, provides
an example of a tame biregular automorphism $F_T$ of
${\mathcal C}^n({\mathbb C})$ corresponding to 
a {\em non-invariant} morphism
$$
T\colon{\mathcal C}^n({\mathbb C})\ni Q\mapsto T(Q)
\in\Aff{\mathbb C}\,, \ \ T(Q)w=cw+\lambda a(Q)\,, \ 
w\in{\mathbb C}\,.
$$
The inverse map is given by
$$
{\mathcal C}^n({\mathbb C})\ni Q\mapsto c'Q+\lambda' a(Q)\cdot\1
\in{\mathcal C}^n({\mathbb C})\,,
$$
where $c'=1/c$ and $\lambda'=-\lambda/[c(c+n\lambda)]$.
\hfill $\bigcirc$
\end{Example}

\begin{Example}\label{Ex: Cn to Cn}
Let $A(z)$,$B(z)$ be homogeneous polynomials in
$z=(z_0,...,z_n)\in{\mathbb C}^{n+1}$ of degrees $2a(n-1)$
and $2b(n-1)$ respectively (one or both of $A,B$ may be $0$).
As in Remark \ref{Rmk: non-degenerate forms}, let \
$Q\in{\mathcal C}^n({\mathbb{CP}}^1)$ \ be the zero set
of a binary projective form
$\phi(x,y;z)=z_0x^n+z_1x_{n-1}y+...+z_ny^n$ with the discriminant
$D_n(z)$. Set $\displaystyle h(z)=\frac{A(z)\exp[B(z)/(D_n(z))^b]}{(D_n(z))^a}$
and define the holomorphic map
$T_{A,B}\colon{\mathcal C}^n({\mathbb{CP}}^1)
\to{\mathbf{SL}}(2,{\mathbb C})$ and the corresponding tame
holomorphic endomorphism $F_{T_{A,B}}
\colon{\mathcal C}^n({\mathbb{CP}}^1)\to{\mathcal C}^n({\mathbb{CP}}^1)$
by
\begin{equation}\label{eq: TAB: Cn to SL(2,C)}
T_{A,B}(Q)=\left(
\begin{array}{ll}
h(z),& h(z)+1\\
h(z)-1, & h(z)
\end{array}
\right)\,, \quad F_{T_{A,B}}(Q)=T_{A,B}(Q)Q\,.
\end{equation}
It is easily seen that $F_{T_{A,B}}$ is homotopic
to the automorphism 
$$
{\mathcal C}^n({\mathbb{CP}}^1)\ni Q=\{q_1,...,q_n\}\mapsto
\{-q_1^{-1},...,-q_n^{-1}\}\Def -Q^{-1}
\in{\mathcal C}^n({\mathbb{CP}}^1)
$$
in the class of holomorphic endomorphisms of
${\mathcal C}^n({\mathbb{CP}}^1)$.
\hfill $\bigcirc$
\end{Example}



\subsection{Automorphisms of certain moduli spaces}
\label{Ss: Automorphisms of certain moduli spaces}
\index{Moduli space!ordered\hfill}
\index{Moduli space!automorphisms of\hfill}
The moduli space $M_o(0,n)$ 
of the Riemann sphere with $n$ {\em ordered} punctures 
may be identified to the orbit space
${\mathcal C}_o^n({\mathbb{CP}}^1)/{\mathbf{PSL}}(2,{\mathbb C})$.
According to Sec. \ref{Ss: Holomorphic universal covering},
the latter orbit space may be viewed as the $(n-3)$-dimensional
ordered configuration space
of the Riemann sphere punctured at $0$, $1$ and $\infty$:
\begin{equation}\label{eq: Con/PSL(2,C)}
\aligned
\hskip-0.25cm M_o&(0,n)
={\mathcal C}_o^n({\mathbb{CP}}^1)/{\mathbf{PSL}}(2,{\mathbb C})
\cong{\mathcal C}_o^{n-3}({\mathbb C}\setminus\{0,1\})\\
&=\{z=(z_1,...,z_{n-3})\in{\mathbb C}^{n-3}|\, 
           z_i\ne 0,1 \ \forall\,i\,, \ \ z_i\ne z_j \
                                                    \forall\, i\ne j\}\\
&=\{q=(q_1,...,q_n)\in{\mathcal C}_o^n({\mathbb{CP}}^1)\,|\ 
                    q_{n-2}=0,\ q_{n-1}=1,\ q_n=\infty\}
={\mathcal D}^{n-3}({\mathbb{CP}}^1).
\endaligned
\end{equation}
Let $g=(g_1,...,g_{n-3})\colon{\mathcal D}^{n-3}({\mathbb{CP}}^1)
\to{\mathcal D}^{n-3}({\mathbb{CP}}^1)$ be a non-constant holomorphic map.
Set $g_{n-2}=0$, $g_{n-1}=1$ and $g_n=\infty$. 
S. Kaliman \cite{Kal76a} discovered the following fact:

\begin{Theorem}[\caps 1st Kaliman theorem]
\label{Thm: 1st Kaliman theorem} 
\index{Endomorphisms of ${\mathcal C}_o^m({\mathbb C}\setminus\{0,1\})$\hfill} 
\index{Kaliman theorem!1st\hfill}
\index{1st Kaliman theorem\hfill}
For $n\ge 4$ every non-constant holomorphic endomorphism 
$g=(g_1,...,g_{n-3})$ of the space 
${\mathcal D}^{n-3}({\mathbb{CP}}^1)
={\mathcal C}_o^{n-3}({\mathbb C}\setminus\{0,1\})$
is an automorphism and its components $g_r$
are of the form
\begin{equation}\label{eq: components of endomorphism D(n-3) to D(n-3)}
g_r(q)=\frac{q_{\sigma(r)}-q_{\sigma(n-2)}}{q_{\sigma(r)}-q_{\sigma(n)}}:
\frac{q_{\sigma(n-2)}-q_{\sigma(n-1)}}{q_{\sigma(n)}-q_{\sigma(n-1)}}\,,
\ \ 1\le r\le n-3\,,
\end{equation}
where $\sigma\in{\mathbf S}(n)$ is a permutation not depending on $r$. 
\end{Theorem}

Since Kaliman's paper \cite{Kal76a} does not contain a proof and exists
only in Russian, I decided to prove this theorem here
(due to the techniques developed in
Sec. \ref{Sec: Holomorphic functions omitting two values},
the proof below is simpler than the original one).  

\begin{proof}
First we prove that each $g_r\ne\const$. Suppose, on the contrary,
that $g_r=c=\const\ne 0,1$ for some $r\in\{1,...,n-3\}$,
whereas $g_p\ne\const$ for some $p\in\{1,...,n-3\}$, $p\ne r$.
Since $g_p(q)\ne g_r(q)$ for all
$q\in{\mathcal D}^{n-3}({\mathbb{CP}}^1)$,
we see that $(*)$ $g_p(q)\ne 0,1,c$ everywhere in
${\mathcal D}^{n-3}({\mathbb{CP}}^1)$. According to Remark
\ref{Rmk: function omitting 0,1 on Com(C**) etc}, $g_p$ must be of the form
$$
g_p(q_1,...,q_{n-3})
=\frac{q_i-q_j}{q_i-q_k}:\frac{q_j-q_l}{q_k-q_l}\,,
$$
where $i,j,k,l\in\{1,...,n\}$ are pairwise distinct.
However the function in the right hand side assumes on 
${\mathcal D}^{n-3}({\mathbb{CP}}^1)$ all values
$a\in{\mathbb C}\setminus\{0,1\}$, which contradicts the property
$(*)$ mentioned above.

Now we consider the composition 
\begin{equation}\label{eq: composition}
G\colon{\mathcal C}_o^n({\mathbb{CP}}^1)
\overset{j_{{\mathbb{CP}}^1,n}}{\cong}{\mathbf{PSL}}(2,{\mathbb C})
\times{\mathcal D}^{n-3}({\mathbb{CP}}^1)
\overset{\nu}{\longrightarrow}{\mathcal D}^{n-3}({\mathbb{CP}}^1)
\overset{g}{\longrightarrow}{\mathcal D}^{n-3}({\mathbb{CP}}^1)\,,
\end{equation}
where $\nu$ is the natural projection. Let $f_r(q)=q_r$
for all $q=(q_1,...,q_{n-3})$ and $r=1,...,n-3$. 
Clearly, the set $\Delta^{n-4}=\{f_1,...,f_{n-3}\}$ is an $(n-4)$-simplex
in the complex $L({\mathcal D}^{n-3}({\mathbb{CP}}^1))$.
Since $g_r\ne\const$ for all $r=1,...,n-3$ and
$g(q)=(g_1(q),...,g_n(q))\in{\mathcal D}^{n-3}({\mathbb{CP}}^1)$
for all $q=(q_1,...,q_{n-3},0,1,\infty)\in{\mathcal D}^{n-3}({\mathbb{CP}}^1)$,
the set $g^*(\Delta^{n-4})=(g_1,...,g_{n-3})$ also is an $(n-4)$-simplex
in the complex $L({\mathcal D}^{n-3}({\mathbb{CP}}^1))$. 
Furthermore, the set $G^*(\Delta^{n-4})
=j_{{\mathbb{CP}}^1,n}^*(\nu^*(g^*(\Delta^{n-4})))$ is an $(n-4)$-simplex
in the complex $L({\mathcal C}_o^n({\mathbb{CP}}^1))$.
By Lemma \ref{Lm: S(n) orbits in SRn and CRn}$(${\bfit b}$)$, the latter
simplex can be carried to the normal form (\ref{eq: 3rd normal form})
(with $m=n-4$) via an appropriate permutation of the coordinates $q_1,...,q_n$
in ${\mathcal C}_o^n({\mathbb{CP}}^1)$. It is easily seen that (up to a choice
of the notation) this is precisely what formulas
(\ref{eq: components of endomorphism D(n-3) to D(n-3)})
say about the components $g_r$ of the map $g$. Finally, a map $g$ whose
components $g_r$ satisfy (\ref{eq: components of endomorphism D(n-3) to D(n-3)})
is certainly an automorphism of ${\mathcal D}^{n-3}({\mathbb{CP}}^1))$. 
\end{proof}

\noindent Clearly, $\sigma Aq=A\sigma q$ for all
$\sigma\in{\mathbf S}(n)$, $A\in{\mathbf{PSL}}(2,{\mathbb C})$ and
$q\in{\mathcal C}_o^n({\mathbb{CP}}^1)$.
Thereby, every $\sigma\in{\mathbf S}(n)$ produces
an automorphism of the moduli space $M_o(0,n)$. The above theorem
implies the following immediate corollary. 

\begin{Corollary}\label{Crl: automorphisms of Mo(0,n)}
Every automorphism of the moduli space
$M_o(0,n)$ of the Riemann sphere with $n$ ordered punctures
is produced by a certain
permutation $\sigma\in{\mathbf S}(n)$ of coordinates
in ${\mathcal C}_o^n({\mathbb{CP}}^1)$.\footnote{I. Dolgachev wrote me
that this was conjectured by W. Fulton, but I did not find the corresponding
reference.} Hence $\Aut M_o(0,n)\cong{\mathbf S}(n)$. 
\end{Corollary}

\subsection{Endomorphisms of balanced confiquration spaces}
\label{Ss: Endomorphisms of balanced confiquration spaces} 
\index{Configuration spaces!balanced\hfill}
\index{Balanced configuration space ${\mathcal C}^{n-1}_b({\mathbb C})$\hfill}
\index{Balanced configuration space ${\mathbf G}_m^0$\hfill}
\index{${\mathcal C}^{n-1}_b({\mathbb C})$\hfill}
\index{${\mathbf G}_m^0$\hfill}
\index{Configuration spaces!balanced!proper endomorphisms of\hfill}
\index{Balanced configuration spaces\hfill}
\index{Balanced configuration spaces!proper endomorphisms of\hfill} 
In Sec. \ref{Ss: Holomorphic maps to the balanced configuration space}
we introduced the balanced configuration spaces
$$
{\mathcal C}^{n-1}_b({\mathbb C})
=\{Q=\{q_1,...,q_n\}\in{\mathcal C}^n({\mathbb C})\,|\ q_1+...+q_n=0\}
$$
and
$$
{\mathcal C}^{n-1}_{ob}({\mathbb C})
=\{q=(q_1,...,q_n)\in{\mathcal C}^n_o({\mathbb C})\,|\ q_1+...+q_n=0\}\,.
$$
As in Remark \ref{Rmk: non-degenerate forms},
we identify ${\mathcal C}^{n-1}_b({\mathbb C})$ to the space
$$
{\mathbf G}_n^0\Def\{z=(0,z_2,...,z_n)\in{\mathbb C}^n\,|\ d_n(z)\ne 0\}
$$
of all complex polynomials of the form
$p_n^0(t;z)=t^n+z_2t^{n-2}+...+z_n$ with simple roots. Clearly,
${\mathbf G}_n\cong{\mathbb C}\times{\mathbf G}_n^0$, so that
$\pi_k({\mathbf G}_n^0)=\pi_k({\mathbf G}_n)=0$ for $k>1$ and
$\pi_1({\mathbf G}_n^0)=\pi_1({\mathbf G}_n)=B_n$; 
we have also the ${\mathbf S}(n)$ Galois covering 
\begin{equation}\label{eq: projection of ordered balanced configurations}
p\colon{\mathcal C}^{n-1}_{ob}({\mathbb C})\ni q=(q_1,...,q_n)\mapsto
\{q_1,...,q_n\}=Q\in{\mathcal C}^{n-1}_b({\mathbb C})={\mathbf G}_n^0
\end{equation}
that carries any ordered balanced configuration $q=(q_1,...,q_n)$
to the unordered one $Q=\{q_1,...,q_n\}$, which, in turn, identifies  
to the polynomial $p_n^0(t;z)=\prod_{j=1}^n (t-q_j)$
with the roots $q_1,...,q_n$.    
We define also the {\em special} balanced configuration spaces 
\index{Configuration spaces!balanced!special\hfill}
\index{Balanced configuration spaces!special\hfill}
\begin{equation}\label{eq: special balanced configuration space}
\aligned
{\mathcal{SC}}^{n-2}_b({\mathbb C})&\Def
\{Q=\{q_1,...,q_n\}\in{\mathcal C}^n({\mathbb C})\,|\ q_1+...+q_n=0\,, \ \
                                    d_n(Q)=1\} \\
&=\{z=(z_1,z_2,...,z_n)\in {\mathbb C}^n\,|\ z_1=0\,, \ \ d_n(z)=1\}
\Def{\mathbf{SG}}_n
\endaligned
\end{equation}
and
\begin{equation}\label{eq: special balanced ordered configuration space}
\aligned
{\mathcal{SC}}^{n-2}_{ob}&({\mathbb C})={\mathbf{SE}}_n^+\\
&\Def\big\{q=(q_1,...,q_n)\in{\mathbb C}^n\,|\ q_1+...+q_n=0\,, \ \
          \prod_{j<k}(q_j-q_k)=i^{n(n-1)/2}\big\}\,.
\endaligned
\end{equation}
The latter space is one of the two connected components of the
space 
$$
{\mathbf{SE}}_n
=\big\{q=(q_1,...,q_n)\in{\mathbb C}^n\,|\ q_1+...+q_n=0\,, \ \ 
D_n(q)=\prod_{j\ne k}(q_j-q_k)=1\big\}\,,
$$
and we have the ${\mathbf A}(n)$ Galois covering 
\begin{equation}\label{eq: A(n) covering}
p'\colon{\mathbf{SE}}_n^+\ni q=(q_1,...,q_n)\mapsto
\prod_{j=1}^n (t-q_j)\in{\mathbf{SG}}_n\,.
\end{equation}
Both spaces ${\mathbf G}_n^0$, ${\mathbf{SG}}_n$ and 
their coverings ${\mathcal C}^{n-1}_{ob}({\mathbb C})$, ${\mathbf{SE}}_n^+$
are smooth irreducible affine algebraic manifolds.
The discriminant map 
\begin{equation}\label{eq: restricted discriminant map}
d_n\colon{\mathbf G}_n^0\ni z\mapsto d_n(z)\in{\mathbb C}^*
\end{equation}
is a holomorphic locally trivial fiber bundle with the fiber
${\mathbf{SG}}_n$. The exact homotopy sequence of this bundle shows that
${\mathbf{SG}}_n$ is an Eilenberg-MacLane $K(\pi,1)$ space for
the commutator subgroup $B'_n$ of the Artin braid group $B_n$.
It is also easely seen that ${\mathbf{SE}}_n^+$ an Eilenberg-MacLane
$K(\pi,1)$ space for the intersection $J_n={PB_n}\cap{B}'_n$
of the pure braid group $PB_n$ with the commutator subgroup $B'_n$.
\vskip0.3cm

\noindent In what follows, we describe automorphisms
${\mathbf G}_n^0$ and ${\mathbf{SG}}_n$; in fact, 
a complete exhibition of endomorphisms of ${\mathbf{SG}}_n$
and proper endomorphisms of ${\mathbf G}_n^0$ will be presented.

Notice that the stabilizer 
$$
St_{{\mathcal C}^{n-1}_b({\mathbb C})}=St_{{\mathbf G}_n^0}
=\{A\in\Aff{\mathbb C}\,|\ A({\mathcal C}^{n-1}_b({\mathbb C}))
\subseteq{\mathcal C}^{n-1}_b({\mathbb C})\}
$$
of the subspace ${\mathcal C}^{n-1}_b({\mathbb C})={\mathbf G}_n^0
\subset{\mathbf G}_n={\mathcal C}^{n-1}({\mathbb C})$ in the
group $\Aut{\mathbb C}=\Aff{\mathbb C}$ is isomorphic to 
${\mathbb C}^*$ and consists of all transformations of the form
\begin{equation}\label{eq: C* action in Cb(n-1)}
A_\zeta\colon{\mathcal C}^{n-1}_b({\mathbb C})\ni Q=\{q_1,...,q_n\}
\mapsto \zeta Q
=\{\zeta q_1,...,\zeta q_n\}\in{\mathcal C}^{n-1}_b({\mathbb C})\,, \ \
\zeta\in{\mathbb C}^*\,.
\end{equation}
Under the identification ${\mathcal C}^{n-1}_b({\mathbb C})
\cong{\mathbf G}_n^0$, the corresponding ${\mathbb C}^*$ action
in ${\mathbf G}_n^0$ looks as follows:
\begin{equation}\label{eq:C* action in Gn0}
\aligned
\forall\,\zeta\in{\mathbb C}^*, \ \
A_\zeta\colon{\mathbf G}_n^0\ni z&=(0, z_2,...,z_{n-1},z_n)
\mapsto (0, \zeta^2 z_2,...,\zeta^{n-1} z_{n-1},\zeta^n z_n)\\
&=t^n+\zeta^2 z_2 t^{n-2}+...+
\zeta^{n-1} z_{n-1} t+\zeta^n z_n\in{\mathbf G}_n^0\,.
\endaligned
\end{equation}
Furthermore, the stabilizer (in the group ${\mathbb C}^*$ acting
in ${\mathbf G}_n^0$ according to (\ref{eq:C* action in Gn0}))
\begin{equation}\label{eq: stabilizer of SGn}
St_{{\mathbf{SG}}_n}
=\{\zeta\in{\mathbb C}^*\,|\ A_\zeta({\mathbf{SG}}_n)
\subseteq{\mathbf{SG}}_n\}
\end{equation}
of the special balanced subspace ${\mathbf{SG}}_n\subset{\mathbf G}_n^0$
is the subgroup
$$
\Gamma_{n(n-1)}=\{1,\omega,...,\omega^{n(n-1)-1}\}
\cong{\mathbb Z}_{n(n-1)}={\mathbb Z}/n(n-1){\mathbb Z}\,,
$$
where $\omega=e^{2\pi i/n(n-1)}$. 

Notice also that every holomorphic function
$h\colon{\mathbf G}_n^0\to{\mathbb C}^*$ is of the form
$h(z)=e^{\varphi(z)}d_n^m(z)$, where $d_n(z)$ is the discriminant,
$m\in{\mathbb Z}$ and $\varphi$ is a holomorphic function
on ${\mathbf G}_n^0$.

\begin{Notation}
Let ${\mathcal A}_n
={\mathcal O}^{{\mathbb C}^*}({\mathbf G}_n^0)$ denote the algebra
of all holomorphic functions $a$ on ${\mathbf G}_n^0$ that are invariant under
the ${\mathbb C}^*$ action $(\ref{eq:C* action in Gn0})$.
In other words, ${\mathcal A}_n$ consists of all holomorphic
functions $a\colon{\mathbf G}_n^0\to{\mathbb C}$ such that
$a(0,\zeta^2 z_2,...,\zeta^n z_n)=a(0,z_2,...,z_n)$ for all
$z=(0,z_2,...,z_n)\in{\mathbf G}_n^0$ and $\zeta\in{\mathbb C}^*$.
For instance, the function
$a(z)=z_2^{n(n-1)/2}/d_n(0,z_2,...,z_n)$ belongs to ${\mathcal A}_n$.
\hfill $\bigcirc$
\end{Notation}

\noindent The following theorem stated in \cite{Lin72b,Lin79}
is an easy consequence of Tame Map Theorem \ref{Thm: Tame Map Thm}
and Theorem \ref{Thm: image of non-cyclic endomorphism of Cn(C)}.

\begin{Theorem}\label{Thm: Endomorphisms of Gn0}
\index{${\mathcal C}^{n-1}_b({\mathbb C})$!endomorphisms of\hfill}
\index{${\mathcal C}^{n-1}_b({\mathbb C})$!proper endomorphisms of\hfill}
\index{${\mathbf G}_m^0$!endomorphisms of\hfill}
\index{${\mathbf G}_m^0$!proper endomorphisms of\hfill}
\index{Configuration spaces!balanced!endomorphisms of\hfill}
\index{Balanced configuration spaces!endomorphisms of\hfill}
\index{Configuration spaces!balanced!proper endomorphisms of\hfill}
\index{Balanced configuration spaces!proper endomorphisms of\hfill}
For $n>4$ every non-cyclic holomorphic map
$F\colon{\mathbf G}_n^0\to{\mathbf G}_n^0$ is ${\mathbb C}^*$-tame,
meaning that it is induced by the ${\mathbb C}^*$ action
$(\ref{eq:C* action in Gn0})$ and a holomorphic function
$h\colon{\mathbf G}_n^0\to{\mathbb C}^*$.
In more details, there exist $m\in{\mathbb Z}$ and a holomorphic function
$a$ on ${\mathbf G}_n^0$ such that $h(z)=e^{a(z)}d_n^m(z)$ and
\begin{equation}\label{eq: non-cyclic F: Gn0 to Gn0}
\aligned
\hskip-0.3cm F&(z)=F_h(z)=A_{h(z)}z
=t^n+z_2 h^2(z)t^{n-2}+\!...\!+z_{n-1} h^{n-1}(z)t+z_n h^n(z)\\
&=t^n+z_2 e^{2a(z)}d_n^{2m}(z)t^{n-2}+\!...\!
+z_{n-1} e^{(n-1)a(z)}d_n^{(n-1)m}(z)t+z_n e^{na(z)}d_n^{nm}(z)\,.
\endaligned
\end{equation}
Moreover, every such $F$ is surjective, and the pre-image $F^{-1}(w)$
of any point $w\in{\mathbf G}_n^0$ is descrete and consists of
all $z$ of the form $z=A_\zeta w$, where $\zeta\in{\mathbb C}^*$ runs over
the set of all roots of the equation
\begin{equation}\label{eq: equation for pre-image}
H_w(\zeta)\Def d_n^m(w)\zeta^{mn(n-1)+1}
\exp a(A_\zeta w)=1\,.\footnote{For each $w\in{\mathbf G}_n^0$
the function $H_w\colon{\mathbb C}^*\to{\mathbb C}^*$
is a well-defined non-constant holomorphic function. By the Picard theorem,
the set $H_w^{-1}(1)$ of all roots $\zeta\in{\mathbb C}^*$ 
of the equation $H_w(\zeta)=1$ is non-empty and descrete;
it is easily seen that distinct $\zeta',\zeta''\in H_w^{-1}(1)$
produce distinct points $A_{\zeta'}w,A_{\zeta''}w\in F^{-1}(w)$.}
\end{equation}
A holomorphic map $F=F_h$ of the form
$(\ref{eq: non-cyclic F: Gn0 to Gn0})$ is proper if and only if the
corresponding holomorphic function $a$ belongs to ${\mathcal A}_n$.
Every proper holomorphic map $F\colon{\mathbf G}_n^0\to{\mathbf G}_n^0$
is everywhere non-degenerate; moreover, it is an unramified finite
${\mathbb Z}/N{\mathbb Z}$ Galois covering of degree
$N\equiv 1\modl n(n-1)$ $($in fact, $N=mn(n-1)+1$ whenever $F=F_h$
with $h=e^a d_n^m$ and some ${\mathbb C}^*$-invariant $a)$. 
The corresponding normal subgroup $\Gamma\vartriangleleft B_n$
of index $N$ consists of all
$g=\sigma_{i_1}^{m_1}\cdots\sigma_{i_q}^{m_q}\in B_n$ such that $N$ divides
$m_1+...+m_q$. Every two such coverings of the same degree are equivalent.
\hfill $\square$
\end{Theorem}

\noindent This theorem implies the following immediate corollary.

\begin{Corollary}\label{Crl: automorphisms of Gn0}
\index{${\mathcal C}^{m-1}_b({\mathbb C})$!automorphisms of\hfill}
\index{${\mathbf G}_m^0$!automorphisms of\hfill}
\index{Configuration spaces!balanced!automorphisms of\hfill}
\index{Balanced configuration spaces!automorphisms of\hfill}
\index{Automorphisms of ${\mathcal C}^{m-1}_b({\mathbb C})$\hfill}
\index{Automorphisms of ${\mathbf G}_m^0$\hfill} 
For $n>4$ every biholomorphic automorphism
$F\colon{\mathbf G}_n^0\to{\mathbf G}_n^0$ is of the form
\begin{equation}\label{eq: automorphism F: Gn0 to Gn0}
F(z)=F_{e^{a(z)}}(z)=A_{e^{a(z)}}z
=t^n+z_2 e^{2a(z)}t^{n-2}+...
+z_{n-1} e^{(n-1)a(z)}t+z_n e^{na(z)}\,,
\end{equation}
where $a\in{\mathcal A}_n$.
The automorphism group $\Aut({\mathbf G}_n^0)$
is isomorphic to the quotient group ${\mathcal A}_n/2\pi i{\mathbb Z}$.
The orbits of the natural $\Aut({\mathbf G}_n^0)$ action in ${\mathbf G}_n^0$
coinside with the orbits of ${\mathbb C}^*$ action
$(\ref{eq:C* action in Gn0})$; each of them is a smooth algebraic
curve isomorphic to ${\mathbb C}^*$.
\hfill $\square$
\end{Corollary}

\subsection{Endomorphisms of ${\mathbf{SG}}_n$}
\label{Ss: Automorphisms of SGn}
The following theorem was stated by S. Kaliman \cite{Kal76a}. 

\begin{Theorem}[\caps 2nd Kaliman theorem]
\label{Thm: 2nd Kaliman theorem}
\index{Endomorphisms of ${\mathbf{SG}}_n$\hfill} 
\index{Kaliman theorem!2st\hfill}
\index{2nd Kaliman theorem\hfill}
For $n\ne 4$ every non-constant holomorphic map 
$F\colon{\mathbf{SG}}_n\to{\mathbf{SG}}_n$ is an automorphism of the form
\begin{equation}\label{eq: automorphism of SGn}
F(z)=F_\omega(z)=t^n+z_2 \omega^{2m} t^{n-2}+...
+z_{n-1} \omega^{(n-1)m}t+z_n \omega^{nm}\,,
\end{equation}
where $\omega=e^{2\pi i/n(n-1)}$ and $m\in\{0,1,...,n-1\}$.
The automorphism group 
$\Aut({\mathbf{SG}}_n)\cong{\mathbb Z}/n(n-1){\mathbb Z}$.
\end{Theorem}

\begin{proof} A complete description of holomorphis endomorphisms of the affine
elliptic curve ${\mathbf{SG}}_3\cong\{(x,y)\in{\mathbb C}^2\,|\ x^2+y^3=1\}=$
{\em a torus with one puncture} is straightforward; therefore we
restrict ourselves to the case $n>4$.
\vskip0.3cm

\noindent According to \cite{Lin74,Lin79}, for $n>4$ the intersection
$J_n=PB_n\cap B'_n$ is a completely characteristic
subgroup of $B'_n$, that is, $\psi(J_n)\subseteq J_n$ for every endomorphism
$\psi$ of the group $B'_n$ (see \cite{Lin96b} or \cite{Lin04b} for a complete
proof). Hence every continuous map $F\colon{\mathbf{SG}}_n\to{\mathbf{SG}}_n$
lifts to a continuous map $f\colon{\mathbf{SE}}_n^+\to{\mathbf{SE}}_n^+$,
which is ${\mathbf A}(n)$ equivariant, meaning that there is an automprphism
$\alpha\in\Aut({\mathbf A}(n))$ such that $f(sq)=\alpha(s)f(q)$ for all
$q\in{\mathbf{SE}}_n^+$ and $s\in{\mathbf A}(n)$. 
In our case $F$ is holomorphic and its lifting $f$ is holomorphic, as well.
\vskip0.3cm

\noindent The cyclic group $\Gamma_{n(n-1)/2}
\cong{\mathbb Z}_{n(n-1)/2}={\mathbb Z}/[n(n-1)/2]{\mathbb Z}$
of $n(n-1)/2$-roots of $1$ acts freely on ${\mathbf{SE}}_n^+$
and on the variety $S$ defined by
$$
S=\big\{(\zeta,z)\in{\mathbb C}\times{\mathcal D}^{n-2}({\mathbb{CP}}^1)\,|\
\zeta^{n(n-1)/2}\prod\limits_{1\le j<k\le n}(z_j-z_k)=i^{n(n-1)/2}\big\}
$$
(here and in what follows $z_{n-1}=0$, $z_n=1$ and $z_{n+1}=\infty$).
In both cases the action is by multiplication:
$(q_1,...,q_n)\mapsto(\gamma q_1,...,\gamma q_n)$
and $(\zeta,z)\mapsto (\gamma\zeta,z)$, where $\gamma\in\Gamma_{n(n-1)/2}$.
These two $\Gamma_{n(n-1)/2}$ spaces are in fact equivalent. This can be seen
from the commutative diagram
$$
\CD
{\hskip-50pt {\mathbf{SE}}_n^+ }\hskip10pt @[1] > {\iota} > {\cong} >\hskip10pt { S } \\
\hskip33pt@[12]/SE/{p\circ\nu}//
    \hskip5pt@/SE/{ \nu }/{\Gamma_{n(n-1)/2}}/
      \hskip90pt@/SW/{\Gamma_{n(n-1)/2}}/{\pi}/@[12]/SW//{p\circ\pi}/ \\ 
\hskip90pt{ {\mathcal D}^{n-2}({\mathbb{CP}}^1)}\\
\hskip45pt@VpVV\\
\hskip90pt{ {\mathcal D}^{n-3}({\mathbb{CP}}^1)}\,,
\endCD
$$
where the isomorphism $\iota$ is defined by
$$
\aligned
\iota\colon &{\mathbf{SE}}_n^+\ni (q_1,...,q_n)\mapsto (\zeta,z)
\in S\subset{\mathbb C}\times{\mathcal D}^{n-2}({\mathbb{CP}}^1)\,,\\
&\zeta=q_n-q_{n-1}\,, \ \ z_j=\frac{q_j-q_{n-1}}{q_n-q_{n-1}} \ \ 
              \text{for} \ j=1,...,n\,, \ \ z_{n+1}=\infty\,,
\endaligned
$$
the quotient map 
$\nu\colon{\mathbf{SE}}_n^+\to{\mathbf{SE}}_n^+/\Gamma_{n(n-1)/2}
={\mathcal D}^{n-2}({\mathbb{CP}}^1)$ is given by 
\begin{equation}\label{eq: nu: SEn to Mo(0,n+1)}
\nu\colon{\mathbf{SE}}_n^+\ni 
(q_1,...,q_n)\mapsto
\left(\frac{q_1-q_{n-1}}{q_n-q_{n-1}},...,
\frac{q_{n-2}-q_{n-1}}{q_n-q_{n-1}},0,1,\infty\right)\in 
{\mathcal D}^{n-2}({\mathbb{CP}}^1)\,,
\end{equation}
the quotient map $\pi\colon S\to S/\Gamma_{n(n-1)/2}
={\mathcal D}^{n-2}({\mathbb{CP}}^1)$ is just the restriction to $S$
of the projection ${\mathbb C}\times{\mathcal D}^{n-2}({\mathbb{CP}}^1)
\to{\mathcal D}^{n-2}({\mathbb{CP}}^1)$, and finally
$$
\aligned
p\colon{\mathcal D}^{n-2}({\mathbb{CP}}^1)\ni 
z&=(z_1,z_2,...,z_{n-2},0,1,\infty)\\
&\mapsto z'=(z_2,...,z_{n-2},0,1,\infty)\in{\mathcal D}^{n-3}({\mathbb{CP}}^1)
\endaligned
$$
is the projection forgetting the first coordinate.
Notice that $p$ is a smooth fibration with fiber $D_{z'}=p^{-1}(z')$
being a Riemann sphere punctured at $n$ points
$z'=(z_2,...,z_{n-2},0,1,\infty)$ and coordinate $z_1$. Moreover,
being coverings of the Kobayashi hyperbolic manifold
${\mathcal D}^{n-2}({\mathbb{CP}}^1)$
(see Remark \ref{Rmk: Dm(X) is Kobayashi hyperbolic})
both ${\mathbf{SE}}_n^+$ and $S$ are Kobayashi hyperbolic as well. 
\vskip0.3cm

\noindent The following claim \cite{Kal75,Kal93}
provides a main technical tool of the proof.
\vskip0.2cm

\noindent {\bfit Claim.}
{\sl Any holomorphic function
$h\colon{\mathbf{SE}}_n^+\cong S\to{\mathbb C}\setminus\{0,1\}$
is $\Gamma_{n(n-1)/2}$ invariant and hence
can be pulled down to ${\mathcal D}^{n-2}({\mathbb{CP}}^1)$.}
\vskip0.2cm

{\bcaps Proof of Claim.}\footnote{A revised version communicated by S. Kaliman.}
The composition
$$
\rho=p\circ\pi\colon S\stackrel{\pi}{\longrightarrow}
{\mathcal D}^{n-2}({\mathbb{CP}}^1)\stackrel{p}{\longrightarrow}
{\mathcal D}^{n-3}({\mathbb{CP}}^1)
$$
is a locally trivial smooth fibration whose fiber
$S_{z'^\circ}=\rho^{-1}(z'^\circ)$
over a point $z'^\circ=(z^\circ_2,...,z^\circ_{n-2},0,1,\infty)
\in{\mathcal D}^{n-3}({\mathbb{CP}}^1)$
is a smooth irreducible algebraic curve in ${\mathbb C}\times D_{z'^\circ}$
given by
$$
S_{z'^\circ}=\Big\{(\zeta,z_1)\in{\mathbb C}\times D_{z'^\circ}|\,
\zeta^{n(n-1)/2}\cdot
\prod\limits_{k=2}^n (z_1-z^\circ_k)\cdot
\prod\limits_{2\le j<k\le n}(z^\circ_j-z^\circ_k)
=i^{n(n-1)/2}\Big\}\,.
$$
The restriction $p|S_{z'^\circ}\colon S_{z'^\circ}\ni (\zeta,z_1)
\to z_1\in D_{z'^\circ}$ of $p$ to $S_{z'^\circ}$ determines
a $\Gamma_{n(n-1)/2}$ covering of the punctured Riemann sphere
$D_{z'^\circ}=\pi^{-1}({z'^\circ})=\overline{\mathbb C}
\setminus\{z^\circ_2,...,z^\circ_{n-2},0,1,\infty\}$.
The curve $S_{z'^\circ}$ have $n-1+N$ punctures. The first $n-1$ of them
denoted by $a(z_2^\circ),...,a(z_{n-2}^\circ),
a(z_{n-1}^\circ)=a(0)$ and $a(z_n^\circ)=a(1)$ are situated, one by one,
above the punctures $z^\circ_2,...,z^\circ_{n-2},0,1$,
whereas the rest $N$ denoded by $b_1(z'^\circ),...,b_N(z'^\circ)$
lie above the puncture at $\infty$ ($N=n-1$ for $n$ even
and $N=(n-1)/2$ for $n$ odd). The $\Gamma_{n(n-1)/2}$ action 
on $S_{z'^\circ}$ coming from $S$ extends uniquelly to the complition
$\overline{S}_{z'^\circ}$ of $S_{z'^\circ}$. This extension keeps each
of the punctures $a$'s at its place; as opposed to this,
it is transitive on the set of the rest $N$ punctures $b$'s. 

Adding to each fiber $S_{z'}$ all its punctured points we get the union
$$
\overline S=\bigcup\limits_{{z'}\in{\mathcal D}^{n-3}({\mathbb{CP}}^1)}
\overline S_{z'}\supset 
\bigcup\limits_{{z'}\in{\mathcal D}^{n-3}({\mathbb{CP}}^1)} S_{z'}=S
$$
of smooth compact curves $\overline{S}_{z'}$. This union
$\overline{S}$ may be viewed as an algebraic variety obtained from
the variety 
$$
\aligned
\Sigma=\Big\{([\zeta:\tau],&[z_1:t]),\,(z_2,...,z_{n-2}))
\in{\mathbb{CP}}^1\times{\mathbb{CP}}^1
\times{\mathcal C}_o^{n-3}({\mathbb C})\,| \\
&\zeta^{n(n-1)/2}
\prod\limits_{k=2}^n (z_1-tz_k)\prod\limits_{2\le j<k\le n}(z_j-z_k)
=i^{n(n-1)/2}\tau^{n(n-1)/2}t^{n-1}\Big\}
\endaligned
$$
by normalization of all fibers ${\widehat S}_{z'}$
of the projection $\Sigma\to{\mathcal C}_o^{n-3}({\mathbb C})$; 
here $[\zeta:\tau]$ and $[z_1:t]$ are homogeneous coordinates 
in two copies of ${\mathbb{CP}}^1$,
$(z_2,...,z_{n-2})\in{\mathcal C}_o^{n-3}({\mathbb C})$, $z_{n-1}=0$ and
$z_n=1$. The projection $\rho$ extends to a projection
${\overline\rho}\colon{\overline S}\to{\mathcal C}_o^{n-3}({\mathbb C})$
with fibers ${\overline S}_{z'}$.

For each $j=2,...,n$ the set $A_j$ of all punctures 
$$
\{a(z_j)\,|\ z'=(z_2,...,z_n)\in{\mathcal D}^{n-3}({\mathbb{CP}}^1)\}
$$
is a hypersurface in $\overline S$, which is a section of the projection
${\overline\rho}\colon{\overline S}\to{\mathcal C}_o^{n-3}({\mathbb C})$
and hence irreducibale.
The $\Gamma_{n(n-1)/2}$ action on $\overline{S}$ coming from $S$
is identical on each $A_j$ and keeps stable the hypersurface 
$$
B=\{b_r(z')\,|z'=(z_2,...,z_n)\in{\mathcal D}^{n-3}({\mathbb{CP}}^1); \ 
r=1,...,N\}
$$
consisting of all punctures above infinity. The hypersurfaces $A_2,...,A_n,B$
do not intersect each other. The projection
$B\to{\mathcal D}^{n-3}({\mathbb{CP}}^1)$ is an unramified Galois covering
of degree $N$ with a transitive Galois group, which is a quotient group
of $\Gamma_{n(n-1)/2}$. To see this one can travel in the base
${\mathcal D}^{n-3}({\mathbb{CP}}^1)$ along a small simple loop
around some hyperplane $z_j=z_k$ with $2\le j<k\le n$;
this trip generates a permutation of the points
$b_1(z'),...,b_N(z')$ that is the same as the
action of a primitive element of $\Gamma_{n(n-1)/2}$.
In particular, $B$ is irreducible.
\vskip0.3cm

\noindent Now take any holomorphic function
$h\colon S\to{\mathbb C}\setminus\{0,1\}$.
If $h=\const$ on each fiber $S_{z'}=\rho^{-1}(z')$
then $h$ is a lift-up of a function on
${\mathcal D}^{n-3}({\mathbb{CP}}^1)$ and, all the more,
it is $\Gamma_{n(n-1)/2}$ invariant.

Suppose that $h\ne\const$ on a certain fiber $S_{z'^\circ}=\rho^{-1}(z'^\circ)$.
Since $S$ is hyperbolic, $h$ extends to a holomorphic map
${\bar h}\colon\overline{S}\to{\mathbb{CP}}^1$. Being non-constant
on $\overline{S}_{z'^\circ}$ the latter extension $\bar h$ takes somewhere
on $\overline{S}_{z'^\circ}$ the values $0$, $1$ and $\infty$;
in fact $\bar h$ must accept each of these values on the whole union
of some of the hypersurfaces $A_2,...,A_n,B$.

Let $\gamma\in\Gamma_{n(n-1)/2}$. As $\gamma (A_j)=A_j$ and $\gamma (B)=B$
we see that the rational function 
$h/(h\circ\gamma)$ has no poles on $\overline{S}$ and hence 
$\bar h/(\bar h\circ\gamma)=\const$ on each fiber $S_{z'}$ of $\rho$.
This constant must be $1$ since $\bar h=1$ at least on one of
the irreducible hypersurfaces $A_2,...,A_n,B$ and on such a hypesurface
$\bar h=\bar h\circ\gamma$. Thus, $h=h\circ\gamma$ on each $S_{z'}$
and hence everywhere on $S$. This is the desired conclusion. 
\hfill $\square$
\vskip0.3cm

\noindent Now we can complete the proof of the theorem.
By Claim, each function
$$
g_j=z_j\circ\nu\circ f\colon{\mathbf{SE}}_n^+\overset{f}{\longrightarrow}
{\mathbf{SE}}_n^+\overset{\nu}{\longrightarrow}
{{\mathcal D}^{n-2}({\mathbb{CP}}^1)}
\overset{z_j}{\longrightarrow}{\mathbb C}\setminus\{0,1\}\,, \ \ 1\le j\le n-2\,,
$$
pulls down to ${\mathcal D}^{n-2}({\mathbb{CP}}^1)$; this
gives a commutative diagram
\begin{equation}\label{eq: push down f}
\CD
{{\mathbf{SE}}_n^+} @ > {f} >> {{\mathbf{SE}}_n^+}\\
@V{\nu}VV @VV{\nu}V\\ 
{{\mathcal D}^{n-2}({\mathbb{CP}}^1)} @ >> {g}>
{{\mathcal D}^{n-2}({\mathbb{CP}}^1)}\,,
\endCD
\end{equation}
with $g=(g_1,...,g_{n-2},0,1,\infty)$ being a non-constant holomorphic map.
According to Theorem \ref{Thm: 1st Kaliman theorem}, there is a permutation
$\sigma\in{\mathbf S}(n+1)$ such that 
$$
g_r(z)=\frac{z_{\sigma(r)}-z_{\sigma(n-1)}}{z_{\sigma(r)}-z_{\sigma(n+1)}}:
\frac{z_{\sigma(n-1)}-z_{\sigma(n)}}{z_{\sigma(n+1)}-z_{\sigma(n)}}\,,
\ \ 1\le r\le n-2\,,
$$
for all
$z=(z_1,...,z_{n-2},0,1,\infty)\in{\mathcal D}^{n-2}({\mathbb{CP}}^1)$.
The existence of diagram (\ref{eq: push down f})
with an ${\mathbf A}(n)$ equivariant $f$ implies that $\sigma(n+1)=n+1$,
for otherwise $g$ cannot have such a lifting $f$ to ${\mathbf{SE}}_n^+$.
Thus $z_{\sigma(n+1)}=z_{n+1}=\infty$ and 
\begin{equation}\label{eq: components of endomorphism D(n-2) to D(n-2) final}
g_r(z)=\frac{z_{\sigma(r)}-z_{\sigma(n-1)}}{z_{\sigma(n)}-z_{\sigma(n-1)}}
\,,
\ \ 1\le r\le n-2\,,
\end{equation}
where $\sigma\in{\mathbf S}(n)$. By (\ref{eq: nu: SEn to Mo(0,n+1)}),
(\ref{eq: push down f})
and (\ref{eq: components of endomorphism D(n-2) to D(n-2) final}),
$$
\aligned
\frac{f_r(q)-f_{n-1}(q)}{f_n(q)-f_{n-1}(q)}&=(\nu\circ f)_r(q)
                =(g\circ\nu)_r(q)=g_r(\nu(q))\\
&\hskip-1.5cm=g_r\left(\frac{q_1-q_{n-1}}{q_n-q_{n-1}},...,
\frac{q_{n-2}-q_{n-1}}{q_n-q_{n-1}},0,1,\infty\right)\\
&\hskip-1.5cm=\left[\frac{q_{\sigma(r)}-q_{n-1}}{q_n-q_{n-1}}
 -\frac{q_{\sigma(n-1)}-q_{n-1}}{q_n-q_{n-1}}\right]:
\left[\frac{q_{\sigma(n)}-q_{n-1}}{q_n-q_{n-1}}
 -\frac{q_{\sigma(n-1)}-q_{n-1}}{q_n-q_{n-1}}\right]\\
&\hskip-1.5cm=\frac{q_{\sigma(r)}-q_{\sigma(n-1)}}{q_{\sigma(n)}-q_{\sigma(n-1)}}
\qquad \text{for} \ 1\le r\le n-2\,,
\endaligned
$$
whereas $(\nu\circ f)_{n-1}(q)=0$, \ 
$(\nu\circ f)_n(q)=1$, \ $(\nu\circ f)_{n+1}(q)=\infty$.
The final relation 
\begin{equation}\label{eq: fr(q)-f(n-1)(q) over fn(q)-f(n-1)(q)}
\frac{f_r(q)-f_{n-1}(q)}{f_n(q)-f_{n-1}(q)}
=\frac{q_{\sigma(r)}-q_{\sigma(n-1)}}{q_{\sigma(n)}-q_{\sigma(n-1)}}
\end{equation}
holds certainly true for all $r=1,...,n$.

In what follows we just solve equations
(\ref{eq: fr(q)-f(n-1)(q) over fn(q)-f(n-1)(q)}).
First, it follows from (\ref{eq: fr(q)-f(n-1)(q) over fn(q)-f(n-1)(q)}) that
\begin{equation}
\label{eq: 2nd simple ratios of f(q) via simple ratios of sigma q}
\frac{f_i(q)-f_j(q)}{f_n(q)-f_{n-1}(q)}
=\frac{q_{\sigma(i)}-q_{\sigma(j)}}{q_{\sigma(n)}-q_{\sigma(n-1)}}
\qquad \text{for} \ 1\le r\le n\,.
\end{equation}
Since $\sigma$ may be odd, for some $q\in{\mathbf{SE}}_n^+$
the point $\sigma^{-1}q$ may be not   
in ${\mathbf{SE}}_n^+$ but in the second connected component
${\mathbf{SE}}_n^-$ of ${\mathbf{SE}}_n$; however in any case
the polynomial $\prod\limits_{j=1}^n (t-q_{\sigma(j)})$ belongs to
${\mathbf{SG}}_n$ so that $q_{\sigma(1)}+...+q_{\sigma(n)}=0$
and $\prod\limits_{1\le j\ne k\le n} (q_{\sigma(j)}-q_{\sigma(k)})=1$.
Therefore, it follows from
(\ref{eq: 2nd simple ratios of f(q) via simple ratios of sigma q})
that
$$
\aligned
\frac{1}{(f_n(q)-f_{n-1}(q))^{n(n-1)}}
&=\frac{\prod\limits_{1\le i\ne j\le n}(f_i(q)-f_j(q))}
                       {(f_n(q)-f_{n-1}(q))^{n(n-1)}}\\
&\\
&=\frac{\prod\limits_{1\le i\ne j\le n}(q_{\sigma(i)}-q_{\sigma(j)})}
                             {(q_{\sigma(n)}-q_{\sigma(n-1)})^{n(n-1)}}
=\frac{1}{(q_{\sigma(n)}-q_{\sigma(n-1)})^{n(n-1)}}
\endaligned
$$
and hence
$$
f_n(q)-f_{n-1}(q)=e^{2\pi im/n(n-1)}(q_{\sigma(n)}-q_{\sigma(n-1)})
\quad \text{for all} \ q\in{\mathbf{SE}}_n^+
$$
with a certain $m\in\{0,1,...,n(n-1)-1\}$.
Together with (\ref{eq: fr(q)-f(n-1)(q) over fn(q)-f(n-1)(q)}),
this shows that
\begin{equation}\label{eq: fr(q)-f(n-1)(q)}
f_r(q)-f_{n-1}(q)\equiv e^{2\pi im/n(n-1)}(q_{\sigma(r)}-q_{\sigma(n-1)})\quad
\text{for all} \ r=1,...,n\,.
\end{equation}
Summing over all $r=1,...,n$ and taking inro account that 
$f_1+...+f_n=0$, we obtain $f_{n-1}(q)=e^{2\pi im/n(n-1)}q_{\sigma(n-1)}$.
This and (\ref{eq: fr(q)-f(n-1)(q)}) show that
\begin{equation}\label{eq: final formula for fr}
f_r(q)=e^{2\pi im/n(n-1)}q_{\sigma(r)}\quad \text{for all} \ \
q\in{\mathbf{SE}}_n^+ \ \ \text{and} \ \ r=1,...,n\,.
\end{equation}
Since $f(q)\in{\mathbf{SE}}_n^+$, we have also 
$e^{2\pi im/n(n-1)}\sigma^{-1}q
=e^{2\pi im/n(n-1)}(q_{\sigma(1)},...,q_{\sigma(n)})\in{\mathbf{SE}}_n^+$,
and (\ref{eq: final formula for fr}) implies that 
$F$ is of the desired form (\ref{eq: automorphism of SGn}).
\end{proof}


\section{Certain coverings of ${\mathcal C}^n(X)$}
\label{Sec: Certain coverings of Cn(X)}

\noindent The main objective of this section 
is describe all unbranched connected $k$-coverings
of the configuration spaces ${\mathcal C}^n(X)$
for $k\le 2n$. We will use the notation introduced in
Remark \ref{Rmk: non-degenerate forms}.

\subsection{Cyclic coverings}
\label{Ss: Cyclic coverings}
The following unbranched holomorphic cyclic coverings will
be referred to as the standard ones:
\begin{eqnarray}
&&\hskip-48pt
E_\infty^n({\mathbb C})=
\{(\zeta,w)\in{\mathbb C}\times{\mathbf G}_n\,|\ e^\zeta=d_n(w)\}, \ 
      p_\infty\colon E_\infty^n({\mathbb C})\ni (\zeta,w)
         \mapsto w\in{\mathbf G}_n; \label{eq: standard Z covering}\\
&&\hskip-48pt
E_k^n({\mathbb C}) \ =
\{(\zeta,w)\in{\mathbb C}^*\times{\mathbf G}_n\,|\ \zeta^k=d_n(w)\}, \
      p_k\colon E_k^n({\mathbb C})\ni (\zeta,w)
         \mapsto w\in{\mathbf G}_n;
                      \label{eq: standard Zk covering of Cn(C)}\\
&&\hskip-50pt
{\aligned
E_{2(n-1)}^n({\mathbb{CP}}^1)=\{z=(z_0,...&,z_n)\in{\mathbb C}^{n+1}\,|\ 
D_n(z)=1\},\label{eq: standard 2(n-1) covering of Cn(CP1)}\\
&p_{2(n-1)}\colon E_{2(n-1)}^n({\mathbb{CP}}^1)
\ni z\mapsto\psi([x:y],z)\in{\mathcal F}^n.
\endaligned}
\end{eqnarray}
For a non-trivial divisor $m$ of $2(n-1)$ set $k=2(n-1)/m$;
the cyclic subgroup 
${\mathbb Z}/m{\mathbb Z}\subset{\mathbb Z}/2(n-1){\mathbb Z}$
acts freely on $E_{2(n-1)}^n({\mathbb{CP}}^1)$ via 
\begin{equation}\label{eq: Zm action in E(2(n-1))}
\aligned
r(z_0,...,z_n)=&(e^{i\pi r/(n-1)}z_0,...,e^{i\pi r/(n-1)}z_n)\\
&\text{for all} \ \ r\in{\mathbb Z}/m{\mathbb Z} \ \ 
\text{and all}\ \ z=(z_0,...,z_n)\in E_{2(n-1)}^n({\mathbb{CP}}^1)\,,
\endaligned
\end{equation}
and the quotient space 
\begin{equation}\label{eq: }
E_k^n({\mathbb{CP}}^1)=[E_{2(n-1)}^n({\mathbb{CP}}^1)]/
({\mathbb Z}/m{\mathbb Z})
\end{equation}
is an unbranched holomorphic ${\mathbb Z}/k{\mathbb Z}$ covering
of ${\mathcal C}^n({\mathbb{CP}}^1)$, which is also referred to
as the standard one.
\vskip0.2cm

\noindent Since the abelianization $B_n(X)/B_n'(X)$ is a cyclic
group, the equivalence classes of cyclic coverings over
${\mathcal C}^n(X)$ are in the one-to-one correspondence
with the subgroups of $B_n(X)/B_n'(X)$. This implies
the following result.  

\begin{Theorem}\label{Thm: Cyclic coverings of Cn(X)}
Every unbranched cyclic covering of ${\mathcal C}^n(X)$ is equivalent to 
one of the standard cyclic coverings exhibited above.
\hfill $\square$
\end{Theorem}

\begin{Remark}\label{Rmk: Z cover of Cn(C)}
{\bfit a}$)$ The ${\mathbb Z}$ cover $E_\infty^n({\mathbb C})$
defined by (\ref{eq: standard Z covering}) is in fact biholomorphic
to the algebraic manifold 
$\widetilde{E}_\infty^n({\mathbb C})
={\mathbb C}\times\{v\in{\mathbb C}^n\,|\ d_n(v)=1\}$.
The isomorphism is given by 
$$
\aligned
E_\infty^n({\mathbb C})\ni(\zeta;\,&w_1,w_2,...,w_n)\\
&\mapsto 
(\zeta;e^{-\zeta/n(n-1)}w_1,e^{-2\zeta/n(n-1)}w_2,...,
e^{-n\zeta/n(n-1)}w_n)
\in\widetilde{E}_\infty^n({\mathbb C})
\endaligned
$$
and the covering map $\widetilde{E}_\infty^n({\mathbb C})
\to{\mathbf G}_n={\mathcal C}^n({\mathbb C})$ is defined by
$$
\aligned
\widetilde{E}_\infty^n({\mathbb C})&\ni (\zeta;v_1,v_2,...,v_n)\\
&\mapsto t^n+e^{\zeta/n(n-1)}v_1t^{n-1}
+e^{2\zeta/n(n-1)}v_2t^{n-2}+...+e^{n\zeta/n(n-1)}v_n
\in{\mathbf G}_n\,.
\endaligned
$$
Moreover, $\pi_1(E_\infty^n({\mathbb C}))
=\pi_1(\widetilde{E}_\infty^n({\mathbb C}))\cong B_n'$ and
$E_\infty^n({\mathbb C})$ is the maximal abelian cover of
${\mathcal C}^n({\mathbb C})$.
\vskip0.2cm

{\bfit b}$)$ If $k\equiv 1\! \mod n(n-1)$ then $E_k^n({\mathbb C})$ is
isomorphic to ${\mathcal C}^n({\mathbb C})$.
The biregular isomorphism ${\mathbf G}_n\ni w=(w_1,...,w_n)
\mapsto (\zeta,v)
=(\zeta;v_1,...,v_n)\in E_k^n({\mathbb C})$
is given by
\begin{equation}\label{eq: explicit isomorphism Cn(C) to E(mn(n-1)+1)}
\zeta=d_n(w), \ \ v_1=(d_n(w))^mw_1, \ \ v_2=(d_n(w))^{2m}w_2,...,
v_n=(d_n(w))^{nm}w_n\,,
\end{equation}
where $m=(k-1)/n(n-1)$. In other words, for $m\in{\mathbb N}$
and $k=mn(n-1)+1$ the regular endomorphism\footnote{Here
$d_n(Q)=d_n(w(Q))$ denotes the discriminant of
the polynomial $p_n(t;w(Q))=t^n+w_1(Q)t^{n-1}+...+w_n(Q)
=\prod_{\zeta\in Q} (t-\zeta)$ corresponding to a point
$Q\in{\mathcal C}^n(\mathbb C)$.}
\begin{equation}\label{eq: mn(n-1)+1 covering of Cn(C)}
\aligned
p\colon{\mathcal C}^n({\mathbb C})\ni Q&=\{q_1,...,q_n\}\\
&\mapsto \{(d_n(Q))^m q_1,...,(d_n(Q))^m q_n\}=(d_n(Q))^mQ
\in{\mathcal C}^n({\mathbb C})
\endaligned
\end{equation}
is an unbranched analytic cyclic covering of order $k$, which is
equivalent to the standard cyclic ${\mathbb Z}/k{\mathbb Z}$
covering defined by (\ref{eq: standard Zk covering of Cn(C)}).
The induced monomorphism $p_*\colon B_n\to B_n$ is given
by $\sigma_i\mapsto\sigma_i\cdot c^m$ ($i=1,...,n-1$), where
$c=a^n=(\sigma_1\cdots \sigma_{n-1})^n$ generates
the center $CB_n\cong{\mathbb Z}$ of the group $B_n$.
The image of $p_*$ consists of
all $g\in B_n$ such that $\chi(g)\in (mn(n-1)+1){\mathbb Z}$,
where $\chi\colon B_n\to{\mathbb Z}$ is the abelianization epimorphism.
\vskip0.2cm

{\bfit c}$)$ $\pi_1(E_{2(n-1)}^n({\mathbb{CP}}^1))\cong B_n'(S^2)$
and $E_{2(n-1)}^n({\mathbb{CP}}^1)$ is the
maximal abelian cover of ${\mathcal C}^n({\mathbb{CP}}^1)$.
\hfill $\bigcirc$
\end{Remark}

\noindent In the rest of this section, assuming that $n>6$ and $k\le 2n$,
we describe all connected unbranched non-cyclic $k$ covers
of ${\mathcal C}^n(X)$. The monodromy concept provides
the one-to-one correspondence between the equivalence classes of such
covers and the conjugacy classes of transitive non-cyclic homomorphisms
$B_n(X)\to{\mathbf S}(k)$, which are described below.

\subsection{Transitive homomorphisms $B_n(X)\to {\mathbf S}(k)$, \
$k\le 2n$}
\label{Ss: Transitive homomorphisms Bn(X) to S(k), k<=2n}

\noindent Artin Theorem quoted in Section \ref{Ss: Proof of Surjectivity Theorem} and 
Theorem \ref{Thm: Bn(X) to S(k) for n>k} give the complete description
of transitive homomorphisms $B_n(X)\to {\mathbf S}(k)$ for $n\ge k$.
Here we present the case $n<k\le 2n$.

\begin{Definition}\label{Def: standard homomorphisms Bn(X) to S(2n)}
For $i=1,...,n-1$ set
$$
\aligned
&\varphi_1(\sigma_i) = \underbrace{(2i-1,2i+2,2i,2i+1)}_{\text {$4$-cycle}};\\
&\varphi_2(\sigma_i)=(1,2)\cdot\cdot\cdot (2i-3,2i-2)
\underbrace{(2i-1,2i+1)(2i,2i+2)}_{\text {two transpositions}}\times\\
&\hskip8.8cm
\times
(2i+3,2i+4)\cdot\cdot\cdot(2n-1,2n);\\
&\varphi_3(\sigma_i)= (1,2)\cdot\cdot\cdot (2i-3,2i-2)
\underbrace{(2i-1,2i+2,2i,2i+1)}_{\text {$4$-cycle}}\times\\
&\hskip8.8cm
\times(2i+3,2i+4)\cdot\cdot\cdot (2n-1,2n).\\
\endaligned
$$
These formulae determine three non-cyclic transitive
homomorphisms $\varphi_j\colon B_n\to{\mathbf S}(2n)$ ($j=1,2,3$).
Since $\varphi_2(\sigma_1\cdots\sigma_{n-1}
\sigma_{n-1}\cdots\sigma_1)=1$, the homomorphism $\varphi_2$
gives also rise to the homomorphism
$\varphi_2\colon B_n(S^2)\to{\mathbf S}(2n)$.\footnote{This is not
the case for $\varphi_1$ and $\varphi_3$.}
All these homomorphisms are referred to as the
{\em standard non-cyclic transitive homomorphisms}
$B_n(X)\to {\mathbf S}(2n)$. 
\hfill $\bigcirc$
\end{Definition}

\noindent For the braid group $B_n$ the next
theorem has been proved in \cite{Lin96b} (part ({\bfit b}) for $n>8$ only;
the cases $n=7$ and $n=8$ were treated by S. Orevkov \cite{Ore98}).
For the sphere braid group $B_n(S^2)$ both ({\bfit a}) and ({\bfit b})
are easy consequences of the corresponding properties of $B_n$ and   
the presentation of $B_n(S^2)$ exhibited in Section
\ref{Ss: Canonical presentations of Bn and Bn(S2)}.
 
\begin{Theorem}\label{Thm: Bn(X) to S(k), 6<n<k<=2n}  
{\bfit a}$)$ Let \ $6<n<k<2n$. \ Then every transitive
homomorphism $B_n(X)\to{\mathbf S}(k)$ is cyclic.
\vskip0.1cm
 
{\bfit b}$)$ For $n>6$ every non-cyclic transitive homomorphism
$B_n(X)\to{\mathbf S}(2n)$ is conjugate to one of the standard
homomorphisms.
\end{Theorem}

\begin{Remark}
\label{Rmk: Non-cyclic homomorphisms B_n to S(rn)} 
For any $n\ge 3$ and any $r\ge 2$ there exist non-cyclic
transitive homomorphisms $B_n\to{\mathbf S}(rn)$; let us present
some of them.
\vskip0.2cm

\noindent For any $l\in{\mathbb Z}$ and $t\in{\mathbb N}$
we denote by $|l|_t$ the image of $l$ along the natural epimorphism
${\mathbb Z}\to{\mathbb Z}/t{\mathbb Z}$.
We identify the set $\boldsymbol \Delta_{n}
=\{1,...,n\}$ to the group ${\mathbb Z}/n{\mathbb Z}$
via $\{1,...,n\}\ni m\leftrightarrow |m|_n-1$
and identify the symmetric groups ${\mathbf S}(n)$ and
${\mathbf S}({\mathbb Z}/n{\mathbb Z})$.
Furthermore, we identify the set $\boldsymbol \Delta_{rn}
=\{1,...,rn\}$ to the direct product ${\mathcal D}(r,n)
=({\mathbb Z}/r{\mathbb Z})\times ({\mathbb Z}/n{\mathbb Z})$
via
$$
\aligned
&\boldsymbol \Delta_{rn}\ni m\mapsto (R(m),N(m))
\in ({\mathbb Z}/r{\mathbb Z})
\times ({\mathbb Z}/n{\mathbb Z})\,,\\
&\text{where} \ R(m)=|m-1|_r\in{\mathbb Z}/r{\mathbb Z}\ \ \text{and} \ \
N(m)=|(m-1-R(m))/r|_n\in{\mathbb Z}/n{\mathbb Z}
\endaligned
$$
and regard ${\mathbf S}(rn)$ as the symmetric group 
${\mathbf S}({\mathcal D}(r,n))$.
\vskip0.2cm

\noindent For each $x,y\in{\mathbb Z}/r{\mathbb Z}$
and any $i=1,...,n-1$ let $s_i=s_{x,y;i}
\in{\mathbf S}({\mathcal D}(r,n))$
denote the permutation that acts
on elements $(R,N)\in{\mathcal D}(r,n)
=({\mathbb Z}/r{\mathbb Z})\times ({\mathbb Z}/n{\mathbb Z})$ as follows:
$$
s_i(R,N) = \left\{
\aligned 
&(R+y,\,N) \hskip30pt {\text{if}} \ \ N\ne i-1,i;\\ 
&(R,\,N+1) \hskip30pt {\text{if}} \ \ N=i-1;\\
&(R+x,\,N-1) \hskip5pt \ {\text{if}} \ \ N=i;
\endaligned
\right.\qquad (1\le i\le n-1)\,.
$$
It is easily shown that the map $\sigma_i\to s_{x,y;i}$, $i=1,...,n-1$,
extends to a uniquely defined non-cyclic homomorphism
$\varphi_{x,y}\colon B_n\to{\mathbf S}({\mathcal D}(r,n))={\mathbf S}(rn)$.
This homomorphism $\varphi_{x,y}$ is transitive if and only
if the elements $x,y\in{\mathbb Z}/r{\mathbb Z}$ generate the whole
group ${\mathbb Z}/r{\mathbb Z}$, or, which is the same,
if and only if $x$,$y$ and $r$ are mutually co-prime. 
However $\varphi_{x,y}$ can never be
primitive. 
For more details see \cite{Lin96b}.
\hfill $\bigcirc$
\end{Remark}

\subsection{Non-cyclic coverings of ${\mathcal C}^n(X)$}
\label{Ss: Non-cyclic coverings of Cn(X)}
We start with some simple examples.

\begin{Example}\label{Ex: standard n covering of Cn(X)}
The unbranched connected analytic covering
\begin{equation}\label{eq: standard n cover of Cn(X)}
{\mathcal E}_n(X)=\{(\zeta,Q)\in X\times{\mathcal C}^n(X)\,|\
\zeta\in Q\}\,,\ \
\nu_n\colon{\mathcal E}_n(X)\ni (\zeta,Q)\mapsto Q\in{\mathcal C}^n(X)
\end{equation}
is referred to as the {\em standard non-cyclic $n$ covering of
${\mathcal C}^n(X)$}. In the case $X={\mathbb C}$ it may also be described as
\begin{equation}
{\mathcal E}_n({\mathbb C})=\{(\zeta,w)\in{\mathbb C}
\times{\mathbf G}_n\,|\
p_n(\zeta,w)=0\}\,,\ \
\nu_n\colon{\mathcal E}_n({\mathbb C})\ni (\zeta,w)\mapsto w
\in{\mathbf G}_n\,, \label{eq: standard n cover of Gn}
\end{equation}
and in the case $X={\mathbb{CP}}^1$ as
\begin{equation}\label{eq: standard n covering of Cn(CP1)}
\aligned
&{\mathcal E}_n({\mathbb{CP}}^1)=\{([x:y],z)\in{\mathbb{CP}}^1
\times{\mathcal F}^n\,|\ \varphi([x:y],z)=0\}\,,\\
&\nu_n\colon{\mathcal E}_n({\mathbb{CP}}^1)\ni ([x:y],z)\mapsto 
z=[z_0:...:z_n]\in{\mathcal F}^n\subset{\mathbb{CP}}^n\,.
\endaligned
\end{equation}
The monodromy representation corresponding to the covering
(\ref{eq: standard n cover of Cn(X)}) coincides with the standard
epimorphism $\mu\colon B_n(X)\to{\mathbf S}(n)$.
\hfill $\bigcirc$
\end{Example}

\begin{Example}\label{Rmk: binomial coverings}
Here is a generalization of the previous example.
For $m\le n$ set
\begin{equation}\label{eq: binom(n,m) cover of Cn(X)}
\aligned
&{\mathcal E}_{\binom{n}{m}}(X)
=\{(O',Q)\in {\mathcal C}^m(X)\times{\mathcal C}^n(X)\,|\
Q'\subseteq Q\}\,,\\
&\nu_{\binom{n}{m}}
\colon{\mathcal E}_{\binom{n}{m}}(X)\ni (Q',Q)\mapsto Q
\in{\mathcal C}^n(X)\,.
\endaligned
\end{equation}
This is an unbranched connected non-cyclic $\binom{n}{m}$-covering
of ${\mathcal C}^n(X)$. Its monodromy representation
is the composition 
$B_n(X)\stackrel{\mu}{\longrightarrow}{\mathbf S}(n)
\stackrel{\tau^n_m}{\longrightarrow}{\mathbf S}(\binom{n}{m})$,
where $\tau^n_m$ stands for the transitive representation
corresponding to the natural ${\mathbf S}(n)$ action on the set
of all $m$ point subsets $M\subset\{1,...,n\}$. 
\hfill $\bigcirc$
\end{Example}

\begin{Example}\label{Ex: standard 2 coverings of En(X)}
Set $p_n'(t,w)=nt^{n-1}+(n-1)w_1t^{n-2}+...+w_{n-1}$.
Define the following non-trivial unbranched analytic $2$ coverings
of the space ${\mathcal E}_n({\mathbb C})$: 
\begin{equation}\label{eq: standard 2 covers of En(C)}
\aligned
&\nu^{(1)}\colon{\mathcal E}_n^{(1)}({\mathbb C})
=\left\{\left.(\xi;(\zeta,w))\in{\mathbb C}^*
\times{\mathcal E}_n({\mathbb C})\right|\,\xi^2=p_n'(\zeta,w)\right\} 
                              \to{\mathcal E}_n({\mathbb C}),\\
&\nu^{(2)}\colon{\mathcal E}_n^{(2)}({\mathbb C})
=\left\{\left.(\xi;(\zeta,w))\in{\mathbb C}^*\times
{\mathcal E}_n({\mathbb C})\right|\, \xi^2=d_n(w) \right\}
                              \to{\mathcal E}_n({\mathbb C}),\\
&\nu^{(3)}\colon{\mathcal E}_n^{(3)}({\mathbb C})
=\left\{\left.(\xi;(\zeta,w))\in{\mathbb C}^*\times
{\mathcal E}_n({\mathbb C})\right|\,\xi^2=d_n(w)p_n'(\zeta,w)\right\}
                              \to{\mathcal E}_n({\mathbb C}),
\endaligned
\end{equation}
all with the same natural projection $(\xi;(\zeta,w))
\mapsto (\zeta,w)$. Any non-trivial $2$ covering of 
${\mathcal E}_n({\mathbb C})$ is equivalent to one of the coverings
$\nu^{(j)}\colon{\mathcal E}_n^{(j)}({\mathbb C})
\to{\mathcal E}_n({\mathbb C})$, $j=1,2,3$ (see \cite{Lin96b}
for the proof; we will not use this fact in what follows).
\vskip0.2cm

\noindent Let us construct a non-trivial unbranched analytic $2$ covering
of the space ${\mathcal E}_n({\mathbb{CP}}^1)$. To this end,
take the standard ${\mathbb Z}/2n(n-1){\mathbb Z}$ covering 
$p_{2(n-1)}\colon E_{2(n-1)}^n({\mathbb{CP}}^1)\to{\mathcal F}^n$
described in (\ref{eq: standard 2(n-1) covering of Cn(CP1)}).
The subgroup ${\mathbb Z}/n(n-1){\mathbb Z}
\subset{\mathbb Z}/2n(n-1){\mathbb Z}$ acts freely on 
$E_{2(n-1)}^n({\mathbb{CP}}^1)$ by the transformations
$r(z_0,...,z_n)=(e^{i\pi r/(n-1)}z_0,...,e^{i\pi r/(n-1)}z_n)$, \
where $r\in{\mathbb Z}/n(n-1){\mathbb Z}$ and
$z=(z_0,...,z_n)\in E_{2(n-1)}^n({\mathbb{CP}}^1)$.
This provides us with the non-trivial unbranched analytic $2$ covering
\begin{equation}\label{eq: 2 cover of Cn(CP1)=Fn}
\lambda\colon E_2^n({\mathbb{CP}}^1)=[E_{2(n-1)}^n({\mathbb{CP}}^1)]/
({\mathbb Z}/n(n-1){\mathbb Z})\to{\mathcal F}^n
={\mathcal C}^n({\mathbb{CP}}^1)\,.
\end{equation}
The desired $2$-covering  
\begin{equation}\label{eq: 2 cover of En(CP1)}
\nu^{(2)}
=\nu_n^*(\lambda)\colon{\mathcal E}_n^{(2)}({\mathbb{CP}}^1)
\to{\mathcal E}_n({\mathbb{CP}}^1)
\end{equation}
is obtained as the pullback $\nu^{(2)}=\nu_n^*(\lambda)$
of the covering (\ref{eq: 2 cover of Cn(CP1)=Fn}) along the mapping 
$\nu_n\colon{\mathcal E}_n({\mathbb{CP}}^1)
\to{\mathcal C}^n({\mathbb{CP}}^1)$
defined by (\ref{eq: standard n covering of Cn(CP1)}).
\vskip0.2cm

\noindent The $2$-coverings  
$\nu^{(j)}\colon{\mathcal E}_n^{(j)}({\mathbb C})
\to{\mathcal E}_n({\mathbb C}) $ ($j=1,2,3$) 
and $\nu^{(2)}\colon{\mathcal E}_n^{(2)}({\mathbb{CP}}^1)
\to{\mathcal E}_n({\mathbb{CP}}^1)$ defined by
(\ref{Ex: standard 2 coverings of En(X)}) and 
(\ref{eq: 2 cover of En(CP1)}) are referred to as the
{\em standard $2$-coverings of} ${\mathcal E}_n(X)$.
\hfill $\bigcirc$
\end{Example}

\begin{Example}\label{Ex: standard non-cyclic 2n coverings of En(X)}
We define $2n$-coverings
$\mu_{2n}^{(j)}(X)\colon{\mathcal E}_n^{(j)}(X)\to{\mathcal C}_n(X)$
as compositions of the $2$-coverings 
${\nu}^{(j)}\colon{\mathcal E}_n^{(j)}(X)\to{\mathcal E}_n(X)$
and the $n$-covering $\nu_n\colon{\mathcal E}_n(X)\to{\mathcal C}_n(X)$
exhibited in Examples \ref{eq: standard 2 covers of En(C)}
and \ref{Ex: standard n covering of Cn(X)} respectively. That is,
\begin{equation}\label{eq: standard 2n coverings of En(X)}
\aligned
&\mu_{2n}^{(j)}({\mathbb C})=\nu_n\circ{\nu}^{(j)}
\colon{\mathcal E}_n^{(j)}({\mathbb C})
\overset{{\nu}^{(j)}}\longrightarrow
{\mathcal E}_n({\mathbb C})\overset{\nu_n}
\longrightarrow{\mathcal C}_n({\mathbb C})\,,\quad j=1,2,3\,,\\
&\mu_{2n}^{(2)}({\mathbb{CP}}^1)=\nu_n\circ\nu^{(2)}
\colon{\mathcal E}_n^{(2)}({\mathbb{CP}}^1)
\overset{\nu^{(2)}}\longrightarrow
{\mathcal E}_n({\mathbb{CP}}^1)\overset{\nu_n}
\longrightarrow{\mathcal C}_n({\mathbb{CP}}^1)\,.
\endaligned
\end{equation}
These covering are referred to as the {\em standard non-cyclic
$2n$-coverings of} ${\mathcal C}_n(X)$.
\vskip0.2cm

\noindent The monodromy representations
$\varphi_j=\mu_{2n}^{(j)*}\colon B_n(X)\to{\mathbf S}(2n)$
corresponding to the above $2n$-coverings
of ${\mathcal C}_n(X)$ are conjugate to the eponymous standard
non-cyclic transitive homomorphisms
$\varphi_j\colon B_n(X)\to{\mathbf S}(2n)$
introduced in Definition \ref{Def: standard homomorphisms Bn(X) to S(2n)}.
\hfill $\bigcirc$
\end{Example}
 
\noindent The description of transitive homomorphisms
$B_n(X)\to{\mathbf S}(k)$ provided by Theorem
\ref{Thm: Bn(X) to S(k) for n>k}, Artin Theorem (quoted in
Section \ref{Ss: Proof of Surjectivity Theorem}) and
Theorem \ref{Thm: Bn(X) to S(k), 6<n<k<=2n}, together 
with Theorem \ref{Thm: Cyclic coverings of Cn(X)}
and the latter example, leads the following complete description
of unbranched connected $k$ coverings of ${\mathcal C}^n(X)$
for $k\le 2n$.

\begin{Theorem}\label{Thm: non-cyclic coverings of Cn(X)}
Let ${\mathcal E}=\{p\colon E\to{\mathcal C}^n(X)\}$
be an unbranched connected $k$-covering.
\vskip0.2cm

{\bfit a}$)$ If $n>\max\{k,4\}$ then ${\mathcal E}$ is equivalent
to one of the cyclic $k$-coverings exhibited
in Section {\rm \ref{Ss: Cyclic coverings}}.
\vskip0.2cm

{\bfit b}$)$ If $6<n\le k\le 2n$ then either   
${\mathcal E}$ is cyclic and then equivalent
to one of the cyclic $k$-coverings exhibited
in Section {\rm \ref{Ss: Cyclic coverings}} or
$k=n$ and ${\mathcal E}$ is equivalent to the
standard non-cyclic $n$-covering 
defined by {\rm(\ref{eq: standard n cover of Cn(X)})}
or $k=2n$ and ${\mathcal E}$ is equivalent to one of the 
standard non-cyclic $2n$-coverings exhibited in
Example {\rm \ref{Ex: standard non-cyclic 2n coverings of En(X)}}.
In the latter case ${\mathcal E}$ splits to the composition
of a $2$-covering and the standard non-cyclic $n$-covering.
\hfill $\square$
\end{Theorem}


\section{Configuration spaces of small dimension and some other
examples}
\label{Sec: Configuration spaces of small dimension and some other
examples}

\noindent Configuration spaces of small dimension admit certain
exceptional morphisms prohibited by results stated above.
In this section we present some of such examples. We will also
describe the Eisenstein automorphism mentioned in
Section \ref{Ss: Some historical remarks}.
 
\subsection{Ferrari map}\label{Ss: Ferrari map}
The Ferrari map (L. Ferrari, 1545), which carries each polynomial
of degree $4$ to its cubic resolvent, provides the oldest and undoubtedly 
the most famous example of a non-cyclic morphism  of configuration spaces.
This morphism $Fe\colon{\mathcal C}^4({\mathbb C})
\to{\mathcal C}^3({\mathbb C})$ 
arises as follows. For every four distinct
$q_1,q_2,q_3,q_4\in{\mathbb C}$
the points $z_1=(q_1-q_2-q_3+q_4)^2/4$, $z_2=(q_1-q_2+q_3-q_4)^2/4$,
$z_3=(q_1+q_2-q_3-q_4)^2/4$ are distinct and a permutation
of $q_1,q_2,q_3,q_4$ leads to a permutation of $z_1,z_2,z_3$.
Thus, we obtain polynomial maps
$fe\colon{\mathcal C}_o^4({\mathbb C})
\to{\mathcal C}_o^3({\mathbb C})$ and
$Fe\colon{\mathcal C}^4({\mathbb C})\to{\mathcal C}^3({\mathbb C})$.
The morphism $Fe$ induces the epimorphism of braid groups
$Fe_*\colon B_4\to B_3$, which carries both $\sigma_1$,$\sigma_3$
to $\sigma_1$ and keeps $\sigma_2$.

\subsection{Disjoint morphisms 
${\mathcal C}^3({\mathbb C})\to{\mathcal C}^6({\mathbb C})$
and ${\mathcal C}^3({\mathbb C})\to{\mathcal C}^9({\mathbb{CP}}^1)$}
\label{Ss: G3 to G6}
For distinct $q_1,q_2,q_3\in{\mathbb C}$ define six points
$\mu_{1}^\pm,\mu_{2}^\pm,\mu_{3}^\pm\in{\mathbb C}$ by the
quadratic equations
$$
(\mu_1^\pm-q_1)^2=(q_1-q_2)(q_1-q_3), \ 
(\mu_2^\pm-q_2)^2=(q_2-q_3)(q_2-q_1), \
(\mu_3^\pm-q_3)^2=(q_3-q_1)(q_3-q_2).
$$
%
All nine points $q_i$, $\mu_j^\pm$, \ $1\le i,j\le 3$,
are distinct and a permutation of $q$'s leads to a permutation
of $\mu$'s. Thus, we obtain a disjoint polynomial map
$$
L\colon{\mathcal C}^3({\mathbb C})\ni \{q_1,q_2,q_3\}
\mapsto\{\mu_1^+,\mu_1^-,\mu_2^+,\mu_2^-,\mu_3^+,\mu_3^-\}
\in{\mathcal C}^6({\mathbb C})\,,
$$
which may also be viewed as a polynomial map 
$L\colon{\mathbf G}_3\to {\mathbf G}_6$ that carries
any polynomial $p_3(t,z)=t^3+z_1t^2+z_2t+z_3\in{\mathbf G}_3$
to the polynomial $p_6(t,L(z))\in{\mathbf G}_6$, which is
co-prime to $p_3(t,z)$ and whose coefficients
$L_1(z),\ldots ,L_6(z)$ are as follows:
$$
\aligned
L_1(z) &= 2z_1;\\
L_2(z) &= 5z_2;\\
L_3(z) &= 20z_3;\endaligned
\qquad
\aligned
L_4(z) &= 20z_1z_3-5z^2_2;\\
L_5(z) &= 8z^2_1z_3 - 2z_1z^2_2 - 4z_2z_3; \\
L_6(z) &= 4z_1z_2z_3 - z^3_2 - 8z^2_3.
\endaligned
$$
The first three of these formulae show that $L$ is an
embedding; one can check that it is non-cyclic
and that the following relation between the discriminants holds:
$$
d_6(L(z)) = -4^9\cdot [d_3(z)]^5.
$$
The following simple geometric interpretation of the map $L$
due to E. A. Gorin (personal communication).
Let $\Delta (Q)\subset{\mathbb C}$ be the
triangle with the vertices at the points $q_1,q_2,q_3$
(it may degenerate to a straightline segment with a distinguished
interior point). For each vertex $A=q_i$,
let $b$ and $c$ denote the lengths of the sides
of $\Delta (Q)$ intersecting at $A$, and let $\ell(A)$
be the bisector line of the inner angle of
$\Delta (Q)$ at $A$.
Take two points on $\ell(A)$ which are 
$\sqrt{bc}$-distant
from $A$. Thus, we obtain six points, which coincide
with the six points $\mu_i^\pm$ defined above
(see the picture below).
\vskip0.2cm

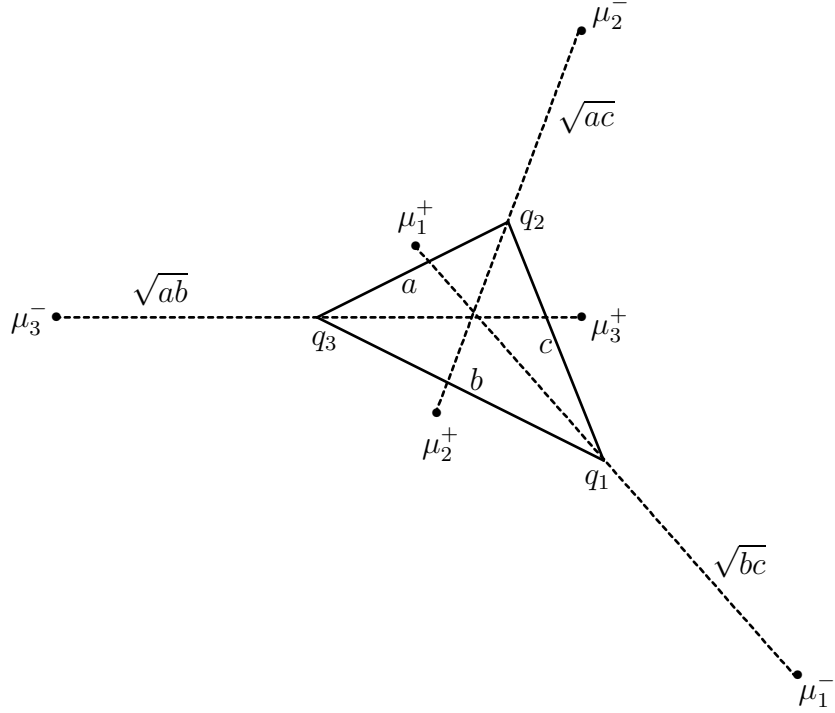
\begin{figure}[tbh]
\mbox{
\hskip1cm\input picture
}
\end{figure}
\vskip0.1cm

\noindent On the other hand, Yoel Feler has recently noticed that the same
morphism $L$ may be defined via the formulae  
$$
\left(\frac{\mu_1^\pm - q_2}{\mu_1^\pm - q_3}\right)^2
=\frac{q_1 - q_2}{q_1 - q_3}\,, \ \
\left(\frac{\mu_2^\pm - q_3}{\mu_2^\pm - q_1}\right)^2
=\frac{q_2 - q_3}{q_2 - q_1}\,, \ \
\left(\frac{\mu_3^\pm - q_1}{\mu_3^\pm - q_2}\right)^2
=\frac{q_3 - q_1}{q_3 - q_2}\,.
$$
This observation let him to construct a disjoint\footnote{Meaning that
each set $Q\in{\mathcal C}^3({\mathbb C})$ is disjoint to its 
image $F(Q)\in{\mathcal C}^9({\mathbb{CP}}^1)$.} morphism
$$
F\colon{\mathcal C}^3({\mathbb C})\ni\{q_1,q_2,q_3\}\mapsto
\{\mu_1^{(1)},\mu_1^{(2)},\mu_1^{(3)},
\mu_2^{(1)},\mu_2^{(2)},\mu_2^{(3)},
\mu_3^{(1)},\mu_3^{(2)},\mu_3^{(3)}\}
\in{\mathcal C}^9({\mathbb{CP}}^1)\,,
$$
where $\mu_i^{(j)}$ ($1\le i,j\le 3$) are defined by the cubic
equations
$$
\aligned
\left(\frac{\mu_1^{(j)} - q_2}{\mu_1^{(j)} - q_3}\right)^3
= \left(\frac{q_1 - q_2}{q_1 - q_3}\right)^2\,,\\
\left(\frac{\mu_2^{(j)} - q_3}{\mu_2^{(j)} - q_1}\right)^3
= \left(\frac{q_2 - q_3}{q_2 - q_1}\right)^2\,,\\
\left(\frac{\mu_3^{(j)} - q_1}{\mu_3^{(j)} - q_2}\right)^3
= \left(\frac{q_3 - q_1}{q_3 - q_2}\right)^2\,.
\endaligned
$$
In terms of polynomials 
$p(t,q)=(t-q_1)(t-q_2)(t-q_3)\in{\mathbf G}_3$
and binary projective forms $\phi([x:y],F(q))\in{\mathcal F}^9$,
the map $F$ may be defined by the formula
$$
\aligned
\phi([x:y],F(q))
   =&((x-yq_1)^3(q_2-q_3)^2-(x-yq_2)^3(q_3-q_1)^2)\\
&\times((x-yq_2)^3(q_3-q_1)^2-(x-yq_3)^3(q_1-q_2)^2)\\
&\times((x-yq_3)^3(q_1-q_2)^2-(x-yq_1)^3(q_2-q_3)^2)\,.
\endaligned
$$
A straightforward computation of the discriminant shows that 
$$
D_\phi=3^{27}[(q_1-q_2)(q_2-q_3)(q_3-q_1)]^{56}\,.
$$

\subsection{Eisenstein automorphism}
\label{Sec: Eisenstein automorphism}
\index{Eisenstein's automorphism!original form\hfill}
Here we present three different constructions
of the Eisenstein automorphism: the original one \cite{Eis1844},
the Cayley construction (1845) and also the new one
in the form of a tame morphism 
${\mathcal C}^3({\mathbb{CP}}^1)\to{\mathcal C}^3({\mathbb{CP}}^1)$
found by Yoel Feler in 2001.
\vskip0.2cm

\noindent The discriminant $D_\phi$ of a binary cubic form 
$\phi(x,y)=z_0x^3 + 3z_1 x^2 y + 3z_2 x y^2 + z_3 y^3$
is given by
$$
D_\phi=z_0^2z_3^2-3z_1^2z_2^2-6z_0z_1z_2z_3+4z_0z_2^3+4z_1^3z_3\,.
$$
For every binary cubic form $\phi$, Eisenstein considers the binary cubic
form $({\mathcal E}\phi)(x,y)=w_0x^3 + 3w_1 x^2 y + 3w_2 x y^2 + w_3 y^3$
whose coefficients $w_i$ are homogeneous cubic polynomials in
$z_0,...,z_3$ defined by\footnote{Eisenstein 
has not pointed out that the discriminant $D_\phi$ plays the part
of a potential function for the coefficients of the form
${\mathcal E}\phi$; I have learned this fact from M. Gizatullin.} 
$$
\aligned
w_0&=\frac{1}{2}\frac{\partial D_\phi}{\partial z_3}
=2z_1^3-3z_0z_1z_2+z_0^2z_3\,,\quad
w_1=\frac{1}{6}\frac{\partial D_\phi}{\partial z_2}
=2z_0z_2^2-z_0z_1z_3-z_1^2z_2\,,\\
w_2&=\frac{1}{6}\frac{\partial D_\phi}{\partial z_1}
=2z_1^2z_3-z_0z_2z_3+z_1z_2^2\,,\quad
w_3=\frac{1}{2}\frac{\partial D_\phi}{\partial z_0}
=2z_2^3-3z_1z_2z_3+z_0z_3^2\,.\\
\endaligned
$$
He proves the identities $\mathcal E({\mathcal E}\phi)
=(D_\phi)^2\cdot \phi$ and
$$
\aligned
D_{{\mathcal E}\phi}&
=w_0^2w_3^2-3w_1^2w_2^2-6w_0w_1w_2w_3+4w_0w_2^3+4w_1^3w_3\\
&=(z_0^2z_3^2-3z_1^2z_2^2-6z_0z_1z_2z_3+4z_0z_2^3+4z_1^3z_3)^3
=(D_\phi)^3\,.
\endaligned
$$
Thus Eisenstein's correspondence $f\mapsto{\mathcal E}\phi$
determines a {\em birational involution} of the space 
${\mathbb{CP}}^3$ of all projective cubic binary forms.  
Moreover, its restriction to the space ${\mathcal F}_3$
of non-degenerate cubic forms is an {\em involutive biregular
automorphism} of the latter space: 
\begin{equation}\label{Eisenstein automorphism}
\CD
{\mathcal F}_3@. @. \\
@V{\mathcal E}V{\cong}V @/SE//{\id}/ @. \\ 
{\mathcal F}_3@>{\cong}>{\mathcal E}
>{\mathcal F}_3\,.
\endCD
\end{equation}
\subsection{Cayley's form of the Eisenstein automorphism}
\label{Ss: Cayley's form of the Eisenstein automorphism}
\index{Eisenstein's automorphism!Cayley's form\hfill} 
Take a complex binary cubic form\footnote{There
is a minor difference between the automorphism
described here and the original Eisenstein form: one has just to
interchange $x$ and $y$ in one of the examples.}
$$
\phi(x,y)=z_0x^3 + 3z_1 x^2 y + 3z_2 x y^2 + z_3 y^3\,. 
$$
Its Hessian 
$$
\psi(x,y)=\det\begin{pmatrix} 
{\partial^2\phi}/{\partial x^2}\,, & 
{\partial^2\phi}/{\partial x\partial y} \\
{\partial^2\phi}/{\partial y\partial x}\,, & {\partial^2\phi}/{\partial y^2}
\end{pmatrix}
$$
is a binary quadratic form. The Jacobian
$$
J(x,y)=\det \begin{pmatrix} 
{\partial\phi}/{\partial x}\,, & {\partial\phi}/{\partial y} \\
{\partial\psi}/{\partial x}\,, & {\partial\psi}/{\partial y}
\end{pmatrix}
$$
is a binary cubic form. 
\medskip

\noindent {\bfit Eisenstein's map:} \qquad 
${\mathcal E}\colon\phi\mapsto{\mathcal E}\phi=J$.

\subsection{Feler's form of the Eisenstein automorphism}
\label{Ss: Feler's form of the Eisenstein automorphism}
\index{Eisenstein's automorphism!Feler's form\hfill}
For any $Q\in{\mathcal C}^3({\mathbb{CP}}^1)$
let $z(Q)=[z_0(Q):z_1(Q):z_2(Q):z_3(Q)]$
denote the coefficients of the projective binary cubic form
$\phi=\phi([x:y],Q)$ whose zero set
$$
Z_\phi
=\{[x:y]\in\mathbb {CP}^1\,|\ \phi([x:y])=0\}
$$
coincides with $Q$, and let $D_\phi$ denote the 
discriminant of $\phi$.
The Eisenstein's automorphism of 
${\mathcal C}^3({\mathbb{CP}}^1)$ coincides with the
tame map
$$
F_T\colon{\mathcal C}^3({\mathbb{CP}}^1)\ni Q
\mapsto T(Q)Q\in{\mathcal C}^3({\mathbb{CP}}^1)
$$
that corresponds to the morphism 
$T\colon{\mathcal C}^3({\mathbb{CP}}^1)
\to{\mathbf{PSL}}(2,{\mathbb C})$, where
\begin{equation}\label{eq: tame Eisenstein}
T(Q)\colon{\mathbb{CP}}^1\ni\zeta\mapsto
\frac{a(Q)\zeta+b(Q)}{c(Q)\zeta+d(Q)}\in{\mathbb{CP}}^1
\end{equation}
is the M{\"o}bius transformation with the coefficients
\begin{equation}\label{eq: coeficients of tame Eisenstein}
\aligned
&a(Q)=(z_1z_2-z_0z_3)/\sqrt{-D_\phi}\,;\\
&b(Q)=2(z_2^2-z_1z_3)/\sqrt{-D_\phi}\,;\\
&c(Q)=2(z_0z_2-z_1^2)/\sqrt{-D_\phi}\,;\\
&d(Q)=-a(Q)=-(z_1z_2-z_0z_3)/\sqrt{-D_\phi}\,.\\
\endaligned
\end{equation}
It is easily seen that $T(Q)\in{\mathbf{PSL}}(2,{\mathbb C})$
is well defined by (\ref{eq: tame Eisenstein}) and
(\ref{eq: coeficients of tame Eisenstein}). 

\newpage
\input biblconf

\input index
\end{document}

%% file: cd.tex
\catcode`\@=11

 \def\AMSTeXfeatures{\Plainheads 
   \let\current@vert=\AMS@vert}

 \def\Plainheads{\sh@ftdiam=0.05em
   \getlabeldims
   \let\vshaftfill=\plnvsolidfill
   \let\hshaftfill=\plnhsolidfill
   \let\th@rhead=\plnrhead
   \let\th@lhead=\plnlhead
   \let\th@dnhead=\plndnhead
   \let\th@uphead=\plnuphead}
 
 \def\glet{\global\let}

 \def\LaTeXfeatures{\catcode`\@=11
   \ifx\@clnwd\undefined \nol@g
      \input ltxcode.tex \dol@g \fi
   \ltxheads \let\current@vert=\new@vert
   \providelto \catcode`\@=\active}

 \def\nol@g{\def\wlog{\edef\garbage}}
 \def\dol@g{\let\wlog=\wl@g} \let\wl@g=\wlog
 \nol@g 

 \newbox\ltobox
 \def\providelto{{\setbox\z@=
   \hbox{$\to$}\minharrlen=\wd\z@
   \global\setbox\ltobox=\hbox{$\activeat>>>$}}
   \def\lto{\mathrel{\copy\ltobox}}}

 \def\ltxheads{\sh@ftdiam=\@wholewidth
   \getlabeldims
   \let\vshaftfill= \ltxvsolidfill
   \let\hshaftfill=\ltxhsolidfill
   \let\th@rhead=\ltxrhead
   \let\th@lhead=\ltxlhead
   \let\th@dnhead=\ltxdnhead
   \let\th@uphead=\ltxuphead}
 {\catcode`\@=\active
   \gdef@#1{\csname #1\string@at\endcsname}
   \glet\activeat=@}
 \def\def@#1{\expandafter\def\csname #1@at\endcsname}

 \def@>#1>#2>{\@rrow R{#1}{#2}}
 \def@<#1<#2<{\@rrow L{#1}{#2}}
 \def@ V#1V#2V{\@rrow V{#1}{#2}}
 \def@ A#1A#2A{\@rrow A{#1}{#2}}
 \def@/#1/#2/#3/{\@rrow{#1}{#2}{#3}}
 \def@.{\ifodd\row\ifmmode\noharrow
     \else\leavevmode.\spacefactor3000 \fi
   \else\novarrow\fi}
 \def@={\ifodd\row\harrow\hequalfill{}{}%
   \else\varrow\vequalfill{}{}\fi}
 \def@:#1{\ifx=#1\harrow\deffill{}{}%
   \else\leavevmode\null:#1\fi}
 \def@|{\current@vert}
  \def\AMS@vert{\varrow\vequalfill{}{}}
  \def\new@vert#1|#2|{\ifodd\row
   \let\nextarrow\vertexvarrow
   \else\let\nextarrow\varrow\fi
   \nextarrow\vshaftfill{#1}{#2}}
 \def@-{\ifmmode\let\next\hl@ne
   \else\let\next\AMSatdash \fi \next}
  \def\hl@ne#1-#2-{\harrow\hshaftfill{#1}{#2}}
  \def\AMSatdash{\let\next\relax\leavevmode
    \def\next@{\ifx\next-%
      \def\next-{\futurelet\next\nextii@}%
     \else\def\next{\hbox{-}}\fi\next}%
    \def\nextii@{\ifx\next-\def\next-{\hbox{---}}%
      \else\def\next{\hbox{--}}\fi\next}%
    \futurelet\next\next@}
 \def@(#1){\tweenarrows{#1}}
 \def@[#1]{\setsp@n#1\relax\activeat}
 \def\fiberbox{\hbox{$\vcenter{\hr@le\hbox{\vr@le
   \kern1ex\vbox{\kern1.2ex}\vr@le}\hr@le}$}}
  \def\hr@le{\hrule height \sh@ftdiam}
  \def\vr@le{\vrule width \sh@ftdiam}
 \def@+#1+#2+#3+{\ifodd\row \harrow{#1}{#2}{#3}%
   \else \varrow{#1}{#2}{#3}\fi}

 \def\Dnarrfill{\vequalfill\Dnhe@d}
 \def\Uparrfill{\Uphe@d\vequalfill}
 
 \def\ontofill{\rtarrfill\kern-0.3em 
   \th@rhead\kern 0.3em} 

 \def\rtarrfill{\hshaftfill\th@rhead}
 \def\ltarrfill{\th@lhead\hshaftfill}
 \def\dnarrfill{\vshaftfill\th@dnhead}
 \def\uparrfill{\th@uphead\vshaftfill}
 \def\hequalfill{\plnhfill=}
 \def\deffill{:\plnhfill=}
 \def\plnvextfill#1{\setbox\z@
   \hbox{\the\textfont3 #1}%
   \dimen@=\dp\z@\advance\dimen@\ht\z@
   \copy\z@ \kern-\dimen@ 
   \cleaders\copy\z@ \vfill
   \kern-\dimen@ 
   \box\z@}
 \def\plnhfill#1{$\m@th\mkern-1.5mu\mathord#1\mkern-6mu
    \cleaders\hbox{$\mkern-2mu\mathord#1\mkern-2mu$}\hfill
    \mkern-6mu\mathord#1\mkern-1.5mu$}
 \def\vequalfill{\plnvextfill{\char'167}}
 \def\plnvsolidfill{\plnvextfill{\char'077}}
 \def\plnhsolidfill{\plnhfill-}
 \def\ltxhsolidfill{\leaders\hrule height\topofshaft depth\botofshaft
   \hfill}
 \def\ltxvsolidfill{\leaders\vrule width\sh@ftdiam\vfill}
 \def\hdashfill{\hd@sh\wd@sh
   \xleaders \hbox{\wd@sh\hd@sh\wd@sh}\hfill
   \wd@sh\hd@sh}
 \def\vdashfill{\vd@sh\wd@sh
   \xleaders \vbox{\wd@sh\vd@sh\wd@sh}\vfill
   \wd@sh\vd@sh}
 \def\dashed{\ifinmeasureCD\else
    \ifodd\row\option{\let\hshaftfill=\hdashfill}%
   \else\option{\let\vshaftfill=\vdashfill}\fi\fi}

 \newdimen\CDstrutht  \newdimen\CDstrutdp
 \newdimen\CDstrutlen \CDstrutlen=\CDstrutht
   \advance\CDstrutlen by \CDstrutdp

 \def\CDstrut{\vrule
   height \ifnum\row=1 \z@\else\CDstrutht \fi
   depth \ifnum\row=\numrows \z@ \else\CDstrutdp \fi
   width\z@}

 \newdimen\CDarrsurr \CDarrsurr=0.375em
 \newdimen\CDdashlen
 \newdimen\CDvarrlen \CDvarrlen=1.5\baselineskip
 \newdimen\minharrlen 
  \setbox\z@\hbox{$\longrightarrow$} \minharrlen=\wd\z@
 \newdimen\minCDharrlen \minCDharrlen=2.5em 
\newdimen \minc@lwd
\def\findminc@lwd{\minc@lwd=2\CDarrsurr
  \advance\minc@lwd\minCDharrlen}

 \newdimen\sh@ftdiam

 \newdimen\labelsurr \labelsurr=1.25 em

\newcount\sp@ncnt \sp@ncnt=\@ne
\newcount\sp@ncnt@ \sp@ncnt@=\@ne
\newdimen\@rrwd \newdimen\@rrdp

 \def\adjustbot#1{\option{\advance\@rrdp#1\relax}}
 
\def\pushvertex#1{\global\p@shlen#1\relax
   \global\let\maybepush=\dopush}

 \newdimen\p@shlen \p@shlen=\z@

 \let\maybepush=\relax
 \def\dopush{\ifinmeasureCD 
   \advance\locdimen by -\p@shlen 
   \else\advance \@rrwd by -\p@shlen \fi 
   \global\let\maybepush=\relax \global\p@shlen=\z@\relax}

 \def\span@ne{\global\sp@ncnt=\@ne\relax}
 \def\setsp@n#1#2{\global\sp@ncnt=#1\relax
   \ifx\relax#2\relax\else\global\sp@ncnt@=#2\relax\fi}

 \def\plnrhead{\llap{$\rightarrow\mkern-1.5mu$}}
 \def\plnlhead{\rlap{$\mkern-1.5mu\leftarrow$}}

 \def\clap#1{\hbox to \z@{\hss #1\hss}}

 \def\plndnhead{\hbox{\the\textfont3 \char'171}}
 \def\plnuphead{\hbox{\the\textfont3 \char'170}}
 \def\Dnhe@d{\hbox{\the\textfont3 \char'177}}
 \def\Uphe@d{\hbox{\the\textfont3 \char'176}}

 \def\ltxrhead{\raise\@xisheight
   \llap{\smash{\@linefnt\@getrarrow(1,0)}}}
 \def\ltxlhead{\raise\@xisheight
   \rlap{\@linefnt\@getlarrow(-1,0)}}
 \def\ltxuphead{\setbox\z@=\rlap{%
   \kern\@halfwidth\@linefnt\char'66}%
   \copy\z@\kern-\ht\z@}
 \def\ltxdnhead{\setbox\z@=\rlap{%
   \kern\@halfwidth\@linefnt\char'77}%
   \ht\z@=\z@\box\z@}

 \def\wd@sh{\kern0.5\CDdashlen}
 \def\hd@sh{\vrule height\topofshaft depth\botofshaft
    width\CDdashlen}
 \def\vd@sh{\hrule height\CDdashlen
   depth\z@ width\sh@ftdiam}

\def\xylist{14{3434}13{2414}12{1723}%
  23{1413}34{1153}11{0867}43{0707}%
  32{0580}21{0414}31{0291}41{0}}
\newcount\tgtcnt@
\def\find@xyargs{\dimen@=\@rrdp
  \advance\dimen@ by \CDstrutlen
  \tgtcnt@=\dimen@ \dimen@=\@rrwd 
  \divide\dimen@ by \@m 
  \divide \tgtcnt@ by \dimen@ 
  \expandafter\testxy\xylist\relax
  \unitlength=\@xarg\@rrdp
  \divide\unitlength by\@yarg\relax}
\def\testxy#1#2#3{\ifnum\tgtcnt@>#3
    \@xarg=#1\relax \@yarg=#2\relax
    \let\next=\ignorerest
  \else\let\next\testxy\fi\next}
\def\ignorerest#1\relax{\relax}

\let\scalefactor=\@ne
\def\SWarrow{\find@xyargs\vector
  (-\@xarg,-\@yarg)\scalefactor\hskip-\wd\@linechar}
\def\NWarrow{\find@xyargs\vector
  (-\@xarg,\@yarg)\scalefactor\hskip-\wd\@linechar}
\def\NEarrow{\find@xyargs\vector
  (\@xarg,\@yarg)\scalefactor}
\def\SEarrow{\find@xyargs\vector
  (\@xarg,-\@yarg)\scalefactor}
\def\rightupline{\find@xyargs\@linelen=\scalefactor
     \unitlength\@sline}
\def\rightdownline{\find@xyargs\@yarg=-\@yarg\relax
     \@linelen=\scalefactor\unitlength\@sline}

\def\Sim{\ifodd\row\setbox\z@=\hbox{$\sim$}\dimen@=\ht\z@
 \advance\dimen@ by -\@xisheight
  \vbox{\box\z@\kern-\@xisheight\kern\dimen@}%
  \else\hbox{$\wr$}\fi}

\def\harrow#1#2#3{\inmeasureCDtrue\findminarrwd
  {#2}{#3}{\sp@ncnt\minharrlen}\inmeasureCDfalse\span@ne
  \mathrel{\hbox{\options\hplace{#1}\ulabel{#2}\dlabel{#3}}}}

\def\noharrow{\harrow\hfill{}{}}
\def\vertexvarrow#1#2#3{\findarrdp \@rrwd=\z@ \setsp@n\@ne\@ne
  \vbox to \z@{\kern-1.2\CDstrutht
  \rlap{\options\vplace{#1}\llabel{#2}\rlabel{#3}}\vss}}

\newif\ifinmeasureCD
\def\measurelabel#1{\setbox\z@
  \hbox{$\scriptstyle#1\kern\labelsurr$}%
  \ifdim\wd\z@>\@rrwd \@rrwd=\wd\z@\fi}
\def\findminarrwd#1#2#3{\@rrwd=#3\relax
   \measurelabel{#1}\measurelabel{#2}}
\def\findCDarrwd#1#2{\@rrwd=\minCDharrlen
   \measurelabel{#1}\measurelabel{#2}%
  }

\newcount\row \row=\@ne \newcount\col \col=\@ne 
 \newcount\numrows 
\numrows=\@ne
 \newcount\numcols
\newcount\arrspan \newdimen\vrtxhalfwd  \newbox\tempbox

\def\DANABUG{\advance\col by \@ne
 \@rrwd=\minCDharrlen
  \advance\@rrwd by \vrtxhalfwd
  \advance\@rrwd by \CDarrsurr
  \ifnum\col>\numcols \numcols=\col
     \newlocdimen{col\the\col}\locdimen=\@rrwd 
  \else \ifdim\@rrwd>\c@l \c@l=\@rrwd\fi\fi}

\def\drop#1\\{
  \findvrtxhalfsum\DANABUG\advance\row by 2 \measureinit}

\def\measureinit{\col=\@ne \vrtxhalfwd=-\CDarrsurr\arrspan=\@ne\@rrwd=\z@
   \setbox\tempbox=\hbox\bgroup$}
\def\measure{
  \let\harrow\measureCDarrow
  \let\CDCR=\measureCR 
   \findminc@lwd 
  \inmeasureCDtrue
  \row=\@ne \numcols=\z@ \measureinit}

\def\endmeasure{\findvrtxhalfsum\DANABUG
  \numrows=\row 
  \inmeasureCDfalse}

\def\newlocdimen#1{\advance\dimenc@unt by \@ne
  \ifnum\dimenc@unt<\insc@unt
     \else\errmessage{No room for the CD}\fi
  \dimendef\locdimen=\dimenc@unt
  \expandafter\dimendef\csname#1\endcsname=\dimenc@unt}

 \def\r@wc@l{\csname row\the\row col\the\col\endcsname}
 \def\c@l{\csname col\the\col\endcsname}

 \def\findvrtxhalfsum{$\egroup
  \newlocdimen{row\the\row col\the\col}
  \locdimen=\vrtxhalfwd 
  \vrtxhalfwd=0.5\wd\tempbox 
  \advance\vrtxhalfwd by \CDarrsurr
  \advance\locdimen by \vrtxhalfwd 
  \advance\@rrwd by \locdimen 
  \maybepush
  \divide\@rrwd by \arrspan\relax
  \ifdim\@rrwd<\minc@lwd
    \ifnum\col>\@ne \@rrwd=\minc@lwd\fi \fi
  \loop 
    \ifnum\col>\numcols \numcols=\col
       \newlocdimen{col\the\col}
       \locdimen=\@rrwd 
    \else \ifdim\@rrwd>\c@l \c@l=\@rrwd\fi \fi
   \ifnum\arrspan>\@ne
      \advance\arrspan by -1 \advance\col by \@ne
  \repeat }

 \def\measureCDarrow#1#2#3{\findvrtxhalfsum
   \arrspan=\sp@ncnt\relax\global\sp@ncnt=1\relax
   \advance\col by \@ne
   \findCDarrwd{#2}{#3}%
   \setbox\tempbox=\hbox\bgroup$}

 \newcount\dr@tn \dr@tn=\z@
 \def\locate#1:#2{\ifinmeasureCD\else
   \count@=-#1
   \multiply\count@ by 2
   \advance\count@ by #2
   \dimen@=\count@\@rrwd
   \ifnum\dr@tn=\@ne\relax \else\dimen@=-\dimen@ \fi
   \dimen@i=\@rrdp
   \ifnum\dr@tn>\z@\advance\dimen@i by \CDstrutlen \fi
   \dimen@i=\count@\dimen@i
   \count@=#2 \multiply\count@ by 2
   \divide\dimen@ by \count@
   \divide\dimen@i by \count@
   \lift\dimen@i\nudge\dimen@\fi}

\def\betweenCDrows{\advance\row by \@ne \col=\@ne
\options}

\def\hbegin{\hbox\bgroup\kern\c@l \kern-\r@wc@l$}
\def\hend{$\glet\maybepush\relax \CDstrut\egroup}
\def\vbegin{\setbox\tempbox=\hbox\bgroup$}
\def\vend{$\egroup\ht\tempbox=\z@\dp\tempbox\CDvarrlen
  \box\tempbox}
\def\setCD{\let\harrow=\setCDarrow
  \let\CDCR=\setCR 
  \row=\@ne \col=\@ne \hbegin}
\let\endsetCD=\hend 

\def\findarrwd{\@rrwd=\z@ \count@=\col \advance\count@ by\sp@ncnt
  \loop\ifnum\count@>\col \advance\count@ by -1
      \advance\@rrwd by\csname col\the\count@\endcsname\repeat}
\def\setCDarrow#1#2#3{\kern\CDarrsurr\advance\col by \@ne
  \findarrwd \advance\@rrwd by -\r@wc@l  
  \@rrdp=\z@ 
  \maybepush
  \advance\col by -\@ne \advance\col by \sp@ncnt \span@ne
  \hbox to \@rrwd{\options
   \@rrwd=\scalefactor\@rrwd\hss
   \hplace{#1}\ulabel{#2}\dlabel{#3}\hss}%
   \kern\CDarrsurr}

\newdimen\labspacei 
\newdimen\labspaceii 

\newdimen\@xisheight
  \@xisheight=\the\fontdimen22\textfont2
\newdimen\labelskip
  \labelskip=\the\fontdimen10\textfont3 
\newdimen\topofshaft
\newdimen\botofshaft
\newdimen\botofulabel
\newdimen\topofdlabel
\def\getlabeldims{
  \topofshaft=0.5\sh@ftdiam
  \botofshaft=\topofshaft
  \advance\topofshaft by \@xisheight  
  \advance\botofshaft by -\@xisheight  
  \botofulabel=\topofshaft
  \advance\botofulabel by \labelskip
  \topofdlabel=\botofshaft
  \advance\topofdlabel by \labelskip}

\def\ulabel{\ifnum\row=\@ne\let\next\ulabeli
   \else\let\next\ulabellap\fi\next}
\def\ulabeli#1{\vbox{
  \clap{\kern-\@rrwd$\scriptstyle#1$}%
  \kern\botofulabel}\maybeoffset}
\def\ulabellap#1{\vbox to \z@{\vss
  \clap{\kern-\@rrwd$\scriptstyle#1$}%
  \kern\botofulabel}\maybeoffset}
\def\dlabel{\ifnum\row=\numrows\let\next\dlabeli
   \else\let\next\dlabellap\fi\next}
\def\dlabeli#1{\vtop{\kern\topofdlabel
  \clap{\kern-\@rrwd$\scriptstyle#1$}%
  }\maybeoffset}
\def\dlabellap#1{\vbox to \z@{\kern\topofdlabel
  \clap{\kern-\@rrwd$\scriptstyle#1$}%
  \vss}\maybeoffset}
\def\rlabel#1{\vbox to \z@{\vss
  \rlap{\kern\labelskip$\scriptstyle#1$}%
  \vss\kern-\@rrdp}\maybeoffset}
\def\llabel#1{\vbox to \z@{\vss
  \llap{$\scriptstyle#1$\kern\labelskip}%
  \vss\kern-\@rrdp}\maybeoffset}
\def\swlabel#1{\vtop{\kern0.5\@rrdp
  \llap{$\scriptstyle#1$\kern\labelskip\kern-0.5\@rrwd}
  }\maybeoffset}
\def\nwlabel#1{\vbox{
  \llap{$\scriptstyle#1$\kern\labelskip\kern-0.5\@rrwd}%
  \kern-0.5\@rrdp}\maybeoffset}
\def\selabel#1{\vtop{\kern0.5\@rrdp
  \rlap{\kern0.5\@rrwd\kern\labelskip$\scriptstyle#1$}%
  }\maybeoffset}
\def\nelabel#1{\vbox{
  \rlap{\kern0.5\@rrwd\kern\labelskip$\scriptstyle#1$}%
  \kern-0.5\@rrdp}\maybeoffset}
\def\cplace#1{\vbox to \z@{\vss
  \clap{$#1$\kern-\@rrwd}%
  \kern-\@rrdp\vss}\maybeoffset}
\def\hplace#1{\hbox to \@rrwd{#1}\maybeoffset}
\def\vplace#1{\clap{\vbox to \z@{#1\kern-\@rrdp}}\maybeoffset}

\newdimen\nudgeamount \nudgeamount=\z@
\newdimen\liftamount \liftamount=\z@
\let\maybeoffset\relax
\newbox\offsetbox \newdimen\lastheight
\def\dooffset{
  \setbox\offsetbox=\lastbox \lastheight=\ht\offsetbox 
  \setbox\offsetbox=\vbox{\kern-\liftamount\box\offsetbox}%
  \ht\offsetbox=\lastheight
  \kern\nudgeamount\box\offsetbox\kern-\nudgeamount
  \global\nudgeamount=\z@ \global\liftamount=\z@
  \glet\maybeoffset=\relax}
\def\nudge#1{\ifinmeasureCD\else
  \global\advance\nudgeamount#1\relax
  \global\let\maybeoffset\dooffset\fi}
\def\lift#1{\ifinmeasureCD\else
  \global\advance\liftamount#1\relax
  \global\let\maybeoffset\dooffset\fi}

\def\findarrdp{\@rrdp=\CDvarrlen
  \ifnum\sp@ncnt@>1
    \advance\@rrdp by \CDstrutlen
    \multiply\@rrdp by \sp@ncnt@
    \advance\@rrdp by -\CDstrutlen \fi
 }

\def\varrow#1#2#3{\ifnum\sp@ncnt>\@ne 
     \sp@ncnt@=\sp@ncnt\relax\fi
  \findarrdp \@rrwd=\z@ 
  \kern\c@l
   \hbox to \z@{\options
   \@rrdp=\scalefactor\@rrdp
    \hss\vplace{#1}\llabel{#2}\rlabel{#3}\hss}%
  \global\advance\col by \@ne \setsp@n\@ne\@ne
  }

\def\novarrow{\varrow\vfill{}{}}

\def\tweenarrows#1{\findarrwd \findarrdp \setsp@n\@ne\@ne
  \rlap{\options\cplace{#1}}}

\def\usarrow #1#2#3{\dr@tn=\@ne
  \findarrwd \findarrdp \setsp@n\@ne\@ne 
  \rlap{\options\cplace{#1}\nwlabel{#2}\selabel{#3}}%
  \dr@tn=\z@}
\def\dsarrow #1#2#3{\dr@tn=\tw@
  \findarrwd \findarrdp \setsp@n\@ne\@ne 
  \rlap{\options\cplace{#1}\swlabel{#2}\nelabel{#3}}%
  \dr@tn=\z@}
 \def\@rrow#1{\csname #1@rrow\endcsname}
 \def\R@rrow{\harrow \rtarrfill}
 \def\L@rrow{\harrow \ltarrfill}
 \def\V@rrow{\varrow \dnarrfill}
 \def\A@rrow{\varrow \uparrfill}
 \def\SE@rrow{\dsarrow \SEarrow}
 \def\NW@rrow{\dsarrow \NWarrow}
 \def\SW@rrow{\usarrow \SWarrow}
 \def\NE@rrow{\usarrow \NEarrow}
 \def\DS@rrow{\dsarrow \dnslope}
 \def\US@rrow{\usarrow \upslope}
 \def\upslope{\find@xyargs
       \@linelen=\unitlength\@sline}
 \def\dnslope{\find@xyargs\@yarg=-\@yarg\relax
       \@linelen=\unitlength\@sline}

\newtoks\optionlist 
\optionlist={}
\let\options\relax
\def\dooptions{\the\optionlist\global\optionlist={}%
  \glet\options=\relax}
\def\option#1{\ifinmeasureCD\else
  \glet\options=\dooptions
  \global\optionlist=\expandafter{\the\optionlist\relax#1}\fi}
\def\wider#1{\ifinmeasureCD\else
   \option{\advance\@rrwd by #1}\fi}
\def\deeper#1{\ifinmeasureCD\else
   \option{\advance\@rrdp by #1}\fi}

{\def\\{\global\let\sptoken= }\\ }

\def\CR{\futurelet\nexttok\testCR}
\def\testCR{\ifx\nexttok\sptoken
   \let\next\eatspaceCR\else\let\next\CDCR\fi\next}
\def\eatspaceCR#1 {\CR}
\def\measureCR{\ifx\nexttok\endmeasure\let\nextCR\relax
    \else\let\nextCR\drop\fi\nextCR}
\def\setCR{\ifodd\row
  \ifx\nexttok\endsetCD\else\hend\betweenCDrows\vbegin\fi
  \else\vend\betweenCDrows\hbegin\fi}

\countdef\dimenc@unt=11
\def\CD#1\endCD{
   \begingroup\let\\=\CR
  \m@th\offinterlineskip
   \measure#1\endmeasure\null\,\vcenter{\setCD#1\endsetCD}\,
   \endgroup
    }

\ifx\@clnwd\undefined \nol@g\else\catcode`\ =14\relax\fi
 \font\@linefnt=line10 
 \newcount\@tempcnta
 \newcount\@tempcntb
 \newdimen\@tempdima
 \newdimen\@tempdimb
 \newdimen\@wholewidth
 \newdimen\@halfwidth
   \@wholewidth\fontdimen8\@linefnt \@halfwidth .5\@wholewidth
 \newdimen\unitlength
 \newcount\@xarg
 \newcount\@yarg
 \newcount\@yyarg
 \newbox\@linechar
 \newdimen\@linelen
 \newdimen\@clnwd
 \newdimen\@clnht
 \newif\if@negarg
 
 \def\@whilenoop#1{}

 \def\@whiledim#1\do #2{\ifdim #1\relax#2\@iwhiledim{#1\relax#2}\fi}

 \def\@iwhiledim#1{\ifdim #1\let\@nextwhile=\@iwhiledim 
         \else\let\@nextwhile=\@whilenoop\fi\@nextwhile{#1}}

 \def\@sline{\ifnum\@xarg< 0 \@negargtrue \@xarg -\@xarg \@yyarg -\@yarg
   \else \@negargfalse \@yyarg \@yarg \fi
 \ifnum \@yyarg >0 \@tempcnta\@yyarg \else \@tempcnta -\@yyarg \fi
 \ifnum\@tempcnta>6 \@badlinearg\@tempcnta0 \fi
 \ifnum\@xarg>6 \@badlinearg\@xarg 1 \fi
 \setbox\@linechar\hbox{\@linefnt\@getlinechar(\@xarg,\@yyarg)}%
 \ifnum \@yarg >0 \let\@upordown\raise \@clnht\z@
    \else\let\@upordown\lower \@clnht \ht\@linechar\fi
 \@clnwd=\wd\@linechar
 \if@negarg \hskip -\wd\@linechar \def\@tempa{\hskip -2\wd\@linechar}\else
      \let\@tempa\relax \fi
 \@whiledim \@clnwd <\@linelen \do
   {\@upordown\@clnht\copy\@linechar
    \@tempa
    \advance\@clnht \ht\@linechar
    \advance\@clnwd \wd\@linechar}%
 \advance\@clnht -\ht\@linechar
 \advance\@clnwd -\wd\@linechar
 \@tempdima\@linelen\advance\@tempdima -\@clnwd
 \@tempdimb\@tempdima\advance\@tempdimb -\wd\@linechar
 \if@negarg \hskip -\@tempdimb \else \hskip \@tempdimb \fi
 \multiply\@tempdima \@m
 \@tempcnta \@tempdima \@tempdima \wd\@linechar \divide\@tempcnta \@tempdima
 \@tempdima \ht\@linechar \multiply\@tempdima \@tempcnta
 \divide\@tempdima \@m
 \advance\@clnht \@tempdima
 \ifdim \@linelen <\wd\@linechar
    \hskip \wd\@linechar
   \else\@upordown\@clnht\copy\@linechar\fi}
 
 \def\@getlinechar(#1,#2){\@tempcnta#1\relax\multiply\@tempcnta 8
 \advance\@tempcnta -9 \ifnum #2>0 \advance\@tempcnta #2\relax\else
 \advance\@tempcnta -#2\relax\advance\@tempcnta 64 \fi
 \char\@tempcnta}
 
 \def\vector(#1,#2)#3{\@xarg #1\relax \@yarg #2\relax
 \@tempcnta \ifnum\@xarg<0 -\@xarg\else\@xarg\fi
 \ifnum\@tempcnta<5\relax
 \@linelen=#3\unitlength
 \ifnum\@xarg =0 \@vvector 
   \else \ifnum\@yarg =0 \@hvector \else \@svector\fi
 \fi
 \else\@badlinearg\fi}
 
 \def\@svector{\@sline
 \@tempcnta\@yarg \ifnum\@tempcnta <0 \@tempcnta=-\@tempcnta\fi
 \ifnum\@tempcnta <5
   \hskip -\wd\@linechar
   \@upordown\@clnht \hbox{\@linefnt  \if@negarg 
   \@getlarrow(\@xarg,\@yyarg) \else \@getrarrow(\@xarg,\@yyarg) \fi}%
 \else\@badlinearg\fi}
 
 \def\@getlarrow(#1,#2){\ifnum #2 =\z@ \@tempcnta='33\else
 \@tempcnta=#1\relax\multiply\@tempcnta \sixt@@n \advance\@tempcnta
 -9 \@tempcntb=#2\relax\multiply\@tempcntb \tw@
 \ifnum \@tempcntb >0 \advance\@tempcnta \@tempcntb\relax
 \else\advance\@tempcnta -\@tempcntb\advance\@tempcnta 64
 \fi\fi\char\@tempcnta}
 
 \def\@getrarrow(#1,#2){\@tempcntb=#2\relax
 \ifnum\@tempcntb < 0 \@tempcntb=-\@tempcntb\relax\fi
 \ifcase \@tempcntb\relax \@tempcnta='55 \or 
 \ifnum #1<3 \@tempcnta=#1\relax\multiply\@tempcnta
 24 \advance\@tempcnta -6 \else \ifnum #1=3 \@tempcnta=49
 \else\@tempcnta=58 \fi\fi\or 
 \ifnum #1<3 \@tempcnta=#1\relax\multiply\@tempcnta
 24 \advance\@tempcnta -3 \else \@tempcnta=51\fi\or 
 \@tempcnta=#1\relax\multiply\@tempcnta
 \sixt@@n \advance\@tempcnta -\tw@ \else
 \@tempcnta=#1\relax\multiply\@tempcnta
 \sixt@@n \advance\@tempcnta 7 \fi\ifnum #2<0 \advance\@tempcnta 64 \fi
 \char\@tempcnta}
\catcode`\ =10

\dol@g 
\catcode`\@=\active
\LaTeXfeatures

%% file: picture.tex
\setlength{\unitlength}{0.075in}
\begin{picture}(53.48,46.35)
\put(36.68,41.00){{\setbox0=\hbox{$\sqrt{ac}$}\lower\ht0\box0}}
\put(47.48,10.){{\setbox0=\hbox{$\sqrt{bc}$}\lower\ht0\box0}}
\put(35.35,24.33){{\setbox0=\hbox{$c$}\lower\ht0\box0}}
\put(6.95,29.){{\setbox0=\hbox{$\sqrt{ab}$}\lower\ht0\box0}}
\put(30.55,22.30){{\setbox0=\hbox{$b$}\lower\ht0\box0}}
\put(25.75,28.50){{\setbox0=\hbox{$a$}\lower\ht0\box0}}
\put(53.48,1.00){{\setbox0=\hbox{$\mu_1^{-}$}\lower\ht0\box0}}
\put(19.48,24.83){{\setbox0=\hbox{$q_3$}\lower\ht0\box0}}
\put(38.55,15.16){{\setbox0=\hbox{$q_1$}\lower\ht0\box0}}
\put(34.01,33.33){{\setbox0=\hbox{$q_2$}\lower\ht0\box0}}
\put(27.21,18.16){{\setbox0=\hbox{$\mu_2^{+}$}\lower\ht0\box0}}
\put(39.,26.5){{\setbox0=\hbox{$\mu_3^{+}$}\lower\ht0\box0}}
\put(25.48,34.){{\setbox0=\hbox{$\mu_1^{+}$}\lower\ht0\box0}}
\put(39.08,48.33){{\setbox0=\hbox{$\mu_2^{-}$}\lower\ht0\box0}}
\put(-1.5,27.){{\setbox0=\hbox{$\mu_3^{-}$}\lower\ht0\box0}}
\put(53.35,0.83){{\setbox0=\hbox{$\scriptstyle\bullet$}\kern-.4\wd0\lower.5\ht0\box0}}
\put(38.28,45.83){{\setbox0=\hbox{$\scriptstyle\bullet$}\kern-.4\wd0\lower.5\ht0\box0}}
\put(26.68,30.83){{\setbox0=\hbox{$\scriptstyle\bullet$}\kern-.4\wd0\lower.5\ht0\box0}}
\put(28.15,19.16){{\setbox0=\hbox{$\scriptstyle\bullet$}\kern-.4\wd0\lower.5\ht0\box0}}
\put(38.28,25.83){{\setbox0=\hbox{$\scriptstyle\bullet$}\kern-.4\wd0\lower.5\ht0\box0}}
\put(1.61,25.83){{\setbox0=\hbox{$\scriptstyle\bullet$}\kern-.4\wd0\lower.5\ht0\box0}}
\special{em:linewidth 0.014in}
\put(1.61,25.83){\special{em:moveto}}
\put(1.95,25.83){\special{em:lineto}}
\put(2.28,25.83){\special{em:moveto}}
\put(2.61,25.83){\special{em:lineto}}
\put(2.95,25.83){\special{em:moveto}}
\put(3.28,25.83){\special{em:lineto}}
\put(3.61,25.83){\special{em:moveto}}
\put(3.95,25.83){\special{em:lineto}}
\put(4.28,25.83){\special{em:moveto}}
\put(4.61,25.83){\special{em:lineto}}
\put(4.95,25.83){\special{em:moveto}}
\put(5.28,25.83){\special{em:lineto}}
\put(5.61,25.83){\special{em:moveto}}
\put(5.95,25.83){\special{em:lineto}}
\put(6.28,25.83){\special{em:moveto}}
\put(6.61,25.83){\special{em:lineto}}
\put(6.95,25.83){\special{em:moveto}}
\put(7.28,25.83){\special{em:lineto}}
\put(7.61,25.83){\special{em:moveto}}
\put(7.95,25.83){\special{em:lineto}}
\put(8.28,25.83){\special{em:moveto}}
\put(8.61,25.83){\special{em:lineto}}
\put(8.95,25.83){\special{em:moveto}}
\put(9.28,25.83){\special{em:lineto}}
\put(9.61,25.83){\special{em:moveto}}
\put(9.95,25.83){\special{em:lineto}}
\put(10.28,25.83){\special{em:moveto}}
\put(10.61,25.83){\special{em:lineto}}
\put(10.95,25.83){\special{em:moveto}}
\put(11.28,25.83){\special{em:lineto}}
\put(11.61,25.83){\special{em:moveto}}
\put(11.95,25.83){\special{em:lineto}}
\put(12.28,25.83){\special{em:moveto}}
\put(12.61,25.83){\special{em:lineto}}
\put(12.95,25.83){\special{em:moveto}}
\put(13.28,25.83){\special{em:lineto}}
\put(13.61,25.83){\special{em:moveto}}
\put(13.95,25.83){\special{em:lineto}}
\put(14.28,25.83){\special{em:moveto}}
\put(14.61,25.83){\special{em:lineto}}
\put(14.95,25.83){\special{em:moveto}}
\put(15.28,25.83){\special{em:lineto}}
\put(15.61,25.83){\special{em:moveto}}
\put(15.95,25.83){\special{em:lineto}}
\put(16.28,25.83){\special{em:moveto}}
\put(16.61,25.83){\special{em:lineto}}
\put(16.95,25.83){\special{em:moveto}}
\put(17.28,25.83){\special{em:lineto}}
\put(17.61,25.83){\special{em:moveto}}
\put(17.95,25.83){\special{em:lineto}}
\put(18.28,25.83){\special{em:moveto}}
\put(18.61,25.83){\special{em:lineto}}
\put(18.95,25.83){\special{em:moveto}}
\put(19.28,25.83){\special{em:lineto}}
\put(19.61,25.83){\special{em:moveto}}
\put(19.95,25.83){\special{em:lineto}}
\put(20.28,25.83){\special{em:moveto}}
\put(20.61,25.83){\special{em:lineto}}
\put(20.95,25.83){\special{em:moveto}}
\put(21.28,25.83){\special{em:lineto}}
\put(21.61,25.83){\special{em:moveto}}
\put(21.95,25.83){\special{em:lineto}}
\put(22.28,25.83){\special{em:moveto}}
\put(22.61,25.83){\special{em:lineto}}
\put(22.95,25.83){\special{em:moveto}}
\put(23.28,25.83){\special{em:lineto}}
\put(23.61,25.83){\special{em:moveto}}
\put(23.95,25.83){\special{em:lineto}}
\put(24.28,25.83){\special{em:moveto}}
\put(24.61,25.83){\special{em:lineto}}
\put(24.95,25.83){\special{em:moveto}}
\put(25.28,25.83){\special{em:lineto}}
\put(25.61,25.83){\special{em:moveto}}
\put(25.95,25.83){\special{em:lineto}}
\put(26.28,25.83){\special{em:moveto}}
\put(26.61,25.83){\special{em:lineto}}
\put(26.95,25.83){\special{em:moveto}}
\put(27.28,25.83){\special{em:lineto}}
\put(27.61,25.83){\special{em:moveto}}
\put(27.95,25.83){\special{em:lineto}}
\put(28.28,25.83){\special{em:moveto}}
\put(28.61,25.83){\special{em:lineto}}
\put(28.95,25.83){\special{em:moveto}}
\put(29.28,25.83){\special{em:lineto}}
\put(29.61,25.83){\special{em:moveto}}
\put(29.95,25.83){\special{em:lineto}}
\put(30.28,25.83){\special{em:moveto}}
\put(30.61,25.83){\special{em:lineto}}
\put(30.95,25.83){\special{em:moveto}}
\put(31.28,25.83){\special{em:lineto}}
\put(31.61,25.83){\special{em:moveto}}
\put(31.95,25.83){\special{em:lineto}}
\put(32.28,25.83){\special{em:moveto}}
\put(32.61,25.83){\special{em:lineto}}
\put(32.95,25.83){\special{em:moveto}}
\put(33.28,25.83){\special{em:lineto}}
\put(33.61,25.83){\special{em:moveto}}
\put(33.95,25.83){\special{em:lineto}}
\put(34.28,25.83){\special{em:moveto}}
\put(34.61,25.83){\special{em:lineto}}
\put(34.95,25.83){\special{em:moveto}}
\put(35.28,25.83){\special{em:lineto}}
\put(35.61,25.83){\special{em:moveto}}
\put(35.95,25.83){\special{em:lineto}}
\put(36.28,25.83){\special{em:moveto}}
\put(36.61,25.83){\special{em:lineto}}
\put(36.95,25.83){\special{em:moveto}}
\put(37.28,25.83){\special{em:lineto}}
\put(37.61,25.83){\special{em:moveto}}
\put(37.95,25.83){\special{em:lineto}}
\put(26.68,30.83){\special{em:moveto}}
\put(26.90,30.60){\special{em:lineto}}
\put(27.11,30.36){\special{em:moveto}}
\put(27.33,30.13){\special{em:lineto}}
\put(27.55,29.90){\special{em:moveto}}
\put(27.76,29.66){\special{em:lineto}}
\put(27.98,29.43){\special{em:moveto}}
\put(28.20,29.18){\special{em:lineto}}
\put(28.41,28.93){\special{em:moveto}}
\put(28.63,28.68){\special{em:lineto}}
\put(28.85,28.43){\special{em:moveto}}
\put(29.06,28.20){\special{em:lineto}}
\put(29.28,27.95){\special{em:moveto}}
\put(29.50,27.71){\special{em:lineto}}
\put(29.71,27.46){\special{em:moveto}}
\put(29.93,27.23){\special{em:lineto}}
\put(30.15,26.98){\special{em:moveto}}
\put(30.36,26.73){\special{em:lineto}}
\put(30.58,26.48){\special{em:moveto}}
\put(30.80,26.23){\special{em:lineto}}
\put(31.01,25.98){\special{em:moveto}}
\put(31.23,25.73){\special{em:lineto}}
\put(31.45,25.48){\special{em:moveto}}
\put(31.66,25.25){\special{em:lineto}}
\put(31.88,25.00){\special{em:moveto}}
\put(32.10,24.76){\special{em:lineto}}
\put(32.31,24.53){\special{em:moveto}}
\put(32.53,24.28){\special{em:lineto}}
\put(32.75,24.03){\special{em:moveto}}
\put(32.96,23.78){\special{em:lineto}}
\put(33.18,23.53){\special{em:moveto}}
\put(33.40,23.28){\special{em:lineto}}
\put(33.61,23.03){\special{em:moveto}}
\put(33.83,22.78){\special{em:lineto}}
\put(34.05,22.55){\special{em:moveto}}
\put(34.26,22.30){\special{em:lineto}}
\put(34.48,22.06){\special{em:moveto}}
\put(34.70,21.81){\special{em:lineto}}
\put(34.91,21.58){\special{em:moveto}}
\put(35.13,21.33){\special{em:lineto}}
\put(35.35,21.08){\special{em:moveto}}
\put(35.56,20.83){\special{em:lineto}}
\put(35.78,20.58){\special{em:moveto}}
\put(36.00,20.33){\special{em:lineto}}
\put(36.21,20.08){\special{em:moveto}}
\put(36.43,19.83){\special{em:lineto}}
\put(36.65,19.60){\special{em:moveto}}
\put(36.86,19.35){\special{em:lineto}}
\put(37.08,19.11){\special{em:moveto}}
\put(37.30,18.88){\special{em:lineto}}
\put(37.51,18.63){\special{em:moveto}}
\put(37.73,18.38){\special{em:lineto}}
\put(37.95,18.13){\special{em:moveto}}
\put(38.16,17.88){\special{em:lineto}}
\put(38.38,17.63){\special{em:moveto}}
\put(38.60,17.38){\special{em:lineto}}
\put(38.81,17.13){\special{em:moveto}}
\put(39.03,16.88){\special{em:lineto}}
\put(39.25,16.65){\special{em:moveto}}
\put(39.46,16.41){\special{em:lineto}}
\put(39.68,16.16){\special{em:moveto}}
\put(39.90,15.93){\special{em:lineto}}
\put(40.11,15.68){\special{em:moveto}}
\put(40.33,15.43){\special{em:lineto}}
\put(40.55,15.18){\special{em:moveto}}
\put(40.76,14.93){\special{em:lineto}}
\put(40.98,14.68){\special{em:moveto}}
\put(41.20,14.43){\special{em:lineto}}
\put(41.41,14.18){\special{em:moveto}}
\put(41.63,13.95){\special{em:lineto}}
\put(41.85,13.70){\special{em:moveto}}
\put(42.06,13.46){\special{em:lineto}}
\put(42.28,13.23){\special{em:moveto}}
\put(42.50,12.98){\special{em:lineto}}
\put(42.71,12.73){\special{em:moveto}}
\put(42.93,12.48){\special{em:lineto}}
\put(43.15,12.23){\special{em:moveto}}
\put(43.36,11.98){\special{em:lineto}}
\put(43.58,11.73){\special{em:moveto}}
\put(43.80,11.48){\special{em:lineto}}
\put(44.01,11.23){\special{em:moveto}}
\put(44.23,11.00){\special{em:lineto}}
\put(44.45,10.76){\special{em:moveto}}
\put(44.66,10.51){\special{em:lineto}}
\put(44.88,10.28){\special{em:moveto}}
\put(45.10,10.03){\special{em:lineto}}
\put(45.31,9.78){\special{em:moveto}}
\put(45.53,9.53){\special{em:lineto}}
\put(45.75,9.28){\special{em:moveto}}
\put(45.96,9.03){\special{em:lineto}}
\put(46.18,8.78){\special{em:moveto}}
\put(46.40,8.53){\special{em:lineto}}
\put(46.61,8.30){\special{em:moveto}}
\put(46.83,8.05){\special{em:lineto}}
\put(47.05,7.81){\special{em:moveto}}
\put(47.26,7.56){\special{em:lineto}}
\put(47.48,7.33){\special{em:moveto}}
\put(47.70,7.08){\special{em:lineto}}
\put(47.91,6.83){\special{em:moveto}}
\put(48.13,6.58){\special{em:lineto}}
\put(48.35,6.33){\special{em:moveto}}
\put(48.56,6.08){\special{em:lineto}}
\put(48.78,5.83){\special{em:moveto}}
\put(49.00,5.58){\special{em:lineto}}
\put(49.21,5.35){\special{em:moveto}}
\put(49.43,5.10){\special{em:lineto}}
\put(49.65,4.86){\special{em:moveto}}
\put(49.86,4.63){\special{em:lineto}}
\put(50.08,4.38){\special{em:moveto}}
\put(50.30,4.13){\special{em:lineto}}
\put(50.51,3.88){\special{em:moveto}}
\put(50.73,3.63){\special{em:lineto}}
\put(50.95,3.38){\special{em:moveto}}
\put(51.16,3.13){\special{em:lineto}}
\put(51.38,2.88){\special{em:moveto}}
\put(51.60,2.65){\special{em:lineto}}
\put(51.81,2.40){\special{em:moveto}}
\put(52.03,2.16){\special{em:lineto}}
\put(52.25,1.91){\special{em:moveto}}
\put(52.46,1.68){\special{em:lineto}}
\put(52.68,1.43){\special{em:moveto}}
\put(52.90,1.18){\special{em:lineto}}
\put(53.11,0.93){\special{em:moveto}}
\put(53.21,0.83){\special{em:lineto}}
\put(28.15,19.16){\special{em:moveto}}
\put(28.26,19.46){\special{em:lineto}}
\put(28.38,19.76){\special{em:moveto}}
\put(28.50,20.06){\special{em:lineto}}
\put(28.61,20.36){\special{em:moveto}}
\put(28.73,20.66){\special{em:lineto}}
\put(28.85,20.96){\special{em:moveto}}
\put(28.96,21.26){\special{em:lineto}}
\put(29.08,21.56){\special{em:moveto}}
\put(29.20,21.86){\special{em:lineto}}
\put(29.31,22.16){\special{em:moveto}}
\put(29.41,22.46){\special{em:lineto}}
\put(29.53,22.76){\special{em:moveto}}
\put(29.65,23.06){\special{em:lineto}}
\put(29.75,23.36){\special{em:moveto}}
\put(29.86,23.66){\special{em:lineto}}
\put(29.98,23.96){\special{em:moveto}}
\put(30.10,24.26){\special{em:lineto}}
\put(30.20,24.56){\special{em:moveto}}
\put(30.31,24.86){\special{em:lineto}}
\put(30.43,25.16){\special{em:moveto}}
\put(30.53,25.46){\special{em:lineto}}
\put(30.65,25.76){\special{em:moveto}}
\put(30.76,26.06){\special{em:lineto}}
\put(30.86,26.36){\special{em:moveto}}
\put(30.98,26.66){\special{em:lineto}}
\put(31.10,26.96){\special{em:moveto}}
\put(31.21,27.26){\special{em:lineto}}
\put(31.31,27.56){\special{em:moveto}}
\put(31.43,27.86){\special{em:lineto}}
\put(31.55,28.16){\special{em:moveto}}
\put(31.66,28.46){\special{em:lineto}}
\put(31.76,28.76){\special{em:moveto}}
\put(31.88,29.06){\special{em:lineto}}
\put(32.00,29.36){\special{em:moveto}}
\put(32.11,29.66){\special{em:lineto}}
\put(32.21,29.96){\special{em:moveto}}
\put(32.33,30.26){\special{em:lineto}}
\put(32.45,30.56){\special{em:moveto}}
\put(32.56,30.86){\special{em:lineto}}
\put(32.66,31.16){\special{em:moveto}}
\put(32.78,31.46){\special{em:lineto}}
\put(32.90,31.76){\special{em:moveto}}
\put(33.01,32.06){\special{em:lineto}}
\put(33.11,32.36){\special{em:moveto}}
\put(33.23,32.66){\special{em:lineto}}
\put(33.33,32.96){\special{em:moveto}}
\put(33.45,33.26){\special{em:lineto}}
\put(33.56,33.56){\special{em:moveto}}
\put(33.68,33.86){\special{em:lineto}}
\put(33.78,34.16){\special{em:moveto}}
\put(33.90,34.46){\special{em:lineto}}
\put(34.01,34.76){\special{em:moveto}}
\put(34.13,35.06){\special{em:lineto}}
\put(34.23,35.36){\special{em:moveto}}
\put(34.35,35.66){\special{em:lineto}}
\put(34.46,35.96){\special{em:moveto}}
\put(34.58,36.26){\special{em:lineto}}
\put(34.68,36.56){\special{em:moveto}}
\put(34.80,36.86){\special{em:lineto}}
\put(34.91,37.16){\special{em:moveto}}
\put(35.01,37.46){\special{em:lineto}}
\put(35.13,37.76){\special{em:moveto}}
\put(35.25,38.06){\special{em:lineto}}
\put(35.36,38.36){\special{em:moveto}}
\put(35.46,38.66){\special{em:lineto}}
\put(35.58,38.96){\special{em:moveto}}
\put(35.70,39.26){\special{em:lineto}}
\put(35.81,39.56){\special{em:moveto}}
\put(35.91,39.86){\special{em:lineto}}
\put(36.03,40.16){\special{em:moveto}}
\put(36.15,40.46){\special{em:lineto}}
\put(36.25,40.76){\special{em:moveto}}
\put(36.36,41.06){\special{em:lineto}}
\put(36.48,41.36){\special{em:moveto}}
\put(36.58,41.66){\special{em:lineto}}
\put(36.70,41.96){\special{em:moveto}}
\put(36.81,42.26){\special{em:lineto}}
\put(36.93,42.56){\special{em:moveto}}
\put(37.03,42.86){\special{em:lineto}}
\put(37.15,43.16){\special{em:moveto}}
\put(37.26,43.46){\special{em:lineto}}
\put(37.38,43.76){\special{em:moveto}}
\put(37.48,44.06){\special{em:lineto}}
\put(37.60,44.36){\special{em:moveto}}
\put(37.71,44.66){\special{em:lineto}}
\put(37.83,44.96){\special{em:moveto}}
\put(37.93,45.26){\special{em:lineto}}
\put(38.05,45.56){\special{em:moveto}}
\put(38.15,45.83){\special{em:lineto}}
\put(19.88,25.83){\special{em:moveto}}
\put(33.21,32.50){\special{em:lineto}}
\put(39.88,15.83){\special{em:lineto}}
\put(19.88,25.83){\special{em:lineto}}
\end{picture}

%% file: biblconf.tex
\bibliographystyle{amsplane}

%% file: index.tex
\begin{theindex}

  \item $\hol B_m$\hfill 55
  \item $t(X)$\hfill 11
  \item ${\mathbb C}^{\bigcirc\hskip-6pt k}$\hfill 21
  \item ${\mathbf G}_m^0$
    \subitem automorphisms of\hfill 69
    \subitem endomorphisms of\hfill 68
    \subitem proper endomorphisms of\hfill 68
  \item ${\mathbf G}_m^0$\hfill 49, 67
  \item ${\mathbf S}(n)$ action in $L({\mathbf q}_n)$\hfill 38
  \item ${\mathbf S}(n)$ orbits of simplices of cross ratios\hfill 39
  \item ${\mathbf S}(n)$ orbits of simplices of simple ratios\hfill 39
  \item ${\mathbf V}^{\prime}(0,l)$\hfill 22
  \item ${\mathcal C}^{m-1}_b({\mathbb C})$
    \subitem automorphisms of\hfill 69
  \item ${\mathcal C}^{m-1}_b({\mathbb C})$\hfill 48, 49, 57
  \item ${\mathcal C}^{n-1}_b({\mathbb C})$
    \subitem endomorphisms of\hfill 68
    \subitem proper endomorphisms of\hfill 68
  \item ${\mathcal C}^{n-1}_b({\mathbb C})$\hfill 67
  \item ${\mathcal C}_o^m({\mathbb C}^{\bigcirc\hskip-6pt 2}\,)$\hfill 
		21
  \item ${\mathcal D}(w)$\hfill 49
  \item 1st Kaliman theorem\hfill 66
  \item 2nd Kaliman theorem\hfill 69

  \indexspace

  \item Algebraic functions with the universal branch locus\hfill 60
  \item Artin Lemma\hfill 28
  \item Automorphisms of ${\mathbf G}_m^0$\hfill 69
  \item Automorphisms of ${\mathcal C}^n(X)$\hfill 64
  \item Automorphisms of ${\mathcal C}^{m-1}_b({\mathbb C})$\hfill 69

  \indexspace

  \item Balanced configuration space ${\mathbf G}_m^0$\hfill 49, 67
  \item Balanced configuration space ${\mathcal C}^{n-1}_b({\mathbb C})$\hfill 
		67
  \item Balanced configuration spaces
    \subitem automorphisms of\hfill 69
    \subitem endomorphisms of\hfill 68
    \subitem proper endomorphisms of\hfill 67, 68
    \subitem special\hfill 67
  \item Balanced configuration spaces\hfill 48, 57, 67
  \item Binary form
    \subitem discriminant of\hfill 10
    \subitem non-degenerate
      \subsubitem projective\hfill 10
    \subitem non-degenerate\hfill 10
  \item Binary form\hfill\phantom\hfill 10
  \item Braid groups
    \subitem Artin braid group
      \subsubitem pure Artin braid group\hfill 7
    \subitem Artin braid group\hfill 7
    \subitem canonical presentations\hfill 25
    \subitem sphere braid group
      \subsubitem pure sphere braid group\hfill 7
    \subitem sphere braid group\hfill 7
  \item Braid groups\hfill \phantom\hfill{7}

  \indexspace

  \item Compositions of multivalued functions\hfill 60
  \item Configuration spaces
    \subitem balanced
      \subsubitem automorphisms of\hfill 69
      \subsubitem endomorphisms of\hfill 68
      \subsubitem proper endomorphisms of\hfill 67, 68
      \subsubitem special\hfill 67
    \subitem balanced\hfill 67
    \subitem holomorphic universal covering of\hfill 21, 23
    \subitem ordered\hfill 7
    \subitem unordered\hfill 7
  \item Configuration spaces\hfill \phantom\hfill{7}
  \item Cyclic morphisms ${\mathcal C}^n({\mathbb C}) \to{\mathcal C}^k({\mathbb C})$\hfill 
		51

  \indexspace

  \item Dimension of $L_{\vartriangle}({\mathcal C}_o^n({\mathbb C}))$\hfill 
		38
  \item Dimension of $L_{\vartriangle}({\mathcal C}_o^n({\mathbb{CP}^1}))$\hfill 
		38

  \indexspace

  \item Eisenstein's automorphism
    \subitem Cayley's form\hfill 81
    \subitem Feler's form\hfill 81
    \subitem original form\hfill 80
  \item Eisenstein\hfill \phantom\hfill{11}
  \item Endomorphism
    \subitem tame\hfill 8
  \item Endomorphisms of ${\mathbf{SG}}_n$\hfill 69
  \item Endomorphisms of ${\mathcal C}_o^m({\mathbb C}\setminus\{0,1\})$\hfill 
		66
  \item Equivariance condition
    \subitem strict\hfill 15
  \item Equivariance condition\hfill 14

  \indexspace

  \item Finite covers of Liouville spaces\hfill 45

  \indexspace

  \item Gorin's relation\hfill 26

  \indexspace

  \item Holomorphic elements of $\pi_i({\mathcal C}^m({\mathbb C}))$\hfill 
		55
  \item Holomorphic elements of homotopy groups\hfill 53, 54
  \item Holomorphic functions omitting two values
    \subitem on ${\mathcal C}_o^n(X)$\hfill 18
    \subitem simplicial complex of\hfill 18, 19, 33
  \item Holomorphic maps
    \subitem of ultra-Picard space to ${\mathcal C}^m(X)$\hfill 45
    \subitem Solvable
      \subsubitem of Liouville space to ${\mathcal C}^m(X)$\hfill 45
  \item Holomorphic maps of ultra-Picard space to ${\mathbf G}_m$\hfill 
		50
  \item Holomorphic part $\hol B_m$ of $B_m$\hfill 55
  \item Holomorphic part $\hol\pi_i(Y)$ of the homotopy group $\pi_i(Y)$\hfill 
		54
  \item Homomorphism
    \subitem abelian\hfill 8
    \subitem cyclic\hfill 8
    \subitem solvable\hfill 8
    \subitem transitive\hfill 27

  \indexspace

  \item Image of non-cyclic endomorphism of ${\mathcal C}^n({\mathbb C})$\hfill 
		57
  \item Invariance property\hfill 13
  \item Invariant maps to $\Aut X$\hfill 64

  \indexspace

  \item Kaliman theorem
    \subitem 1st\hfill 66
    \subitem 2st\hfill 69

  \indexspace

  \item Lifting property\hfill 12
  \item Liouville and Picard spaces\hfill 45

  \indexspace

  \item Map
    \subitem abelian\hfill 8
    \subitem cyclic\hfill 8
    \subitem disjoint\hfill 11
    \subitem linked\hfill 11
    \subitem solvable\hfill 8
  \item Mixed simplices\hfill 38
  \item Model holomorphic maps ${\mathbb C}^*\to{\mathcal C}^m({\mathbb C})$\hfill 
		51
  \item Moduli space
    \subitem automorphisms of\hfill 65
    \subitem ordered\hfill 65
  \item Moduli space\hfill \phantom\hfill{9}
  \item Morphism
    \subitem degenerate tame\hfill 8
    \subitem invariant\hfill 16
    \subitem orbit-like\hfill 8

  \indexspace

  \item Normal form of simplices\hfill 39

  \indexspace

  \item Parameterization of ${\mathbb C}^*$ orbits in ${\mathbf G}_m^0$\hfill 
		49
  \item Proper divisors\hfill 35

  \indexspace

  \item Simple and cross ratios\hfill 34
  \item Simplices of cross ratios\hfill 36
  \item Simplices of simple ratios\hfill 36
  \item Simplicial complex of holomorphic functions omitting two values\hfill 
		33
  \item Simplicial complexes of simple and cross ratios\hfill 35
  \item Solvable holomorphic maps of Liouville space to ${\mathbf G}_m$\hfill 
		50
  \item Special homomorphisms $B_m\to B_k$\hfill 53
  \item Special presentations of $B_m$\hfill 53
  \item Special system of generators in $B_m$\hfill 53
  \item Sphere braid group $B_n(S^2)$\hfill 26
  \item Standard normal series for $PB_n(X)$\hfill 23
  \item Standard system of generators in $B_m$\hfill 53

  \indexspace

  \item Tame algebraic ensembles\hfill 60
  \item Tame endomorphism
    \subitem of ${\mathcal C}_o^n(X)$
      \subsubitem degenerate\hfill 16
    \subitem of ${\mathcal C}_o^n(X)$\hfill 16
  \item Tame endomorphism\hfill 8
  \item Teichm{\" u}ller space\hfill 21
  \item Theorem
    \subitem Artin Theorem\hfill 27
    \subitem Artin-Markov-Fadell-Neuwirth\hfill 24
    \subitem Coherence Theorem
      \subsubitem explanation of\hfill 17
      \subsubitem proof of\hfill 41, 42
    \subitem Coherence Theorem\hfill 17
    \subitem Commutator Theorem
      \subsubitem proof of\hfill 26
    \subitem Commutator Theorem\hfill 13
    \subitem Composition Theorem
      \subsubitem proof of\hfill 28
    \subitem Composition Theorem\hfill 13
    \subitem Equivariance Theorem
      \subsubitem explanation of\hfill 15
      \subsubitem proof of\hfill 29
    \subitem Equivariance Theorem\hfill 15
    \subitem Equivariant Map Theorem
      \subsubitem proof of\hfill 17
    \subitem Equivariant Map Theorem\hfill 16
    \subitem Homomorphisms $B_n\to B_k$
      \subsubitem induced by morphisms\hfill 53, 55
    \subitem Hubbard-Earl-Kra Theorem\hfill 42
    \subitem Invariance Theorem
      \subsubitem explanation of\hfill 13
      \subsubitem proof of\hfill 29
    \subitem Invariance Theorem\hfill 13
    \subitem Kernel Theorem
      \subsubitem proof of\hfill 29
    \subitem Kernel Theorem\hfill 13
    \subitem Lifting Theorem
      \subsubitem proof of\hfill 29
    \subitem Lifting Theorem\hfill 13
    \subitem Linked Map Theorem
      \subsubitem for ${\mathcal C}_o^m({\mathbb C}\setminus\{0,1\})$\hfill 
		43, 44
      \subsubitem for ${\mathcal C}_o^n(X)$\hfill 44, 45
    \subitem Linked Map Theorem\hfill 11, 42
    \subitem Murasugi Theorem\hfill 52
    \subitem on $n$-valued algebraic functions with the universal branch locus\hfill 
		59
    \subitem on algebraic functions containing the universal one\hfill 
		60
    \subitem on compositions without superfluous branch points\hfill 
		62, 63
    \subitem on invariant morphisms to $\Aut X$ and automorphisms of ${\mathcal C}^n(X)$\hfill 
		64
    \subitem on non-cyclic special homomorphisms $B_n\to B_k$\hfill 55
    \subitem Surjectivity Theorem
      \subsubitem proof of\hfill 27
    \subitem Surjectivity Theorem\hfill 14
    \subitem Tame Map Theorem
      \subsubitem proof of\hfill 42
      \subsubitem strategy of the proof\hfill 12
    \subitem Tame Map Theorem\hfill 9
  \item Torsion of $B_m/CB_m$ and $B_m(S^2)$\hfill 52

  \indexspace

  \item Ultra-Picard spaces\hfill 45
  \item Universal Teichm{\" u}ller family\hfill 21

\end{theindex}

%% file: confspace.bbl
\begin{thebibliography}{FadBusk62} 
\baselineskip12pt

\bibitem[Arn70a]{Arn70a} 
V. I. Arnold,
{\em Cohomology classes of algebraic functions
invariant under Tschirnhausen transformation},
Funct. Anal. Appl. {\bf 4}, n. 1 (1970), 74--75.

\bibitem[Arn70b]{Arn70b} 
V. I. Arnold,
{\em On certain topological invariants of algebraic functions},
Trudy Mosk. Mat. Obshch. {\bf 21} (1970), 27--46
(in Russian). {\it English translation:}

\bibitem[Arn70c]{Arn70c} 
V. I. Arnold,
{\em Topological invariants of algebraic functions, {\rm II}},
Funkt. Analiz i Priloz. {\bf 4}, no. 2 (1970), 1--9
(in Russian). {\it English translation:}
Funct. Anal. Appl. {\bf 4}, no. 2 (1970), 91--98.

\bibitem[Art25]{Art25} 
E. Artin,
{\em Theorie der Z\"opfe},
Abh. Math. Sem. Univ. Hamburg {\bf 4} (1925), 47--72.

\bibitem[Art47a]{Art47a} 
E. Artin,
{\em Theory of braids},
Annals of Math. {\bf 48}, no. 1 (1947), 101--126.

\bibitem[Art47b]{Art47b} 
E. Artin,
{\em Braids and permutations},
Ann. of Math. {\bf 48}, no. 3 (1947), 643--649.

\bibitem[BeRo86]{BeRo86}
L. Bers, H. L. Royden
{\em Holomorphic families of injections},
Acta Math. {\bf 157} (1986), no. 3-4, 259--286

\bibitem[Bir74]{Bir74} 
J. S. Birman,
Braids, links and mapping class groups.
Annals of Math. Studies {\bf 82}.
Princeton Univercity Press, Princeton, N. J., 1974.

\bibitem[BorNar67]{BorNar67} 
A. Borel, R. Narasimhan,
{\em Uniqueness conditions for certain holomorphic mappings},
Invent. Math. {\bf 2} (1967), 247--255.

\bibitem[Bohn47]{Bohn47} 
F. Bohnenblust,
{\em The algebraic braid group},
Annals of Math. {\bf 48}, no. 1 (1947), 127--136.

\bibitem[Bro76]{Bro76}
F. E. Browder,
{\em Problems of Present Day Mathematics},
Mathematical Developments Arising from Hilbert Problems,
Proc. Symp. Pure Math. {\bf 28} (F. E. Browder, ed.)
Amer. Math. Soc., 1976, 35--80.

\bibitem[Bur32]{Bur32} 
W. Burau,
{\em \"Uber Zopfinvarianten},
Abh. Math. Sem. Univ. Hamburg {\bf 9} (1932), 117--124.


\bibitem[Car40]{Car40} 
H. Cartan,
{\em Sur les matrices holomorphes de $n$ variables complexes},
J. Math. Pures et Appl. {\bf 19} (1940), 1--26.

\bibitem[Chow48]{Chow48} 
W.-L. Chow,
{\em On the algebraic braid group},
Annals of Math. {\bf 49}, no. 3 (1948), 654--658.

\bibitem[DanGiz72]{DanGiz72}
V. I. Danilov and M. H. Gizatullin,
{\em Automorphisms of affine surfaces} II,
Math. USSR Izvestija {\bf 41}, no. 2 (1977), 

\bibitem[Del72]{Del72} 
P. Deligne,
{\em Les immeubles des groupes de tresses g\'en\'eralis\'es},
Invent. Math. {\bf 17}, no. 4 (1972),273--302.

\bibitem[DyeGro81]{DyeGro81} 
J. Dyer and E. Grossman,
{\em The automorphism groups of the braid groups},
Amer. J. Math. {\bf 103}, no. 6 (1981), 1151--1169.

\bibitem[Eis1844]{Eis1844} 
G. Eisenstein,
{\em \"Uber eine merkw\"urdige identische
Gleichung},
Journal f\"ur die reine und angewandte Mathematik {\bf 27} (1844),
105--106 (see also: Gotthold Eisenstein, Mathematische Werke,
Band I, Chelsea, New York, 1975,1989, pp. 26--27).  

\bibitem[EarKra74]{EarKra74}
C. J. Earle and I. Kra
{\em On holomorphic mappings between Teichm\"uller spaces},
in book: Contributions to Analysis. A collection of papers dedicated
to L. Bers. L. V. Ahlfors etc., ed.,
Academic Press, New-York 1974, 107--124.

\bibitem[EarKra76]{EarKra76}
C. J. Earle and I. Kra
{\em On selections of some families of closed Riemann surfaces},
Acta Math. {\bf 137} (1976), 49--79

\bibitem[Enr23]{Enr23}
{\em Sulla costrucione delle funcioni algebriche di due variabili
possedenti una data curva di diramazione}, Annali di matematica pure ed
applicata, Ser. 4, {\bf 1} (1923), 185--198.

\bibitem[FadBusk61]{FadBusk61} 
E. Fadell and J. Van Buskirk,
{\em On the braid groups of $E^2$ and $S^2$},
Bull. Amer. Math. Soc. {\bf 67}, no. 2 (1961), 211--213.

\bibitem[FadBusk62]{FadBusk62} 
E. Fadell and J. Van Buskirk,
{\em The braid groups of $E^2$ and $S^2$},
Duke Math. J. {\bf 29}, no. 2 (1962), 243--257.

\bibitem[Fad62]{Fad62}
E. Fadell,
{\em Homotopy groups of configuration spaces and the string problem
of Dirac}, Duke Math. J. {\bf 29}, no. 2 (1962), 231--242.

\bibitem[FadNeu62]{FadNeu62} 
E. Fadell and L. Neuwirth,
{\em Configuration spaces},
Math. Scand. {\bf 10}, no. 1 (1962), 111--118.

\bibitem[FadHus01]{FadHus01}
E. R. Fadell and S. Y. Husseini, 
Geometry and topology of configuration spaces.
Springer, 2001.

\bibitem[FeiZie99]{FeiZie99}
E. M. Feichtner, M. Ziegler,
{\em The integral cohomology algebras of ordered configuration spaces
of spheres}, Documenta Math. {\bf 5} (2000), 115-139.

\bibitem[Fin45]{Fin45} 
P. Finsler,
{\"U}ber die Primzahlen zwischen $n$ und $2n$.
Speiser-Festschrift,
Orell--F\"ussli, Z\"urich, 1945.

\bibitem[Fo51]{Fo51}
L. R. Ford, Automorphic Functions, Chelsea, 1951.

\bibitem[FoxNeu62]{FoxNeu62} 
R. H. Fox and L. Neuwirth,
{\em The braid groups},
Math. Scand. {\bf 10}, no. 1 (1962), 119--126.

\bibitem[Fuks70]{Fuks70} 
D. B. Fuks,
{\em Cohomologies of a braid group mod 2},
Funct. Anal. Appl. {\bf 4}, no. 2 (1970), 91--98.

\bibitem[Gar]{Gar}
F. Garside,
{\em The braid group and other groups}, 
Quart. J. Math. Oxford {\bf 20}, n. 78 (1969), 235--254.  

\bibitem[Gol74]{Gol74} 
D. L. Goldsmith,
{\em Homotopy of braids in answer to a question of E. Artin},
Topology Conference,
Lect. Notes Math. {\bf 375} (R. F. Dickman, Jr. and
P. Fletcher, eds.),
Springer, 1974, 91--96.

\bibitem[GorLin69]{GorLin69} 
E. A. Gorin and V. Ya. Lin,
{\em Algebraic equations with continuous coefficients
and some problems of the algebraic theory of braids},
Matem. Sbornik {\bf 78\,(120)}, no. 4 (1969), 579--610
(in Russian). {\em English translation:}
Math. USSR Sbornik {\bf 7} (1969), 569--596.

\bibitem[GorLin74]{GorLin74} 
E. A. Gorin and V. Ya. Lin,
{\em On separable polynomials over commutative Banach algebras},
Dokl. Akad. Nauk SSSR {\bf 218}, no. 3 (1974), 505--508
(in Russian). {\em English translation:}
Soviet Math. Dokl. {\bf 15}, no. 5 (1974), 1357--1361.

\bibitem[Gra58]{Gra58} 
H. Grauert,
{\em Analytische Faserungen {\" u}ber holomorphvollst{\" a}ndigen
R{\" a}umen},
Math. Ann. {\bf 135} (1958), 263--273.

\bibitem[Hall]{Hall} 
M. Hall, 
The theory of groups.
McMillan Co., New York, 1959.

\bibitem[Hub72]{Hub72} 
J. H. Hubbard,
{\em Sur la non-existence de sections analytiques {\` a} la
courbe universelle de Teichm{\" uller}},
C.R. Acad. Sci. Paris S{\' e}r. A-B {\bf 274} (1972),
A978--A979.

\bibitem[Hub76]{Hub76} 
J. H. Hubbard,
{\em Sur les sections analytiques de la
courbe universelle de Teichm{\" uller}},
Mem. Amer. Math. Soc. {\bf 4} (1976), no.166, ix+137pp.

\bibitem[Kal75]{Kal75} 
S. I. Kaliman,
{\em Holomorphic universal covering of the space of
polynomials without multiple roots},
Funkt. Analiz i Priloz. {\bf 9}, no. 1 (1975), 71
(in Russian). {\em English translation:}
Funct. Anal. Appl. {\bf 9} (1975), no. 1, 67--68.

\bibitem[Kal76a]{Kal76a} 
S. I. Kaliman,
{\em Holomorphic endomorphisms of the manifold of complex
polynomials with discriminant $1$},
Uspehi Mat. Nauk {\bf 31} (1976), no. 1, 251--252 (in Russian).

\bibitem[Kal76b]{Kal76b} 
S. I. Kaliman,
{\em A holomorphic universal covering of
the space of polynomials without multiple roots},
Teor. Funkc., Funkcional. Anal. i Prilozen. Vyp. 28, (1977), 25--35
(in Russian). MR0590056 (58 \#28654)

\bibitem[Kal93]{Kal93}
Sh. I. Kaliman
{\em The holomorphic universal covers of the spaces of
polynomials without multiple roots}, 
Selected translations. Selecta Math. Soviet. {\bf 12} (1993), no. 4,
395--405. 

%
\bibitem[Kle]{Kle}
F. Klein,
Vorlesungen {\" u}ber das Ikosaeder
und die Aufl{" o}sung der Gleichungen vom f{" u}nften Grade,
Birkh{\" a}user Verlag \& B. G. Teubner Verlagsgesellschaft,
1993.             


\bibitem[Koj02]{Koj02}
S. Kojima
{\em Complex hyperbolic cone structures on the configuration spaces},
Rend. Istit. Mat. Univ. Trieste 32 (2001), suppl. 1, 149--163 (2002).
(Electronic version: arXiv:math.GT/9907147)

\bibitem[Kurth97]{Kurth97} 
A. Kurth,
{\em $SL_2$-equivariant polynomial automorphisms of binary forms},
Ann. Inst. Fourier, Grenoble {\bf 47}, no. 2 (1997), 
585--597.

\bibitem[Lin71]{Lin71} 
V. Ya. Lin,
{\em Algebroidal functions and holomorphic elements
of homotopy groups of complex manifolds},
Dokl. Akad. Nauk  SSSR {\bf 201:1} (1971), 28--31
(in Russian). {\em English translation:}
Soviet Math. Dokl. {\bf 12}, n. 6 (1971), 1608--1612.

\bibitem[Lin72a]{Lin72a} 
V. Ya. Lin,
{\em On representations of the braid group by permutations},
Uspehi Matem. Nauk {\bf 27}, n. 3 (1972), 192.

\bibitem[Lin72b]{Lin72b} 
V. Ya. Lin,
{\em Algebraic functions with universal discriminant
manifolds}, Funkt. Analiz i Priloz. {\bf 6}, no. 1 (1972), 81--82
(in Russian). {\em English translation:}
Funct. Anal. Appl. {\bf 6}, n. 1 (1972), 73--75.

\bibitem[Lin72c]{Lin72c} 
V. Ya. Lin,
{\em On superpositions of algebraic functions},
Funkt. Analiz i Priloz. {\bf 6}, no. 3 (1972), 77--78
(in Russian). {\em English translation:}
Funct. Anal. Appl. {\bf 6}, n. 3 (1972), 240--241.

\bibitem[Lin74]{Lin74} 
V. Ya. Lin,
{\em Representations of braids by permutations},
Uspehi Mat. Nauk {\bf 29}, n. 1 (1974), 173--174
(in Russian).

\bibitem[Lin76]{Lin76} 
V. Ya. Lin,
{\em Superpositions of algebraic functions},
Funkt. Analiz i Priloz. {\bf 10}, no. 1 (1976), 37--45
(in Russian). {\em English translation:}
Funct. Anal. Appl. {\bf 10}, n. 1 (1976), 32--38.

\bibitem[Lin79]{Lin79} 
V. Ya. Lin,
{\em Artin braids and the groups and spaces connected
with them},
Itogi Nauki i Tekhniki, Algebra, Topologiya,
Geometriya {\bf 17}, VINITI, Moscow, 1979, pp. 159--227
(in Russian). {\em English translation:}
Journal of Soviet Math. {\bf 18} (1982), 736--788.

\bibitem[Lin88]{Lin88} 
V. Ya. Lin,
{\em Liouville coverings of complex spaces and amenable groups},
Matem. Sbornik {\bf 132\,(174)}, no. 2 (1987), 202--224
(in Russian). {\em English translation:}
Math. USSR Sbornik {\bf 60}, n. 1 (1988), 197--216.

\bibitem[Lin96a]{Lin96a} 
V. Ya. Lin,
{\em Around the $13$th Hilbert problem for algebraic functions},
Israel Math. Conference Proc. {\bf 9} (1996), 307--327.

\bibitem[Lin96b]{Lin96b} 
V. Ya. Lin,
{\em Braids, permutations, polynomials}--I,  
Preprint MPI 96-118, Max-Planck-Institut f\"ur Mathematik in Bonn, 
August 1996, 112pp.

\bibitem[LinZai98]{LinZai98}
V. Ya. Lin and M. Zaidenberg, 
{\em Liouville and Carath{\'e}odory coverings
in Riemannian and complex geometry}
in book: Voronezh Winter Mathematical School.
P. Kuchment \& V. Lin eds.,
Amer. Math. Soc. Transl. (2), {\bf 184} (1998), 111--130.

\bibitem[Lin03]{Lin03} 
V. Ya. Lin, 
{\em Configuration spaces of $\mathbb C$ and ${\mathbb{CP}}^1$:
some analytic properties}, Max-Planck-Institut f{\" u}r
Mathematik Preprint Series 2003 (98), Bonn 2003, 80pp.

\bibitem[Lin04a]{Lin04a} 
V. Ya. Lin, 
{\em Configuration spaces of $\mathbb C$ and ${\mathbb{CP}}^1$:
some analytic properties}, Electronic version arXiv:math.AG/0403120.

\bibitem[Lin04b]{Lin04b} 
V. Ya. Lin,
{\em Braids and permutations}, Electronic version arXiv:math

\bibitem[Mar45]{Mar45} 
A. A. Markov,
{\em Foundations of the algebraic theory of tresses},
Proc. Steklov Math. Institute {\bf 16}, 1945 
(in Russian).

\bibitem[Mur82]{Mur82} 
K.
Murasugi,
{\em Seifert fiber spaces and braid groups},
Proc. London Math. Soc. {\bf 44} (1982), 71--84.

\bibitem[MurKur99]{MurKur99} 
K. Murasugi, B. I. Kurpita,
A study of braids. Kluwer, 1999.
 
\bibitem[HNeum]{HNeum} 
H. Neumann,
Varieties of groups.
Springer, 1967.

\bibitem[Ore98]{Ore98}
S. Orevkov, {\em Ttransitive homomorphisms $B_7\to{\mathbf S}(14)$
and $B_8\to{\mathbf S}(16)$}, private letter, October 4, 1998.

\bibitem[Rams65]{Rams65} 
K. J. Ramspott,
{\em Stetige und holomorphe Schnitte in B\"undeln mit homogener Faser},
Math. Z. {\bf 89}, no. 3 (1965), 234--246.

\bibitem[Trost]{Trost} 
E. Trost,
Primzahlen.
Birkh\"auser Verlag, Basel--Stuttgart, 1968.

\bibitem[Vas88]{Vas88} 
V. A. Vassiliev,
{\em Braid group cohomologies and algorithm complexity},
Funct. Anal. Appl. {\bf 22}, no. 3 (1988), 182--189.

\bibitem[Vas98]{Vas98} 
V. A. Vassiliev,
Complements of Discriminants of Smooth Maps:
Topology and Applications.
Translations of Mathematical Monographs {\bf 98},
Amer. Math. Soc., Providence, Rhode Island, 1992.

\bibitem[Wi64]{Wi64} 
H. Wielandt,
Finite permutation groups.
Academic Press, 1964.

\bibitem[Zar29]{Zar29}
O. Zariski, {\em On the problem of existence of algebraic functions
of two variables possessing a given branch curve},
Amer. J. Math. {\bf 51} (1929).

\bibitem[Zar36]{Zar36}
O. Zariski, {\em On the Poincar{\' e} group of rational
plane curves}, Amer. J. Math. {\bf 58} (1936), 607--619. 

\bibitem[Zar37a]{Zar37a}
O. Zariski, {\em The topological discriminant group of a Riemann
surface of genus $p$}, Amer. J. Math. {\bf 59} (1937), 335--358.

\bibitem[Zar37b]{Zar37b}
O. Zariski, {\em A theorem on the Poincar{\' e} group of an algebraic
hypersurface}, Ann. Math. {\bf 38} (1937). 

\bibitem[Zin77]{Zin77}
V. Zinde,
{\em Holomorphic maps of regular orbit spaces of Coxeter
groups of series $B$ and $D$},
Sibirsk. Math. J. {\bf 18}, no. 5 (1977), 1016--1026 (in Russian).
{\em English translation:}


\end{thebibliography}
